\documentclass[
12pt, 
a4paper, 
oneside, 
headinclude,footinclude, 
BCOR5mm, 
]{book}

\usepackage[
nochapters, 
beramono, 
eulermath,
pdfspacing, 
dottedtoc 
]{classicthesis} 

\usepackage{arsclassica} 

\usepackage[T1]{fontenc} 

\usepackage[utf8]{inputenc} 

\usepackage{graphicx} 

\usepackage{enumitem} 

\usepackage{amsmath,amssymb,amsthm} 

\usepackage{varioref} 

\hypersetup{
colorlinks=true, breaklinks=true, bookmarks=true,bookmarksnumbered,
urlcolor=webbrown, linkcolor=RoyalBlue, citecolor=webgreen, 
pdftitle={}, 
pdfauthor={\textcopyright}, 
pdfsubject={}, 
pdfkeywords={}, 
pdfcreator={pdfLaTeX}, 
pdfproducer={LaTeX with hyperref and ClassicThesis} 
} 

\usepackage{bm}						
\usepackage{a4wide} 
\usepackage{multirow}
\usepackage{mathabx} 

\newtheorem{thm}{Theorem}
\newtheorem{lem}[thm]{Lemma}
\newtheorem{prop}[thm]{Proposition}
\newtheorem{conj}[thm]{Conjecture}
\newtheorem{coro}[thm]{Corollary}

\theoremstyle{definition}				

\newcommand{\R}{\mathbb{R}}
\newcommand{\C}{\mathbb{C}}
\newcommand{\N}{\mathbb{N}}
\newcommand{\Z}{\mathbb{Z}}
\newcommand{\Q}{\mathbb{Q}}
\newcommand{\K}{\mathbb{K}}
\newcommand{\Sph}{\mathbb{S}}
\newcommand{\V}{\detokenize{Vol}}
\newcommand{\e}{\detokenize{e}}

\numberwithin{equation}{chapter}		
\numberwithin{thm}{chapter}

\title{\normalfont\spacedallcaps{Homogeneous Forms Inequalities}} 

\author{\spacedlowsmallcaps{Faustin ADICEAM* \& Oscar MARMON\textsuperscript{1}}} 

\date{} 

\begin{document}


\renewcommand{\sectionmark}[1]{\markright{\spacedlowsmallcaps{#1}}} 
\lehead{\mbox{\llap{\small\thepage\kern1em\color{halfgray} \vline}\color{halfgray}\hspace{0.5em}\rightmark\hfil}} 

\renewcommand{\bibname}{\centerline{References}} 

\pagestyle{scrheadings} 


\maketitle 

\setcounter{tocdepth}{3} 

\newpage
\thispagestyle{empty}
$\quad$
\newpage

{
\cleardoublepage
\let\clearpage\relax
\vspace{3mm}
\flushright{\textit{To all our students}}
\setcounter{page}{1}
\section*{Abstract} 
}

Let $F_1(\bm{x}), \dots, F_p(\bm{x})$ be a set of $p\ge 1$ real homogeneous forms in $n\ge 2$ variables, all of the same degree. Given a compact set $K\subset\R^n$ and a nondecreasing map $T\mapsto b(T)\ge 1$, consider the region of the space determined by the inequalities $$  \left|F_i(\bm{x})\right|\le b(T)\quad \textrm{ for all }\quad i=1,\dots, p\quad \textrm{ when }\quad \bm{x}\in T\cdot K.$$ The following intertwined results are established~: \\

\begin{itemize}
\item The asymptotic behavior as $T$ tends to infinity  of the volume of this region is determined  under various regularity assumptions on the set $K$. The growth rate of the volume is furthermore related to the largest root of the Sato-Bernstein  polynomial associated to a homogeneous form derived from $F_1(\bm{x}), \dots, F_p(\bm{x})$, and to its multiplicity;

\item In relation with this result, a classical open problem in the theory of $\mathcal{D}$-modules is settled~: it is known that the largest root of the Sato-Bernstein polynomial of a multivariate polynomial coincides with the opposite of the smallest pole of the complex meromorphic distribution attached to it. The cor\-res\-ponding multiplicity and order are proved to also coincide;

\item In the case that the $p$ homogeneous forms under consideration are twisted by composing them with unimodular matrices, the solubility of the  above set of in\-e\-qua\-li\-ties in integer lattice points with bounded height  is established in the form of a metric statement depending on the convergence of a volume sum. This ef\-fec\-tive, metric and uniform generalisation of the Oppenheim Conjecture to the case of a system of homogeneous forms answers a question raised by Athreya and Margulis (2018). So does the related statement determining the asymptotic behavior, generic  in a metric sense, of the number of integer solutions to this system of twisted inequalities;

\item In the deterministic case where the forms are fixed, an upper bound is established for the function counting the number of integer lattice points satisfying the system of inequalities. The main term in this bound is related to the volume of the region of the space under consideration and is sharp under suitable assumptions. The error term is furthermore shown to admit a power saving provided that a quantitative measure of flatness emerging from geometric tomography is large enough. This settles a conjecture stated by Sarnak (1997).
\end{itemize}

\newpage
$\quad$
\newpage

\tableofcontents 




\let\thefootnote\relax\footnotetext{*\textit{Laboratoire d’analyse et de mathématiques appliquées, Université Paris-Est Créteil, Créteil, France.} \texttt{faustin.adiceam@u-pec.fr}}

\let\thefootnote\relax\footnotetext{\textsuperscript{1} \textit{Centre for Mathematical Sciences, Lund University, Box 118, 221 00 Lund, Sweden.} \texttt{oscar.marmon@math.lu.se}}



\chapter*{Notations}

\addcontentsline{toc}{chapter}{\protect\numberline{}Notation}

\paragraph{\textbf{Notations in Euclidean spaces~:}}
\begin{itemize}
\item Vectors are always denoted with bold characters (e.g., $\bm{x}, \bm{y}, \bm{z}$...). 
\item The Euclidean norm is denoted by $\left\|\,\cdot\,\right\|$. 
\item Given $r>0$, $B_n(r)$ stands for the closed Euclidean ball centered at the origin with radius $r$ in dimension $n\ge 1$, and  $B_n$ for $B_n(1)$.  
\end{itemize}

\paragraph{\textbf{Algebraic Notations~:}} 
\begin{itemize}
\item The set $\K$ stands for either the field of reals $\R$ or the field of complex numbers $\C$. 
\item The notation $\N$ refers to the set of (strictly) positive integers and $\N_0$ to the set of nonnegative integers (i.e.~$\N_0=\N\cup\{0\}$~).
\item Given a set of $p\ge 1$ homogeneous forms $\bm{F}(\bm{x})=\left(F_1(\bm{x}), \dots, F_p(\bm{x})\right)$, their common set of zeros over $\K$ is denoted by $\mathcal{Z}_\K(\bm{F})$; that is, $$\mathcal{Z}_\K(\bm{F})\;=\; \left\{\bm{y}\in \K^n\; :\; \forall i= 1,\dots,  p, \; \; F_i(\bm{x})=0\right\}.$$ 
\item Given a subset $A\subset \K^n$, the above notation is extended by setting $\mathcal{Z}_A(\bm{F})\;=\; \mathcal{Z}_\K(\bm{F})\cap A.$
\end{itemize}

\paragraph{\textbf{Analytic Notations~:}} 
\begin{itemize}
\item The partial derivative with respect to the $i^{th}$ coordinate in dimension $n$ (where $1\le i\le n$) is denoted by $\partial_i$.   The gradient operator is then  $\nabla=\left(\partial_1, \dots, \partial_n\right)$. 
\item When $\bm{k}=\left(k_1, \dots, k_n\right)$ is an $n$-tuple of nonnegative integers, define its size by $\left|\bm{k}\right|\;=\;\sum_{i=1}^nk_i$ and the differential operator $\bm{\partial}^{\bm{k}}$  by $\bm{\partial}^{\bm{k}}=\prod_{i=1}^{n}\partial_i^{k_i}$. 
\item Similarly, given $\bm{x}\in\R^n$, $\bm{x}^{\bm{k}}$ is shorthand notation for $\prod_{i=1}^{n}x_i^{k_i}$.
\item Given a subset $A\subset\R^n$,  its characteristic function is denoted by $\chi_{A}$ . Similarly, given a property $\mathcal{P}$ depending on a vector $\bm{y}\in\R^n$,  $\chi_{\left\{\mathcal{P}(\bm{y})\right\}}$ is the boolean function equal to 1 when $\mathcal{P}(\bm{y})$ holds and 0 otherwise.
\item When $k\in\N_0\cup\{\infty\}$,  $\mathcal{C}^{k}(\K^n)$ stands for the set of $k$ times continuously dif\-fe\-ren\-tiable functions over $\K^n$ and $\mathcal{C}_{c}^{k}(\K^n)$  for the subset of such functions with compact support; the support of an element $\psi$ in $\mathcal{C}_{c}^{k}(\K^n)$ is then denoted by $\textrm{Supp }\psi$. 
\item Given two functions $f$ and $g$ depending on a variable $\bm{y}\in\K^n$, the notation \mbox{$f(\bm{y})\ll g(\bm{y})$} is used to mean the existence of a real $c>0$, referred to as the \emph{implicit constant}, such that $\left|f(\bm{y})\right|\le c\cdot \left|g(\bm{y})\right|$ for all admissible values of the variable $\bm{y}$. This is equivalent to the notation $f(\bm{y})=O\left(g(\bm{y})\right)$. Also, the relation $f(\bm{y})\asymp g(\bm{y})$ means that  the two conditions $f(\bm{y})\ll g(\bm{y})$ and $g(\bm{y})\ll f(\bm{y})$ are simultaneously met.
\item The Lebesgue measure in $\R^n$ is denoted by $\V_n$.
\item Given $x\in \R$, set $\e(x)=\exp(2i\pi x)$.
\end{itemize}

\paragraph{\textbf{Miscellaneous Notations~:}} 
\begin{itemize}
\item The following total order over pairs of nonnegative reals is used in Chapters~\ref{secvolestim} and~\ref{lctrootpole}: given $a,b,c,d\ge 0$,
\begin{equation*}\label{ordering}
(a,b)\preccurlyeq (c,d)\quad \Longleftrightarrow\quad \left(a<c\right)\;\vee\;\left(\left(a=c\right)\wedge\left(b\ge d\right)\right).
\end{equation*}
Explicitly, this corresponds to the order induced by the asymptotic growth of the map $x\mapsto x^{-a}\cdot\left(\log x\right)^b$ at infinity. Set also $(a,b)\succcurlyeq (c,d)$ to mean $(c,d)\preccurlyeq (a,b)$.

\item The notational use of the hat symbol has two different meanings which shall not cause confusion as they are employed in distinct places~: in Chapter~\ref{lctrootpole} and in the corresponding introductory Section~\ref{intrologcanthr} in Chapter~\ref{intro}, it refers to the complex analogue of a real quantity. In Chapters~\ref{secdeterlarge} and~\ref{secdeternarrow} however, it refers to the Fourier transform of a function.

\item In the final Chapter~\ref{secdeternarrow}, $\overline{X}$ stands for the topological closure of a set $X\subset\R^n$.
\end{itemize}

\chapter{Introduction}\label{intro}

\section{The General Setup}

Let $p\ge 1, n\ge 2$ and $d\ge 2$ be integers. Fix a $p$--tuple $\bm{F}(\bm{x})=\left(F_1(\bm{x}), \dots, F_p(\bm{x})\right)$ of non--constant real homogeneous forms $F_i(\bm{x})$ of degree $d$ in $n$ variables. Let $K$ be a subset of $\R^n$. Throughout, it is assumed that
\begin{quotation}
$(\mathcal{H}_1)$~: \emph{the set $K\subset\R^n$ is  compact and is equal to the closure of its interior}.
\end{quotation}

\noindent Given reals $a, b>0$ and a unimodular matrix $\mathfrak{g}\in SL_n(\R)$, define the set 
\begin{equation}\label{ensalg0}
\mathcal{S}_{\bm{F}}(K, \mathfrak{g}, a, b)\; =\; \left\{\bm{x}\in a\cdot K\; :\; \left\|\left(\bm{F}\circ \mathfrak{g}\right)\left(\bm{x}\right)\right\|\le b\right\}.
\end{equation}

\noindent  Here, $a\cdot K =\left\{a\bm{x}\; :\; \bm{x}\in K\right\}$ and $\left\|\;\cdot\;\right\|$ denotes the Euclidean norm. Thus,
\begin{equation*}
\left\|\left(\bm{F}\circ \mathfrak{g}\right)\left(\bm{x}\right)\right\|\;=\; \sqrt{\sum_{i=1}^{p}F_i^2\left(\mathfrak{g}\bm{x}\right)}.
\end{equation*}
The overarching undertaken goal is the study of the behaviour of the lattice point counting function
\begin{equation}\label{countfct}
\mathcal{N}_{\bm{F}}(K, \mathfrak{g}, a, b)\; =\; \#\left(\mathcal{S}_{\bm{F}}(K, \mathfrak{g}, a, b)\cap\Z^n\right)
\end{equation}
under various assumptions on the parameters $a$ and $b$, on the set $K$ and on the unimodular matrix  $\mathfrak{g}$. \\

Let $\bm{x}\in K$. The trivial inequalities $$a^d\cdot \alpha_K\left(\bm{F}\circ \mathfrak{g}\right)\;\le\; \left\|\left(\bm{F}\circ \mathfrak{g}\right)\left(a\cdot\bm{x}\right)\right\|\;\le\;  a^d\cdot \beta_K\left( \bm{F}\circ \mathfrak{g}\right),$$  where $$\alpha_K\left( \bm{F}\circ \mathfrak{g}\right)\;=\; \min_{\bm{y}\in K}\; \left\| \left(\bm{F}\circ \mathfrak{g}\right)\left(\bm{y}\right)\right\|\qquad\textrm{and}\qquad \beta_K\left( \bm{F}\circ \mathfrak{g}\right)\;=\; \max_{\bm{y}\in K}\; \left\| \left(\bm{F}\circ \mathfrak{g}\right)\left(\bm{y}\right)\right\|,$$ show that the set $\mathcal{S}_{\bm{F}}(K, \mathfrak{g}, a, b)$ is of interest only when 
\begin{equation}\label{ineqrestriabc}
a^d\cdot \alpha_K\left(\bm{F}\circ \mathfrak{g}\right)\;\le\; b\;<\; a^d\cdot \beta_K\left(\bm{F}\circ \mathfrak{g}\right).
\end{equation} 

Of singular interest is the case when $a=T$ is a large parameter and when $b=b(T)$ is taken as a function of this parameter. Then, for the sake of simplicity of notations, set
\begin{equation*}\label{ensalg0fafadde}
\mathcal{S}'_{\bm{F}}(K, \mathfrak{g}, b(T))\; =\; \mathcal{S}_{\bm{F}}(K, \mathfrak{g}, T, b(T))\qquad \textrm{and}\qquad \mathcal{N}'_{\bm{F}}(K, \mathfrak{g}, b(T))\; =\; \mathcal{N}_{\bm{F}}(K, \mathfrak{g}, T, b(T)),
\end{equation*}
where, in view of the right--hand side inequality in~\eqref{ineqrestriabc}, it is always assumed that the limit condition 
\begin{equation}\label{limiassump}
c(T)\;:=\; \frac{b(T)}{T^d}\;\underset{T\rightarrow\infty}{\longrightarrow}\; 0
\end{equation}
holds. The left--hand side inequality in~\eqref{ineqrestriabc} then shows that the problem is nontrivial only when 
\begin{quotation}
$(\mathcal{H}_2)$~: \emph{the set $K\subset\R^n$ intersects non trivially the algebraic variety $\mathcal{Z}_\R(\bm{F})$}.
\end{quotation}
Here, 
\begin{equation}\label{defalgvar}
\mathcal{Z}_\R(\bm{F})=\left\{\bm{x}\in \R^n\; :\; F_i(\bm{x})=0\quad \textrm{for all}\quad 1\le i\le p\right\}.
\end{equation}

Conventionally, the dependence on the matrix $\mathfrak{g}$ and on the set $K$ in the above va\-ri\-ous sets and quantities is dropped when $\mathfrak{g}=I_n$ is the identity matrix and/or when $K=B_n$ is the closed Euclidean ball centered at the origin.\\

The benchmark for most of the results is the generic case where the set of homogeneous forms $\bm{F}(\bm{x})$ has smooth complete intersection over $K$. This is understood in the sense that for all $\bm{x}\in K$, the gradient vectors $\nabla F_1(\bm{x}), \;\dots\;, \nabla F_p(\bm{x})$ are linearly independent. Equivalently, this is saying that the map 
\begin{equation}\label{mapsmoothcominter}
\bm{x}\in K\;\mapsto\; \bigwedge_{i=1}^p \nabla F_i(\bm{x})
\end{equation} 
does not vanish.

\section[Volume, Homogeneous Forms and Globally Semianalytic Domains]{Volume determined by Homogeneous Forms over Globally Semianalytic Domains} \label{sec0.1intro} 

The matrix $\mathfrak{g}\in SL_n(\R)$ is first seen as  a fixed parameter. Without loss of generality, upon redefining the forms under consideration, it can thus be taken as the identity. With this in mind, define the set 
\begin{equation}\label{ensalg0fafaddebis}
\underline{\mathcal{S}}_{\bm{F}}(K, b(T))\; =\; \mathcal{S}'_{\bm{F}}(K, I_n, b(T))
\end{equation}
and the counting function
\begin{equation}\label{ensalg0fafadis}  
\underline{\mathcal{N}}_{\bm{F}}(K,  b(T))\; =\; \mathcal{N}'_{\bm{F}}(K, I_n, b(T)).
\end{equation}
Recall also that both the topological conditions $(\mathcal{H}_1)$ and $(\mathcal{H}_2)$ and the limit condition~\eqref{limiassump} are assumed to hold. \\

It should be expected that, in favorable circumstances made precise in the following sections, the counting function $\underline{\mathcal{N}}_{\bm{F}}(K,  b(T))$ should be comparable to the volume of the domain $\underline{\mathcal{S}}_{\bm{F}}(K, b(T))$. Denoting by $\V_n$ the $n$--dimensional Lebesgue volume, the goal in this first section is thus to determine the behavior of the quantity $\V_n\left(\underline{\mathcal{S}}_{\bm{F}}(K, b(T))\right)$ as the parameter $T$ tends to infinity. For the stated results to be effective upon re\-sol\-ving singularities, it is convenient to make  the further general assumption that the set $K$ should be a \emph{globally semianalytic domain} (in the sense that it is defined by a finite number of inequalities involving analytic maps --- see Section~\ref{locglobbs} for details). The growth of the quantity $\V_n\left(\underline{\mathcal{S}}_{\bm{F}}(K, b(T))\right)$ is then estimated in three cases realising together a trade--off between, on the one hand, the accuracy of the  bounds that can be obtained and, on the other, the generality of the assumptions under which they hold. Specifically, the asymptotic growth of $\V_n\left(\underline{\mathcal{S}}_{\bm{F}}(K, b(T))\right)$  is  expressed as a function of the properties of the zeta distribution $\zeta_{\bm{F}}$ attached to the $p$--tuple $\bm{F}(\bm{x})$.\\ 

The distribution $\zeta_{\bm{F}}$  is defined for any complex number $s$ with non--positive real part and any element $\psi$ lying in the space of compactly supported smooth functions $\mathcal{C}_{c}^{\infty}(\R^n)$ by setting 
\begin{equation}\label{distrizetaR}
\left\langle\zeta_{\bm{F}}(s), \psi\right\rangle\;=\; \int_{\R^n}\frac{\psi(\bm{x})}{\left\|\bm{F}(\bm{x})\right\|^{s}}\cdot\textrm{d}\bm{x}.
\end{equation}
As recalled in Sections~\ref{tauberian} and~\ref{locglobrlct}, it can be extended meromorphically to the entire complex plane. \\

In Case (1), the most general considered, an upper bound is obtained for the quantity $\V_n\left(\underline{\mathcal{S}}_{\bm{F}}(K, b(T))\right)$ as a function of the smallest pole $r_{\bm{F}}(\psi)$ of the function $s\mapsto \left\langle\zeta_{\bm{F}}(s), \psi\right\rangle$, and of its order $m_{\bm{F}}(\psi)$. This is assuming only that the set $K$ meets the above conditions $(\mathcal{H}_1)$ and $(\mathcal{H}_2)$ and that the smooth test function $\psi$  is nonnegative and bounds from above the characteristic function of the set $K$.\\

In Case (2), the set $K$ is furthermore assumed to be generic enough so that the singularities of the algebraic variety $\mathcal{Z}_\R(\bm{F})$ should  not be located on its boundary. This is formalised in the following assumption~:
\begin{quote}
$(\mathcal{H}_{3})$~: \emph{there exist a real number $r_{\bm{F}}(K)\ge 0$, an integer $m_{\bm{F}}(K)\ge 1$, a compact set $C$ contained in the interior of $K$ and an open set $U$ containing $K$ satisfying the following property~:  the smallest poles and the corresponding orders of the functions $s\mapsto \left\langle\zeta_{\bm{F}}\left(s\right), \psi\right\rangle$ remain constant, respectively equal to $r_{\bm{F}}(K)$ and  $m_{\bm{F}}(K)$, for any nonnegative map $\psi$  in $\mathcal{C}_{c}^{\infty}(\R^n)$ positive over $C$ with support contained in $U$.}
\end{quote}
In particular, for any smooth test function $\psi$ meeting this condition, 
 \begin{equation}\label{suppproper}
 C\;\subset \; \textrm{Supp}\, \psi\;\subset \; U.
 \end{equation}
 
As a matter of fact, Assumption $(\mathcal{H}_3)$ turns out to be always verified in the generic case where the homogeneous forms defining $\bm{F}(\bm{x})$ have smooth complete intersection over any given  set $K$ (see Proposition~\ref{smoothcase} in Chapter~\ref{secvolestim}). \\

Under the three conditions $(\mathcal{H}_1)-(\mathcal{H}_3)$ assumed to hold in this second case, Theorem~\ref{volestim} below provides the exact asymptotic order of $\V_n\left(\underline{\mathcal{S}}_{\bm{F}}(K, b(T))\right)$  as a function of the pair $\left(r_P(K), m_P(K)\right)$; that is, its sharp asymptotic growth up to multiplicative constants.\\

In Case (3), which is perhaps the most natural in view of the assumption of the homogeneity of the forms, the set $K$ is assumed to be star-shaped with respect to the origin and to  contain the origin in its interior. Then, the precise asymptotic behavior of the volume of the domain $\underline{\mathcal{S}}_{\bm{F}}(K, b(T))$ is determined in Theorem~\ref{volestim} below. It is expressed as a function of the smallest pole $r_P>0$ and of the corresponding order $m_P\ge 1$ of the above-defined zeta distribution $\zeta_{\bm{F}}$ attached to the $p$--tuple $\bm{F}(\bm{x})$. \\

From the above discussion, the assumptions of Case (2) clearly contain those of Case (1). Also, from the homogeneity of the polynomials under consideration, the condition $(\mathcal{H}_2)$ required in the first two cases is  immediately verified in Case (3) since the origin then  lies in the set $K$. As established in Lemma~\ref{egpolorder} of Chapter~\ref{secvolestim}, it also turns out that the condition $(\mathcal{H}_3)$ assumed in Case (2) always hold under the assumptions of Case (3). This succession of observations show that the various cases considered provide an increasing degree of refinement in the following sense~: $$ \textrm{Case (1)}\qquad \Longleftarrow \qquad \textrm{Case (2)}\qquad \Longleftarrow \qquad \textrm{Case (3)}.$$  Correspondingly, the volume estimates stated in Theorem~\ref{volestim} below become sharper as the cases are more refined.\\

Theorem~\ref{volestim} below also shows that the locations and orders of the poles $r_{\bm{F}}$ and $r_{\bm{F}}(K)$ are closely related to the Sato--Bernstein polynomial $B_{\bm{F}}(s)$ of  the homogeneous form $\left\|\bm{F}(\bm{x})\right\|^2$. The Sato--Bernstein theory and the related theory of $\mathcal{D}$--modules are reviewed in Section~\ref{locglobbs}. It suffices here to mention that the Sato--Bernstein polynomial is a univariate, monic, split polynomial with negative rational roots that is uniquely associated to any given non--constant (real or complex) multivariate polynomial. It can furthermore be computed explicitly. 

\begin{thm}[Volume growth of domains determined by homogeneous forms ine\-qua\-li\-ties over globally semianalytic sets]\label{volestim}
Let the set $K$ be globally  semianalytic and let it satisfy the assumptions $(\mathcal{H}_{1})$ and $(\mathcal{H}_{2})$. Let  the limit condition~\eqref{limiassump}, whereby the quantity $c(T)$ is defined, hold. \\

Then, the volume of the set $\underline{\mathcal{S}}_{\bm{F}}(K, b(T))$ can be estimated with an increasing degree of precision as follows~: 

\begin{itemize}

\item[\emph{(1)\,--}] assume that $\psi\ge 0$ is a smooth test function bounding from above the characteristic function of the set $K$. Denote by $r_{\bm{F}}(\psi)$ the smallest real pole of the meromorphic map $s\mapsto \left\langle\zeta_{\bm{F}}(s), \psi\right\rangle$ defined in~\eqref{distrizetaR} and by $m_{\bm{F}}(\psi)\ge 1$ its order. Then, the pair $\left(r_{\bm{F}}(\psi), m_{\bm{F}}(\psi)\right)$ is well-defined and for $T\ge 1$,
\begin{equation}\label{formulevol0}
\V_n\left(\underline{\mathcal{S}}_{\bm{F}}(K, b(T))\right)\; \ll\;  T^n\cdot c(T)^{r_{\bm{F}}(\psi)}\cdot\left|\log\left(c(T)\right)\right|^{m_{\bm{F}}(\psi)-1}.
\end{equation}

\item[(2)\,--] assume that \emph{$(\mathcal{H}_{3})$ is verified}, and let $\left(r_{\bm{F}}(K), m_{\bm{F}}(K)\right)$ be the pair introduced in $(\mathcal{H}_{3})$. Then, $r_{\bm{F}}(K)$ is a strictly positive rational number. Moreover, as the parameter $T$ tends to infinity, it holds that 
\begin{equation}\label{formulevolbis}
\V_n\left(\underline{\mathcal{S}}_{\bm{F}}(K, b(T))\right)\; \asymp\;  T^n\cdot c(T)^{r_{\bm{F}}(K)}\cdot\left|\log\left(c(T)\right)\right|^{m_{\bm{F}}(K)-1}.
\end{equation}
Here, one has that $(r_{\bm{F}}(K),m_{\bm{F}}(K))=\left(p, 1\right)$ when the set of homogeneous forms $\bm{F}(\bm{x})$ has smooth complete intersection over $K$. Furthermore, under this smoothness assumption, the condition $(\mathcal{H}_{3})$ can be guaranteed under $(\mathcal{H}_{1})$ provided that $(\mathcal{H}_{2})$ is strenghtened as follows~: the algebraic variety $\mathcal{Z}_{\R}(\bm{F})$ should intersect the \emph{interior} of the set $K$.

\item[\emph{(3)\,--}] assume that \emph{the set $K$ is star-shaped with respect to the origin and that it contains the origin in its interior}. Denote by $r_{\bm{F}}$ the smallest real pole of the meromorphic distribution $\zeta_{\bm{F}}$ and by $m_{\bm{F}}\ge 1$  its order. Then, $r_{\bm{F}}$ is a well-defined rational number lying in the interval $\left(0, n/d\right]$. Furthermore, there exists a constant $\gamma_{\bm{F}}(K)>0$ such that, as the parameter $T$ tends to infinity, 
\begin{equation}\label{formulevol}
\V_n\left(\mathcal{S}_{\bm{F}}(K, b(T))\right)\; =\; \left(\gamma_{\bm{F}}(K)+o(1)\right)\cdot T^n\cdot c(T)^{r_{\bm{F}}}\cdot\left|\log\left(c(T)\right)\right|^{m_{\bm{F}}-1}.
\end{equation}
\end{itemize}

The pairs $(r_{\bm{F}}(K),m_{\bm{F}}(K))$ and $(r_{\bm{F}},m_{\bm{F}})$ involved in the above estimates can moreover be determined by resolution of singularities in a finite number of steps. Also, the rationals  $-r_{\bm{F}}(K)/2$ and $-r_{\bm{F}}/2$   are  roots of the Sato--Bernstein polynomial $B_{\bm{F}}(s)$ with multiplicity at most $m_{\bm{F}}(K)$ and $m_{\bm{F}}$, respectively.
\end{thm}

It must be emphasised that the error term in  the  volume estimate~\eqref{formulevol} can be made explicit depending on invariants attached to the homogeneous form $\left\|\bm{F}(\bm{x})\right\|^2$ --- see Section~\ref{tauberian} in Chapter~\ref{secvolestim} for details.\\

Specialising Case (3)  to the situation where $b(T)$ is constant equal to 1 and where $K=\left[-1, 1\right]^n$ is the closed cube centered at the origin with sidelength 2, one obtains a statement interesting on its own~:

\begin{coro}[Finiteness of the volume of domains determined by homogeneous forms inequalities]
The set of points $\bm{x}\in\R^n$ determined by the system of homogeneous forms inequalities $$\left|F_1(\bm{x})\right|\le 1, \quad \left|F_2(\bm{x})\right|\le 1, \quad \dots \quad, \quad \left|F_p(\bm{x})\right|\le 1$$has a finite $n$--dimensional volume if, and only if, the distribution $\zeta_{\bm{F}}$ has a simple pole at $n/d$ which is, furthermore, its smallest real pole.
\end{coro}

\section[Order and Multiplicity of the Log--Canonical Threshold]{Order \; and\; Multiplicity\; of\; the\; Log--Canonical \; Threshold}\label{intrologcanthr}

The conclusion of Theorem~\ref{volestim} prompts a finer analysis of the relation between the roots of the Sato--Bernstein polynomial of a multivariate polynomial and the poles of the corresponding zeta function. Let then $P(\bm{z})\in\C\left[\bm{z}\right]$ be a  (non--necessarily homogeneous) polynomial in $n$ variables $\bm{z}=(z_1, \dots, z_n)$ with associated Sato--Bernstein polynomial $B_P(s)\in\Q[s]$. The roots of $B_P(s)$ happen to be closely related to the poles of the \emph{complex} zeta distribution $\widehat{\zeta}_{P}$ associated to $P(\bm{z})$, which is defined for all complex numbers $s$ with non--positive real parts and all elements $\psi$  in the space of complex-valued compactly supported smooth functions $\mathcal{C}_{c}^{\infty}(\C^n)$ by the formula 
\begin{equation}\label{zetadistribcompl}
\left\langle\widehat{\zeta}_{P}(s), \psi\right\rangle\;=\; \int_{\C^n}\left|P(\bm{z})\right|^{-2s}\cdot \psi(\bm{z}, \overline{\bm{z}}) \cdot\textrm{d}\bm{z}\wedge\textrm{d}\overline{\bm{z}}.
\end{equation}
As justified in Chapter~\ref{lctrootpole}, the distribution  $\widehat{\zeta}_{P}$ also admits a meromorphic extension to the entire complex plane.\\

Malgrange~\cite{malg1}  proved that the roots of $B_P(s)$ are determined by the geometry of the algebraic variety $\mathcal{Z}_\C(P)=\left\{\bm{z}\in\C^n\; :\; P(\bm{z})=0\right\}$ in the following sense~: if $\rho\in\Q$ is a root of $B_P(s)$, then $\exp(2i\pi\rho)$ is an eigenvalue of the local monodromy of $P(\bm{z})$ at some point of $\mathcal{Z}_\C(P)$; conversely,  all eigenvalues are obtained this way. It can then be shown that, if $\exp(2i\pi\rho)$ is an eigenvalue of monodromy, then $\rho$ is a pole of the complex zeta distribution $\widehat{\zeta}_{P}$ associated to $P(\bm{z})$; here again, each eigenvalue arises this way --- see~\cite{barlet} and~\cite{malg2}. For detailed studies on the location of the poles of $\widehat{\zeta}_{P}$, the reader is referred to~\cite{loes1} and~\cite{loes2}. \\

The \emph{multiplicity} of the roots of the Sato--Bernstein polynomial and their relation to the zeta distribution are not as well understood. The only known general result seems to be due to Saito and is concerned with the roots of $\tilde{B}_P(s)=B_P(s)/(s+1)$, which is always a non--constant rational polynomial when $P(\bm{z})$ is singular over $\C$ (see~Section~\ref{locglobbs} in Chapter~\ref{secvolestim} for a justification of this claim). Denoting by $\hat{\rho}_P$ the largest root of $\tilde{B}_P(s)$, it is indeed proved in~\cite{saitoprems} that any root $\hat{\rho}$ of this polynomial lies in the interval $\left[\hat{\rho}_P-n, \hat{\rho}_P\right]\subset [-n-1, 0]$ and has multiplicity at most $n+\hat{\rho}_P+\hat{\rho}+1\le n+2\hat{\rho}_P+1\le n+1$. In some cases where the polynomial $P(\bm{z})$ is not ``too singular'' (in a suitable sense), the multiplicity of the roots of its Sato--Bernstein polynomial can also be deduced from an explicit formula for  $B_P(s)$ --- see~\cite{brigrma} and~\cite{brigrmabis} for further details.\\

The statement below characterises the  multiplicity of the largest root $\hat{r}_P>0$ of $B_P(-s)$ (recall that all roots of $B_P(s)$ are rational and negative). It is known,  see~\cite[\S 10.6]{kollar}, that this root coincides with the smallest real pole of the distribution $\widehat{\zeta}_{P}$. It is referred to as the \emph{(complex) log--canonical threshold} of the polynomial $P(\bm{z})$.  More precisely, the following statement elucidates the relationship between the multiplicity of the log--canonical threshold as a root of the polynomial $B_P(-s)$ and its order as a pole of the distribution $\widehat{\zeta}_{P}$. It answers a classical problem in the theory of $\mathcal{D}$--modules.

\begin{thm}[Order and multiplicity of the log--canonical threshold]\label{lctordermultplicity}
Let $P(\bm{z})\in \C[\bm{z}]$ be a non--constant polynomial and let $\hat{r}_P$ be its  log--canonical threshold. Then, the multiplicity of $-\hat{r}_P$ as the largest root of the Sato--Bernstein polynomial $B_p(s)$ is also the order of $\hat{r}_P$ as the smallest real pole of the meromorphic distribution $\widehat{\zeta}_{P}$.
\end{thm}

\section{Generic Unimodular Distorsions of a Family of Homogeneous Forms}

Beyond its instrinsic interest, some of the main applications of the volume estimate obtained in Theorem~\ref{volestim} are to determine, on the one hand, the generic existence of solutions to Diophantine inequalities defined by the system of homogeneous forms $\bm{F}(\bm{x})$ and, on the other, the generic asymptotic behavior of variants, natural in this context, of the counting function~\eqref{countfct}. A property is here said to be generic if its holds for almost all unimodular transformations  $\mathfrak{g}$ when $\textrm{SL}_n(\R)$ is equipped with its Haar measure.\\

More precisely, given  $\mathfrak{g}\in SL_n(\R)$, $T\ge 1$ and reals $a<b$, set  
\begin{align}
\widetilde{\mathcal{S}}_{\bm{F}}(\mathfrak{g}, a, b, T)\; &=\; \left(\left\|\left(\bm{F}\circ \mathfrak{g}\right)\right\|^{-1}\left([a,b]\right)\right)\;\cap\; B_n(T)\nonumber\\
&=\;\left\{\bm{x}\in B_n(T)\; :\; a< \left\|\left(\bm{F}\circ \mathfrak{g}\right)(\bm{x})\right\|\le b\right\} \label{defslast}\\
&\underset{\eqref{ensalg0}}{=}\; \mathcal{S}_{\bm{F}}\left(B_n, \mathfrak{g}, T, b\right)\backslash \;\mathcal{S}_{\bm{F}}\left(B_n, \mathfrak{g}, T, a\right), \label{defslastbis}
\end{align}
where $B_n(T)$ denotes the closed Euclidean ball centered at the origin with radius $T\ge 1$, and where $B_n=B_n(1)$. (The reader should not be confused by the fact that the rôle of the variable $a$  in the above relation~\eqref{defslastbis} differs from the one it plays in the definition~\eqref{ensalg0}.)\\

Athreya and Margulis~\cite{athrmargbib} gave impetus to the study of the Diophantine corollaries that can be derived from  the volume growth of the set $\widetilde{\mathcal{S}}_{\bm{F}}(\mathfrak{g}, a, b, T)$. Specifically, their focus was on the case where $d=2$ and $p=1$, and where $\bm{F}(\bm{x})=F_1(\bm{x})$ is a quadratic form in $n\ge 3$ variables --- denote it by $Q(\bm{x})$ for convenience -- such that for some constant $c_Q>0$ and all reals $a<b$, it holds that 
\begin{equation}\label{volsqgtab}
\V_n\left(\widetilde{\mathcal{S}}_{Q}(\mathfrak{g}, a, b, T)\right)\; =\; c_Q\cdot\left(b-a\right)\cdot T^{n-2}+o\left(T^{n-2}\right)
\end{equation}
as $T$ tends to infinity. From~\cite[Lemma~3.8]{EMM}, this is a full measure condition on the set of quadratic forms. Athreya and Margulis  show in~\cite[Theorem~1.1]{athrmargbib} that, under these assumptions, for every $\delta>0$ and for almost all $\mathfrak{g}\in SL_n(\R)$, there are constants $\kappa_Q(\mathfrak{g})\ge 1$ and $\varepsilon_{Q}(\mathfrak{g})>0$ such that for all $\varepsilon\in\left(0, \varepsilon_{Q}(\mathfrak{g})\right)$, there exists a nonzero integer $\bm{x}\in\Z^n$ satisfying the inequalities
\begin{equation}\label{unifquadra}
\left|\left(Q\circ\mathfrak{g}\right)(\bm{x})\right|\;<\;\varepsilon\qquad \textrm{and}\qquad \left\|\bm{x}\right\|\;\le\;\kappa_Q(\mathfrak{g})\cdot\varepsilon^{-(1/(n-2)+\delta)}.
\end{equation}
This can be rephrased as follows~: almost all unimodular twists of a quadratic form for which the volume estimate~\eqref{formulevol} holds with parameters $(r, m)=(1, 1)$ when $K=B_n$ satisfies  a uniform Diophantine approximation property with exponent $1/(n-2)+\delta$ for any $\delta>0$ (in the sense that the system~\eqref{unifquadra} admits a nontrivial solution for all $\varepsilon>0$ small enough). The last section of~\cite{athrmargbib} raises the problem of generalising this uniform Diophantine approximation property to any homogeneous form.\\

Athreya and Margulis' result is motivated by the problem of proving quantitative forms of the Oppenheim Conjecture. Recall that in its simplest form, this conjecture, solved by Margulis in~\cite{margop}, claims that for any indefinite quadratic form $Q(\bm{x})$ in $n\ge 3$ variables which is not a real multiple of a form with rational coefficients, the origin is  an accumulation point of the set $Q\left(\Z^n\right)\subset \R$. From a metric point of view, the problem of providing a quantitative estimate for this statement in the case of a \emph{generic} indefinite quadratic form (in a metric sense which varies depending on the context) has been undertaken by many authors. The existing literature in this area is therefore vast. For the purpose of the present work, it is relevant to mention the contribution by Bourgain~\cite{bourg} which is concerned with  ternary diagonal forms which are generic in a suitable sense. In~\cite{GGN}, Ghosh, Gorodnik and Nevo  use spectral methods to obtain general results of this type which imply the case $n=3$ of the above uniform approximation result~\eqref{unifquadra}. 
They also provide a heuristic to show that a homogeneous form of degree $d$ in $n>d$ variables, considered as generic in a suitable sense, should admit an exponent of uniform approximation of the form $1/(n-d)+\delta$ for any $\delta>0$.  For more background on such metric and effective statements in the context of the Oppenheim Conjecture, the reader is referred to the works by Gosh, Kelmer and Yu~\cite{GKY}, 
by Buterus, Götze and Hille~\cite{BGH0} and to the references within.\\

As far as the problem raised by Athreya and Margulis of generalising the uniform approximation result~\eqref{unifquadra} to the case of  unimodular twists of any homogeneous form is concerned, two notable contributions should be mentioned. The first one, due to Kelmer and Yu~\cite{dubyu}, proves the sought approximation property in the case of the in\-de\-fi\-ni\-te diagonal  homogeneous form 
\begin{equation}\label{diagkelyu}
F_1(\bm{x})\;=\;\sum_{i=1}^{p}x_i^{2d}- \sum_{i=p+1}^{n}x_i^{2d},
\end{equation} 
where $1\le p \le n-1$. More precisely, the particular case of~\cite[Corollary~2]{dubyu} relevant to this discussion establishes the required uniform Diophantine approximation property with exponent $1/(n-2d)+\delta$ for any $\delta>0$ provided that $n>2d$. The proof relies on the fact that for this homogeneous form,  the volume estimate~\eqref{formulevol} holds with exponents $(r, m)=(1, 1)$ when $K=B_n$. \\

The second contribution that should be mentioned is that of Kleinbock and Sken\-deri~\cite[Theorem~1.3]{kleinshk} which, when specialised to the case of homogeneous forms, establishes the following claim~: a generic unimodular twists of a system of homogeneous forms $\bm{F}(\bm{x})$ satisfies a uniform Diophantine approximation property at a given scale (i.e.~for a certain function of $\varepsilon>0$ on the right--hand side of the second inequality in~\eqref{unifquadra}) provided that a series converges. This series depends on the scaling function and is expressed as the volume of sets of the form $\widetilde{\mathcal{S}}_{\bm{F}}(\mathfrak{g}, a, b, T)$, where $a$ and $b$ are taken as functions of $T$. However, Kleinbock and Skenderi do not determine the vo\-lu\-me of such sets. \\

The sharp volume estimate obtained in Case (3) of Theorem~\ref{volestim} yields the following solution to the problem raised by Athreya and Margulis in the case of any system of homogeneous forms. This solution contains as particular cases the above--mentioned results by Kelmer and Yu~\cite{dubyu} and by Kleinbock and Skenderi~\cite{kleinshk}.

\begin{thm}[Uniform Diophantine approximation property of generic unimodular twists of systems of homogeneous forms]\label{athmarg}
In the notation of Theorem~\ref{volestim}, let  $$(r, m)\;=\;\left(r_{\bm{F}}, m_{\bm{F}}\right)\in\Q_+\times\N$$ be the pair such that the volume estimate~\eqref{formulevol} holds for the $p$--tuple $\bm{F}(\bm{x})=\left(F_1(\bm{x}), \dots, F_p(\bm{x})\right)$. Let $f : \R_+\rightarrow (1, \infty)$ be a non--decreasing function tending to infinity at infinity and satisfying the growth condition 
\begin{equation}\label{growthfcth}
\limsup_{j\rightarrow +\infty}\frac{f\left(2^{j+1}\right)}{f\left(2^{j}\right)}\;<\; \infty
\end{equation} 
and also the limit condition
\begin{equation}\label{growthfcthbis}
\frac{2^{j}}{f\left(2^{j}\right)^d}\;\underset{j\rightarrow\infty}{\longrightarrow}\; 0.
\end{equation} 
Assume furthermore that $$\sum_{j=0}^{\infty}\frac{2^{jr}}{f\left(2^j\right)^{n-rd}\cdot\left|\log \left(2^{j}\cdot f\left(2^{j}\right)^{-d}\right)\right|^{m-1}}\;<\;\infty.$$
Then, for almost every $\mathfrak{g}\in SL_n(\R)$,  there exists constants $\kappa_{\bm{F}}(\mathfrak{g}, f)\ge 1$ and $\varepsilon_{\bm{F}}(\mathfrak{g}, f)\in (0,1)$ such that for any $0<\varepsilon<\varepsilon_{\bm{F}}(\mathfrak{g}, f)$, the system of Diophantine inequalities
\begin{equation*}
\left\|\left(\bm{F}\circ \mathfrak{g}\right)\left(\bm{m}\right)\right\|< \varepsilon \qquad \textrm{and}\qquad 1\le \left\|\bm{m}\right\|\le \kappa_{\bm{F}}(\mathfrak{g}, f)\cdot f\left(\varepsilon^{-1}\right)
\end{equation*}
admits a solution in $\bm{m}\in\Z^n$.\\
\end{thm}

The growth condition~\eqref{growthfcth} means that the function $f$ is required not to admit "abrupt" variations. As for the limit condition~\eqref{growthfcthbis}, it is not a restrictive one inasmuch as it is the analogue, in this context, of the limit condition~\eqref{limiassump} meant to avoid the cases where the considered inequalities are trivially satisfied.\\

\sloppy Theorem~\ref{athmarg} supersedes the above--mentioned heuristic by Gosh, Gorodnik and Nevo~\cite{GGN} providing an exponent of uniform approximation of the form $1/(n-d)+\delta$ for any $\delta>0$ when $n>d$ in the following sense~: it shows that the underlying approximation property relies on the values of the pair $(r, m)$ for which the volume estimate~\eqref{formulevol}  holds. The particular case when $(r, m)=(1,1)$ then recovers the heuristic bound proposed by Gosh, Gorodnik and Nevo~\cite{GGN}.\\

Fixing once again  $\mathfrak{g}\in SL_n(\R)$, and reals $T\ge 1$ and $a<b$, define now the counting function
\begin{align}\label{defsnblast}
\widetilde{\mathcal{N}}_{\bm{F}}(\mathfrak{g}, a, b, T)\;  &=\; \#\left(\left(\widetilde{\mathcal{S}}_{\bm{F}}(\mathfrak{g}, a, b, T)\right)\cap\Z^n\right)\nonumber\\
&=\; \#\left\{\bm{m}\in\Z^n \cap B_n(T) : a<\left\|\left(\bm{F}\circ \mathfrak{g}\right)\left(\bm{m}\right)\right\|\le b\right\}.
\end{align}
Thus, $\widetilde{\mathcal{N}}_{\bm{F}}(\mathfrak{g}, a, b, T)$ is the lattice points counting function associated to the difference set $\mathcal{S}_{\bm{F}}\left(B_n, \mathfrak{g}, T, b\right)\backslash \;\mathcal{S}_{\bm{F}}\left(B_n, \mathfrak{g}, T, a\right)$ defined from~\eqref{ensalg0}. \\

Another problem raised by Athreya and Margulis in~\cite{athrmargbib} is to determine the generic asymptotic behavior of the function $\widetilde{\mathcal{N}}_{\bm{F}}(\mathfrak{g}, a, b, T)$ as the parameter $T$ tends to infinity. In the case of a quadratic form $Q(\bm{x})$ in $n\ge 3$ variables satisfying the volume growth~\eqref{volsqgtab}, Theorem~1.2 in~\cite{athrmargbib} shows that for all $\delta>0$ and almost every $\mathfrak{g}\in SL_n(\R)$, there exists a constant $c_{Q, \mathfrak{g}}>0$ such that 
\begin{equation}\label{estimcountam}
\widetilde{\mathcal{N}}_{\bm{F}}(\mathfrak{g}, a, b, T)\;  =\; \left(b-a\right)\cdot\left(c_{Q, \mathfrak{g}}+O\left(T^{-(n-5)/2+\delta}\right)\right)\cdot T^{n-2}
\end{equation} as $T$ tends to infinity.  This has been extended in~\cite{GKY0} to the case of shifts of quadratic forms. More generally, the reader is referred to~\cite{BGHM, GKY} and to the references within for the state of the art in metric counting problems related to quadratic forms.\\

When the homogenous form has degree $d\ge 3$, the only result towards a solution to this counting problem the authors are aware of is  due to Kelmer and Yu~\cite{dubyu} (and its slight extension to a system of closely related forms by Bandi, Ghosh and Han~\cite{BGH}). It deals with the case of the homogeneous form $F_1(\bm{x})$ defined in~\eqref{diagkelyu} when $n>2d$. Specifically, Kelmer and Yu are concerned with the related but slightly different pro\-blem when the bounds $a=a(T)$ and $b=b(T)$ of the real intervals depend on the parameter $T$ and satisfy the shrinking target property that $b(T)=a(T)+c\cdot T^{-\gamma}$ for some $c, \gamma>0$. Then, \cite[Theorem~1]{dubyu} states that for every $\beta\in\left(0, 2\left(n-d-\gamma\right)/\left(n^2+n+4\right)\right)$ and for almost all $\mathfrak{g}\in SL_n(\R)$, 
\begin{equation}\label{estimcountdy}
\widetilde{\mathcal{N}}_{F_1}\left(\mathfrak{g}, T, a(T), a(T)+c\cdot T^{-\gamma}\right)\;  =\; c\cdot T^{n-2d-\gamma}\cdot\left(c_{F_1, \mathfrak{g}}+O\left(T^{-\beta}\right)\right)\cdot T^{n-2},
\end{equation}
where $c_{F_1, \mathfrak{g}}>0$.\\

The proof of both counting estimates~\eqref{estimcountam} and~\eqref{estimcountdy} share two common features~: on the one hand, they crucially rely on the fact that the corresponding algebraic domains $\widetilde{\mathcal{S}}_{Q}(\mathfrak{g}, a, b, T)$
 and $\widetilde{\mathcal{S}}_{F_1}(\mathfrak{g}, a, b, T)$ are ``large enough'' in the sense that their volumes grow at least polynomially. Equivalently, when considering the volume estimate~\eqref{formulevol}  in the case that $K=B_n$ and $b(T)=1$, this is saying that the  pairs $(r, m)$ corresponding to each of these forms are both such that the exponent $r$ does not take the extremal value $n/d$, where $d$ denotes the degree of either forms $Q(\bm{x})$ or $F_1(\bm{x})$ (in other words, one requires that $0<r<n/d$). On the other hand, the power savings in the error terms in~\eqref{estimcountam} and~\eqref{estimcountdy} strongly rely on the fact that, in both cases, it holds that $m=1$ (in fact, as mentioned above, $(r, m)=(1,1)$ for the two forms $Q(\bm{x})$ and $F_1(\bm{x})$ under consideration).\\

The following statement generalises these observations and settles the counting problem raised by Athreya and Margulis in the more general case of systems of homogeneous forms in $n\ge 3$ variables.

\begin{thm}[Generic asymptotic behavior of the counting function of unimodular twists of systems of homogeneous forms]\label{athmarg2} Assume that $n\ge 3$ and fix real numbers $0\le a< b$. In the notation of Theorem~\ref{volestim}, let  $(r, m)=(r_{\bm{F}}, m_{\bm{F}})\in\Q_+\times\N$ be the pair such that the volume estimate~\eqref{formulevol} holds for the $p$--tuple $\bm{F}(\bm{x})=\left(F_1(\bm{x}), \dots, F_p(\bm{x})\right)$.  Then, under the assumption that $r<n/d$, there exists, for almost all $\mathfrak{g}\in SL_n(\R)$,   a constant $\gamma_{\bm{F}, \mathfrak{g}}>0$ such that 
\begin{equation}\label{coungenthù}
\widetilde{\mathcal{N}}_{\bm{F}}(\mathfrak{g}, a, b, T)= \left(\gamma_{\bm{F}, \mathfrak{g}}+o(1)\right)\cdot T^{n-rd}\cdot \left( b^r\cdot\left|\log\left(\frac{b}{T^d}\right)\right|^{m-1}- a^r\cdot\left|\log\left(\frac{a}{T^d}\right)\right|^{m-1}\right)
\end{equation}
as the parameter $T$ tends to infinity. 
\end{thm}

As a matter of fact, this theorem is deduced in Chapter~\ref{ranunidisfamhomfor} upon specialising a much more general statement to the case of the homogeneous form $P_{\bm{F}}(\bm{x})=\left\|\bm{F}(\bm{x})\right\|^2$. This more general statement deals with  any real homogeneous form $P(\bm{x})$ and provides an explicit form for the error term in~\eqref{coungenthù} as a function of various invariants attached to $P(\bm{x})$. It suffices to say here that in the case $m=1$, it enables one to recover the power saving estimates obtained in~\eqref{estimcountam} and~\eqref{estimcountdy} with a comparable strength. When $m\ge 2$ however, the presence of logarithmic factors in the volume estimate~\eqref{formulevol} imposes that the induced saving can only be of the order of a power of a logarithm.\\

Various other problems related to metric twists of algebraic functions have been consi\-de\-red in the literature. To cite but a few relevant to this work, Vanderkam~\cite{vdk}  studies 
the metric behavior of a counting function similar to~\eqref{defsnblast} when considering all coefficients of a polynomial as variables behaving independently. The paper~\cite{ho} also surveys the state of the art in the related problem of counting integral points \emph{on} random algebraic curves. Finally, in~\cite{adi} is considered the problem of approximating values taken by a polynomial twisted by a random nonlinear  transformation, namely a vertical shift.

\section{The Number of Solutions to Homogeneous Forms Inequalities}

Let $K\subset\R^n$ be a set meeting the assumptions $(\mathcal{H}_1)$, $(\mathcal{H}_2)$ and $(\mathcal{H}_3)$ (recall here that $(\mathcal{H}_3)$ is necessarily satisfied when $K$ is star-shaped with respect to the origin and that it contains the origin in its interior). Given a parameter $\alpha>0$, consider the case where 
\begin{equation}\label{counboundalpha}
b(T)=T^{d-\alpha}
\end{equation} 
in the definitions of the set $\underline{\mathcal{S}}_{\bm{F}}\left(K, b(T)\right)$  in~\eqref{ensalg0fafaddebis} and of the corresponding counting function $\underline{\mathcal{N}}_{\bm{F}}\left(K, b(T)\right)$ in~\eqref{ensalg0fafadis}. Note that the limit condition~\eqref{limiassump}  then holds. In order to single out the rôle of the parameter $\alpha$ in the forthcoming statements, define
\begin{align}
\mathcal{S}^{\dag}_{\bm{F}}\left(K,T, \alpha\right)\;\underset{\eqref{ensalg0fafaddebis}}{=}\; \underline{\mathcal{S}}_{\bm{F}}\left(K, T^{d-\alpha}\right) \;=\; \left\{\bm{x}\in T\cdot K \; :\; \left\|\bm{F}(\bm{x})\right\|\le T^{d-\alpha}\right\}\label{colsfsar}
\end{align}
and 
\begin{equation}\label{nbsfsar} 
\mathcal{N}^{\dag}_{\bm{F}}\left(K,T, \alpha\right)\;\underset{\eqref{ensalg0fafadis}}{=}\; \underline{\mathcal{N}}_{\bm{F}}\left(K, T^{d-\alpha}\right) \;=\; \#\left(\mathcal{S}^{\dag}_{\bm{F}}\left(K,T, \alpha\right)\cap \Z^n\right).
\end{equation}

The focus on a power function such as in~\eqref{counboundalpha} is motivated by the following deterministic counting result due to Sarnak~: assuming that the zero set $\mathcal{Z}_K(\bm{F}):=\left\{\bm{x}\in K\; :\; \bm{F}(\bm{x})=\bm{0}\right\}$ does not lie  in a linear subspace of dimension $n-p$ and that the smooth complete intersection condition determined by the non-vanishing of the map~\eqref{mapsmoothcominter} over the set $K$ holds, it is established in~\cite[Appendix~ 1]{sarnakprinci} that the estimate 
\begin{equation}\label{estimsarnakun}
\mathcal{N}^{\dag}_{\bm{F}}\left(K,T, \alpha\right)\;\ll\; T^{n-p\alpha}+T^{n-p-\delta} 
\end{equation} 
is verified for some $\delta>0$ provided that the set $K$ is "nice". Sarnak does not make  explicit the meaning of the latter condition; his proof nevertheless certainly works if, say, $K$ meets the assumptions $(\mathcal{H}_1)$ and $(\mathcal{H}_2)$ and if, furthermore, it is convex and has piecewise smooth boundary. Moreover, the arguments of the proof show that one can then take the value $\delta=1/(nd)$. \\

The fact that there should exist a power saving factor $\delta>0$ in the error term of the estimate~\eqref{estimsarnakun} has  a particularly interesting consequence~: as shown in~\cite[p.~177]{se} (with the choice of $\alpha=1/2$ in the above notation), the right-hand side of~\eqref{estimsarnakun} then bounds nontrivially the number of solutions to the system of Diophantine \emph{equations} $\bm{F}(\bm{m})=\bm{0}$ when $\bm{m}\in T\cdot K$. Inspired by this implication, Sarnak conjectures at the end of his proof that an estimate of the form~\eqref{estimsarnakun} should hold without the assumption of smoothness in the following sense~: assuming that the zero set $\mathcal{Z}_{K}(\bm{F})$ is non empty, not contained in a linear space of dimension $n-p$ and that the set $K$ is again "nice", there should exist $\delta>0$ such that 
\begin{equation}\label{estimsarnakunbis}
\mathcal{N}^{\dag}_{\bm{F}}\left(K,T, \alpha\right)\;\ll\; Vol_n\left(\mathcal{S}^{\dag}_{\bm{F}}\left(K,T, \alpha\right)\right)+T^{n-p-\delta}.
\end{equation} 
Interpreting the assumption that the set $K$ is "nice" in the sense that the conditions $(\mathcal{H}_1)$ and $(\mathcal{H}_2)$ are met, Case (2) in  Theorem~\ref{volestim} confirms that the above upper bound indeed reduces to the inequality~\eqref{estimsarnakun} in the case of smooth complete intersection. \\

It nevertheless turns out that a power saving in the error term such as in~\eqref{estimsarnakunbis} is not unconditionnally true  in the general case just assuming that the zero set $\mathcal{Z}_{K}(\bm{F})$ is not empty, not contained in a linear subspace of dimension $n-p$ and that the set $K$ is "nice" (or regular) enough. To see this, consider the case where $p=1$ and where the homogeneous form under consideration is 
\begin{equation}    \label{deff_n}
F_n(\bm{x})\; :=\;\prod_{i=1}^{n}x_i, 
\end{equation}
which has degree $d=n$. The zero set of $F_n(\bm{x})$ is here the union of the coordinate hyperplanes. Taking $K=K_n$ to be the unit cube $[0,1]^n$ and choosing any $\alpha> 0$, the counting function $\mathcal{N}^{\dag}_{F_n}\left(K,T, \alpha\right)$ is closely related to the generalised divisor summatory function $$\Delta_n(t)\; =\; \sum_{1\le F_n(\bm{x})\le t} 1$$ in the sense that  
\begin{align}
0&\; \le\; \Delta_{n}\left(T^{n-\alpha}\right)+n\cdot\left(T+1\right)^{n-1} - \mathcal{N}^{\dag}_{F_n}\left(K,T, \alpha\right)\nonumber \\
&\;\qquad\qquad\qquad\qquad\qquad\qquad\le\; \frac{n(n+1)}{2}\cdot (T+1)^{n-2}+n\cdot\Delta_{n-1}\left(T^{n-1-\alpha}\right).\label{contreexsarnakml}
\end{align}

Indeed, this follows upon considering the contributions to $\mathcal{N}^{\dag}_{F_n}\left(K,T, \alpha\right)$ of the coordinate hyperplanes separately and, as far as the upper bound is concerned, upon furthermore  noticing this elementary fact~: the number of integer solutions to the ine\-qua\-li\-ty $1\le F_n(\bm{x})\le T^{n-\alpha}$ when $x_i> T$ for some coordinate $1\le i\le n$ is bounded above by $n$ times the number of integer solutions to the inequality $1\le F_{n-1}(\bm{y})\le T^{n-1-\alpha}$, where $\bm{y}\in\R^{n-1}$. \\

Now, it is known that there exists a polynomial $Q_{n}(t)$ of degree $n-1$ such that $$\Delta_n(t)\; =\; t\cdot Q_n(\log t)+O\left(t^{1-\varepsilon}\right)$$for some $\varepsilon>0$. (Determining an optimal exponent for the error term is the Piltz Divisor Problem, which reduces to the better known Dirichlet Divisor Problem when $n=2$.) More precisely, the coefficients of the polynomials $Q_n(t)$ can be evaluated from the formula $$Q_n(\log t)\;=\; \textrm{Res}_{s=1} \left(t^{s-1}\cdot \zeta^n(s)\cdot s^{-1}\right)$$ and from the Laurent series expansion of the Riemann zeta function $$\zeta(s)\; =\; \frac{1}{s-1}+\sum_{k=0}^\infty \frac{(-1)}{k!}\cdot\gamma_k\cdot (s-1)^k$$valid in a neighbourhood of $s=1$. Here, $\left(\gamma_k\right)_{k\ge 0}$ denotes the sequence of the Stieltjes constants. Explicit calculations of the coefficients of $Q_n(t)$ are carried out, for instance, in~\cite[Theorem~ 1]{lavrik}. In particular, one obtains that $$Q_2(t)\;=\; t+ (2\gamma -1)\qquad\textrm{ and }\qquad Q_3(t)\;=\; \frac{t^2}{2}+(3\gamma -1)\cdot t +\left(3\gamma^2-3\gamma+3\gamma_1+1\right),$$where $\gamma \approx 0.577$ is the Euler–Mascheroni constant and where $\gamma_1 \approx -0.073$. Since elementary calculations show that the volume of the region $\mathcal{S}^{\dag}_{F_n}\left(K,T, \alpha\right)$ grows, up to a multiplicative constant, as $T^{n-\alpha}\cdot \left(\log T\right)^{n-1}$ (this is saying that $\left(r_{F_n}, m_{F_n}\right)=\left(1, n\right)$ in the notations of Theorem~\ref{volestim}), the estimates~\eqref{contreexsarnakml} imply that there can be no power saving in the error term when $n=3$ and, more generally, in the generic case where the coefficients do not cancel out suitably. \\

This counterexample makes clear the nature of the obstacle to overcome for Sarnak's estimate~\eqref{estimsarnakunbis} to hold~: in the presence of singularities (such as the origin in the case of the homogeneous form~\eqref{deff_n}), logarithmic terms appearing in the volume term need not be compensated by a power saving in the error term. The presence of such logarithmic terms should be seen as an indication of the fact that the zero set  $\mathcal{Z}_{K}(\bm{F})$ is in some sense "too flat" (whithout necessarily being contained in a linear space as illustrated by the case of the homogeneous form~\eqref{deff_n}). In the particular case of a planar curve, this link between flatness and behavior of the counting function of Diophantine inequalities has been observed and formalised (in a related but slightly different context) by Vaughan and Velani~\cite{VV} and then by Huang~\cite{H1}. \\

With this in mind, the final goal of this paper is to determine, under the assumptions $(\mathcal{H}_1)$, $(\mathcal{H}_2)$ and $(\mathcal{H}_3)$, conditions related to the "the level of flatness" of the zero set $\mathcal{Z}_{K}(\bm{F})$ under which Sarnak's claim holds. As this goal is achieved working in the context of semialgebraic geometry, it is natural to assume furthermore that the set $K$ is semialgebraic, meaning that it is defined as a finite union of sets that can be represented as finitely many equalities and inequalities involving polynomials. The statement of the results furthermore depends on the choice of  a compact semialgebraic set $\mathcal{K}$ containing $K$ in its interior and being contained in the open set $U$ introduced in~\eqref{suppproper}. This is because to make the values of the power savings effective whenever they hold, one indeed needs to consider  the behaviour of the set of homogeneous forms $\bm{F}(\bm{x})$ over a compact neighbourhood  of $K$. \\

Define the dimension $\dim\left(\mathcal{Z}_\mathcal{K}(\bm{F})\right)$ of the algebraic variety $\mathcal{Z}_\R(\bm{F})$ over the set $\mathcal{K}$  as the dimension of the tangent vector space at any non singular point of $\mathcal{Z}_\R(\bm{F})\cap \mathcal{K}$ (this is well--defined). Letting $\textrm{codim}\left(\mathcal{Z}_\mathcal{K}(\bm{F})\right) = n-dim\left(\mathcal{Z}_\mathcal{K}(\bm{F})\right)$, set for the sake of simplicity of notation
\begin{equation*}\label{defdimcodim0}
\tau_{\bm{F}}(\mathcal{K})\;=\; \dim\left(\mathcal{Z}_\mathcal{K}(\bm{F})\right)\qquad \textrm{and}\qquad \widehat{\tau}_{\bm{F}}(\mathcal{K})\;=\; \textrm{codim}\left(\mathcal{Z}_\mathcal{K}(\bm{F})\right).
\end{equation*}
For instance, if the variety $\mathcal{Z}_{\R}(\bm{F})$ has smooth  complete intersection over $K$ and that the compact neighbourhood $\mathcal{K}$ is chosen small enough around $K$ so that this property is preserved over $\mathcal{K}$, then $\widehat{\tau}_{\bm{F}}(\mathcal{K})=p$. The exponent in the error term in~\eqref{estimsarnakun}  valid in the smooth case can then  be interpreted as $\tau_{\bm{F}}(\mathcal{K})-\delta$. This is this form of the exponent that is generalised in what follows to the non--smooth case. \\

The existence of a power saving error term is expressed as a condition depending on a measure of "the level of flatness" of the zero set $\mathcal{Z}_{\mathcal{K}}(\bm{F})$ defined as follows~:  given a unit vector $\bm{v}\in\Sph^{n-1}$ and a real $\sigma$, set first $$\bm{v}^{\perp}(\sigma)\;:=\;\left\{\bm{x}\in\R^n\; :\; \bm{v\cdot x}=\sigma\right\}.$$The flatness of  the zero set $\mathcal{Z}_{\mathcal{K}}(\bm{F})$ along such an affine space is measured by the volume section
\begin{equation*}\label{defmupcvte0}
\mu_{\bm{F}, \mathcal{K}}\left(\bm{v}, \sigma, \varepsilon\right)\;=\; \V_{n-1}\left(\left\{\bm{x}\in \mathcal{K}\cap \bm{v}^{\perp}(\sigma)\; :\; \left\|\bm{F}(\bm{x})\right\|\le\varepsilon\right\}\right),
\end{equation*}
where $\varepsilon>0$. The "biggest" such volume section is then defined as 
\begin{equation*}\label{defmpueps0} 
M_{\bm{F}}\left(\mathcal{K}, \varepsilon\right)\;=\; \sup_{\bm{v}\in\Sph^{n-1}}\;\sup_{\sigma\in\R}\;\mu_{\bm{F},\mathcal{K}}\left(\bm{v}, \sigma, \varepsilon\right).
\end{equation*}
From this, the measure of flatness relevant to the problem under con\-si\-de\-ration is defined as
\begin{equation}\label{conditionforall0} 
q_{\bm{F}}(\mathcal{K})\;=\; \liminf_{\varepsilon\rightarrow 0^+} \left(\frac{\log M_{\bm{F}}\left(\mathcal{K}, \varepsilon\right)}{\log \varepsilon}\right).
\end{equation}
This is a measure of the "biggest level of flatness" that can be achieved when intersecting the sublevel set $\left\{\bm{x}\in \mathcal{K}\; :\; \left\|\bm{F}(\bm{x})\right\|\le\varepsilon\right\}$  with affine hyperplanes in the following sense~: the smaller it is, the more flat the intersection of the sublevel set in some direction. This measure of flatness is a well--defined real number under the assumption that $\mathcal{Z}_K(\bm{F})\neq\emptyset$. It is furthermore shown in Chapter~\ref{secdeternarrow} (see the last section therein) that the real $q_{\bm{F}}(\mathcal{K})$ can then be estimated effectively in a finite number of steps (in a suitable sense). The statement of the main result in this section (Theorem~\ref{sarnak} below) relies on a comparison between the value of $q_{\bm{F}}(\mathcal{K})$ and the volume growth exponent $r_{\bm{F}}(K)$ present in the statement of Theorem~\ref{volestim} and part of the assumption $(\mathcal{H}_3)$.\\

Specifically, the following provides a solution to Sarnak's problem generalising the smooth case~: 

\begin{thm}[Counting solutions to homogeneous forms inequalities]\label{sarnak} Recall that the conditions $(\mathcal{H}_1)$, $(\mathcal{H}_2)$ and $(\mathcal{H}_3)$ are assumed to hold and that the above defined sets $K\subset \mathcal{K}$ are assumed to be semialgebraic. Fix a parameter $\alpha>0$ and assume that the measure of flatness~\eqref{conditionforall0} is big enough in the sense that 
\begin{equation}\label{cndiestimvolsa1reutilise}
q_{\bm{F}}(\mathcal{K})\;>\; \max\left\{r_{\bm{F}}(K)-1,  n-1-r_{\bm{F}}(K)\cdot\frac{\tau_{\bm{F}}(\mathcal{K})}{\widehat{\tau}_{\bm{F}}(\mathcal{K})}\right\}.
\end{equation}
Then, there exists a real $\delta(K, \mathcal{K})>0$ depending on the semialgebraic sets $K$ and $\mathcal{K}$ such that the following estimate with power saving error term holds~: 
\begin{equation}\label{estimvolsa0}
\mathcal{N}^{\dag}_{\bm{F}}\left(K,T, \alpha\right)\;\ll\; \V_n\left(\mathcal{S}^{\dag}_{\bm{F}}\left(K,T, \alpha\right)\right)+T^{\tau_{\bm{F}}(\mathcal{K})-\delta_{\bm{F}}(\mathcal{K}, K)}.
\end{equation}
Furthermore, the inequality~\eqref{cndiestimvolsa1reutilise} is satisfied when the following two conditions are simultaneously met~: the set of polynomials $\bm{F}(\bm{x})$ has smooth complete intersection over $\mathcal{K}$ and the zero set $\mathcal{Z}_K(\bm{F})$ is not contained in a linear subspace of dimension $n-p$, where $n-p=\tau_{\bm{F}}(\mathcal{K})$.
\end{thm}

It should be noted that in the case of the counterexample~\eqref{deff_n} to Sarnak's claim, the inequality~\eqref{cndiestimvolsa1reutilise} is indeed not met since in any semialgebraic compact neighbourhood $\mathcal{K}_n^+$ of $K_n^+=[0,1]^n$, one easily checks that $\widehat{\tau}_{F_n}(\mathcal{K}_n^+)=1$, that $r_{F_n}(\mathcal{K}_n^+)=n-1$ and that $q_{F_n}(\mathcal{K}_n^+)=0$ (this last relation follows from the fact that the quantity $M_{F_n}(\mathcal{K}_n^+, \varepsilon)$ is constant independent of $\varepsilon>0$ as can be seen upon considering the volume of the sections of the sublevel set $\left|F_n(\bm{x})\right|<\varepsilon$ with coordinate hyperplanes).\\

In fact, Theorem~\ref{sarnak} as stated above is an extremely abridged version of a class of statements providing a much more general picture of the behaviour of the coun\-ting function $\mathcal{N}^{\dag}_{\bm{F}}\left(K,T, \alpha\right)$. Thus, it is established in Chapter~\ref{secdeterlarge} that the volume $\V_n\left(\mathcal{S}^{\dag}_{\bm{F}}\left(K,T, \alpha\right)\right)$ determines the \emph{exact} asymptotic order of this coun\-ting function as $T$ tends to infinity in the case that $\alpha\in (0,1)$ (i.e.~when the algebraic set $\mathcal{S}^{\dag}_{\bm{F}}\left(K,T, \alpha\right)$ is "large enough") and also that the counting function is, up to a multiplicative constant, bounded above by this volume in the case $\alpha=1$. These estimates are further refined when $\alpha\in (0,1)$ under the stronger assumptions of Case (3) stated in the above  Section~\ref{sec0.1intro}~: an asymptotic expansion with power saving error term is then established for the quantity $\mathcal{N}^{\dag}_{\bm{F}}\left(K,T, \alpha\right)$ as $T$ tends to infinity. \\

In the complementary case when $\alpha>1$ (i.e.~when the algebraic set $\mathcal{S}^{\dag}_{\bm{F}}\left(K,T, \alpha\right)$ is "small"), it is proved in Chapter~\ref{secdeternarrow} that an estimate of the form~\eqref{estimvolsa0} can still hold in the degenerate case when  $q_{\bm{F}}(\mathcal{K})\le r_{\bm{F}}(K)-1$ under an explicit condition (which is, of course,  still not satisfied by the counterexample~\eqref{deff_n}). Also, in all instances where an estimate with power saving error term holds, a range of admissible values of the exponent $\delta_{\bm{F}}(K, \mathcal{K})>0$ is explicitly determined. \\
 
Theorem~\ref{sarnak} complements and extends a number of previously known results. Thus, in the case of a single \emph{positive} homogeneous function $F(\bm{x})$ (not necessarily assumed to be a polynomial), Randol~\cite{randol} establishes a counting bound with power saving error term depending on the decay of a Fourier operator acting on the boun\-da\-ry surface $\left\{F=1\right\}$, which is supposed to be regular enough (in a suitable sense). This has been slightly generalised by Colin de Verdière~\cite{colin}  who  parallelly considers the behavior of the counting bound when rotating randomly the domain $\left\{F\le 1\right\}$.  \\

In the general case where one considers the problem of counting rational points \emph{near} manifolds (Theorem~\ref{sarnak} being concerned with the particular case where the manifold is a homogeneous algebraic variety), the most powerful results have been obtained under various assumptions of non degeneracy (which often amount to imposing constraints on the behavior of the local curvature). In this respect, Beresnevich \emph{et al.}~\cite{beres} show that, provided that a certain matrix of second order derivatives of a prametrisation of a twice continuously differentiable manifold admits a determinant bounded below by a strictly positive constant, one can recover a counting bound with a power saving error term. Under the stronger assumption that the Gaussian curvature remains strictly positive, Huang~\cite{huang} bounds from above the counting function of the number of rational points near a manifold by the volume of a neighbourhood of this manifold (which takes a particularly simple form under the considered assumption). This establishes the heuristics that the number of rational points in a "nice enough domain" should be given by the volume of this domain provided that it is "large enough". In the case that one considers a domain with analytic boundary, Greenblatt~\cite{greenblatt} proves a counting bound related to Theorem~\ref{sarnak} by the use of suitable resolutions of singularities. This nevertheless comes with an additional restriction  not present in the above Theorem~\ref{sarnak}~: in his work, the boundary of the domain must be locally the graph of an analytic function which separates the interior of the domain from its complement in exactly two pieces.\\

Finally, it should be noted that although Theorem~\ref{sarnak} holds in the binary case of $n=2$ variables, explicit and sharper results are then known if one considers a single homogeneous form with integer coefficents (one is then dealing with a Thue inequality). To see recent progress in this well--studied topic, see the works by Fouvry and Waldschmidt~\cite{fw} , by Stewart \& Xiao~\cite{sx} and the references within.

\section{Structure of the Memoir} 

The following table summarises the sections where the statements presented in the introduction are established. As almost all of these have not been stated in full generality for the sake of the coherence of the exposition, the table also provides the statements either generalising them or relating them to other aspects of the problem they are concerned with.\\

\begin{center}
\begin{tabular}{||c||c||c||}
\hline
\hline
Statement & Proof & Related extensions \\
\hline
\hline
\multirow{4}{*}{Theorem~\ref{volestim}}& \multirow{4}{*}{ Sections~\ref{tauberian} \&~\ref{locglobrlct}}& Corollary~\ref{proptmh1.1}  \\
& & Proposition~\ref{locsmalpole}  \\
& & Proposition~\ref{smoothcase}  \\
& &Theorem~\ref{lctvolcomplex}\\
\hline
Theorem~\ref{lctordermultplicity} & Chapter~\ref{lctrootpole}&\\
\hline
Theorem~\ref{athmarg} & Chapter~\ref{ranunidisfamhomfor} & \\
\hline
Theorem~\ref{athmarg2} & Chapter~\ref{ranunidisfamhomfor} & Theorem~\ref{thmprobacountvol}\\
\hline
\multirow{2}{*}{Theorem~\ref{sarnak}} & Chapter~\ref{secdeterlarge} when $0<\alpha\le 1$& Theorem~\ref{countinglargedomains}\\
& \& Chapter~\ref{secdeternarrow} when $\alpha> 1$&  Theorem~\ref{thmgeneralsarnak}\\
\hline
\hline
\end{tabular}
\end{center}

\paragraph{Acknowledgments} This memoir is the conclusion of a long project during which the authors exchanged ideas with many people. In this respect, they would like to thank Profs.~Victor Beresnevich, Raf Cluckers, Pietro Corvaja, Roger Heath--Brown, Tobias Kaiser, Mircea Mustaţă and Mihnea Popa, and  also Quentin Guignard and Gareth O.~Jones.\\

The first-named author would especially like to thank Dr. Detta Dickinson who quite unexpectedly initiated him to this area of research in 2011 when he started his PhD under her supervision~: as far as he is concerned, the problem she suggested him to work on somehow gave impetus to this project.

\chapter{Volume Estimates}\label{secvolestim}

The goal of this chapter is to establish Theorem~\ref{volestim}. The set $\K$ denotes either the field of reals $\R$ or the field of complex numbers $\C$. Furthermore, $P(\bm{y})$ is a non--constant polynomial of degree $q\ge 1$ with coefficients in $\K$ which will be specialised when needed to the homogeneous form $\left\|\bm{F}(\bm{y})\right\|^2$.\\

The following form of Hironaka's Theorem on simultaneous resolution of sin\-gu\-la\-ri\-ties  stated in~\cite[p.147]{atiyah} (when $l=1$) and in~\cite[Theorem 2.8]{watanabe} (when $l\ge2$) will be  used repeatedly~:

\begin{thm}[Simultaneous Resolution of Singularities, Hironaka, 1964]\label{hironaka}
 Let $l\ge 1$ be an integer and let $f_1, \dots, f_l$ be analytic maps defined over an open set $\mathcal{U}\subset\K^n$ containing the origin. Assume that given $1\le i\le l$, the map $f_i$ is non--constant and satisfies $f_i(\bm{0})=0$. Then, there exists an open set $\mathcal{W}\subset \mathcal{U}$ containing the origin, an $n$--dimensional $\K$--analytic manifold $\mathcal{M}$ and a $\K$--analytic map $g~: \mathcal{M}\rightarrow \mathcal{W}$ such that the following hold~: 
\begin{itemize}
\item the map $g$ is proper (that is, the preimage of a compact set is compact) and surjective;
\item the map $g$ is an analytic isomorphism from $\mathcal{M}\backslash\mathcal{M}_0$ to $\mathcal{W}\backslash \mathcal{W}_0$, where $$\mathcal{M}_0=\bigcup_{i=1}^{l}\left\{\bm{y}\in\mathcal{M}\; :\; f_i(g\left(\bm{y})\right)=0\right\}\qquad  \textrm{and}\qquad \mathcal{W}_0=\bigcup_{i=1}^{l}\left\{\bm{x}\in \mathcal{W}\; :\; f_i(\bm{x})=0\right\};$$
\item given any point $\bm{y}_0\in\mathcal{M}_0$, there exists local coordinates $\bm{y}=\left(y_1, \dots, y_n\right)\in\mathcal{M}$ centered at $\bm{y}_0$ such that for all $1\le i\le l$, 
\begin{equation}\label{monoform} 
f_i(g(\bm{y}))\;=\; \alpha_i(\bm{y})\cdot \bm{y}^{\bm{k}_i} .
\end{equation} 
Furthermore, the Jacobian of the change of variables $\bm{x}=g(\bm{y})$ is such that 
\begin{equation}\label{monoformbis}
\textrm{Jac}_g(\bm{y})=\beta(\bm{y})\cdot \bm{y}^{\bm{h}}.
\end{equation} 
In these relations, $\alpha_1$ is a constant map equal to $\pm 1$; $\alpha_2, \dots , \alpha_n$ and $\beta$ are non--vanishing anaytic maps, and $\bm{k}_i$ and $\bm{h}$ are $n$--tuples of nonnegative integers.
\end{itemize}
\end{thm}

\section{Local and Global Sato--Bernstein Polynomials} \label{locglobbs} 

This section is devoted to recalls on and extensions of the elements of the Sato--Bernstein theory needed for the purpose of this work. Classical references where the known claims stated hereafter are proved include~\cite{igusa} and~\cite{kashiwara}.\\

Given a non--constant polynomial $P(\bm{y})\in\K\left[\bm{y}\right]$, there exists a univariate nonzero polynomial $b(s)\in\K[s]$ and a differential operator $\mathcal{D}\left(\bm{y}, s, \bm{\partial}\right)$ with coefficients in $\K$  such that the following identity holds formally~:
\begin{equation}\label{sbrelation}
\mathcal{D}\left(\bm{y}, s, \bm{\partial}\right)\cdot P(\bm{y})^{s+1}\; =\;  b(s) P(\bm{y})^s.
\end{equation}

Furthermore, the operator $\mathcal{D}\left(\bm{y}, s, \bm{\partial}\right)$ can be chosen as a (non--cummutative) polynomial in the variables $\bm{y}=\left(y_1, \dots, y_n\right)$, in an additional variable $s$ and in the set of partial derivatives $\bm{\partial}=\left(\partial_1, \dots, \partial_n\right)$. Clearly, the set of all those polynomials $b(s)\in\K[s]$ for which a differential operator satisfying~\eqref{sbrelation} exists  forms an ideal in the ring $\K[s]$. The \emph{ (global) Sato--Bernstein polynomial} $B_P(s)$ of $P(\bm{y})$ is then defined as the monic generator of this ideal; in other words, it is the monic polynomial of smallest degree satisfying the relation~\eqref{sbrelation} for some differential operator $\mathcal{D}$. The following fundamental result due to Kashiwara~\cite{kashi} implies that $B_P(s)$ has rational coefficients.

\begin{thm}[Kashiwara, 1976]\label{kashiroot}
The Sato--Bernstein polynomial $B_P(s)$ of a polynomial $P(\bm{y})\in\K[\bm{y}]$ is a split polynomial all of whose roots are rational and negative.
\end{thm}

For instance, in the case that $P(\bm{y})=P_1(\bm{y})=\sum_{i=1}^{n}y_i^2$,  one obtains when taking the Laplacian as  the differential operator that
$$\sum_{i=1}^n\partial_i^2P_1(\bm{y})^{s+1}\;=\; 4\left(s+1\right)\cdot\left(s+\frac{n}{2}\right)\cdot P_1(\bm{y})^{s}.$$ It is then easy to deduce from this relation that the corresponding Sato--Bernstein polynomial is $B_{P_1}(s)=\left(s+1\right)\cdot\left(s+\frac{n}{2}\right)$.\\

The proof of Kashiwara's Theorem~\ref{kashiroot} shows that computing a Sato--Bernstein polynomial can be achieved by resolving singularities (although more efficient algorithms are known --- see, e.g., \cite{berley}). In view of this, the case where $$P(\bm{y})\;=\;P_2(\bm{y})\;=\; \bm{y}^{\bm{k}}$$is a monomial with integral powers $\bm{k}=\left(k_1, \dots, k_n\right)\in\N_0^n\backslash\left\{\bm{0}\right\}$ is of particular interest. As established in~\cite[p.49]{igusa}, the relation $$\bm{\partial}^{\bm{k}} P_2(\bm{y})^{s+1}\;=\; \left(\prod_{i=1}^{n}\prod_{j=1}^{k_i}\left(k_is+j\right)\right) \cdot P_2(\bm{y})^{s}$$ implies that the Sato--Bernstein polynomial of $P_2(\bm{y})$ is $$B_{P_2}(s)\;=\; \prod_{i=1}^{n}\prod_{j=1}^{k_i}\left(s+\frac{j}{k_i}\right).$$

In the general case, specialising the identity~\eqref{sbrelation} to the case where $b(s)=B_P(s)$ and $s=-1$ reveals that $B_P(s)$ is always divisible by the factor $s+1$. Furthermore, $B_P(s)=s+1$ if and only if $\mathcal{Z}_\C\left(\nabla P\right)=\emptyset$; that is, if and only if the polynomial $P(\bm{y})$ admits no singular point in $\C$ (see~\cite[p.48 (v)]{igusa} for a proof of the direct implication and~\cite{brima} for the converse). 

\paragraph{} The proof of Theorem~\ref{volestim} (and also that of Theorem~\ref{lctordermultplicity} in Chapter~\ref{lctrootpole}) requires the concept of a local Sato--Bernstein polynomial $B_{P,\bm{y}_0}^{\K}(Y, s)$  at a point $\bm{y}_0$ lying in the topological closure of  a domain $Y$ which, when $\K=\C$, will always be chosen as $Y=\C^n$. Then, for the sake of simplicity of notation, set $B_{P,\bm{y}_0}^{\C}(\C^n, s)= B_{P,\bm{y}_0}^{\C}( s).$ In the case where $\K=\R$, the domain $Y\subset\R^n$ shall be assumed to be contained in the topological closure of its interior and also to be  globally semianalytic in the sense that it is defined as a finite union of sets of the form  
\begin{equation}\label{setomega}
\left\{\bm{x}\in\R^n\; :\; f_i(\bm{x})\;\square_i\; 0\; \textrm{ for all }\; 1\le i\le l\right\}.
\end{equation} 
Here, for each $1\le i \le l$ (with $l\ge 1$ an integer), the symbol $\square_i\in\{\ge, >\}$ stands for an inequality and  $f_i~: \R^n\rightarrow\R^n$ is an analytic map.\\

The definition of the local polynomial $B_{P,\bm{y}_0}^{\K}(Y, s)$ relies on Lemma~\ref{lemmeintersb} below. Before stating it, given an open set $\mathcal{U}\subset\K^n$ intersecting non--trivially the interior of $Y$, denote  by $\underline{B}_{P}^{\K}(\mathcal{U}\cap Y, s)$ an auxiliary Sato--Bernstein polynomial. It is defined from the restriction of $P(\bm{y})$ to the set $\mathcal{U}\cap Y$ as the monic polynomial of least degree satisfying the identity~\eqref{sbrelation} for all $\bm{y}\in \mathcal{U}\cap Y$ when the differential operator $\mathcal{D}\left(\bm{y}, s, \bm{\partial}\right)$ is allowed to be a polynomial in $s$ and $\bm{\partial}$ with coefficients which are power series  in $\bm{y}\in\K^n$  convergent in $\mathcal{U}\cap Y$. Since the global Bernstein--Sato polynomial $B_P$ satisfies such a relation, $\underline{B}_{P}^{\K}(\mathcal{U}\cap Y, s)$ divides $B_{P}(s)$. In particular, the conclusions of Kashiwara's Theorem~\ref{kashiroot} remain true for the polynomial $\underline{B}_{P}^{\K}(\mathcal{U}\cap Y, s)$.

\begin{lem}\label{lemmeintersb}
Keep the above notation and let $\bm{y}_0$ be a point in the topological closure of the set $Y$. Then, there exists an open neighbourhood $\mathcal{U}_{\bm{y}_0}$ of $\bm{y}_0$ such that for any other open neighbourhood $\mathcal{V}_{\bm{y}_0}$, the polynomial $\underline{B}_{P}^{\K}(\mathcal{U}_{\bm{y}_0}\cap Y, s)$ divides the polynomial $\underline{B}_{P}^{\K}(\mathcal{V}_{\bm{y}_0}\cap Y, s)$. 
\end{lem}

\begin{proof}
Consider a neighbourhood $\mathcal{U}_1$ of $\bm{y}_0$ and assume that it does not satisfies the conclusion of the lemma.  \sloppy Then, there exists another neighbourhood $\mathcal{U}_2$ of $\bm{y}_0$ such that $\underline{B}_{P}^{\K}\left(\mathcal{U}_1\cap Y, s\right)$ does not divide $\underline{B}_{P}^{\K}\left(\mathcal{U}_2\cap Y, s\right)$. Since the Sato--Bernstein identities on  $\mathcal{U}_1\cap Y$ and $\mathcal{U}_2\cap Y$ can be restricted to the intersection $\mathcal{U}_1\cap\mathcal{U}_2\cap Y$, the polynomial $\underline{B}_{P}^{\K}\left(\mathcal{U}_1\cap\mathcal{U}_2\cap Y, s\right)$ is a proper factor of the polynomial $\underline{B}_{P}^{\K}\left(\mathcal{U}_1\cap Y, s\right)$.  \\

If the neigbourhood $\mathcal{U}_1\cap\mathcal{U}_2$ does not satisfy the property stated in the lemma, the process can be iterated~: there exists a neighbourhood $\mathcal{U}_3$ of $\bm{y}_0$ such that the polynomial $\underline{B}_{P}^{\K}(\mathcal{U}_1\cap\mathcal{U}_2\cap\mathcal{U}_3\cap Y, s)$ is a proper factor of $\underline{B}_{P}^{\K}(\mathcal{U}_1\cap Y, s)$.\\

The number of possible iterations in this process is finite since the number of proper factors of the polynomial $\underline{B}_{P}^{\K}(\mathcal{U}_1\cap Y, s)$ is finite. If the process stabilises after $r\ge 1$ steps, then the neighbourhood of $\bm{y}_0$ defined as $\mathcal{U}_{\bm{y}_0}=\cap_{i=1}^{r}\mathcal{U}_i$ meets the required property.
\end{proof}

Given $\bm{y}_0$ lying in the closure of $Y$, let $\mathcal{U}_{\bm{y}_0}$ be a neighbourhood of $\bm{y}_0$ satisfying the conclusion of Lemma~\ref{lemmeintersb}. Recall that when $\K=\R$, the set $Y\subset\R^n$ is assumed to be globally semianalytic and to be contained in the topological closure of its interior. The \emph{local Sato--Bernstein polynomial $B_{P,\bm{y}_0}^{\K}(Y, s)$ with respect to the domain $Y\subset\K^n$  at the point $\bm{y}_0$} is then 
\begin{equation*}
B_{P,\bm{y}_0}^{\K}(Y, s)\;=\; \underline{B}_{P}^{\K}(\mathcal{U}_{\bm{y}_0}\cap Y, s)\;\in\;\Q[s].
\end{equation*}
It is thus also a split polynomial with negative rational roots dividing the global Sato--Bernstein polynomial $B_P(s)$. The \emph{(global) Sato--Bernstein polynomial $B_{P}^{\K}(Y, s)$ with respect to the domain $Y\subset\K^n$} is 
\begin{equation*}
B_{P}^{\K}(Y, s)\;=\; \underset{\bm{y}_0\in Y}{\textrm{lcm}}\;
B_{P,\bm{y}_0}^{\K}(Y, s)\;\in\;\Q[s].
\end{equation*}
In view of Lemma~\ref{lemmeintersb}, it is the lowest common multiple of finitely many local Sato--Bernstein polynomials when the closure of $Y$ is compact.

\paragraph{} In the case when $\K=\C$, it is known, see~\cite[Lemma~2.5.2]{gyoja},  that the local and global Sato--Bernstein polynomials are related by the formula 
\begin{equation}\label{locglobsb}
B_P(s)\; =\; \underset{\bm{z}_0\in\C^n}{\textrm{lcm}}\; B_{P,\bm{z}_0}^{\C}(s).
\end{equation}
This relation does not necessarily hold anymore when $\K=\R$ and $Y=\R^n$. \\

Given $\bm{z}_0\in\C^n$ and $\bm{x}_0\in\R^n$, let  from now on and for the sake of simplicity of notation
\begin{equation}\label{globalsbomega}
\widehat{B}_{P, \bm{z}_0}(s)\;=\; B_{P, \bm{z}_0}^{\C}(s)\qquad \textrm{and}\qquad B_{P, \bm{x}_0}(Y, s)\;=\; B_{P, \bm{x}_0}^{\R}(Y, s)
\end{equation}
be the local Sato--Bernstein polynomials at $\bm{z}_0\in\C^n$ and $\bm{x}_0\in\R^n$,  respectively. Let also
\begin{equation}\label{globalsbomegabis}
B_{P}(Y, s)\;=\; B_{P}^{\R}(Y, s)
\end{equation}
be the global Sato--Bernstein polynomial with respect to the domain $Y$.

\section{A Consequence of a Tauberian Theorem}\label{tauberian}

Take $\K=\R$ and let $K\subset\R^n$ be a set satisfying the conditions $(\mathcal{H}_1)$ and $(\mathcal{H}_2)$. Assume in this section that the polynomial $P(\bm{x})\in\R[\bm{x}]$ is homogeneous of degree $q\ge 1$. Theorem~\ref{volestim} is established in the case of the polynomial $P(\bm{x})$ upon estimating the volume of the set 
\begin{equation}\label{ensalgbis}
\underline{\mathcal{S}}_{P}(K, b(T))\; =\; \left\{\bm{x}\in T\cdot K\; :\; \left|P\left(\bm{x}\right)\right|\le b(T)\right\}
\end{equation}
under the successive assumptions of Cases (1), (2) and (3). In all cases, the limit condition~\eqref{limiassump} is assumed to be verified with the exponent $q$ instead of  $d$; namely, it is assumed throughout that 
\begin{equation}\label{limiassumpbis}
c(T)\;:=\; \frac{b(T)}{T^q}\;\underset{T\rightarrow\infty}{\longrightarrow}\; 0.
\end{equation}

The zeta distribution associated to the polynomial $P(\bm{x})$ is here defined by setting 
\begin{equation}\label{zetaP}
\left\langle \zeta_{P}\left(s\right), \psi\right\rangle\;=\; \int_{\R^n} \left|P(\bm{x})\right|^{-s}\cdot\psi(\bm{x})\cdot\textrm{d}\bm{x}
\end{equation}
for any complex number $s$ with non--positive real part and any $\psi$ in $C_c^{\infty}\left(\R^n\right)$.\\

Denote by $\mathfrak{O}^+$ the collection of all those non--empty open sets obtained by taking the interiors of the connected components of the region 
\begin{equation}\label{defregion}
\left\{\bm{x}\in\R^n\; : \; P(\bm{x})>0\right\}.
\end{equation}
It is then known that $\mathfrak{O}^+$ has a finite cardinality (see~\cite[Theorem 2.23]{coste}) and that each element of $\mathfrak{O}^+$ is a semialgebraic domain (i.e.~a domain which is a finite union of sets  such as~\eqref{setomega} when the maps $f_i$ are polynomials for all  $1\le i\le l$) --- see~\cite[Proposition 3.1]{basu} for a proof of this claim. 
These properties also hold for the collection of sets $\mathfrak{O}^-$ defined as above when replacing $P(\bm{x})$ with $-P(\bm{x})$ in~\eqref{defregion}, and therefore also for their union 
\begin{equation}\label{defoobig}
\mathfrak{O}=\mathfrak{O}^+\cup\mathfrak{O}^-.
\end{equation}

From the homogeneity of the polynomial $P(\bm{x})$,  the volume of the  domain~\eqref{ensalgbis} can then be expanded as 
\begin{align}
\V_n\left(\underline{\mathcal{S}}_{P}(K, b(T))\right)&\;=\; \int_{\R^n}\chi_{\left\{\left|P(\bm{x})\right|\le b(T)\right\}}\cdot\chi_{\left\{\bm{x}\in T\cdot K\right\}}\cdot\textrm{d}\bm{x}\nonumber\\
&\;\underset{\eqref{limiassumpbis}}{=}\; T^n\int_{K} \chi_{\left\{\left|P(\bm{x})\right|\le c(T)\right\}}\cdot\textrm{d}\bm{x}  \label{decompvol0}\\
&\;=\; T^n\sum_{\Omega\in\mathfrak{O}} \int_{K}\chi_{\Omega}(\bm{x})\cdot \chi_{\left\{\left|P(\bm{x})\right|\le c(T)\right\}}\cdot\textrm{d}\bm{x}.\label{decompvol}
\end{align}

The explicit form of the volume in~\eqref{decompvol0} will be used to determine its asymptotic order of growth in Case (2). The decomposition~\eqref{decompvol} will be required to establish the precise asymptotic  behavior of the volume in  Case (3) at a level of generality needed to also prove related statements in Chapter~\ref{ranunidisfamhomfor}. This level of generality essentially requires the determination of the asymptotic behavior of each of the integrals appearing in the sum~\eqref{decompvol} when the parameter  $T$ tends to infinity. This is achieved by introducing the zeta distributions induced by the corresponding decomposition.\\

To do so and to fix the notations, given $\Omega\in\mathfrak{O}$, assume  that  
\begin{equation}\label{defPomega}
\forall \bm{x}\in\Omega, \; \; \; P_{\Omega}(\bm{x})\;>\; 0,
\end{equation}
where $P_{\Omega}(\bm{x})$ is either the polynomial $P(\bm{x})$ or its opposite. Extend the map $\bm{x}\mapsto P_{\Omega}(\bm{x})$ to $\R^n$ by setting $P_{\Omega}(\bm{x})=0$ for all $\bm{x}\not\in\Omega$. The restricted zeta distribution $\zeta_{P}(\Omega, \,\cdot\,)$ is then defined by setting for all $\psi$  in $\mathcal{C}_{c}^{\infty}(\R^n)$ and all $s\in\C$ with $Re(s)\le 0$
\begin{equation}\label{restriczeta}
\left\langle \zeta_{P}\left(\Omega, s\right), \psi\right\rangle\;=\; \int_{\R^n} P_{\Omega}(\bm{x})^{-s}\cdot\psi(\bm{x})\cdot\textrm{d}\bm{x}.
\end{equation}
Clearly, it is analytic in the domain of complex numbers with non-positive real parts. Furthermore, it holds that 
\begin{equation}\label{zetaPrel}
\zeta_{P}\;=\; \sum_{\Omega\in\mathfrak{O}} \zeta_{P}(\Omega, \,\cdot\,).
\end{equation}
Given $\Omega\in\mathfrak{O}$, whenever well-defined, let
\begin{equation}\label{romegmomeg}
\left(r_P, m_P\right)\qquad \textrm{and}\qquad \left(r_P(\Omega), m_P(\Omega)\right)
\end{equation}
be the pairs made from the smallest real poles $r_P$  and $r_P(\Omega)$ of the distributions $\zeta_{P}$  and $\zeta_{P}(\Omega, \,\cdot\,)$, respectively, and from the corresponding respective orders $m_P\ge 1$ and $m_P(\Omega)\ge 1$. If either of these quantities does not exist, set conventionally  $\left(r_P, m_P\right)=\left(\infty, -\right)$ and $\left(r_P(\Omega), m_P(\Omega)\right)=\left(\infty, -\right)$ (in particular, the integers determining the orders are then left undefined).\\

Recall here that a meromorphic distribution, say $\zeta$, has a pole at $r\in\C$ of order $m\ge 1$ if for all test functions $\psi$ in $\mathcal{C}_{c}^{\infty}(\R^n)$, the map $s\mapsto \left(s-r\right)^m\cdot\left\langle \zeta(s), \psi\right\rangle$ is holomorphic and if $m$ is the smallest integer with this property.

\begin{lem}\label{propzetamerogro}
Let $\Omega\in\mathfrak{O}$. The distribution $\zeta_{P}\left(\Omega, \, \cdot\,\right)$ defined by~\eqref{restriczeta} extends meromorphically to the entire complex plane. Furthermore, its poles, whenever they exist, lie in a finite number of arithmetic progressions contained in the interval $(0, \infty)$. The same therefore holds for the distribution $\zeta_{P}$.
\end{lem}

Let $B_P(s)\in\Q[s]$ be the (global) Sato--Bernstein polynomial of the homogeneous form $P(\bm{x})$. Before justifying the claims made in Lemma~\ref{propzetamerogro}, denote by 
\begin{equation}\label{expdiffop}
\mathcal{D}\left(\bm{x}, s, \bm{\partial}\right)\;=\;\sum_{\bm{\alpha}\in\N_0^n} a_{\bm{\alpha}}(\bm{x},s)\bm{\partial}^{\bm{\alpha}}
\end{equation}
a differential operator associated to $B_P(s)$, where the sum has finite support and where $a_{\bm{\alpha}}(\bm{x},s)$ is a polynomial. Define also the \emph{dual operator} 
\begin{equation}\label{opdual}
\mathcal{D}^*\left(\bm{x}, s, \bm{\partial}\right)\;=\;\sum_{\bm{\alpha}\in\N_0^n} (-1)^{|\bm{\alpha}|}\bm{\partial}^{\bm{\alpha}}\left(a_{\bm{\alpha}}(\bm{x},s)\;\cdot\;\right).
\end{equation} 
Given $\psi$ in $\mathcal{C}_{c}^{\infty}(\R^n)$ and  $s\in\C$ such that $\textrm{Re}(s)\ge -1$, it satisfies the integration by parts formula 
\begin{equation*}\label{intbypartsdualop}
\int_{\R^n} \mathcal{D}\left(\bm{x}, s, \bm{\partial}\right)P_{\Omega}(\bm{x})^{s+1}\cdot\psi(\bm{x})\cdot\textrm{d}\bm{x}\;=\;  \int_{\R^n} P_{\Omega}(\bm{x})^{s+1}\cdot\mathcal{D}^*\left(\bm{x}, s, \bm{\partial}\right)\psi(\bm{x})\cdot\textrm{d}\bm{x}.
\end{equation*}
If $\textrm{Re}(s)\le 0$, the defining identity~\eqref{sbrelation} satisfied by the polynomial $B_P(s)$ then implies that 
\begin{align}\label{elemdual}
B_P(-s)\cdot \left\langle \zeta_{P}\left(\Omega, s\right), \psi\right\rangle\;& =\;  \left\langle \zeta_{P}\left(\Omega, s-1\right), \mathcal{D}^*\left(\bm{x}, -s, \bm{\partial}\right)\cdot\psi(\bm{x})\right\rangle.
\end{align}

\begin{proof}[Proof of Lemma~\ref{propzetamerogro}.] The argument is well--known and is reproduced here for the sake of completeness as it is short and as it will be referred to often in what follows. Let  $s$ be a complex number such that $\textrm{Re}(s)\le 0$ and let $\psi$ be in $\mathcal{C}_{c}^{\infty}(\R^n)$~: since in the equation~\eqref{elemdual}, the right--hand side is well--defined whenever $\textrm{Re}(s)\le 1$, dividing this relation by $B_P(-s)$ and proceeding recursively shows the existence of the meromorphic continuation over the entire complex plane. This relation also implies that the real poles of the distribution $\zeta_{P}\left(\Omega, \,\cdot\,\right)$, whenever they exist, are roots of the polynomial $B_P(-s)$ modulo one. This completes the proof of  Lemma~\ref{propzetamerogro} since the last claim is a direct consequence of the decomposition~\eqref{zetaPrel}.\\

Alternative proofs of this (well--known) result based on resolutions of singularities can also be found in the classical works by Atiyah~\cite{atiyah} and by Bernstein and Gel'fand~\cite{berngel}.
\end{proof}

The following lemma relates the properties of the zeta distributions under consideration to the volume of the domains appearing in the expansion~\eqref{decompvol0} and in the decomposition~\eqref{decompvol}. The first two claims in the statement are particular cases adapted to the present situation of  a Tauberian Theorem due to Chambert--Loir and Tschinkel, see~\cite[Appendix A, Theorem A.1]{chambloitshcik}.\\ 

In order to state the lemma, given $\psi$  in $\mathcal{C}_{c}^{\infty}(\R^n)$ , denote by
\begin{equation}\label{deforderpole} 
\left(r_P(\psi),\; m_P(\psi)\right)\qquad \textrm{and}\qquad\left(r_P(\Omega, \psi),\; m_P(\Omega, \psi)\right)
\end{equation}
the pairs made from the smallest real poles  $r_P(\psi)\ge 0$ and $r_P(\Omega, \psi)\ge 0$ of the meromorphic functions $s\in\C \mapsto \left\langle \zeta_{P}\left(s\right), \psi\right\rangle$  and $s\in\C \mapsto \left\langle \zeta_{P}\left(\Omega, s\right), \psi\right\rangle$, respectively, and from the corresponding orders $m_P(\psi)\ge 1$ and $m_P(\Omega, \psi)\ge 1$, respectively. When either of these quantities is not well-defined, set here again conventionally $\left(r_P(\psi), m_P(\psi)\right)=\left(\infty, -\right)$ and $\left(r_P(\Omega, \psi), m_P(\Omega, \psi)\right)=\left(\infty, -\right)$. Also, when $c$ is a real lying in the interval  $(0,1)$, $\psi$ a nonnegative map in $\mathcal{C}_{c}^{\infty}(\R^n)$ and $\Omega$ a set in $\mathfrak{O}$, let 
\begin{equation}\label{defmuvolweighted}
\mu_P(\psi, c)\;=\; \int_{\R^n}\chi_{\left\{\left|P(\bm{x})\right|\le c\right\}}\cdot\psi(\bm{x})\cdot\textrm{d}\bm{x}\quad\textrm{and}\quad\mu_P(\Omega, \psi, c)\;=\; \int_{\R^n}\chi_{\left\{P_{\Omega}(\bm{x})\le c\right\}}\cdot\psi(\bm{x})\cdot\textrm{d}\bm{x}.
\end{equation}

\begin{lem}\label{lemloir} 
Let $c\in (0,1)$ be a real and let $\psi$ be a non--negative map  in $\mathcal{C}_{c}^{\infty}(\R^n)$.\\

Assume first that $r_P(\psi) \in (0, \infty)$ and thus that $m_P(\psi)\ge 1$. Then, there exists a univariate monic polynomial $R_{\psi}(x)\in\R[x]$ of degree $m_P(\psi)-1$ and some $\varepsilon>0$ such that, as $c$ tends to zero, 
\begin{equation}\label{asympmu0}
\mu_P(\psi, c)\;=\;  \Theta_P(\psi)\cdot c^{r_P(\psi)}\cdot R_{\psi}\left(\left|\log c\right|\right) + O\left(c^{r_P(\psi)+\varepsilon}\right).
\end{equation} 
Here, $$ \Theta_P(\psi)\; =\; \frac{1}{m_P(\psi)!}\lim_{\sigma\rightarrow r_P(\psi)} (\sigma-r_P(\psi))^{m_P(\psi)}\left\langle \zeta_{P}\left(s\right), \psi\right\rangle\;>\; 0.$$\\

Fix from now on $\Omega\in\mathfrak{O}$.  Assume that $r_P(\Omega, \psi)\in (0, \infty)$  and thus that  $m_P(\Omega, \psi)\ge 1$. Then, there exists a univariate monic polynomial $R_{\psi}^{\Omega}(x)\in\R[x]$ of degree $m_P(\Omega, \psi)-1$ and some $\varepsilon>0$ such that, as $c$ tends to zero, 
\begin{equation}\label{asympmu}
\mu_P(\Omega, \psi, c)\;=\;  \Theta_P(\Omega, \psi)\cdot c^{r_P(\Omega, \psi)}\cdot R_{\psi}^{\Omega}\left(\left|\log c\right|\right) + O\left(c^{r_P(\Omega, \psi)+\varepsilon}\right).
\end{equation} 
Here, $$ \Theta_P(\Omega, \psi)\; =\; \frac{1}{m_P(\Omega, \psi)!}\lim_{\sigma\rightarrow r_P(\Omega, \psi)} (\sigma-r_P(\Omega, \psi))^{m_P(\Omega, \psi)}\left\langle \zeta_{P}\left(\Omega, s\right), \psi\right\rangle\;>\; 0.$$\\

Assume finally that $\mathcal{F}$  is a family of nonnegative maps in $\mathcal{C}_{c}^{\infty}(\R^n)$ satisfying the following two properties~: 
\begin{itemize}
\item[(A1)] [uniform boundedness of the support] there exists a compact set $\mathcal{K}\subset\R^n$ such that  $\textrm{Supp } \psi\subset \mathcal{K}$ for all $\psi\in\mathcal{F}$;
\item[(A2)] [uniform boundedness of the range] there exists $\varphi_{\mathcal{F}}$ in $\mathcal{C}_{c}^{\infty}(\R^n)$ such that for all $\psi\in\mathcal{F}$, it holds that $0\le \psi\le \varphi_{\mathcal{F}}$ and, with the above notation, that $\left(r_P\left(\Omega, \psi\right), m_P\left(\Omega, \psi\right)\right) = \left(r_P\left(\Omega, \varphi_{\mathcal{F}}\right), m_P\left(\Omega,  \varphi_{\mathcal{F}}\right)\right)$  (in particular, the pair $\left(r_P\left(\Omega, \psi\right), m_P\left(\Omega, \psi\right)\right)$ remains constant over the choice of $\psi\in\mathcal{F}$).
\end{itemize}
Then, the following uniformity claims on the parameters appearing in the expansion~\eqref{asympmu} hold when the map $\psi$ varies in $\mathcal{F}$~: 
\begin{itemize}
\item[(V1)] the leading coefficient $\Theta_P(\Omega, \psi)$ is uniformly bounded above by a constant depending on the choice of the smooth map  $\varphi_{\mathcal{F}}$;
\item[(V2)] there exist integers $M, N\ge 1$ depending on the polynomial $P(\bm{x})$ and on the compact set $\mathcal{K}$  and, given a value of $\varepsilon>0$ and a choice of $\varphi_{\mathcal{F}}$ as above, there exists a constant $\Gamma(\mathcal{F}, \varphi_F, \varepsilon)>0$  such that the following holds~: the  implicit constant in the error term can be taken as $2\max\left\{\Gamma(\mathcal{F}, \varphi_{\mathcal{F}}, \varepsilon), \left(C_{N}(\mathcal{F}, \mathcal{K})\right)^M\right\}$ provided that the quantity 
\begin{equation}\label{defcnfk}
C_{N}(\mathcal{F}, \mathcal{K})\; : =\; \sup_{\psi\in\mathcal{F}}\; \max_{0\le \left|\bm{k}\right|\le N} \; \max_{\bm{x}\in\mathcal{K}}\; \left|\frac{\partial^{\left|\bm{k}\right|}\psi}{\partial\bm{x}^{\bm{k}}}(\bm{x})\right|
\end{equation}
is finite. In particular, as soon as $\max_{\bm{x}\in\mathcal{K}}\, \psi(\bm{x}) = 1$ for some $\psi\in\mathcal{F}$, this implicit constant can be taken as  $\left(C_{N}(\mathcal{F}, \mathcal{K})\right)^M$ up to a multiplicative constant depending on $P(\bm{x})$, on $\mathcal{F}$, on $\mathcal{K}$, on a given a value of $\varepsilon>0$ and on the choice of $\varphi_{\mathcal{F}}$;
\item[(V3)] there exist integers $M, N\ge 1$ depending on the polynomial $P(\bm{x})$ and on the compact set $\mathcal{K}$ such that the coefficients of the univariate polynomial  $\Theta_P(\Omega, \psi)\cdot  R_{\psi}^{\Omega}\left(x\right) $ are, up to a multiplicative constant depending on $\mathcal{K}$, upper bounded  by the above defined quantity $\left(C_{N}(\mathcal{F}, \mathcal{K})\right)^M$, provided it is finite.
\end{itemize}

With obvious modifications (consisting in dropping all notational dependencies on the set $\Omega$), the uniformity claims (V1)--(V3) also hold under the assumptions (A1) and (A2) for the expansion~\eqref{asympmu0}. 
\end{lem}

All the coefficients of the polynomials $R_{\psi}(x), R_{\psi}^{\Omega}(x)\in\R[x]$ can be made explicit --- see, for instance, \cite[Chap.~II.5]{tenen}. Also, the main interest of the above statement is that the growths of the functions $c\mapsto\mu_P(\psi, c)$ and $c\mapsto\mu_P(\Omega, \psi, c)$ when $c$ tends to zero are expressed as functions of the smallest poles of the zeta functions $\left\langle \zeta_{P}, \psi\right\rangle$  and $\left\langle \zeta_{P}\left(\Omega, \;\cdot\;\right), \psi\right\rangle$, respectively, and of their orders. Such  growth rates can indeed otherwise be obtained from resolutions of singularities not relying on the properties of the zeta functions under consideration (albeit the error terms are then not so explicit) --- see, e.g., \cite[Proposition 7.2]{benoitoh} for details.

\begin{proof}[Proof of Lemma~\ref{lemloir}] Let $\Omega\in\mathfrak{O}$. It is enough to establish the relation~\eqref{asympmu} as~\eqref{asympmu0} follows upon reproducing \emph{verbatim} the same proof after removing the references to the set $\Omega$ in the successive relations. \\

For the sake of the simplicity of the notation, set in this proof $\left(r, m\right)=\left(r_P(\Omega, \psi), m_P(\Omega, \psi)\right)$ and recall that it is assumed that $r\in (0, \infty)$.\\
 
Then, with evident notational modifications, the volume estimate~\eqref{asympmu} is Theorem A.1 in~\cite[Appendix A]{chambloitshcik} applied to the case where the map $f$ therein is taken as $\bm{x}\mapsto P_{\Omega}(\bm{x})$  and where the mesure $\mu$ therein is defined for any Borel set $A\subset\R^n$ by setting $\mu(A)=\int_{A}\psi(\bm{x})\cdot\textrm{d}\bm{x}.$ According to this statement, the relation~\eqref{asympmu} holds provided that two assumptions are met~: 
\begin{itemize}
\item[(a)] the meromorphic function $s\in\C \mapsto \left\langle \zeta_{P}\left(\Omega, s\right), \psi\right\rangle$ has no pole in the half--plane $Re(s)\le r-\delta$ for some value of $\delta\in (0,1)$ (which determines the value of $\varepsilon$ in the statement);
\item[(b)] this function has moderate growth rate on the vertical strip $r-\delta+i\R$ in the sense that there exists $\kappa>0$  such that for any $\tau\in\R$, 
\begin{equation}\label{constantkappa}
\left|\left\langle \zeta_{P}\left(\Omega, r-\delta+i\tau\right), \psi\right\rangle\right|\;\ll\; \left(1+\left|\tau\right|\right)^\kappa.
\end{equation}
\end{itemize}
Assumption (a) follows immediately from Lemma~\ref{propzetamerogro}, which gives the location of the poles of the meromorphic function $\left\langle \zeta_{P}\left(\Omega, \;\cdot\;\right), \psi\right\rangle$  as a subset of a finite number of arithmetic progressions. As for (b), fix $\delta>0$ satisfying the first assumption and take it smaller if needed to ensure that the Sato--Bernstein polynomial $B_P(s)$ does not vanish at $-r+\delta$ (and therefore on the vertical strip $-r+\delta+i\R$, since its roots are rational). Then, denoting by $\mathcal{D}^*\left(\bm{x}, s, \bm{\partial}\right)$ the dual operator associated  to $B_P(s)$ and expanding it as in~\eqref{opdual}, one obtains from the equation~\eqref{elemdual} that
\begin{align*}
\left|B_P( -r+\delta-i\tau)\right|\cdot & \left|\left\langle \zeta_{P}\left(\Omega, r-\delta+i\tau\right), \psi\right\rangle\right|\\ 
&=\;\left|\left\langle \zeta_{P}\left(\Omega, r-\delta-1+i\tau\right), \mathcal{D}^*\left(\bm{x}, -r+\delta-i\tau, \bm{\partial}\right)\cdot\psi(\bm{x})\right\rangle\right|\\
&\le\; \sum_{\bm{\alpha}\in\N_0^n} \left| \left\langle\zeta_{P}^\R\left(\Omega, r-\delta-1+i\tau\right), \bm{\partial}^{\bm{\alpha}}\left(a_{\bm{\alpha}}(\bm{x}, -r+\delta-i\tau)\cdot\psi(\bm{x})\right)   \right\rangle   \right|\\
&\le\; \sum_{\bm{\alpha}\in\N_0^n} \left\langle \zeta_{P}\left(\Omega, r-\delta-1\right),\left|\bm{\partial}^{\bm{\alpha}}\left(a_{\bm{\alpha}}(\bm{x}, -r+\delta-i\tau)\cdot\psi(\bm{x})\right)\right| \right\rangle,
\end{align*}
where the sum has a finite support. The moderate growth assumption is then easily seen to be implied by the properties that $\inf_{\tau\in\R} \left|B_P(-r+\delta-i\tau)\right|>0$, that $a_{\bm{\alpha}}(\bm{x}, s)$ is a polynomial in $s$ and that $r-\delta-1$ is not a pole of the distribution $\zeta_{P}\left(\Omega, \;\cdot\;\right)$.\\

As for the uniformity claims, as above, it is here again enough to establish them in the case of the expansion~\eqref{asympmu}. Note then first that the assumption (A2) implies that $\mu_P\left(\Omega, \psi, c\right)\le \mu_P\left(\Omega, \varphi_{\mathcal{F}}, c\right)$. Letting $c$ tend to zero, one infers from the same assumption that 
\begin{equation*}
\Theta_P\left(\Omega, \psi\right)\;\le\; \Theta_P\left(\Omega, \varphi_{\mathcal{F}}\right),
\end{equation*}
whence (V1). \\

It is easier to establish the claim (V3) before (V2). To this end, the proof of~\cite[Appendix A, Theorem A.1]{chambloitshcik} shows that the coefficients of the polynomial $\Theta_P(\Omega, \psi)\cdot  R_{\psi}^{\Omega}(x)$ depend polynomially on the con\-stants $r$ and $\kappa$ and also on the coefficients of the (well--defined) univariate polynomial $\widetilde{R}_{\psi}^{\Omega}(x)$ of degree $m-1$ determined by the residue relation 
\begin{align}
B^{-r}\cdot \widetilde{R}_{\psi}^{\Omega}\left(\log B\right)\; &=\; \textrm{Res}_{s=r}\left(\frac{B^{-s}\cdot \left\langle \zeta_{P}\left(\Omega, s\right), \psi\right\rangle}{s^{k+1}}\right)\nonumber\\
&=\; \frac{1}{(m-1)!}\cdot \lim_{s\rightarrow r}\left[\frac{\textrm{d}^{m-1}}{\textrm{d}s^{m-1}}\left((s-r)^m\cdot s^{-k-1}\cdot B^{-s}\cdot \left\langle \zeta_{P}\left(\Omega, s\right), \psi\right\rangle\right)\right].\label{residueeqt}
\end{align} 
Here, $k$ is again an integer larger than $\kappa>0$. To isolate the coefficients of $\widetilde{R}_{\psi}^{\Omega}(x)$, expand $\left\langle \zeta_{P}\left(\Omega, s\right), \psi\right\rangle$ as a Laurent series around $s=r$~: 
\begin{equation}\label{laurentexpanszeta}
\left\langle \zeta_{P}\left(\Omega, s\right), \psi\right\rangle\;=\; \sum_{l=1}^{m}\frac{a_{r, l}\left(\Omega, \psi\right)}{\left(s-r\right)^{l}}\;+\; \left\langle h_{P, r}\left(\Omega, s\right), \psi\right\rangle,
\end{equation} 
where $h_{P, r}\left(\Omega, \, \cdot\,\right)$ is a  distribution holomorphic in a neighbourhood of $s=r$. Cauchy's Integral Formula implies that the maps $\psi\in \mathcal{C}_{c}^{\infty}(\R^n) \mapsto a_{r,l}\left(\Omega, \psi\right)$ ($1\le l \le m$) 
determining the coefficients of the Laurent expansion are well--defined distributions (see~\cite[Appendix A.2.3]{gelshi} for details). This is saying that these linear functionals are continuous, thereby guaranteeing the existence of an integer $N\ge 1$ and of an implicit constant both depending only on the compact set $\mathcal{K}$ such that for all $\psi$ in $\mathcal{C}_{c}^{\infty}(\R^n)$ with support contained in $\mathcal{K}$, 
\begin{equation}\label{ineqlaurent}
\left|a_{r,l}\left(\Omega, \psi\right)\right|\;\ll\; \max_{0\le \left|\bm{k}\right|\le N}\; \max_{\bm{x}\in\mathcal{K}}\; \left|\frac{\partial^{\left|\bm{k}\right|}\psi}{\partial\bm{x}^{\bm{k}}}(\bm{x})\right|.
\end{equation} 
Inserting the expansion~\eqref{laurentexpanszeta} into the equation~\eqref{residueeqt},  the coefficients of the polynomial $\widetilde{R}_{\psi}^{\Omega}(x)$ are seen to be linear combinations of elements in the set $\left\{1, a_{r,1}\left(\Omega, \psi\right), \dots, a_{r,m}\left(\Omega, \psi\right)\right\}.$ The inequality~\eqref{ineqlaurent} then establishes (V3).\\

As for the claim (V2), it follows from the proof of~\cite[Appendix A, Theorem A.1]{chambloitshcik}  that the implicit constant in the error term in~\eqref{asympmu} can be taken as the sum of two quantities~: 

\begin{itemize}
\item a multiple, depending on the exponent $r= r_P\left(\Omega, \psi\right)$, on the order $m= m_P\left(\Omega, \psi\right)$ and on a small enough parameter $\delta>0$ (determining the value of $\epsilon>0$ in~\eqref{asympmu}) of the factor $$\left|\int_{r-\delta+i\R}\frac{\textrm{d}s}{s^{k+1}}\cdot\int_{\R^n}\frac{\psi(\bm{x})}{P_{\Omega}(\bm{x})^s}\cdot\textrm{d}\bm{x}\right|.$$Here, $k$ is any integer larger than the exponent $\kappa>0$ appearing in~\eqref{constantkappa}. This factor is clearly bounded above by a multiple (depending on $r, \delta$ and $\kappa$) of the nonnegative quantity $\left\langle \zeta_{P}\left(\Omega, r-\delta\right), \psi\right\rangle$  which, in turn, is bounded above by the constant $\left\langle \zeta_{P}\left(\Omega, r-\delta\right), \varphi_{\mathcal{F}}\right\rangle$ uniformly in $\psi\in\mathcal{F}$. This constant is well-defined from the assumption (A2);

\item the coefficients of the polynomials $\Theta_P(\Omega, \psi)\cdot  R_{\psi}^{\Omega}(x)$ (this dependency is captured in the error term denoted by $\widetilde{O}\left(B^{a-\varepsilon}\right)$ in the proof of~\cite[Appendix A, Lemma A.5]{chambloitshcik}). As established in (V3) and with the notation therein, these coefficients are, up to a multiplicative constant depending on $\mathcal{K}$, bounded above by the quantity $\left(C_{N}(\mathcal{F}, \mathcal{K})\right)^M$, whenever it is finite.
\end{itemize}
These two observations prove the first part of the claim (V2). As far as the second one is concerned, it is enough to note that, under the assumption that $\max_{\bm{x}\in\mathcal{K}} \psi(\bm{x})=1$ for some $\psi\in\mathcal{F}$, it is immediate that $C_{N}(\mathcal{F}, \mathcal{K})\ge 1$. As a consequence, $$2\max\left\{\Gamma(\mathcal{F}, \varphi_F, \varepsilon), \left(C_{N}(\mathcal{F}, \mathcal{K})\right)^M\right\}\;\le 2\max\left\{1, \Gamma(\mathcal{F}, \varphi_F, \varepsilon)\right\}\cdot \left(C_{N}(\mathcal{F}, \mathcal{K})\right)^M.$$

This completes the proof of Lemma~\ref{lemloir}.
\end{proof}

The goal is now to deduce from Lemma~\ref{lemloir} sharp volume estimates for the set~\eqref{ensalgbis} defined from the polynomial $P(\bm{x})$ in the three cases considered in Theorem~\ref{volestim}. These estimates are established under the assumption that the largest poles of the various meromorphic distributions under consideration are well-defined. The proofs of the existence and of the determination of these poles are the subject of the next Section~\ref{locglobrlct} (see in particular Proposition~\ref{proprlctgene} therein).\\

Consider first Case (1). If $\psi$ is a map in $\mathcal{C}_{c}^{\infty}(\R^n)$  satisfying the inequalities
\begin{equation*}
0\;\le\; \chi_K\;\le\; \psi
\end{equation*} 
which is such that the pair $\left(r_P(\psi), m_P(\psi)\right)$ introduced in~\eqref{deforderpole} is well-defined with $r_P(\psi)>0$, then an immediate consequence of the expansion~\eqref{decompvol0} and of the relation~\eqref{asympmu0} in the above Lemma~\ref{lemloir} is that 
\begin{equation*}\label{case0thmvol}
\V_n\left(\underline{\mathcal{S}}_{P}(K, b(T))\right)\;\ll\; T^n\cdot c(T)^{r_{P}(\psi)}\cdot\left|\log\left(c(T)\right)\right|^{m_{P}(\psi)-1}.
\end{equation*}

This is enough to establish the volume estimate~\eqref{formulevol0} in Theorem~\ref{volestim} upon specialising this inequality to the case where $P(\bm{x})$ is the homogeneous polynomial $\left\|\bm{F}(\bm{x})\right\|^2$ of degree $q=2d$  and  upon noticing that the poles and orders of the distribution $\zeta_{\bm{F}}$ defined in~\eqref{distrizetaR} then coincide with those of the distribution $\zeta_P(\,\cdot\, / 2)$ defined by~\eqref{zetaP}.\\

Consider now Case (2), thus assuming that the condition $(\mathcal{H}_3)$ additionally holds. Denote by $\left(r_P(K), m_P(K)\right)$ the pair provided by $(\mathcal{H}_3)$ in the case of the polynomial $P(\bm{x})$ (recall that the assumption $(\mathcal{H}_3)$ relates the analytic properties of the distribution $\zeta_P$ to the topological properties of the set $K$). The following result, which turns out to be an easy consequence of Lemma~\ref{lemloir}, then provides the sharp order of asymptotic growth of the volume of the set $\underline{\mathcal{S}}_{P}(K, b(T))$. Upon specialising it as above to the case of the polynomial $P(\bm{x})=\left\|\bm{F}(\bm{x})\right\|^2$ , it establishes the relation~\eqref{formulevolbis} in Theorem~\ref{volestim}.

\begin{thm}[Volume Estimate in Case (2)]\label{mainthmvolomega0}
Let  the set $K\subset\R^n$  satisfy the three assumptions  $(\mathcal{H}_1)-(\mathcal{H}_3)$. Assume that the pair $\left(r_P(K), m_P(K)\right)$   is such that $r_P(K)>0$ and that the limit condition~\eqref{limiassumpbis}, where the quantity $c(T)$ is defined, holds. Then for $T\ge 1$, 
\begin{equation*}
\V_n\left(\underline{\mathcal{S}}_{P}(K, b(T))\right)\;\asymp\; T^{n}\cdot c(T)^{r_P\left(K\right)}\cdot\left|\log c(T)\right|^{m_P\left(K\right)-1}.
\end{equation*} 
\end{thm}

For the sake of simplicity of the terminology, say from now on that a smooth test function $\psi$ whose support satisfies the inclusions~\eqref{suppproper} while being strictly positive over the set $C$ meets the \emph{support restriction condition~\eqref{suppproper}}.

\begin{proof}[Proof of Theorem~\ref{mainthmvolomega0}]
Let $\psi_1, \psi_2$ be two maps in $\mathcal{C}_{c}^{\infty}(\R^n)$ meeting the support restriction condition~\eqref{suppproper} such that 
\begin{equation*}
0\;\le\; \psi_1\;\le\; \chi_K\;\le\; \psi_2.
\end{equation*} 
Under the assumption $(\mathcal{H}_3)$, it holds that 
\begin{equation}\label{rpkrppsi1rppsi2}
\left(r_P(K), m_P(K)\right)= \left(r_P(\psi_1), m_P(\psi_1)\right)= \left(r_P(\psi_2), m_P(\psi_2)\right),
\end{equation}  
where, given $\psi$ in $\mathcal{C}_{c}^{\infty}(\R^n)$,  the pair $\left(r_P(\psi),  m_P(\psi)  \right)$ is defined in~\eqref{deforderpole}. Furthermore, the expansion~\eqref{decompvol0} implies that 
\begin{equation*}
\mu_P(\psi_1, c(T))\;\le\; \frac{\V_n\left(\underline{\mathcal{S}}_{P}(\alpha, b(T))\right)}{T^n}\;\le\;  \mu_P(\psi_2, c(T)),
\end{equation*}
where the quantities $\mu_P(\psi_1, c(T))$ and $\mu_P(\psi_2, c(T))$ are explicitly given in~\eqref{defmuvolweighted}. The conclusion of the statement is then a consequence of the relation~\eqref{rpkrppsi1rppsi2} and of Equation~\eqref{asympmu0} in Lemma~\ref{lemloir}.
\end{proof}

Consider finally Case (3); that is, assume from now on that $K$ is a body star-shaped with respect to the origin containing the origin in its interior. With a view towards the metric counting results proved in Chapter~\ref{ranunidisfamhomfor}, the volume estimate~\eqref{formulevol}  is established at a level of generality bigger than what is needed in Case (3) of Theorem~\ref{volestim}. To this end, fix a non--empty subcollection of sets $\mathfrak{O}'$ contained in the collection $\mathfrak{O}$ defined in~\eqref{defoobig}. While it is obvious that the condition $(\mathcal{H}_2)$  is verified in  Case (3) since the set $K$ then contains the origin in its interior (recall that the polynomial $P(\bm{x})$ is assumed to be homogeneous), Lemma~\ref{egpolorder} below makes it clear that $(\mathcal{H}_3)$ also holds under this same assumption on $K$. With this in mind, define, with a slight abuse of notations, 
\begin{equation}\label{def123}
\zeta_{P}\left(\mathfrak{O}', \,\cdot\,\right)\;=\; {\sum_{\Omega\in\mathfrak{O}'}} \;\zeta_{P}(\Omega, \,\cdot\,),
\end{equation} 
where the restricted distributions $\zeta_{P}(\Omega, \,\cdot\,)$ are defined in~\eqref{restriczeta}. \\ 

It should be clear that  Lemma~\ref{propzetamerogro} also applies to the distribution $\zeta_{P}\left(\mathfrak{O}', \,\cdot\,\right)$ which therefore extends to a meromorphic distribution defined over the complex plane. Let then $r_P\left(\mathfrak{O}'\right)\ge 0$ be its smallest real pole (assuming it exists) and let $m_P\left(\mathfrak{O}'\right)\ge 1$  denote the corresponding order. \\

If $\psi$ is a nonnegative map in $\mathcal{C}_c^\infty(\R^n)$ which does not vanish at the origin, it is not hard to see that, given $\Omega\in\mathfrak{O}'$,  the pair $\left(r_P(\Omega, \psi), m_P(\Omega, \psi)\right)$ introduced in~\eqref{deforderpole} is well--defined under the assumption of the homogeneity of the polynomial $P(\bm{x})$, and that the same holds for the pair $\left(r_P( \psi), m_P(\psi)\right)$  also introduced in~\eqref{deforderpole} (this is also established in much more generality in Proposition~\ref{proprlctgene} in Section~\ref{locglobrlct} below). An immediate consequence of the following lemma  is that the assumption $(\mathcal{H}_3)$ is ne\-ces\-sarily met when the set $K$ contains the origin in its interior (and thus under the assumptions of Case (3)).

\begin{lem}\label{egpolorder}
Let $\psi\ge 0$ be in $\mathcal{C}_c^\infty(\R^n)$ such that $\psi(\bm{0})>0$. Then, under the assumption of the homogeneity of the polynomial $P(\bm{x})$, it holds on the one hand that
\begin{equation}\label{dernia}
r_P(\Omega)\;=\; r_P(\Omega, \psi) \qquad\qquad \textrm{and} \qquad \qquad m_P(\Omega)\;=\; m_P(\Omega, \psi),
\end{equation}
and on the other that 
\begin{equation}\label{derniabis}
r_P\;=\; r_P( \psi) \qquad\qquad \textrm{and} \qquad \qquad m_P\;=\; m_P(\psi).
\end{equation}
In these relations, the restricted pair $\left(r_P(\Omega), m_P(\Omega)\right)$ and the global pair $\left(r_P, m_P\right)$ are  both defined in~\eqref{romegmomeg}.
\end{lem}

To justify the claim that $(\mathcal{H}_3)$ is satisfied when the set $K$ contains the origin in its interior, it is enough to define, in the notation of the support restriction condition~\eqref{suppproper}, the set $C$ to be any compact neighbourhood of the origin contained in $K$ and the set $U$ to be the full space $\R^n$. In view of the above lemma, it then suffices to set the pair $\left(r_P(K), m_P(K)\right)$ as $\left(r_P, m_P\right)$ independently of the choice of such a set $K$.

\begin{proof} 
It suffices to establish the relations~\eqref{dernia} as those in~\eqref{derniabis} follow from this first case upon reproducing  \emph{verbatim} the proof after removing the notational dependence on $\Omega$ of the various quantities  under consideration. \\

By definition,  one has that $r_P(\Omega, \psi)\ge r_P(\Omega)$ and, in the case where equality holds, that $m_P(\Omega, \psi)\le m_P(\Omega)$. \\

As for the reverse inequalities, let  $\varphi$ be in $\mathcal{C}_c^\infty(\R^n)$ such that the meromorphic map 
\begin{equation}\label{mapzetacettepreuve}
s\in\C \mapsto \left\langle \zeta_{P}\left(\Omega, s\right), \varphi\right\rangle
\end{equation} 
achieves its smallest real pole at $r_P(\Omega)$ with an order $m_P(\Omega)$. Clearly, when $\sigma<r(\Omega)$ and when $\xi$ is a map in $\mathcal{C}_c^\infty(\R^n)$ which has the same support as $\varphi$ and  which bounds from above $\left|\varphi\right|$, it holds that 
\begin{equation}\label{majpoleordre}
\left|\left\langle \zeta_{P}\left(\Omega, \sigma\right), \varphi\right\rangle\right|\;\le\; \left\langle \zeta_{P}\left(\Omega, \sigma\right), \xi\right\rangle .
\end{equation} 
 
Letting $\sigma$ tend to $r_P(\Omega)^-$, the meromorphic distribution on the right--hand side appears to have a smallest pole at most equal to $r_P(\Omega)$, and therefore equal to this quantity from its definition. Furthermore, one then also obtains that its order at this pole is at least  $m_P(\Omega)$, and therefore also equal to this integer from its definition. This shows that one may assume without loss of generality that $\varphi\ge 0$ and, even if it means renormalising it, that $0\le \varphi <1$.\\

Fix again $\sigma<r_P(\Omega)$ and let $\lambda>0$ be a parameter large enough so  that \mbox{$\lambda^{-1}\cdot\textrm{Supp}(\varphi) :=\left\{\bm{x}/\lambda\in\R^n :\bm{x}\in \textrm{Supp}(\varphi)\right\}$} is contained in a small enough neighbourhood of the origin where $\psi$ remains strictly positive. Upon setting $\varphi_{\lambda}~:\bm{y}\in\R^n\mapsto \varphi(\lambda\bm{y})$, the change of variables $\bm{y}=\bm{x}/\lambda$ in the integral defining the map~\eqref{mapzetacettepreuve} yields
\begin{align}\label{inhomgchgevar} 
 \left\langle \zeta_{P}\left(\Omega, \sigma\right), \varphi\right\rangle  \; &=\;  \lambda^{n-q\sigma}\cdot \left\langle \zeta_{P}\left(\Omega, \sigma\right), \varphi_{\lambda}\right\rangle, 
\end{align} 
where it is recalled that $q\ge 1$ is the degree of the polynomial $P(\bm{x})$. This shows that the meromorphic distribution $s\mapsto \left\langle \zeta_{P}\left(\Omega, \sigma\right), \varphi_{\lambda}\right\rangle$ is such that 
\begin{equation}\label{locrmpl}
\left(r_P(\Omega, \varphi_{\lambda}), m_P(\Omega, \varphi_{\lambda})\right)\;=\;\left(r_P(\Omega, \varphi), m_P(\Omega, \varphi)\right)\;=\; \left(r_P(\Omega), m_P(\Omega)\right).
\end{equation}

Also, 
\begin{align}\label{zetadecompphilamb}
 \left\langle \zeta_{P}\left(\Omega, \sigma\right), \psi\right\rangle  \; &=\; \left\langle \zeta_{P}\left(\Omega, \sigma\right), \psi\cdot \varphi_{\lambda}\right\rangle \; +\; \left\langle \zeta_{P}\left(\Omega, \sigma\right), \psi\cdot\left(1- \varphi_{\lambda}\right)\right\rangle,
\end{align}
where all the  terms remain positive for $\sigma<r_P(\Omega)$ from the relation~\eqref{restriczeta} defining the distribution $\zeta_P\left(\Omega, \;\cdot\;\right)$. Since the test function $\psi\cdot\varphi_{\lambda}$ is, up to strictly positive multiplicative constants, bounded from above and below by $\varphi_{\lambda}$  on its support (this follows from the assumption made on $\textrm{Supp}(\varphi_{\lambda})\;=\; \lambda^{-1}\cdot\textrm{Supp}(\varphi)$), the smallest pole and the corresponding order of the function $\left\langle \zeta_{P}\left(\Omega, \:\cdot\:\right), \psi\cdot \varphi_{\lambda}\right\rangle $ are also given by~\eqref{locrmpl}. Upon letting $\sigma$ tend to $r_P(\Omega)^-$, it is then easily deduced that the decomposition~\eqref{zetadecompphilamb} implies the reverse inequalities $r_P(\Omega)\ge r_P(\Omega, \psi)$ and $m_P(\Omega)\le m_P(\Omega, \psi)$. This establishes the lemma. 
\end{proof}

The following notational reductions shall be convenient to either adopt or keep in mind~:
\begin{itemize}
\item when $\mathfrak{O}' = \mathfrak{O} $, the definition of $\zeta_{P}\left(\mathfrak{O}, \,\cdot\,\right)$ in~\eqref{def123} coincide with that of $\zeta_P$ in~\eqref{zetaP};

\item when $\mathfrak{O}' = \mathfrak{O}^{\pm}$
, set $\zeta_P^{\pm}=\zeta_{P}\left(\mathfrak{O}^{\pm}, \,\cdot\,\right)$;

\item when $\mathfrak{O}' = \mathfrak{O} $ and   when $\mathfrak{O}' = \mathfrak{O}^{\pm}$, let respectively
\begin{equation}\label{latedef}
\left(r_P, m_P\right) = \left(r_P\left(\mathfrak{O}\right), m_P\left(\mathfrak{O}\right)\right)\quad \textrm{and}\quad  \left(r_P^{\pm}, m_P^{\pm}\right)= \left(r_P\left(\mathfrak{O}^{\pm}\right), m_P\left(\mathfrak{O}^{\pm}\right)\right). 
\end{equation}
\end{itemize}

The precise asymptotic behavior of the quantity $\V_n\left(\underline{\mathcal{S}}_{P}(K, b(T))\right)$  in Case (3)  is a corollary of the following general theorem. To state it, note that if $\mathfrak{g}\in SL_n(\R)$, then the map $\Omega\in\mathfrak{O}'\mapsto \mathfrak{g}^{-1}\cdot\Omega\in \mathfrak{g}^{-1}\cdot \mathfrak{O}'$ induces a bijection between connected sets where the polynomial $P(\bm{x})$ keeps a constant sign and connected sets where the polynomial $\left(P\circ\mathfrak{g}\right)(\bm{x})$ satisfies the same property. 

\begin{thm}[Volume estimate in Case (3) under a unimodular distorsion]\label{mainthmvolomega}
Assume that the set $K$ is star-shaped with respect to the origin, that it meets the condition $(\mathcal{H}_1)$ and that it contains the origin in its interior. Let furthermore  the limit  condition~\eqref{limiassumpbis}, where the quantity $c(T)$ is defined, hold. Let $\mathfrak{g}\in SL_n(\R)$ and let $\mathfrak{O}'$ be a nonempty subcollection of sets in $\mathfrak{O}$ such that  there exists $\Omega\in\mathfrak{O}'$ for which $\left(\mathfrak{g}^{-1}\cdot \Omega\right)\cap K\neq\emptyset$. \\

Then, there exist  strictly positive constants $\delta_P\left(\mathfrak{O}', K\right)$ and $ \gamma_{P, \mathfrak{g}}\left(\mathfrak{O}', K\right)$ such that, as the parameter $T$ tends to infinity, 
\begin{align}
&\V_n\left(\bigcup_{\Omega\in\mathfrak{O}'}\left\{\bm{x}\in T \cdot K\; : \; 0< \left(P\circ\mathfrak{g}\right)_{\left(\mathfrak{g}^{-1}\cdot\Omega\right)}(\bm{x})<b(T)\right\}\right)\;=\; \qquad\qquad\qquad\qquad\qquad\qquad\qquad\qquad\qquad\nonumber\\
&\left(\gamma_{P, \mathfrak{g}}\!\left(\mathfrak{O}', K\right)+O\!\left(\left(f_{P}\!\left(\mathfrak{O}', c(T)^{-1}\right)\right)^{-\delta_P\left(\mathfrak{O}', K\right)}\right)\right)\cdot T^{n}\cdot c(T)^{r_P\left(\mathfrak{O}'\right)}\cdot\left|\log c(T)\right|^{m_P\left(\mathfrak{O}'\right)-1}. \label{volestimformsign0}
\end{align}
In this relation, the pair $\left(r_P\left(\mathfrak{O}'\right), m_P\left(\mathfrak{O}'\right)\right)$ is well--defined and is such that $r_P\left(\mathfrak{O}'\right)>0$; the restricted polynomials $\left(P\circ\mathfrak{g}\right)_{\left(\mathfrak{g}^{-1}\cdot\Omega\right)}(\bm{x})$ are defined as in~\eqref{defPomega} and, given $c\in (0,1)$, 
$$
f_{P}\left(\mathfrak{O}', c\right)\;=\;
\begin{cases}
c & \textrm{if for all } \Omega\in \mathfrak{O}', \textrm{it holds that } m_P\left(\Omega\right)=1 \textrm{ whenever } \\
&\qquad r_P\left(\Omega\right)= r_P\left(\mathfrak{O}'\right);\\
\log(c) & \textrm{otherwise}.
\end{cases}
$$
Finally, given any compact set $B\subset\textrm{SL}_n(\R)$, the implicit constant as well as the leading coefficient in the volume estimate~\eqref{volestimformsign0} can be chosen uniformly over all polynomials in the family $\left\{\left(P\circ\mathfrak{g}\right)(\bm{x})\right\}_{\mathfrak{g}\in B}$.
\end{thm}

Recall that  from Lemma~\ref{egpolorder}, $(\mathcal{H}_3)$ is verified  under the assumptions of Case~(3). The above statement can then be slightly refined in the following sense~: given $\mathfrak{g}\in SL_n(\R)$, one can take $f_{P}\left(\mathfrak{O}', c\right)=c$ in the asymptotic relation~\eqref{volestimformsign0} whenever for any nonnegative map $\psi$ in  $\mathcal{C}_{c}^{\infty}(\R^n)$ satisfying the support condition~\eqref{suppproper} present in $(\mathcal{H}_3)$, the non--leading coefficients of the monic univariate polynomials  $R_{\psi}^{\Omega}(x)=R_{\mathfrak{g}, \psi}^{\Omega}(x)\in\R[x]$ appearing in the statement of Lemma~\ref{lemloir}  vanish for all $\Omega$ in the collection $\mathfrak{O}'$ (there is one such univariate polynomial associated to each restricted polynomial $\left(P\circ\mathfrak{g}\right)_{\left(\mathfrak{g}^{-1}\cdot\Omega\right)}(\bm{x})$ as $\Omega$ varies in $\mathfrak{O}'$). In other words, this is requiring that $R_{\mathfrak{g}, \psi}^{\Omega}(x)=x^{m-1}$ for  all  such sets $\Omega$  and all such maps $\psi\ge 0$. See the proof of Lemma~\ref{estimerreurunif} below for a justification of this claim. \\

Taking  $B\subset\textrm{SL}_n(\R)$ to be  reduced to the identity element and $\mathfrak{O}'=\mathfrak{O}$, Theorem~\ref{mainthmvolomega} admits the following immediate corollary. As previously, upon specialising this corollary to the case where $P(\bm{x})$ is the homogeneous polynomial $\left\|\bm{F}(\bm{x})\right\|^2$ of degree $q=2d$  and  upon noticing that the poles and orders of the distribution $\zeta_{\bm{F}}$ defined in~\eqref{distrizetaR} then coincide with those of the distribution $\zeta_P(\,\cdot\, / 2)$, it implies (in a stronger form) the asymptotic expansion~\eqref{formulevol} part of Theorem~\ref{volestim}.

\begin{coro}[Volume estimate in Case (3)]\label{proptmh1.1}
Assume that the set $K$ meets the condition $(\mathcal{H}_1)$, that it is star-shaped with respect to the origin  and that it contains the origin in its interior. Let furthermore  the limit  condition~\eqref{limiassumpbis}, where the quantity $c(T)$ is defined, hold.\\

Then, there exists constants $\delta_P(K), \gamma_P(K)>0$ such that 
$$\V_n\left(\underline{\mathcal{S}}_{P}(K, b(T))\right)\;=\;  \left(\gamma_P(K)+O\left(\left(f_P\left(c(T)^{-1}\right)\right)^{-\delta_P(K)}\right)\right)\cdot T^{n}\cdot c(T)^{r_P}\cdot\left|\log c(T)\right|^{m_P-1}$$ as the parameter $T$ tends to infinity. Here, $r_P>0$ and, given $c\in (0,1)$,
$$
f_{P}\left(c\right)\;=\;
\begin{cases}
c & \textrm{if for all } \Omega\in \mathfrak{O}, \textrm{it holds that } m_P\left(\Omega\right)=1 \textrm{ whenever } \\
&\qquad r_P\left(\Omega\right)= r_P;\\
\log(c) & \textrm{otherwise}.
\end{cases}
$$
\end{coro}

Before beginning the proof of Theorem~\ref{mainthmvolomega},  another corollary is derived from this statement upon specialising it to the case where $\mathfrak{O}'$ is either of the collections of sets $\mathfrak{O}^{\pm}$ and where $K=B_n$ is the closed unit Euclidean ball. This corollary plays a crucial role in the proofs of the metric counting results established in Chapter~\ref{ranunidisfamhomfor}.

\begin{coro}\label{zetapm}
Let $B$ be a compact subset of $\textrm{SL}_n(\R)$ and let  $\mathfrak{g}\in B$. Assume that the limit condition~\eqref{limiassumpbis}, where the quantity $c(T)$ is defined, holds and assume that the pair $\left(r_P^{\pm}, m_P^{\pm}\right)$ is well--defined.\\

Then, there exist constants $\delta_{P}^{\pm}, \gamma_{P, \mathfrak{g}}^{\pm}> 0$ such that, as the parameter $T$ tends to infinity, 
\begin{align}
&\V_n\left(\left\{\bm{x}\in B_n(T) : 0<\pm \left(P\circ\mathfrak{g}\right)(\bm{x})<b(T)\right\}\right)\qquad \qquad \qquad \qquad \qquad \qquad \qquad \qquad \qquad  \nonumber\\ 
&\qquad\qquad\qquad\qquad=\; \left(\gamma_{P, \mathfrak{g}}^{\pm}+O\left(\left(f_{P}^{\pm}\left(c(T)^{-1}\right)\right)^{-\delta_{P}^{\pm}}\right)\right)\cdot T^n\cdot c(T)^{r_P^{\pm}}\cdot
\left|\log c(T)\right|^{m_P^{\pm}-1}.\label{volestimformsign}
\end{align}
Here, given $c\in (0,1)$,
$$
f_{P}^{\pm}\left(c\right)\;=\;
\begin{cases}
c & \textrm{if for all } \Omega\in \mathfrak{O}^{\pm}, \textrm{ it holds that } m_P\left(\Omega\right)=1 \textrm{ whenever } r_P\left(\Omega\right)= r_P^{\pm};\\
\log(c) & \textrm{otherwise}.
\end{cases}
$$Furthermore, the leading coefficient $\gamma_{P, \mathfrak{g}}^{\pm}$ and the implicit constant  in the relation~\eqref{volestimformsign} are bounded above uniformly over all polynomials in the family $\left\{\left(P\circ\mathfrak{g}\right)(\bm{x})\right\}_{\mathfrak{g}\in B}$.
\end{coro}

The first step in the proof of Theorems~\ref{mainthmvolomega} is to give an explicit expression for the pair $\left(r_P\left(\mathfrak{O}'\right), m_P\left(\mathfrak{O}'\right)\right)$ and to relate it to the corresponding pair $\left(r_{\left(P\circ\mathfrak{g}\right)}\left(\mathfrak{g}^{-1}\cdot\mathfrak{O}'\right), m_{\left(P\circ\mathfrak{g}\right)}\left(\mathfrak{g}^{-1}\cdot\mathfrak{O}'\right)\right)$ associated to the polynomial $\left(P\circ\mathfrak{g}\right)(\bm{x})$. This is achieved in the following lemma.

\begin{lem}\label{lempolemulti}
Let $\mathfrak{g}$ be in $\textrm{SL}_n(\R)$ and let  $\mathfrak{O}'$ be a nonempty subcollection of $\mathfrak{O}$. \\  

Then, the pair $\left(r_{\left(P\circ\mathfrak{g}\right)}\left(\mathfrak{g}^{-1}\cdot\mathfrak{O}'\right), m_{\left(P\circ\mathfrak{g}\right)}\left(\mathfrak{g}^{-1}\cdot\mathfrak{O}'\right)\right)$ is well-defined if, and only if, so is the pair $\left(r_P\left(\mathfrak{O}'\right), m_P\left(\mathfrak{O}'\right)\right)$, in which case 
\begin{equation}\label{ordremultdecomposum1}
r_{\left(P\circ\mathfrak{g}\right)}\left(\mathfrak{g}^{-1}\cdot\mathfrak{O}'\right)\;=\; r_P\left(\mathfrak{O}'\right)\;=\; \min_{\Omega\in\mathfrak{O}'} r_P(\Omega)
\end{equation} 
and
\begin{equation}\label{ordremultdecomposum2}
m_{\left(P\circ\mathfrak{g}\right)}\left(\mathfrak{g}^{-1}\cdot\mathfrak{O}'\right)\;=\; m_P\left(\mathfrak{O}'\right)\;=\; \max_{\underset{r_P(\Omega) = r_P\left(\mathfrak{O}'\right)}{\Omega\in\mathfrak{O}'}} m_P(\Omega).
\end{equation} 
\end{lem}

\begin{proof}
\sloppy Given a map $\psi$ in $\mathcal{C}_c^\infty(\R^n)$ and a set $\Omega$ in $\mathfrak{O}'$,  an obvious change a variables shows  that $\left\langle  \zeta_{(P\circ\mathfrak{g})}\left(\mathfrak{g}^{-1}\cdot\left(\Omega\right), \,\cdot\,\right), \psi \right\rangle = \left\langle  \zeta_{P}(\Omega, \,\cdot\,), \psi\left(\mathfrak{g}^{-1}\,\cdot\,\right) \right\rangle$. As a consequence, $\left(r_{(P\circ\mathfrak{g})}\left(\mathfrak{g}^{-1}\cdot\left(\Omega\right)\right), m_{(P\circ\mathfrak{g})}\left(\mathfrak{g}^{-1}\cdot\left(\Omega\right)\right)\right) = \left(r_{P}(\Omega), m_{P}(\Omega)\right)$ whenever either of these pairs exists. It is therefore enough to establish the second equations, in~\eqref{ordremultdecomposum1} and in~\eqref{ordremultdecomposum2} respectively, to prove the lemma.\\

To this end, note that the definition~\eqref{def123} of the distribution $\zeta_{P}\left(\mathfrak{O}', \,\cdot\,\right)$ immediately implies the inequalities 
\begin{equation}\label{ineqordremultp1}
r_P\left(\mathfrak{O}'\right)\;\ge\; \min_{\Omega\in\mathfrak{O}'} r_P(\Omega)\qquad \textrm{ and }\qquad m_P\left(\mathfrak{O}'\right)\;\le\; \max_{\underset{r_P(\Omega) = r_P\left(\mathfrak{O}'\right)}{\Omega\in\mathfrak{O}'}} m_P(\Omega). 
\end{equation}
The reverse inequalities are obtained upon noticing that the smallest poles of the distributions $\zeta_{P}\left(\Omega, \,\cdot\,\right)$ and their respective orders cannot cancel out as $\Omega$ varies in $\mathfrak{O}'$. To see this, fix $\Omega\in\mathfrak{O}'$ and consider $\psi_{\Omega}$ in $\mathcal{C}_c^\infty(\R^n)$ for which the meromorphic map $s\in\C \mapsto \left\langle \zeta_{P}\left(\Omega, s\right), \psi_{\Omega}\right\rangle$ has a pole at $r_P(\Omega)$ of order $m_P(\Omega)\ge 1$. Let $\varphi\ge 0$ be a map in $\mathcal{C}_c^\infty(\R^n)$ bounding from above $\left|\psi_{\Omega}\right|$. Assume that $s=\sigma<r_P(\Omega)$ is real. From the Triangle Inequality, 
\begin{align}
&\left(r_P(\Omega)-\sigma\right)^{m_P(\Omega)}\cdot \left|\left\langle \zeta_{P}\left(\Omega, \sigma\right), \psi_{\Omega}\right\rangle\right| \; \le\; \left(r_P(\Omega)-\sigma\right)^{m_P(\Omega)}\cdot \left\langle \zeta_{P}\left(\Omega, \sigma\right), \varphi\right\rangle. \label{ineqordremultp}
\end{align}
Upon taking the limit $\sigma\rightarrow r_P(\Omega)^-$, it follows that $ r_P(\Omega)\ge r_P(\Omega, \varphi)$ and that, in the case where equality holds, $m_P(\Omega)\le m_P(\Omega, \varphi)$. Since, by the definition of the pole $ r_P(\Omega)$ and of its order $m_P(\Omega)$, it holds that $  r_P(\Omega, \varphi)\ge r_P(\Omega)$ and that,  in the case of equality, $m_P(\Omega, \varphi)\le m_P(\Omega)$, one obtains that $\left(r_P(\Omega, \varphi), m_P(\Omega, \varphi)\right)= \left(r_P(\Omega), m_P(\Omega)\right)$. (The argument is similar to the one used in~\eqref{majpoleordre}.)\\

Also, the coefficient of $(r_P(\Omega)-s)^{-m_P(\Omega)}$ in the Laurent expansion of the map $s\in\C \mapsto \left\langle \zeta_{P}\left(\Omega, s\right), \varphi\right\rangle$ at $r_P(\Omega)$ is positive whenever this map is well-defined (since it is then clearly positive when $s=\sigma<r_P(\Omega)$ is real). Choose $\varphi\ge 0$ bounding from above the absolute values of all smooth functions $\psi_{\Omega}$ when $\Omega$ varies over all those sets achieving the minimum in~\eqref{ordremultdecomposum1}. Then, the sum over these sets $\Omega$ of the  Laurent coefficients of $(r_P(\Omega)-s)^{-m_P(\Omega)}$ of the meromorphic maps $\left\langle \zeta_{P}\left(\Omega, \,\cdot\,\right), \varphi\right\rangle$ does not va\-nish, which shows that  the converse inequalities in~\eqref{ineqordremultp1} also hold. This completes the proof.
\end{proof}

 Lemma~\ref{lempolemulti} reduces the calculations of the pairs  $\left(r_P\left(\mathfrak{O}'\right), m_P\left(\mathfrak{O}'\right)\right)$  and $\left(r_{\left(P\circ\mathfrak{g}\right)}\left(\mathfrak{g}^{-1}\cdot\mathfrak{O}'\right), m_{\left(P\circ\mathfrak{g}\right)}\left(\mathfrak{g}^{-1}\cdot\mathfrak{O}'\right)\right)$ to that of a finite number or pairs $\left( r_P(\Omega),  m_P(\Omega)\right)$ when $\Omega\in \mathfrak{O}'$. The determination of the latter quantities can be achieved thanks to the above Lemma~\ref{egpolorder}.\\

\begin{proof}[Proof of Theorem~\ref{mainthmvolomega}] The proof of the claim that  the pair $\left(r_P\left(\mathfrak{O}'\right), m_P\left(\mathfrak{O}'\right)\right)$ is well-defined and that its first component is strictly positive is an immediate consequence of the above Lemma~\ref{lempolemulti} and of Proposition~\ref{proprlctgene} from Section~\ref{locglobrlct} below. The goal is here to establish the remaining claims in the statement of Theorem~\ref{mainthmvolomega}.\\

In view of the result to establish, even if it means discarding elements from the collection $\mathfrak{O}'$, it may be assumed without loss of generality that  for each $\Omega\in \mathfrak{O}'$, it holds that $\left(\mathfrak{g}^{-1}\cdot \Omega\right) \cap K\;\neq\; \emptyset$. It is then enough to prove Theorem~\ref{mainthmvolomega} in the case that $\mathfrak{O}'$ contains a unique element $\Omega$ for which the pair $\left( r_P\left(\Omega\right), m_P\left(\Omega\right)\right)$ is well--defined. Indeed, assume that  for some $\delta_P\left(\Omega, K\right)>0$,
\begin{align}
&\V_n\left(\left\{\bm{x}\in T\cdot K \; :\; 0< \left(P\circ\mathfrak{g}\right)_{\left(\mathfrak{g}^{-1}\cdot \Omega\right)}(\bm{x})<b(T)\right\}\right)\;=\;\qquad\qquad\qquad\qquad\qquad\qquad\qquad\qquad\qquad\nonumber\\
&  \left(\gamma_{P, \mathfrak{g}}\left(\Omega, K\right)+O\left(\left(f_P\left(\Omega, c(T)^{-1}\right)\right)^{-\delta_P\left(\Omega, K\right)}\right)\right)\cdot T^{n}\cdot c(T)^{r_P\left(\Omega\right)}\cdot\left|\log c(T)\right|^{m_P\left(\Omega\right)-1} \label{volreunelmt}
\end{align}
with the required uniformity claims, where $\gamma_{P, \mathfrak{g}}\left(\Omega, K\right)>0$ and where for $c\in (0,1)$,
\begin{equation}\label{deffpomegt}
f_{P}\left(\Omega, c\right)\;=\;
\begin{cases}
c & \textrm{if } m_P\left(\Omega\right)=1;\\
\log(c) & \textrm{otherwise}.
\end{cases}
\end{equation}
From the disjointness of the supports of the restricted polynomials $\left(P\circ\mathfrak{g}\right)_{\left(\mathfrak{g}^{-1}\cdot\Omega\right)}(\bm{x})$, 
\begin{align*}
&\V_n\left(\bigcup_{\Omega\in\mathfrak{O}'}\left\{\bm{x}\in T\cdot K\;:\; 0< \left(P\circ\mathfrak{g}\right)_{\left(\mathfrak{g}^{-1}\cdot\Omega\right)}(\bm{x})<b(T)\right\}\right) \qquad\qquad\qquad\qquad\qquad\qquad\qquad\qquad\qquad\qquad \\
& \quad\qquad\qquad \qquad\qquad\qquad\;\;\;=\; \sum_{\Omega\in\mathfrak{O}'} \V_n\left(\left\{\bm{x}\in T\cdot K \;:\; 0< \left(P\circ\mathfrak{g}\right)_{\left(\mathfrak{g}^{-1}\cdot\Omega\right)}(\bm{x})<b(T)\right\}\right).
\end{align*}
Combining the above two relations, it follows from Lemma~\ref{lempolemulti} that Theorem~\ref{mainthmvolomega} is obtained upon setting $$\gamma_{P, \mathfrak{g}}\left(\mathfrak{O}', K\right)\;=\; \sum_{\stackrel{\Omega\in\mathfrak{O}'}{\left(r_P(\Omega), m_P(\Omega)\right)=\left(r_P\left(\mathfrak{O}'\right), m_P\left(\mathfrak{O}'\right)\right)}} \gamma_{P, \mathfrak{g}}(\Omega, K),$$ $$ \delta_P\left(\mathfrak{O}', K\right)\;=\;\min_{\stackrel{\Omega\in\mathfrak{O}'}{r_P\left(\Omega\right)=r_P\left(\mathfrak{O}'\right)}}\delta_P\left(\Omega, K\right)$$ and $$f_{P}\left(\Omega, c\right)\;=\;\min_{\stackrel{\Omega\in\mathfrak{O}'}{\delta_P\left(\Omega, K\right)=\delta_P\left(\mathfrak{O}', K\right)}}\;f_{P}\left(\Omega, c\right).$$ \\

In order to establish the relation~\eqref{volreunelmt}, consider two sequences $\left(\psi_k^-\right)_{k\ge 1}$ and $\left(\psi_k^+\right)_{k\ge 1}$ of elements in  $\mathcal{C}_c^\infty(\R^n)$ satisfying the following two properties~:
\begin{itemize}
\item[(i)] the sequence $\left(\psi_k^-\right)_{k\ge 1}$ is pointwise increasing and the sequence $\left(\psi_k^+\right)_{k\ge 1}$ pointwise decreasing;
\item[(ii)]  for all $k\ge 1$, it holds that $\chi_{K(1-1/k)}\;\le\; \psi_k^- \;\le\; \chi_{K}\;\le\; \psi_k^+\le \chi_{K(1+1/k)}$, where, given $\lambda>0$, $K(\lambda)$ denotes the dilate of the set $K$ by the factor $\lambda$; in other words, $K(\lambda)=\left\{\lambda\cdot\bm{x}\;:\; \bm{x}\in K\right\}$.
\end{itemize}

In the above, for any given integer $k\ge 1$, the inequalities $\chi_{K(1-1/k)}\le \chi_{K}\le \chi_{K(1+1/k)}$ hold under the assumption that the set $K$ is star-shaped with respect to the origin. When $\mathfrak{g}\in SL_n(\R)$, $k\ge 1$ and $\bm{x}\in\R^n$, set furthermore
\begin{equation}\label{defpsiplgk}
\psi_{\mathfrak{g}, k}^{\pm}\left(\bm{x}\right)=\psi_{k}^{\pm}\left(\mathfrak{g}^{-1}\cdot\bm{x}\right),\qquad\textrm{where}\qquad \textrm{Supp }\psi_{\mathfrak{g}, k}^{\pm}\subset K\left(2\cdot \kappa_{\mathfrak{g}}\right).
\end{equation}
Here and throughout, $\kappa_{\mathfrak{g}}>0$ denotes the operator norm of the linear map induced by $\mathfrak{g}$ (with respect to the Euclidean norm). In the notations introduced in the defining relations~\eqref{defmuvolweighted} and in Lemma~\ref{lemloir}, given $c\in (0,1)$ and $k\ge 1$, set first
\begin{align}
\nu^\pm_{P, \mathfrak{g}}\left(\Omega, k, c\right)\; &=\;\mu_{\left(P\circ\mathfrak{g}\right)}\left(\mathfrak{g}^{-1}\cdot\Omega, \psi_{k}^\pm , c\right)\;=\; \int_{\mathfrak{g}^{-1}\cdot\Omega} \chi_{\left\{\left(P\circ\mathfrak{g}\right)(\bm{x})\le c\right\}}\cdot\psi_{k}^\pm(\bm{x})\cdot\textrm{d}\bm{x} \nonumber \\
&=\; \mu_P\left(\Omega, \, \psi_{\mathfrak{g}, k}^\pm\, , c\right)\;=\; \int_{\R^n} \chi_{\left\{P_{\Omega}(\bm{x})\le c\right\}}\cdot\psi_{\mathfrak{g},k}^\pm(\bm{x})\cdot\textrm{d}\bm{x},\label{defnupgomegkT}
\end{align}  
then
$$\nu_{P, \mathfrak{g}}\left(\Omega,  K, c\right)\;=\; \int_{K\cap\left(\mathfrak{g}^{-1}\cdot\Omega\right)} \chi_{\left\{\left(P\circ\mathfrak{g}\right)(\bm{x})\le c\right\}}\cdot\textrm{d}\bm{x}$$ and finally $$\Theta_{P, \mathfrak{g}}^\pm\left(\Omega, k\right)\;=\; \Theta_{\left(P\circ\mathfrak{g}\right)}\left(\mathfrak{g}^{-1}\cdot\Omega, \psi^\pm_k\right)\;=\; \Theta_{P}\left(\Omega,  \psi_{\mathfrak{g}, k}^\pm\right).$$ Thus, by an elementary change of variables relying on the homogeneity of the polynomial $P(\bm{x})$, the volume appearing in~\eqref{volreunelmt} can be expressed as
\begin{equation}\label{volnnupgomegT}
\V_n\left(\left\{\bm{x}\in T\cdot K\; :\; 0< \left(P\circ\mathfrak{g}\right)_{\left(\mathfrak{g}^{-1}\cdot\Omega\right)}(\bm{x})<b(T)\right\}\right)\;=\; T^n\cdot \nu_{P, \mathfrak{g}}\left(\Omega, K, c(T)\right)
\end{equation}
with the quantity $c(T)$ defined in~\eqref{limiassumpbis}. \\

Clearly, given $k\ge 1$, the above  assumptions (i) and (ii) imply that 
\begin{align}
\nu^-_{P, \mathfrak{g}}\left(\Omega, k, c\right)\;\le\; \nu^-_{P, \mathfrak{g}}\left(\Omega, k+1, c\right)\; &\le\; \nu_{P, \mathfrak{g}}\left(\Omega, K, c\right)\nonumber \\
&\le\; \nu^+_{P, \mathfrak{g}}\left(\Omega, k+1, c\right)\;\le\;  \nu^+_{P, \mathfrak{g}}\left(\Omega, k, c\right).\label{inegnu}
\end{align}

Since from  Lemma~\ref{lemloir} and Lemma~\ref{lempolemulti} (applied to the case when $\mathfrak{O}'=\left\{\Omega\right\}$), it holds that
\begin{equation}\label{limnutheta}
\lim_{c\rightarrow 0^{+}} \left(\frac{\nu_{P, \mathfrak{g}}^{\pm}\left(\Omega, k, c\right)}{c^{r_P(\Omega)}\cdot\left|\log c\right|^{m_P(\Omega)-1}}\right)\;=\; \Theta_{P, \mathfrak{g}}^\pm\left(\Omega, k\right),
\end{equation}
one deduces from the inequalities~\eqref{inegnu} that the two sequences $\left(\Theta_{P, \mathfrak{g}}^\pm\left(\Omega, k\right)\right)_{k\ge 1}$ are both monotonic and bounded. They therefore  converge to respective values $\Theta_{P, \mathfrak{g}}^\pm\left(\Omega, K\right)$ which satisfy the following relations for any $k\ge 1$~: 
\begin{align}
0\;<\; \Theta_{P, \mathfrak{g}}^-\left(\Omega, 1\right)\;\le\;  \Theta_{P, \mathfrak{g}}^-\left(\Omega, k\right) &\le\; \Theta_{P, \mathfrak{g}}^-\left(\Omega, K\right) \nonumber\\
&\le\;  \Theta_{P, \mathfrak{g}}^+\left(\Omega, K\right)  \;\le\; \Theta_{P, \mathfrak{g}}^+\left(\Omega, k\right)   \;\le\; \Theta_{P, \mathfrak{g}}^+\left(\Omega, 1\right). \label{ineaqcst}
\end{align} 
 
\begin{lem}\label{estimcoefdomomeggT}
Let $\mathfrak{g}\in SL_n(\R)$ with operator norm $\kappa_{\mathfrak{g}}>0$. With the above notations, there exists $\varepsilon>0$ such that when $c\in (0,1)$,
\begin{equation}\label{lastvpomkT}
\left| \nu_{P, \mathfrak{g}}^{+}\left(\Omega, k, c\right)-\nu_{P, \mathfrak{g}}^{-}\left(\Omega, k, c\right)\right|\;\ll\; \lambda_{P}\left(\Omega, \varepsilon, \mathfrak{g}\right) \cdot \frac{\nu_{P}^{+}\left(\Omega, 1, c\right)}{k},
\end{equation}
where 
\begin{equation}\label{defgammaprems}
\lambda_{P}\left(\Omega, \varepsilon, \mathfrak{g}\right)\;=\;\kappa_{\mathfrak{g}}^{n-q\cdot r_P(\Omega)}\cdot \max\left\{1, \kappa_{\mathfrak{g}}^{-\epsilon}\right\}\cdot \left(1+\left|\log \kappa_{\mathfrak{g}}\right|\right)^{m_P(\Omega)-1} 
\end{equation}
and where $\nu_{P}^{+}\left(\Omega, 1, c\right)$ is shorthand notation for $\nu_{P, \mathfrak{g}}^{+}\left(\Omega, 1, c\right)$ when $\mathfrak{g}$ is the identity. 
In par\-ti\-cu\-lar, one obtains that
\begin{equation}\label{ineqen1k}
\left| \Theta_{P, \mathfrak{g}}^+\left(\Omega, k\right) - \Theta_{P, \mathfrak{g}}^-\left(\Omega, k\right)\right|\; \ll\; \frac{\lambda_{P}\left(\Omega, \varepsilon, \mathfrak{g}\right)}{k}
\end{equation}
and also that
\begin{equation}\label{smevulethta}
\Theta_{P, \mathfrak{g}}^-\left(\Omega, K\right)\; =\; \Theta_{P, \mathfrak{g}}^+\left(\Omega, K\right) \; \; \;(=\;\Theta_{P, \mathfrak{g}}\left(\Omega, K\right), \textrm{say}).
\end{equation} 
In the above relations, the implicit constants depend only on the sets $\Omega$ and $K$ and on the polynomial $P(\bm{x})$.
\end{lem}

\begin{proof} From  the inequalities in (ii) above  and from the assumption of the homogeneity of the polynomial $P(\bm{x})$, it follows from an elementary change of variables that 
\begin{align*}
\left| \nu_{P, \mathfrak{g}}^{+}\left(\Omega, k, c\right)-\nu_{P, \mathfrak{g}}^{-}\left(\Omega, k, c\right)\right|\; &\le\; \int_{\left(K(1+1/k)\backslash K(1-1/k)\right)} \chi_{\left\{\left(P\circ\mathfrak{g}\right)_{\left(\mathfrak{g}^{-1}\cdot\Omega\right)}(\bm{x})\le c\right\}}\cdot\textrm{d}\bm{x}\\
&=\; \left(1+\frac{1}{k}\right)^{n}\cdot\int_{K} \chi_{\left\{\left(P\circ\mathfrak{g}\right)_{\left(\mathfrak{g}^{-1}\cdot\Omega\right)}(\bm{y})\le c\cdot (1+1/k)^{-q}\right\}}\cdot\textrm{d}\bm{y}  \\
&\quad\; - \left(1-\frac{1}{k}\right)^{n}\cdot\int_{K} \chi_{\left\{\left(P\circ\mathfrak{g}\right)_{\left(\mathfrak{g}^{-1}\cdot\Omega\right)}(\bm{y})\le c\cdot (1-1/k)^{-q}\right\}}\cdot\textrm{d}\bm{y}.
\end{align*}

Expanding the $n^{\textrm{th}}$ powers multiplying the integrals, rewriting the whole  expression as a polynomial in $1/k$ and isolating the constant term in this polynomial, one obtains that 
\begin{align*}
\left| \nu_{P, \mathfrak{g}}^{+}\left(\Omega, k, c\right)-\nu_{P, \mathfrak{g}}^{-}\left(\Omega, k, c\right)\right|\; &\ll\;  \int_{K} \chi_{\left\{c\cdot (1+1/k)^{-q}\;<\;\left(P\circ\mathfrak{g}\right)_{\left(\mathfrak{g}^{-1}\cdot\Omega\right)}(\bm{x})\;\le\; c\cdot (1-1/k)^{-q}\right\}}\cdot\textrm{d}\bm{x}\nonumber\\
&\qquad\qquad+ \frac{1}{k}\cdot \int_{K} \chi_{\left\{\left(P\circ\mathfrak{g}\right)_{\left(\mathfrak{g}^{-1}\cdot\Omega\right)}(\bm{x})\;\le\; 2^q\cdot c\right\}}\cdot\textrm{d}\bm{x}\\ 
&\ll  \int_{\mathfrak{g}\cdot K} \chi_{\left\{c\cdot (1+1/k)^{-q}\;<\;P_{\Omega}(\bm{y})\;\le\; c\cdot (1-1/k)^{-q}\right\}}\cdot\textrm{d}\bm{y}\nonumber\\
&\qquad\qquad + \frac{1}{k}\cdot \int_{\mathfrak{g}\cdot K} \chi_{\left\{P_{\Omega}(\bm{y})\;\le \;2^q\cdot c\right\}}\cdot\textrm{d}\bm{y}.
\end{align*}
Since $\mathfrak{g}\cdot K\subset K\left(\kappa_{\mathfrak{g}}\right)$, this yields
\begin{align*}
\left| \nu_{P, \mathfrak{g}}^{+}\left(\Omega, k, c\right)-\nu_{P, \mathfrak{g}}^{-}\left(\Omega, k, c\right)\right|\; &\ll\; \int_{ K\left(\kappa_{\mathfrak{g}}\right)} \chi_{ \left\{c\cdot (1+1/k)^{-q}\;<\;P_{\Omega}(\bm{y})\;\le\; c\cdot (1-1/k)^{-q}\right\}}\cdot\textrm{d}\bm{y}\nonumber\\
&\qquad\qquad+ \frac{1}{k}\cdot \int_{ K\left(\kappa_{\mathfrak{g}}\right)} \chi_{\left\{P_{\Omega}(\bm{y})\;\le \; 2^q\cdot c\right\}}\cdot\textrm{d}\bm{y}\\
&\ll \kappa_{\mathfrak{g}}^n\cdot \left( \int_{K} \chi_{\left\{c\cdot \left(\kappa_{\mathfrak{g}} (1+1/k)\right)^{-q}\;<\;P_{\Omega}(\bm{y})\;\le\; c\cdot \left(\kappa_{\mathfrak{g}}(1+1/k)\right)^{-q}\right\}}\cdot\textrm{d}\bm{y}\right.\nonumber\\
&\left.\qquad\qquad+ \frac{1}{k}\cdot \int_{K} \chi_{\left\{P_{\Omega}(\bm{y})\;\le\; c\cdot (\kappa_{\mathfrak{g}}/2)^{-q}\right\}}\cdot\textrm{d}\bm{y}\right).
\end{align*}
To bound from above this quantity, note that the inequality $\chi_{K}\le \psi_1^{+}$ implies one the one hand that 
\begin{align*}
&\int_{K} \chi_{\left\{c\cdot \left(\kappa_{\mathfrak{g}} (1+1/k)\right)^{-q}\;<\;P_{\Omega}(\bm{y})\;\le\; c\cdot \left(\kappa_{\mathfrak{g}}(1-1/k)\right)^{-q}\right\}}\textrm{d}\bm{y}\\
&\qquad\qquad\qquad\qquad\qquad\;\;\;\le \;  \int_{ \R^n} \chi_{\left\{c\cdot \left(\kappa_{\mathfrak{g}} (1+1/k)\right)^{-q}\;<\;P_{\Omega}(\bm{y})\;\le\; c\cdot \left(\kappa_{\mathfrak{g}}(1-1/k)\right)^{-q}\right\}}\cdot\psi_1^{+}(\bm{y})\cdot\textrm{d}\bm{y}
\end{align*} 
and, on the other, that $$\int_{K} \chi_{\left\{P_{\Omega}(\bm{y})\;\le\; c\cdot (\kappa_{\mathfrak{g}}/2)^{-q}\right\}}\cdot\textrm{d}\bm{y}\;\le\;  \int_{ \R^n} \chi_{\left\{P_{\Omega}(\bm{y})\;\le\; c\cdot (\kappa_{\mathfrak{g}}/2)^{-q}\right\}}\cdot\psi_1^{+}(\bm{y})\cdot\textrm{d}\bm{y}.$$  Consequently, it follows from the definition of $\nu_{P}^{+}\left(\Omega, 1, c\right)$ that
\begin{align*}
\left| \nu_{P, \mathfrak{g}}^{+}\left(\Omega, k, c\right)-\nu_{P, \mathfrak{g}}^{-}\left(\Omega, k, c\right)\right|\; &\ll \kappa_{\mathfrak{g}}^n \cdot\left(\nu_{P}^{+}\left(\Omega, 1, c\cdot \left(\kappa_{\mathfrak{g}} \left(1-\frac{1}{k}\right)\right)^{-q}\right)- \right. \nonumber\\
&\left.  \nu_{P}^{+}\left(\Omega, 1, c\cdot \left(\kappa_{\mathfrak{g}} \left(1+\frac{1}{k}\right)\right)^{-q}\right) +\frac{\nu_{P}^{+}\left(\Omega, 1, c\cdot (\kappa_{\mathfrak{g}}/2)^{-q}\right)}{k}\right).\nonumber
\end{align*}
Since the set $K$ contains the origin in its interior and that $\chi_{K}\le \psi_1^+$, one infers from Lemma~\ref{egpolorder} that $\left(r_P\left(\Omega, \psi\right), m_P\left(\Omega, \psi\right)  \right)=\left(r_P\left(\Omega\right), m_P\left(\Omega\right)\right)$. One then deduces from the above inequality and from the asymptotic expansion~\eqref{asympmu} that there exists $\varepsilon>0$ such that
\begin{align*}
\left| \nu_{P, \mathfrak{g}}^{+}\left(\Omega, k, c\right)-\nu_{P, \mathfrak{g}}^{-}\left(\Omega, k, c\right)\right|\; &\ll\; \kappa_{\mathfrak{g}}^{n-q\cdot r_P(\Omega)}\cdot \max\left\{1, \kappa_{\mathfrak{g}}^{-\varepsilon}\right\}\cdot \left(1+\left|\log \kappa_{\mathfrak{g}}\right|\right)^{m_P(\Omega)-1}\quad \times\\ 
&\qquad \nu_{P}^{+}\left(\Omega, 1, c\right)\left(\left(1-\frac{1}{k}\right)^{-q r_P(\Omega)}-  \left(1+\frac{1}{k}\right)^{ -q r_P(\Omega)}+\frac{1}{k}\right)\nonumber\\
&\underset{\eqref{defgammaprems}}{\ll}\; \lambda_{P}\left(\Omega, \varepsilon, \mathfrak{g}\right)\cdot \frac{\nu_{P}^{+}\left(\Omega, 1, c\right)}{k}\cdotp
\end{align*}
This establishes~\eqref{lastvpomkT}. Dividing through by $c^{r_P(\Omega)}\cdot \left|\log c\right|^{m_P(\Omega)-1}$ and calling on the limit relation~\eqref{limnutheta}, one obtains the inequality~\eqref{ineqen1k}. Letting $k$ tend to infinity then proves the identity~\eqref{smevulethta}. 
\end{proof}

Set from now on
\begin{equation}\label{cptheap}
\gamma_{P, \mathfrak{g}}(\Omega, K)= \Theta_{P, \mathfrak{g}}\left(\Omega, K\right),
\end{equation} 
where the quantity $\Theta_{P, \mathfrak{g}}\left(\Omega, K\right)$ is defined in~\eqref{smevulethta}. Furthermore, given integers $k, L\ge 1$, let
\begin{equation*}\label{cstpkikpm}
C_L(k)\;=\; \max_{0\le \left|\bm{l}\right|\le L}\; \max_{\bm{x}\in\R^n} \left|\frac{\partial^{\left|\bm{l}\right|}\Psi_k^{+}}{\partial\bm{x}^{\bm{l}}}(\bm{x})\right|.
\end{equation*}

\begin{lem}\label{estimerreurunif}
Let $\mathfrak{g}$ belong to a compact subset $B$ of $\textrm{SL}_n(\R)$. Then, there exist a constant $\delta_P\left(\Omega, K\right)>0$ and integers $M,N\ge 1$ such that for all $k\ge 1$
\begin{equation*}\label{estimerreurunif1}
\left|\frac{\nu_{P, \mathfrak{g}}\left(\Omega, K,  c\right)}{c^{r_P(\Omega)}\cdot \left|\log c\right|^{m_P(\Omega)-1}} -  \gamma_{P, \mathfrak{g}}\left(\Omega, K\right)\right| \ll  \frac{C_N(k)^M}{\left(f_{P}\left(\Omega, c^{-1}\right)\right)^{\delta_P\left(\Omega, K\right)}}+ \frac{1}{k},
\end{equation*}
where the map $f_{P}\left(\Omega, \,\cdot\,\right)$ is defined in~\eqref{deffpomegt}. Furthermore, the implicit constant depends only on the polynomial $P(\bm{x})$, on the set $\Omega$ and on the choice of the compact subset $B$.
\end{lem}

\begin{proof}
Let $k\ge 1$ and $c\in (0,1)$. Then, 
\begin{align*}
\left|\frac{\nu_{P, \mathfrak{g}}\left(\Omega, K, c\right)}{c^{r_P(\Omega)}\cdot \left|\log c\right|^{m_P(\Omega)-1}} -  \gamma_{P, \mathfrak{g}}\left(\Omega, K\right)\right|\;
&\underset{\eqref{cptheap}}{\le}\; \frac{\left|\nu_{P}\left(\Omega,  K, c\right)- \nu_{P, \mathfrak{g}}^{+}\left(\Omega, k, c\right)\right| }{c^{r_P(\Omega)}\cdot \left|\log c\right|^{m_P(\Omega)-1}}\\
&\quad  + \left| \frac{\nu_{P, \mathfrak{g}}^{+}\left(\Omega, k, c\right)}{c^{r_P(\Omega)}\cdot \left|\log c\right|^{m_P(\Omega)-1}} -  \Theta_{P, \mathfrak{g}}^{+}\left(\Omega, k\right)\right| \\
&\quad+\left| \Theta_{P, \mathfrak{g}}^+\left(\Omega, k\right) -  \Theta_{P, \mathfrak{g}}\left(\Omega, K\right)\right|.
\end{align*}
Each of the terms on the right--hand side of this inequality is estimated separately.\\

First, note that
\begin{align}\label{estim1unif}
\left| \Theta_{P, \mathfrak{g}}^+\left(\Omega, k\right) -  \Theta_{P, \mathfrak{g}}\left(\Omega, K\right)\right|\;\underset{\eqref{ineaqcst}}{\le}\; \left| \Theta_{P, \mathfrak{g}}^+\left(\Omega, k\right) -  \Theta_{P, \mathfrak{g}}^-\left(\Omega, k\right)\right|\; \underset{\eqref{ineqen1k}}{\ll}\; 
\frac{1}{k},
\end{align}
where the implicit constant can be chosen uniformly over the compact subset $B$ of $\textrm{SL}_n(\R)$. \\

Then, it follows from the relations~\eqref{inegnu} that 
\begin{align*}
\left|\nu_{P, \mathfrak{g}}\left(\Omega, K,  c\right)- \nu_{P, \mathfrak{g}}^{+}\left(\Omega, k, c\right)\right| \;\le\; \left|\nu_{P, \mathfrak{g}}^{+}\left(\Omega, k, c\right)- \nu_{P, \mathfrak{g}}^{-}\left(\Omega, k, c\right)\right| 
\end{align*}
in such a way that
\begin{align}
\frac{\left|\nu_{P, \mathfrak{g}}\left(\Omega, K,  c\right)- \nu_{P, \mathfrak{g}}^{+}\left(\Omega, k, c\right)\right|}{c^{r_P(\Omega)}\cdot \left|\log c\right|^{m_P(\Omega)-1}}\;&\underset{\eqref{lastvpomkT}}{\ll}\; \frac{\nu_{P}^{+}\left(\Omega, 1, c\right)}{c^{r_P(\Omega)}\cdot \left|\log c\right|^{m_P(\Omega)-1}}  \cdot \frac{1}{k}\nonumber \\
&\ll \; 
\frac{1}{k}\cdotp \label{estim2unif}
\end{align}
This last inequality is a consequence of Lemma~\ref{lemloir} applied to the smooth map $\psi_1^{+}$. All implicit constants are, again, uniform over $B\subset \textrm{SL}_n(\R)$.\\

Finally, apply the uniform estimates stated in the second part of Lemma~\ref{lemloir} with the following assumptions~: 
\begin{itemize}
\item the family $\mathcal{F}$ is taken as $\mathcal{F}\;=\; \left\{\Psi_{\mathfrak{g}, k}^{+}\right\}_{\mathfrak{g}\in B}$;
\item the compact set $\mathcal{K}\subset\R^n$ is defined as $K\left(2\kappa\right)$, where $\kappa>0$ is an upper bound for the operator norms of all $\mathfrak{g}\in B$.
\item the smooth map $\varphi_{\mathcal{F}}$ is defined by setting $\varphi_{\mathcal{F}}(\bm{x})=\Psi_1^{+}\!\left(\bm{x}/(2\kappa)\right)$ for all $\bm{x}\in\R^n$.
\end{itemize} 
It is then a consequence of the assumption (ii) above, of the relations~\eqref{defpsiplgk} and of Lemma~\ref{egpolorder} that the smooth map $\varphi_{\mathcal{F}}$ satisfies the assumption (A2) in Lemma~\ref{lemloir}. Furthermore, from repeated applications of the chain rule, it should be clear that for any multi--index $\bm{l}\in\N^n$, any $\mathfrak{g}\in B$ and any $\bm{x}\in\R^n$,
\begin{equation*}
\left|\frac{\partial^{\left|\bm{l}\right|}\Psi_{\mathfrak{g}, k}^{+}}{\partial\bm{x}^{\bm{l}}}(\bm{x})\right|\;\ll\; \max_{0\le\left|\bm{j}\right|\le \left|\bm{l}\right|}\;\max_{\bm{x}\in\R^n}\;\left|\frac{\partial^{\left|\bm{j}\right|}\Psi_{k}^{+}}{\partial\bm{x}^{\bm{j}}}(\bm{x})\right|
\end{equation*} with an implicit constant depending on the compact set $B$ only. Recalling the explicit definition of the quantity $\nu_{P, \mathfrak{g}}^{+}\left(\Omega, k, T\right)$ in~\eqref{defnupgomegkT},  Lemma~\ref{lemloir} then implies the existence of  integers $M, N\ge 1$ such that for some implicit constant depending on  $B$,
\begin{equation}\label{estim3unif}
\left| \frac{\nu_{P, \mathfrak{g}}^{+}\left(\Omega, k, c\right) }{c^{r_P(\Omega)}\cdot \left|\log c\right|^{m_P(\Omega)-1}}-  \Theta_{P, \mathfrak{g}}^{+}\left(\Omega, k\right)\right| \;\ll\; \frac{C_N(k)^M}{\left(f_{P}\left(\Omega, c^{-1}\right)\right)^{\delta_P\left(\Omega, K\right)}}\cdotp
\end{equation}
To define precisely the exponent $\delta_P\left(\Omega, K\right)$ and the choice of the function $f_{P}\left(\Omega, \;\cdot\;\right)$ in this relation, let $0\le j_{\mathfrak{g}, k}^{+}\left(\Omega\right)\le m_P(\Omega)-2$ denote the exponent associated to the second largest nonzero coefficient of the real, univariate and monic polynomial $R_{\psi_{\mathfrak{g}, k}^+}^{\Omega}[x]$ introduced in Lemma~\ref{lemloir}. If there exists an integer $k_0\ge 1$ such that this exponent is well-defined  for all $k\ge k_0$, set $f_{P}\left(\Omega, c\right)=\log(c)$ and $$\delta_P\left(\Omega, K\right)\;=\; \min_{k\ge k_0}\;\left\{m_P(\Omega)-1- j_{\mathfrak{g}, k}^{+}\left(\Omega\right)\right\}\;\ge\; 1.$$Otherwise (this is the case when $m_P(\Omega)=1$), set $f_{P}\left(\Omega, c\right) = c$ and $\delta_P\left(\Omega, K\right)=\varepsilon$, where $\varepsilon>0$ is the quantity introduced in the asymptotic relation~\eqref{asympmu}.\\

Adding up the three inequalities~\eqref{estim1unif}, \eqref{estim2unif} and~\eqref{estim3unif} then concludes the proof of the lemma.
\end{proof}

The completion of the proof of Theorem~\ref{mainthmvolomega} relies on a suitable choice of the sequences of smooth maps $\left(\psi_k^{\pm}\right)_{k\ge 1}$  in  $\mathcal{C}_c^\infty(\R^n)$. A classical convolution argument detailed in~\cite[Theorem~1.4.1, p.25]{horm} shows that there exist such  sequences 
 satisfying the above as\-sump\-tions (i) and (ii) with the additional property that for all $k\ge 1$, $$\max_{\bm{x}\in\R^n}\;\left|\frac{\partial^{\left|\bm{l}\right|}\psi_k^{\pm}}{\partial\bm{x}^{\bm{l}}}(\bm{x})\right|\;\ll\; k^{\left|\bm{l}\right|}.$$Here, the implicit constant depends only on the multi--index $\bm{l}\in\N^n$ and on the dimension $n\ge 1$. With the notation of Lemma~\ref{estimerreurunif1}, this yields the existence of integers  $M,N\ge 1$ such that for all $k\ge 1$, 
\begin{equation*}
\left|\frac{\nu_{P, \mathfrak{g}}\left(\Omega, K, c\right) }{c^{r_P(\Omega)}\cdot \left|\log c\right|^{m_P(\Omega)-1}} -  \gamma_{P, \mathfrak{g}}\left(\Omega, K\right)\right| \;\ll\; \frac{k^{NM}}{\left(f_{P}\left(\Omega, K, c^{-1}\right)\right)^{\delta_P\left(\Omega, K\right)}}+ \frac{1}{k},
\end{equation*}
where $\delta_P\left(\Omega, K\right)$ is strictly positive and where the multiplicative constant is uniform over all compact sets containing the matrix $\mathfrak{g}\in SL_n(\R)$. Let $c=c(T)$ be as in the relation~\eqref{limiassumpbis}. Choose then $k=k(T)$ as the integer part of $\left(f_{P}\left(\Omega, c(T)^{-1}\right)\right)^{\delta_P\left(\Omega, K\right)/(1+MN)}$ in such a way that $$\frac{\nu_{P, \mathfrak{g}}\left(\Omega, K, c(T)\right)}{c(T)^{r_P(\Omega)}\cdot \left|\log c(T)\right|^{m_P(\Omega)-1}}\;=\; \gamma_{P, \mathfrak{g}}\left(\Omega, K\right)+O\left(\left(f_{P}\left(\Omega, c(T)^{-1}\right)\right)^{-\delta_P\left(\Omega, K\right)/(1+MN)}\right)$$with the error term satisfying the uniformity property  required in the statement of Theorem~\ref{mainthmvolomega}. From the equation~\eqref{volnnupgomegT}, this completes the proof of Theorem~\ref{mainthmvolomega} upon redefining the value of the exponent $\delta_P\left(\Omega, K\right)$ as $\delta_P\left(\Omega, K\right)/\left(1+MN\right)$.
\end{proof}

\section{Local and Global Real Log--Canonical Thresholds}\label{locglobrlct}  

Throughout this section, the set $Y\subset\R^n$ refers either to a globally semianalytic set $K$ meeting the assumptions $ (\mathcal{H}_1)-(\mathcal{H}_3)$ (with respect to the polynomial $P(\bm{x})\in\R[\bm{x}]$ as far as $(\mathcal{H}_3)$ is concerned) or else to a set $\Omega\subset\R^n$ belonging to the collection  $\mathfrak{O}$ defined in~\eqref{defoobig}. Then, $Y$ is a finite union of sets such as~\eqref{setomega}, where the maps $f_i$ ($1\le  i\le l$) are all analytic (and even polynomials when $Y=\Omega$). Set also
\begin{equation*}
\zeta_{P}(Y, \, \cdot\,)\;=\;
\begin{cases}
\zeta_{P}(\Omega, \, \cdot\,) &\textrm{when } Y=\Omega;\\ 
\zeta_{P} &\textrm{when } Y=K,\\
\end{cases}
\end{equation*} 
where the meromorphic distribution $\zeta_{P}(\Omega, \, \cdot\,)$ is defined at all complex numbers with non--positive real parts by~\eqref{restriczeta}. When $Y=K$, the dependency on $K$ in the above notation is expressed by the fact that the distribution $\zeta_{P}$ will be tested against maps $\psi$ in $\mathcal{C}_c^\infty(\R^n)$ meeting the support restriction condition~\eqref{suppproper} part of the assumption $(\mathcal{H}_3)$. Let then 
\begin{equation*}
\left(r_P(Y), m_P(Y)\right)\;=\;
\begin{cases}
\left(r_P(K), m_P(K)\right)&\textrm{when } Y=K;\\ 
\left(r_P(\Omega), m_P(\Omega)\right)  &\textrm{when } Y=\Omega.\\
\end{cases}
\end{equation*} 
Here, the quantities $r_P(K)$ and $m_P(K)$ are defined as part of $(\mathcal{H}_3)$, and $r_P(\Omega)$, whenever well-defined,  denotes the smallest real pole of the meromorphic distribution $\zeta_{P}(\Omega, \, \cdot\,)$, the integer $m_P(\Omega)$ being its order. Let finally $B_P(Y, s)$ denote the Sato--Bernstein polynomial defined in~\eqref{globalsbomega}. \\

The goal of this section is to establish two claims~: on the one hand, that the quantity $r_P(Y)$  is a strictly positive rational number that can be determined by re\-so\-lu\-tions of singularities in a finite number of steps under the assumption $(\mathcal{H}_2)$; on the other, that $r_P(Y)$ is a root of the polynomial $B_P(Y, -s)$   whose multiplicity is at most its order $m_P(Y)$ as a pole of $\zeta_{P}(Y , \, \cdot\,)$. 
The same claim is proved to hold for the global meromorphic distribution $\zeta_{P}$ upon considering its smallest pole $r_P$, the cor\-res\-pon\-ding multiplicity $m_P$ and the global Sato-Bernstein polynomial $B_P(s)$. Finally, the specialisation to the case where the polynomial $P(\bm{x})$ is the homogeneous form $\left\|\bm{F}(\bm{x})\right\|^2$ will complete the proof of Theorem~\ref{volestim}.\\

To establish the aforementioned claims, define, whenever it makes sense, the \emph{(global) real log--canonical threshold of the map determined by $P(\bm{x})$ over the domain $\Omega$ with respect to the test function $\psi$} as the pair 
\begin{equation}\label{globalrlct}
\R\textrm{-}LCT_{P}(\Omega, \psi)=\left(r_P(\Omega,\psi), m_P(\Omega,\psi)\right),
\end{equation} 
where, here again,  $r_P(\Omega,\psi)$ denotes the smallest real pole of the map $\left\langle\zeta_P\left(\Omega, \,\cdot\, \right), \psi\right\rangle$  and where $m_P(\Omega,\psi)$ is its order.  If not well--defined, let $RLCT_{P}(\Omega, \psi)=\left(\infty, -\right)$; the quantity $m_P(\Omega,\psi)$ is in particular left undefined (this happens, for instance, when the variety $\mathcal{Z}_\R(P)$ does not intersect the support of the test function $\psi$). Similarly, given a test function $\psi$, let   
\begin{equation}\label{globalrlctbis}
\R\textrm{-}LCT_{P}(K, \psi)=\left(r_P(\psi), m_P(\psi)\right),
\end{equation} 
where $r_P(\psi)$ denotes the smallest real pole of the map $\left\langle\zeta_P, \psi\right\rangle$  and $m_P(\psi)$ its order (adopting the same convention as above when these quantities  do not exist). When $\psi$ meets the support restriction~\eqref{suppproper}, $\R\textrm{-}LCT_{P}(K, \psi)$ equals the pair $\left(r_P(K), m_P(K)\right)$, which is well-defined by assumption.\\

Given a  subset $Z\subset\R^n$ with positive measure, the notation $\R\textrm{-}LCT_{P}(Z|\psi)$ is used, provided it makes sense, to refer to the pair given by the smallest real pole of the zeta function defined for $\textrm{Re}(s)\le 0$ by $s\mapsto \int_Z\left|P(\bm{x})\right|^{-s}\cdot\psi(\bm{x})\textrm{d}\bm{x}$, and by its order. It is here again set to be $(\infty, -)$ if not well-defined.\\

The following lemma is concerned with the determination of the real log--canonical threshold in  the case where the polynomial $P(\bm{x})$ is a monomial.

\begin{lem}\label{lemmonorlct}
Assume that $\psi$ is an element in $\mathcal{C}_{c}^{\infty}(\R^n)$ such that $\psi(\bm{0})>0$. Let $Y^+_n=(0,\infty)^n$ be the nonnegative orthant in $\R^n$ and let $\bm{k}=(k_1, \dots, k_n)\in\N_0^n$ and  $\bm{h}=(h_1, \dots, h_n)\in\N_0^n$ be $n$--tuples of non--negative integers. Then, $$\R\textrm{-}LCT_{\bm{x}^{\bm{k}}}\left(Y^+_n|\bm{x}^{\bm{h}}\cdot\psi\right)\;=\; \left(r, m\right),$$ where $$ r\;=\; \min_{1\le i\le n}\frac{h_i+1}{k_i} \qquad \qquad  \textrm{and}\qquad \qquad  m\;=\; \#\left\{1\le i\le n\; :\; \frac{h_i+1}{k_i}\;=\; r\right\}.$$ In these relations, a ratio is conventionally taken to be  infinite if the denominator  vanishes, and  the quantity $m$ is left undefined whenever $r=\infty$.
\end{lem}

\begin{proof}
The claim regarding the value of the rational $r$  is the specific calculation worked out in the proof of~\cite[Lemma 7.3]{AGZV}. The value of the integer $m$ follows from~\cite[Lemma 7.5, (1)]{AGZV}.
\end{proof}

Lemma~\ref{lemmonorlct} extends to the case of any polynomial $P(\bm{x})$ by resolving singularities. This is the main content of the following statement.
 
\begin{lem}\label{locrlct}
Let $\bm{x}_0$ be a point in the closure 
of the set $Y\subset\R^n$. Then, there exists a neighbourhood $Y_{\bm{x}_0}$ of $\bm{x}_0$ in $\R^n$ such that the following properties hold~: 
\begin{itemize}
\item[(W1)]  for any  two functions $\varphi_{\bm{x}_0}, \varphi'_{\bm{x}_0}\ge 0$ in $\mathcal{C}_{c}^{\infty}(\R^n)$ with supports contained in $Y_{\bm{x}_0}$ and taking  strictly positive values at $\bm{x}_0$, one has that $\R\textrm{-}LCT_{P}\left(Y, \varphi_{\bm{x}_0}\right) = \R\textrm{-}LCT_{P}\left(Y, \varphi'_{\bm{x}_0}\right)$. Let  
\begin{equation}\label{ocalrlct}
\R\textrm{-}LCT_{P, \bm{x}_0}(Y)=\left(r_{P, \bm{x}_0}(Y), m_{P, \bm{x}_0}(Y)\right)
\end{equation}
denote this common value.

\item[(W2)] When $P(\bm{x}_0)=0$, the pair $\R\textrm{-}LCT_{P, \bm{x}_0}\left(Y\right)$ can be determined by resolving sin\-gu\-la\-ri\-ties and the quantity $r_{P, \bm{x}_0}(Y)$ is then a strictly positive rational number. If $P(\bm{x}_0)\neq 0$, the  quantity $r_{P, \bm{x}_0}(Y)$ is undefined (and is thus conventionally set to equal $\infty$).

\item[(W3)] let $\zeta_{P, \bm{x}_0}(Y, \, \cdot\,)$ be the restriction of the distribution $\zeta_{P}(Y, \, \cdot\,)$ to the open neighbourhood $Y_{\bm{x}_0}$~: it is  defined for all $\varphi_{\bm{x}_0}$ in $\mathcal{C}_{c}^{\infty}(\R^n)$ with support contained in $Y_{\bm{x}_0}$ and all $s\in\C$ such that $\textrm{Re}(s)\le 0$ by setting $\left\langle \zeta_{P, \bm{x}_0}\left(Y, s\right), \varphi_{\bm{x}_0}\right\rangle= \left\langle \zeta_{P}\left(Y, s\right), \varphi_{\bm{x}_0}\right\rangle$. Then, assuming that $P(\bm{x}_0)=0$, the rational $r_{P, \bm{x}_0}(Y)> 0$ introduced in~\eqref{ocalrlct} is the smallest real pole of the local distribution $\zeta_{P, \bm{x}_0}(Y, \, \cdot\,)$ and $m_{P, \bm{x}_0}(Y)$ is its order. 

\item[(W4)] Whenever finite, $r_{P, \bm{x}_0}(Y)$ is also a root of the local Sato--Bernstein polynomial $B_{P,\bm{x}_0}(Y, -s)$ defined in~\eqref{globalsbomega} with multiplicity at most  $m_{P, \bm{x}_0}(Y)$.
\end{itemize}
\end{lem}

The quantity $\R\textrm{-}LCT_{P, \bm{x}_0}(Y)$ introduced in~\eqref{ocalrlct} defines the \emph{(local) real log--canonical threshold of the map $\bm{x}\mapsto P(\bm{x})$ at the point $\bm{x}_0$ of the closure of the domain $Y$}. One of the main features of the proof of the lemma is to provide an explicit formula for each of the components of this pair. Related expressions are known for the first component $r_{P, \bm{x}_0}(Y)$, whenever well--defined --- see, e.g., ~\cite[\S 8.5 \& \S 10.7]{kollar} and~\cite{saito}.

\begin{proof}
If $P(\bm{x}_0)\neq 0$, there exists a neighbourhood $Y_{\bm{x}_0}$ of $\bm{x}_0$  such that $P(\bm{x})$ does not vanish on $Y_{\bm{x}_0}$. It is then immediate that $\R\textrm{-}LCT_{P, \bm{x}_0}(Y)=(\infty, -)$ and the lemma is trivially true in this case. Assume therefore that $P(\bm{x}_0) =  0$. \\

Recall that $Y$ is a finite union of sets of the form~\eqref{setomega}, where the maps $f_i$ ($1\le  i\le l$) are  analytic. Apply Hironaka's Theorem~\ref{hironaka} on simultaneous resolution of singularities to the set of all maps $\bm{x}\mapsto f_i(\bm{x}_0+\bm{x})$ ($1\le  i\le l$) and $\bm{x}\mapsto P(\bm{x}_0+\bm{x})$ vanishing at the origin. It provides the existence of an analytic morphism $g$ defined on a real manifold $\mathcal{M}$ and taking values in a neighbourhood $\mathcal{W}$ of $\bm{x_0}$ in $\R^n$ such that $P(g(\bm{y}))$ and the relevant $f_i(g(\bm{y}))$'s take the monomial form~\eqref{monoform} on some chart $\mathcal{M}_{\bm{y}_0}\subset \mathcal{M}$. These charts $\mathcal{M}_{\bm{y}_0}$ are furthermore centered around points $\bm{y}_0$ of the preimage $g^{-1}(\{\bm{x}_0\})$. Since $g$ is  proper, this preimage is compact and therefore admits a finite open subcover $\left\{\mathcal{M}_{\bm{y}_0}\right\}_{\bm{y}_0\in I(\bm{x}_0)}$, where $I(\bm{x}_0)$ is a finite index set.\\

The first step is  to show that $g\left(\bigcup_{\bm{y}_0\in I(\bm{x}_0)}\mathcal{M}_{\bm{y}_0}\right)$ contains an open neighbourhood $Y^{(1)}_{\bm{x}_0}$ of $\bm{x}_0$ in $\R^n$. Indeed, assume this is not the case. Then, there exists a bounded sequence of points $\left\{\bm{x}_0^{(j)}\right\}_{j\ge 1}$ in $\mathcal{W}\backslash g\left(\bigcup_{\bm{y}_0\in I(\bm{x}_0)}\mathcal{M}_{\bm{y}_0}\right)$ with limit $\bm{x}_0$. Since $g$ is surjective, there exists a sequence $\left\{\bm{y}_0^{(j)}\right\}_{j\ge 1}$ in $\mathcal{M}$ such that $g\left(\bm{y}_0^{(j)}\right)=\bm{x}_0^{(j)}$ for all $j\ge 1$. The sequence $\left\{\bm{x}_0^{(j)}\right\}_{j\ge 1}$ being bounded, the properness of the map $g$ implies that the sequence  $\left\{\bm{y}_0^{(j)}\right\}_{j\ge 1}$ lies in a compact set. Upon considering a subsequence, it can  be assumed to converge to a point $\bm{y}_0^{*}$, which remains outside the set $\bigcup_{\bm{y}_0\in I(\bm{x}_0)}\mathcal{M}_{\bm{y}_0}$ as it is open. However, the relation $g\left(\bm{y}_0^{*}\right)=\bm{x}_0$ shows that $\bm{y}_0^{*}$ lies in the preimage $g^{-1}(\{\bm{x}_0\})$, which is covered by $\left\{\mathcal{M}_{\bm{y}_0}\right\}_{\bm{y}_0\in I(\bm{x}_0)}$. This contradiction establishes the claim. \\

Consider now a maximal subset $J(\bm{x}_0)$ of the index set $I(\bm{x}_0)$ such that for all $\bm{y}_0\in J(\bm{x}_0)$, the chart $\mathcal{M}'_{\bm{y}_0}=\mathcal{M}_{\bm{y}_0}\cap g^{-1}\left(Y^{(1)}_{\bm{x}_0}\cap Y\right)$ has positive measure. Let $\left\{\pi_{\bm{y}_0}\right\}_{\bm{y}_0\in J(\bm{x}_0)}$ be a partition of unity subordinate to the collection of open sets $\left\{\mathcal{M}'_{\bm{y}_0}\right\}_{\bm{y}_0\in J(\bm{x}_0)}$ such that each $\pi_{\bm{y}_0}$ takes a positive value at $\bm{y}_0$. Let  $r_{P, \bm{x}_0}\left(Y, \varphi_{\bm{x}_0}\right)$ denote the smallest real pole of the map $\sigma\mapsto \left\langle \zeta_{P, \bm{x}_0}\left(Y, \sigma\right), \varphi_{\bm{x}_0}\right\rangle$, assuming for the time being that this quantity  is well--defined. Then, given a smooth test function $\varphi_{\bm{x}_0}$ with support in $Y_{\bm{x}_0}^{(1)}$, assuming that the considered partition of unity has been chosen so that it sums to the constant function 1 on its domain of definition intersected with $g^{-1}\left(\textrm{Supp}\; \varphi_{\bm{x}_0}\right)$, it holds that for any $\sigma<r_{P, \bm{x}_0}\left(Y, \varphi_{\bm{x}_0}\right)$,
\begin{align}
\left\langle \zeta_{P, \bm{x}_0}\left(Y, \sigma\right), \varphi_{\bm{x}_0}\right\rangle\;=\; \sum_{\bm{y}_0\in J(\bm{x}_0)} \int_{\mathcal{M}'_{\bm{y}_0}} \left|P(g\left(\bm{y}\right))\right|^{-\sigma}\cdot\varphi_{\bm{x}_0}\left(g\left(\bm{y}\right)\right)\cdot \left|\textrm{Jac}_g(\bm{y})\right|\cdot \pi_{\bm{y}_0}\left(\bm{y}\right)\cdot\textrm{d}\bm{y}.\label{eqtdeczeat}
\end{align}
For each $\bm{y}_0\in J(\bm{x}_0)$, the boundary constraints $f_i(g(\bm{y}))\ge 0$ and $P(g(\bm{y}))\ge 0$ are monomial inequalities. The set $\mathcal{M}'_{\bm{y}_0}$ is thus a finite union of closed orthants which intersect any neighbourhood of $\bm{y}_0$. As a consequence, each of the above integrals can be reduced to a finite sum of, say, $K(\bm{y}_0)\ge 1$ integrals of the form 
\begin{equation}\label{zetasimplifieddec}
\int_{Y_n^+}\bm{y}^{\bm{h}-\sigma\bm{k}}\cdot\psi(\bm{y})\cdot\textrm{d}\bm{y},
\end{equation} 
where $Y_n^+=(0, \infty)^n$, where $\psi(\bm{0})>0$ and where for $1\le j\le K(\bm{y}_0)$, $\bm{k} = \bm{k}_j(\bm{y_0}) \in\N_0^n$ and $\bm{h}=\bm{h}_j(\bm{y_0})\in\N_0^n$. Here, the exponents  $\bm{k}$ and $ \bm{h}$   depend neither on the smooth function $\psi$ nor on the choice of the orthant indexing the finite sum. The real log--canonical threshold related to each of the integrals~\eqref{zetasimplifieddec} is then determined from Lemma~\ref{lemmonorlct}, and so is $\R\textrm{-}LCT_{P}\left(Y, \varphi_{\bm{x}_0}\right)$ in view of the decomposition~\eqref{eqtdeczeat}. Since the resulting  formula is independent of the choice of the test function $\varphi_{\bm{x}_0}$, the claim (W1) in the statement follows. Explicitly, if, in the above notations, $\bm{k}_j(\bm{y_0})=\left(k_j^{(1)}(\bm{y_0}), \dots, k_j^{(n)}(\bm{y_0})\right)$ and $\bm{h}_j(\bm{y_0})=\left(h_j^{(1)}(\bm{y_0}), \dots, h_j^{(n)}(\bm{y_0})\right)$, then $\R\textrm{-}LCT_{P, \bm{x}_0}(Y)=\left(r_{P, \bm{x}_0}(Y), \, m_{P, \bm{x}_0}(Y)\right)$, where $$r_{P, \bm{x}_0}(Y)\;=\; \min_{\bm{y}_0\in J(\bm{x}_0)}\;\min_{1\le j\le K(\bm{y}_0)}\;\min_{1\le l\le n}\; \frac{h_j^{(l)}(\bm{y_0})+1}{k_j^{(l)}}$$ and $$m_{P, \bm{x}_0}(Y)\;=\;\max_{\bm{y}_0\in J(\bm{x}_0)}\;\max_{1\le j\le K(\bm{y}_0)}\;\#\left\{1\le l\le n\; :\; \frac{h_j^{(l)}(\bm{y_0})+1}{k_j^{(l)}}\;=\; r_{P, \bm{x}_0}(Y, \psi)\right\}.$$\\

These formulae, and the fact that they do not depend on the choice of the resolution of singularities, hold provided that the existence of $\R\textrm{-}LCT_{P, \bm{x}_0}(Y, \varphi_{\bm{x}_0})$ is justified when  $\varphi_{\bm{x}_0}$ is a smooth test function with support in any small enough neighbourhood of $\bm{x}_0$. This follows from the fact that the local zeta distribution $\zeta_{P, \bm{x}_0}(Y, \, \cdot\,)$ admits a meromorphic continuation to the whole complex plane. This, in turn, is a consequence of  the following local version of the identity~\eqref{elemdual}, which is deduced from the existence of the local Sato--Bernstein polynomial $B_{P,\bm{x}_0}(Y, s)$~: keeping the same notation as in the aforementioned identity, there exis\-ts a neighbourhood $Y_{\bm{x}_0}^{(2)}$ of $\bm{x}_0$ in $\R^n$ such that for any smooth test function $\varphi_{\bm{x}_0}$  supported in $Y_{\bm{x}_0}^{(2)}$, it holds that 
\begin{align}\label{elemdualbis}
B_{P,\bm{x}_0}(Y, -s)\cdot \left\langle \zeta_{P, \bm{x}_0}\left(Y, s\right), \varphi_{\bm{x}_0}\right\rangle\;& =\;  \left\langle \zeta_{P, \bm{x}_0}\left(Y, s-1\right), \mathcal{D}^*\left(\bm{x}, -s, \bm{\partial}\right)\varphi_{\bm{x}_0}(\bm{x})\right\rangle.
\end{align}
To complete the proof of (W1) and to establish (W2), it is therefore enough to set  $$Y_{\bm{x}_0}= Y_{\bm{x}_0}^{(1)}\cap Y_{\bm{x}_0}^{(2)}.$$ 

As for (W3) and (W4), they are also consequences of~\eqref{elemdualbis}. Indeed, to establish first the claim (W3), let $\tilde{\varphi}_{\bm{x}_0}$ be any \sloppy smooth test function  supported in $Y_{ \bm{x}_0}$ such that  the meromorphic  map $s\mapsto \left\langle \zeta_{P, \bm{x}_0}\left(Y, s\right), \tilde{\varphi}_{\bm{x}_0}\right\rangle$ achieves the smallest real pole, say $r>0$, of the zeta distribution $\zeta_{P, \bm{x}_0}\left(Y, \, \cdot\, \right)$. Assume that the smooth test function $\varphi_{\bm{x}_0}$ supported in $Y_{\bm{x}_0}$ bounds from above $\left|\tilde{\varphi}_{\bm{x}_0}\right|$. Then, given $\sigma<r$, the inequality 
\begin{equation}\label{ineqlast-}
\left|\left\langle \zeta_{P, \bm{x}_0}\left(Y, \sigma\right), \tilde{\varphi}_{\bm{x}_0}\right\rangle\right|\;\le\; \left\langle \zeta_{P, \bm{x}_0}\left(Y, \sigma\right), \varphi_{\bm{x}_0}\right\rangle
\end{equation} 
shows that $r$ is also the smallest pole of the meromorphic map defined on the right--hand side, and that the orders coincide (the argument is similar to the one used in the relation~\eqref{ineqordremultp}). From (W1), this establishes that the smallest pole of the distribution $\zeta_{P, \bm{x}_0}\left(Y, \, \cdot\, \right)$ is $r_{P, \bm{x}_0}(Y)$ and that it has order $m_{P, \bm{x}_0}(Y)$, whence (W3).\\

Finally, (W4) is a consequence of the following observation~: since the right--hand side of~\eqref{elemdualbis} is holomorphic at $r_{P, \bm{x}_0}(Y)$ (this is a consequence the definition of $r_{P, \bm{x}_0}(Y)$ as the smallest real pole of the distribution $ \zeta_{P, \bm{x}_0}\left(Y, \,\cdot\,\right)$), for any given smooth test function $\varphi_{\bm{x}_0}$ supported in $Y_{\bm{x}_0}$, the order of the pole $r_{P, \bm{x}_0}(Y)$ of the meromorphic map $\sigma\mapsto \left\langle \zeta_{P, \bm{x}_0}\left(Y, \sigma\right), \varphi_{\bm{x}_0}\right\rangle$ cannot be larger than the multiplicity of $r_{P, \bm{x}_0}(Y)$ as a root of the polynomial  $B_{P,\bm{x}_0}(Y, -s)$. \\

This completes the proof of the Lemma~\ref{locrlct}.
\end{proof}

The relation between the global real log--canonical thresholds defined in~\eqref{globalrlct} and in~\eqref{globalrlctbis} and the local ones introduced in Lemma~\ref{locrlct} is made clear in the following proposition. In order to state it, introduce an order over pairs of nonnegative reals by setting, given $a,b,c,d\ge 0$, 
\begin{equation}\label{ordering}
(a,b)\preccurlyeq (c,d)\quad \Longleftrightarrow\quad \left(a<c\right)\;\vee\;\left(\left(a=c\right)\wedge\left(b\ge d\right)\right).
\end{equation}
Explicitly, this corresponds to the ordering induced by the asymptotic growth of the map $x\mapsto x^{-a}\cdot\left(\log x\right)^b$ at infinity. Set also $(a,b)\succcurlyeq (c,d)$ to mean $(c,d)\preccurlyeq (a,b)$.\\

Recall that when a smooth test function $\psi$ meets the support restriction conditon~\eqref{suppproper}, it takes, by assumption,  \emph{strictly} positive values on the set $C$ introduced  in this condition.

\begin{prop}\label{proprlctgene}
Keep the notation of Lemma~\ref{locrlct}. When $Y=\Omega$, recall that the real polynomial $P(\bm{x})$ is assumed to be homogeneous and when $Y=K$, recall that the set $K$ is assumed to be  semianalytic and to meet the conditions $\left(\mathcal{H}_1\right)-\left(\mathcal{H}_3\right)  $.  The following claims then hold true~: 
\begin{itemize}
\item[(X1)] Let $\psi\ge 0$ be a smooth test function such that $\psi(\bm{0})>0$ when $Y=\Omega$ and such that support restriction condition~\eqref{suppproper} is met when $Y=K$. Then, 
\begin{equation}\label{rlctsemiglob}
\R\textrm{-}LCT_{P}(Y, \psi)\;=\;
\min_{\bm{x}_0\in\R^n} \R\textrm{-}LCT_{P, \bm{x}_0}(Y)\;=\;\min_{\bm{x}_0\in\mathcal{Z}_{\R}(P)}\R\textrm{-}LCT_{P, \bm{x}_0}(Y),
\end{equation}
where the minimum is taken with respect to the order induced by~\eqref{ordering}.  This pair can furthermore be computed by resolving singularities in a finite number of steps and its first component is a strictly positive rational number.

\item[(X2)] The pair $\R\textrm{-}LCT_{P}(Y, \psi)$ is independent of the choice of $\psi\ge 0$ provided that $\psi(\bm{0})>0$ when $Y=\Omega$ and provided that support restriction condition~\eqref{suppproper} is met when $Y=K$. Denote the common value by $\R\textrm{-}LCT_{P}(Y)= \left(r_P(Y), m_P(Y)\right)$. Then, $r_P(Y)$  is the smallest real pole of the distribution $\zeta_{P}(Y, \, \cdot\,)$ and $m_P(Y)$ is its order,  under the following additional requirement when $Y=K$~:  only those poles obtained when  testing the distribution against smooth maps satisfying the support restriction condition~\eqref{suppproper} must be considered.

\item[(X3)] The pole $r_P(Y)$ is a root of the global Sato--Bernstein polynomial $B_P\left(Y, -s\right)$ with respect to the domain $Y$ as defined in~\eqref{globalsbomegabis}. Its multiplicity as a root is, furthermore, at most its order $m_P(Y)$ as a pole. 
\end{itemize}
\end{prop}

The proof of the statement makes it clear that the above defined pair $\R\textrm{-}LCT_{P}(Y, \psi)$ belongs to the set $\Q_{>0}\times \N_{\ge 1}$ if, and only if, the support of $\psi$ intersects non trivially the algebraic variety $\mathcal{Z}_{\R}(P)$. This is what is guaranteed by the assumption of homogeneity when $Y=\Omega$ and $\psi(\bm{0})>0$ and by the condition $(\mathcal{H}_2)$ when $Y=K$. In the latter case, it is immediate that the definition of the pair $\R\textrm{-}LCT_{P}(K)$ ensuing from the above point (X2) coincides with the pair $\left(r_P(K), m_P(K)\right)$ introduced in the assumption $(\mathcal{H}_3)$.\\

Let $\mathfrak{O}$ be the finite collection of sets defined in~\eqref{defoobig}. Define the \emph{real log--canonical threshold of the polynomial $P(\bm{x})$} as 
\begin{equation}\label{rlctpolezetaprim}
\mathbb{R}\!\!-\!\!LCT_P\;=\; \min_{\Omega\in\mathfrak{O}} \; \R\textrm{-}LCT_{P}(\Omega).
\end{equation} 
From Lemma~\ref{lempolemulti} and from the above proposition,  
\begin{equation}\label{rlctpolezeta}
\R\textrm{-}LCT_P\;=\; \left(r_P, m_P\right)\;\in\; \Q_{> 0}\times \N,
\end{equation}
where $r_P$ is the smallest real pole of the zeta distribution $\zeta_{P}$ and $m_P$ its order.

\begin{proof}
The second relation in~\eqref{rlctsemiglob} follows from the fact that the local real-log canonical threshold is $(\infty, -)$ at a point where the polynomial $P(\bm{x})$ does not vanish. To prove the first equation, assign to each $\bm{x}_0$ lying in the support of the test function $\psi\ge 0$ an open neighbourhood $Y_{\bm{x}_0}$ in $\R^n$ satisfying all the conclusions of Lemma~\ref{locrlct}. Extract a cover $\left\{Y_{\bm{x}_0}\right\}_{\bm{x}_0\in S(\psi)}$ of the support of $\psi$ , where the index set $S(\psi)$ has a finite cardinality. Consider a smooth partition of unity $\left\{\pi_{\bm{x}_0}\right\}_{\bm{x}_0\in S(\psi)}$ subordinate to this cover such that $\pi_{\bm{x}_0}(\bm{x}_0)>0$ for each $\bm{x}_0\in S(\psi)$ and such that this partition of unity sums to the constant function 1 in a small enough neighbourhood of $\textrm{Supp}\; \psi$. Then, 
\begin{equation}\label{dzfzgeh}
\left\langle \zeta_{P}\left(Y, s\right), \psi\right\rangle\;=\; \sum_{\bm{x}_0\in S(\psi)}\left\langle \zeta_{P}\left(Y, s\right), \psi\cdot \pi_{\bm{x}_0}\right\rangle 
\end{equation} 
whenever $s\in\C$ is not a pole of any of the maps in the sum of the right--hand side. As in the proof of Lemma~\ref{lempolemulti}, this decomposition implies that 
\begin{equation}\label{rlctminpsi}
\R\textrm{-}LCT_{P}(Y, \psi)\;=\;\min_{\bm{x}_0\in S(\psi)} \R\textrm{-}LCT_{P, \bm{x}_0}\left(Y, \psi\cdot \pi_{\bm{x}_0}\right). 
\end{equation} 
In this relation, the claim (W1) in Lemma~\ref{locrlct} implies that $\R\textrm{-}LCT_{P, \bm{x}_0}\left(Y, \psi\cdot \pi_{\bm{x}_0}\right) = \R\textrm{-}LCT_{P, \bm{x}_0}\left(Y\right)$ whenever $\psi\left(\bm{x}_0\right)>0$. Dealing with the case where $\psi\left(\bm{x}_0\right)=0$ requires more work. Assume therefore that for a given $\bm{x}_0\in S(\psi)$, it holds that $\psi\left(\bm{x}_0\right)=0$. The argument differs depending on whether $Y=\Omega$ or $Y=K$.$$\quad$$

Consider first the case where $Y=K$. Let $ \widetilde{\psi}$ be a nonnegative map in $\mathcal{C}_{c}^{\infty}(\R^n)$ such that~:
\begin{itemize} 
\item[(a)] with the notations of the support restriction condition~\eqref{suppproper}, the support of $\widetilde{\psi}$ is contained both in the open neighbourhood $U$ of $K$ and in the small enough neighbourhood of $\textrm{Supp}\; \psi$ where the partition of unity $\left\{\pi_{\bm{x}_0}\right\}_{\bm{x}_0\in S(\psi)}$ sums to the constant function 1;
\item[(b)] it holds that $0\le \psi\le \widetilde{\psi}$ and that $\widetilde{\psi}$ is strictly positive in the interior of its support.
\end{itemize}

It is then clear that for any $\bm{x}_0\in S(\psi)$, $$ 0\;\le\; \left\langle \zeta_{P}\left(K, \sigma\right), \psi\cdot \pi_{\bm{x}_0}\right\rangle\;\le \; \left\langle \zeta_{P}\left(K, \sigma\right), \widetilde{\psi}\cdot \pi_{\bm{x}_0}\right\rangle$$as long as the real $\sigma$  remains smaller than the smallest real pole of the left-hand side. As a consequence, 
\begin{equation}\label{comparlctco0}
\R\textrm{-}LCT_{P}\left(K, \psi\cdot \pi_{\bm{x}_0}\right)  \;\succcurlyeq\; \R\textrm{-}LCT_{P}\left(K, \widetilde{\psi}\cdot \pi_{\bm{x}_0}\right).
\end{equation} 
Since $\widetilde{\psi}(\bm{x}_0)>0$ from the above assumption (b),  the claim (W1) in Lemma~\ref{locrlct} implies that the right-hand side is precisely $\R\textrm{-}LCT_{P, \bm{x}_0}\left(K\right)$. Given that the map $\psi$ meets the support restriction condition~\eqref{suppproper}, the assumption $(\mathcal{H}_3)$ yields 
\begin{align}\label{kjh}
\R\textrm{-}LCT_{P}\left(K\right) \;=\; \R\textrm{-}LCT_{P}(K, \psi) \; &\underset{\eqref{rlctminpsi}}{=}\; \min_{\bm{x}_0\in S(\psi)} \R\textrm{-}LCT_{P, \bm{x}_0}\left(K, \psi\cdot \pi_{\bm{x}_0}\right)\nonumber \\
& \underset{\eqref{comparlctco0}}{\succcurlyeq}\; 
\min_{\bm{x}_0\in S(\psi)} \R\textrm{-}LCT_{P, \bm{x}_0}\left(K, \widetilde{\psi}\cdot \pi_{\bm{x}_0}\right)\nonumber \\
&=\; \min_{\bm{x}_0\in S(\psi)} \R\textrm{-}LCT_{P, \bm{x}_0}\left(K\right).
\end{align}
From the above assumption (a), upon testing the distribution in~\eqref{dzfzgeh} against $\widetilde{\psi}$ instead of $\psi$, the equation~\eqref{rlctminpsi} remains true with $\widetilde{\psi}$ inside the real log-canonical thresholds. From the assumption $(\mathcal{H}_3)$ again, one thus has that $$ \R\textrm{-}LCT_{P}\left(K\right)\;=\; \R\textrm{-}LCT_{P}(K, \widetilde{\psi})  \;=\; \min_{\bm{x}_0\in S(\psi)} \R\textrm{-}LCT_{P, \bm{x}_0}\left(K, \widetilde{\psi}\cdot \pi_{\bm{x}_0}\right).$$Combined with the equation~\eqref{kjh}, this finally establishes that when  $Y=K$, it indeed holds that $$ \R\textrm{-}LCT_{P}(K, \psi) \;=\; \min_{\bm{x}_0\in S(\psi)} \R\textrm{-}LCT_{P, \bm{x}_0}\left(K\right).$$\\

Showing that the same relation is valid in the case where $Y=\Omega$ requires a different and more elaborate argument. To this end, the first goal is to show that there exists $\bm{x}$ such that  $\psi\left(\bm{x}\right)>0$  satisfying the property that
\begin{equation}\label{comparlctco}
\R\textrm{-}LCT_{P, \bm{x}_0}\left(\Omega, \psi\cdot \pi_{\bm{x}_0}\right)\;\succcurlyeq\; \R\textrm{-}LCT_{P, \bm{x}}\left(\Omega, \varphi_{\bm{x}}\right).
\end{equation} 
In this relation, it is required that $\varphi_{\bm{x}}\ge 0$ should be a test function supported in a small enough neighbourhood of $\bm{x}$ such that $\varphi_{\bm{x}}(\bm{x})>0$. This can be obtained from the homogeneity of the polynomial $P(\bm{x})$. Indeed, under this assumption,   for any real $\sigma<r_{P, \bm{x}_0}(\Omega)$, 
\begin{equation}\label{relzetaloc}
\left\langle \zeta_{P}\left(\Omega, \sigma\right), \psi\cdot \pi_{\bm{x}_0}\right\rangle \;\ll\; \left\langle \zeta_{P}\left(\Omega, \sigma\right), \pi_{\bm{x}_0}\right\rangle \; \ll\;  \left\langle \zeta_{P}\left(\Omega, \sigma\right), \pi_{\bm{x}_0, \lambda}\right\rangle.
\end{equation}
The second inequality is  here a consequence of the change of variables $\bm{y}= \bm{x}/\lambda$  upon setting $\pi_{\bm{x}_0, \lambda}~:\bm{x}\in\R^n\mapsto \pi_{\bm{x}_0}\left(\lambda\bm{x}\right)$. The parameter $\lambda$ is chosen  to be large enough  so that $\psi(\bm{x}_0/\lambda)>0$. 
Setting  $\bm{x}=\bm{x}_0/\lambda$  and $\varphi_{\bm{x}} = \pi_{\bm{x}_0, \lambda}$, the inequality between the first and the last terms in~\eqref{relzetaloc} then implies the one in~\eqref{comparlctco}. In fact, one may even assume that $\varphi_{\bm{x}}$ is supported in an arbitrarily small neighbourhood of $\bm{x}$. To see this, it is enough to note that if $\pi'_{\bm{x}_0}\ge 0$ is any smooth test function taking a  strictly positive value at $\bm{x}_0$, then the claim (W1) in Lemma~\ref{locrlct} implies that the smallest real poles of the meromorphic maps $s\mapsto \left\langle \zeta_{P}\left(\Omega, \sigma\right), \pi_{\bm{x}_0}\right\rangle$ and $s\mapsto \left\langle \zeta_{P}\left(\Omega, \sigma\right), \pi_{\bm{x}_0}\cdot \pi'_{\bm{x}_0}\right\rangle$ coincide, and so do their respective multiplicities. Therefore, for any real $\sigma$ strictly less than the common value $r_{P, \bm{x}_0}(\Omega)$ taken by these smallest poles, $$ \left\langle \zeta_{P}\left(\Omega, \sigma\right), \pi_{\bm{x}_0}\right\rangle \; \ll\; \left\langle \zeta_{P}\left(\Omega, \sigma\right), \pi_{\bm{x}_0}\cdot \pi'_{\bm{x}_0}\right\rangle \; \ll\; \left\langle \zeta_{P}\left(\Omega, \sigma\right), \pi_{\bm{x}_0, \lambda}\cdot \pi'_{\bm{x}_0, \lambda}\right\rangle,$$where $\pi'_{\bm{x}_0, \lambda}$ is defined in the same way as $\pi_{\bm{x}_0, \lambda}$. One may then set $\varphi_{\bm{x}} = \pi_{\bm{x}_0, \lambda}\cdot \pi'_{\bm{x}_0, \lambda}$ upon choosing adequately the support of $\pi'_{\bm{x}_0}$. \\

To conclude the proof of the relation~\eqref{rlctsemiglob} when $Y=\Omega$, consider a small enough neighbourhood $\Omega_{\bm{x}}$ of $\bm{x}$ on which $\psi$ remains positive and where the conclusions of Lemma~\ref{locrlct} hold. Working with a partition of unity $\left\{\pi_{\bm{x}_0}\right\}_{\bm{x}_0\in S(\psi)}\cup\left\{\pi_{\bm{x}}\right\}$ subordinate to the cover $\left\{\Omega_{\bm{x}_0}\right\}_{\bm{x}_0\in S(\psi)}\cup\left\{\Omega_{\bm{x}}\right\}$ of the support of $\psi$, the same argument as the one used to obtain~\eqref{rlctminpsi} implies that $$ \R\textrm{-}LCT_{P}(\Omega, \psi)\;=\; \min\left\{\R\textrm{-}LCT_{P, \bm{x}}\left(\Omega, \psi\cdot \pi_{\bm{x}}\right), \min_{\bm{x}_0\in S(\psi)} \R\textrm{-}LCT_{P, \bm{x}_0}\left(\Omega, \psi\cdot \pi_{\bm{x}_0}\right)\right\}.$$ Combining this with~\eqref{rlctminpsi} thus yields that  
\begin{equation}\label{rlctxpipsi}
\R\textrm{-}LCT_{P, \bm{x}}\left(\Omega, \psi\cdot \pi_{\bm{x}}\right)\;\succcurlyeq \; \R\textrm{-}LCT_{P}(\Omega, \psi).
\end{equation}

Since, by assumption,  $\pi_{\bm{x}}(\bm{x})>0, \psi(\bm{x})>0$ and $\varphi_{\bm{x}}(\bm{x})>0$, one infers from the claim (W1) in Lemma~\ref{locrlct} and from the relations~\eqref{comparlctco} and \eqref{rlctxpipsi}  that, at any point $\bm{x}_0$ where $\psi(\bm{x}_0)=0$, $$\R\textrm{-}LCT_{P, \bm{x}_0}\left(\Omega, \psi\cdot \pi_{\bm{x}_0}\right)\;\succcurlyeq \; \R\textrm{-}LCT_{P, \bm{x}}\left(\Omega\right) \;\succcurlyeq \; \R\textrm{-}LCT_{P}\left(\Omega, \psi\right).$$The minimum in the equation~\eqref{rlctminpsi} is thus achieved at those points where $\psi$ takes a positive value. At such points, the local log--canonical threshold can be determined by resolution of singularities from Lemma~\ref{locrlct}. This establishes (X1). $$\quad$$

As for (X2), that the pair $\R\textrm{-}LCT_{P}(Y, \psi)$ should be independent of the choice of the test function $\psi\ge 0$ under the assumptions of the statement is immediate from the relations~\eqref{rlctsemiglob}. \\

Assume that for such a map $\psi$, it holds that $\R\textrm{-}LCT_{P}(Y, \psi)=  \R\textrm{-}LCT_{P}(Y)$ and denote this pair by $\R\textrm{-}LCT_{P}(Y)=\left(r_P(Y), m_P(Y)\right)$.  Recall that the argument based on the Triangle Inequality in~\eqref{ineqlast-} shows that one may assume without loss of generality that the log-canonical thresholds of the various zeta distributions under consideration are achieved when tested against a nonnegative map in $\mathcal{C}_{c}^{\infty}(\R^n)$. Then, when $Y=K$,  it is clear that the assumption $(\mathcal{H}_3)$ implies that  $r_P(K)$ is the smallest real pole of the distribution $\zeta_{P}\left(Y, \,\cdot\,\right)$ when it is tested against all nonnegative smooth test functions satisfying the support restriction condition~\eqref{suppproper}, and that $m_P(K)$ is the corresponding order.\\

When $Y=\Omega$, it follows from the homogeneity of the polynomial $P(\bm{x})$ that the pair $\R\textrm{-}LCT_{P}(\Omega)$ is also the log-canonical threshold of the restriction of the distribution $\zeta_P\left(\Omega, \;\cdot\;\right) $ to any neighbourhood of the origin (this is, here again, a simple consequence of the relation~\eqref{inhomgchgevar}). That  $r_P(\Omega)$ should be the smallest real pole of the distribution $\zeta_{P}\left(\Omega, \,\cdot\,\right)$  and that $m_P(K)$ should be the corresponding order is then easily deduced from the claim (W3) in Lemma~\ref{locrlct}. This establishes (X2).$$\quad$$

To prove (X3), fix $\psi\ge 0$ in $\mathcal{C}_c^\infty(\R^n)$ such that $\psi(\bm{0})>0$ when $Y=\Omega$ and such that the support restriction~\eqref{suppproper} is met when $Y=K$. From the claim (X2), the meromorphic map $s\mapsto \left\langle \zeta_P(Y, s),\psi \right\rangle$  achieves the log-canonical threshold $\R\textrm{-}LCT_{P}(Y)=\left(r_P(Y), m_P(Y)\right)$. Note then that there exists a differential operator, say $\mathcal{D}_{Y}\left(\bm{x}, s, \bm{\partial}\right)$, such that 
\begin{equation}\label{globallocalop}
B_{P}(Y, -s)\cdot \left\langle \zeta_{P}\left(Y, s\right), \psi\right\rangle\; =\;  \left\langle \zeta_{P}\left(Y, s-1\right), \mathcal{D}_{Y}^*\left(\bm{x}, -s, \bm{\partial}\right)\psi\right\rangle,
\end{equation}
where the polynomial $B_P\left(Y, s\right)$ is defined in~\eqref{globalsbomegabis} and where $\mathcal{D}_{Y}^*\left(\bm{x}, s, \bm{\partial}\right)$ denotes here also the dual operator. To see this, let $\left\{\sigma_{\bm{x}_0}\right\}_{\bm{x}_0\in I\left(\overline{Y}\right)}$ be a locally finite partition of unity subordinate to the open cover $\left\{\mathcal{U}_{\bm{x}_0}\right\}_{\bm{x}_0\in I\left(\overline{Y}\right)}$ of  the topological closure $\overline{Y}$ of $Y$, where the sets $\mathcal{U}_{\bm{x}_0}$ result from Lemma~\ref{lemmeintersb}. For each $\bm{x}_0\in I\left(\overline{Y}\right)$, there exists a differential operator $\mathcal{D}_{Y, \bm{x}_0}\left(\bm{x}, s, \bm{\partial}\right)$ acting on $Y\cap\mathcal{U}_{\bm{x}_0}$ such that 
\begin{equation}\label{equasblocome}
B_{P, \bm{x}_0}\left(Y, s\right) \cdot P(\bm{x})^s\;=\; \mathcal{D}_{Y, \bm{x}_0}\left(\bm{x}, s, \bm{\partial}\right)  P(\bm{x})^{s+1}.
\end{equation}
Here,  the local polynomial $B_{P, \bm{x}_0}\left(Y, s\right)$ defined in~\eqref{globalsbomega} divides the global one $B_{P}(Y, s)$; consequently, there exists some polynomial $T_{P, \bm{x}_0}\left(Y, s\right)$ with rational coefficients such that $B_{P}(Y, s)=T_{P, \bm{x}_0}\left(Y, s\right)\cdot B_{P, \bm{x}_0}\left(Y, s\right)$. Choose the partition of unity $\left\{\sigma_{\bm{x}_0}\right\}_{\bm{x}_0\in I\left(\overline{Y}\right)}$ in such a way that it sums to the constant function 1 in a small enough neighbourhood of $Y$, and assume without loss of generality that the support of $\psi$ is contained in this neighbourhood. It then follows from the relation~\eqref{equasblocome} that 
\begin{align*}
B_{P}(Y, -s)\cdot \left\langle \zeta_{P}\left(Y, s\right), \psi\right\rangle\; &=\; \sum_{\bm{x}_0\in I\left(\overline{Y}\right)} T_{P, \bm{x}_0}\left(Y, -s\right)\cdot B_{P, \bm{x}_0}\left(Y, -s\right)\cdot \left\langle \zeta_{P}\left(Y, s\right), \sigma_{\bm{x}_0}\cdot \psi\right\rangle\\
&= \sum_{\bm{x}_0\in I\left(\overline{Y}\right)} T_{P, \bm{x}_0}\left(Y, -s\right)\left\langle \zeta_{P}\left(Y, s-1\right), \mathcal{D}_{Y, \bm{x}_0}^*\left(\bm{x}, -s, \bm{\partial}\right) \left( \sigma_{\bm{x}_0}\cdot \psi\right)\right\rangle. 
\end{align*}
In view of the formula defining the dual operator in~\eqref{opdual}, it is thus enough to let $$\mathcal{D}_{Y}\left(\bm{x}, s, \bm{\partial}\right)\;=\;  \sum_{\bm{x}_0\in I\left(\overline{Y}\right)} T_{P, \bm{x}_0}\left(Y, s\right)\cdot\sigma_{\bm{x}_0}\left(\bm{x}\right)\cdot \mathcal{D}_{Y, \bm{x}_0}\left(\bm{x}, s, \bm{\partial}\right)$$ so that the identity~\eqref{globallocalop} holds.\\

Set then $s=\sigma$ to be a real number and let $\sigma<r_P(Y)$ tend to $r_P(Y)$ in the identity~\eqref{globallocalop}~: since the right--hand side  is holomorphic at $\sigma=r_P(Y)$, the claim (X3) follows.\\

This completes the proof of the Proposition~\ref{proprlctgene}.
\end{proof}

\section[Parameters of the Volume Estimate]{Parameters of the Volume Estimate under Additional Regularity and Homogeneity Assumptions}

The goal in this concluding section of the chapter is to establish the few remaining claims to complete the proof of the main Theorem~\ref{volestim}. To begin with, the following proposition establishes under more general assumptions the claim on the localisation of the smallest real pole of the zeta distribution $ \zeta_{\bm{F}}$ made in Theorem~\ref{volestim}. 

\begin{prop}\label{locsmalpole} 
Let $\mathfrak{O}'$ be a subcollection of the collection of sets $\mathfrak{O}$ defined in~\eqref{defoobig}. Assume that the corresponding distribution $\zeta_P\left(\mathfrak{O}', \,\cdot\,\right)$ introduced in~\eqref{def123} is well--defined. Then, under the assumption of the homogeneity of the polynomial $P(\bm{x})$, it holds that the smallest real pole  $r_P\left(\mathfrak{O}'\right)$ of $\zeta_P\left(\mathfrak{O}', \,\cdot\,\right)$ lies in the interval $\left(0, n/q\right]$, where $q$ denotes the degree of $P(\bm{x})$. 
\end{prop}

\begin{proof}
From Lemma~\ref{lempolemulti}, it is enough to establish the statement for the quantity $r_P(\Omega)$ for any fixed $\Omega\in\mathfrak{O}'$ such that the distribution $\zeta_P\left(\Omega, \,\cdot\,\right)$ is well-defined. To this end, consider a nonnegative, smooth and compactly supported function $\xi$ defined over the real line and taking a positive value at the origin. Define a radial function $\psi$ in $\mathcal{C}_{c}^{\infty}(\R)$ by setting $\psi\left(\bm{x}\right)=\xi\left(\left\|\bm{x}\right\|\right)$ for all $\bm{x}\in\R^n$. Then, denoting by $\tau_{n-1}$ the induced Lebesgue measure on the $n$--dimensional sphere $\Sph^{n-1}\subset\R^n$, it holds that for any real $\sigma<r_P(\Omega)$, 
\begin{align*}
\left\langle \zeta_{P}\left(\Omega, \sigma\right), \psi\right\rangle \; &=\;  \int_{\R^n}\frac{\psi(\bm{x})}{P_{\Omega}(\bm{x})^{\sigma}}\cdot\textrm{d}\bm{x}\\
&=\; \left( \int_{\R}\xi(r)\cdot r^{n-1-q\sigma}\cdot \textrm{d}r\right)\cdot\left(\int_{\Sph^{n-1}} P_{\Omega}(\bm{u})^{-\sigma}\cdot \textrm{d}\tau_{n-1}\left(\bm{u}\right)\right). 
\end{align*}
Let $\sigma_0>0$ be the smallest real pole of the second factor on the right--hand side of this relation. Then, clearly, $r_P(\Omega)\le \min\left\{n/q, \sigma_0\right\}\le n/q$ under the assumption that $\xi(0)>0$. The sought conclusion then follows from the point (X2) in Proposition~\ref{proprlctgene}.
\end{proof}

Recall that $\bm{F}(\bm{x})=\left(F_1(\bm{x}), \cdots, F_p(\bm{x})\right)$ is assumed to be a set of $p\ge 1$ homogeneous forms in $n\ge 2$ variables with common degree $d\ge 2$. Set \begin{equation}\label{normpolynsq}
P_{\bm{F}}(\bm{x})\;=\;\left\|\bm{F}(\bm{x})\right\|^2\;=\; \sum_{i=1}^{p} F_i(\bm{x})^2.
\end{equation}

\begin{prop}[The real log-canonical threshold in the case of a smooth complete intersection]\label{smoothcase}
With the above notations, assume that the set of homogeneous forms $\bm{F}(\bm{x})$  has a smooth complete intersection over a domain $K\subset\R^n$ satisfying the conditions $(\mathcal{H}_1)$ and $(\mathcal{H}_2)$. Suppose furthermore that, if the condition $(\mathcal{H}_3)$ is not assumed, then the algebraic variety $\mathcal{Z}_{\R}(\bm{F})$ defined in~\eqref{defalgvar} intersects non trivially the interior of $K$. 

Then, the condition $(\mathcal{H}_3)$ holds in all cases. Furthermore, the real log-canonical threshold of the polynomial $P_{\bm{F}}(\bm{x})$ over the domain $K$ is well-defined and equals $$\R\textrm{-}LCT_{P_{\bm{F}}}(K)=\left(\frac{p}{2}, 1\right).$$
\end{prop}

\begin{proof} In the case that it is just assumed that $\mathcal{Z}_{\R}(\bm{F})$ intersects non trivially the interior of $K$, define in the notation of the support restriction condition~\eqref{suppproper} the open set $U$ as an open neighbourhood of $K$ where the map~\eqref{mapsmoothcominter} does not vanish. Let also the set $C$ be any compact subset of the interior of $K$ with an interior intersecting $\mathcal{Z}_{\R}(\bm{F})$ non trivially. Fix then a smooth test function $\psi\ge 0$ meeting the support restriction condition~\eqref{suppproper} for this choice of the sets $C$ and $U$. The goal is to show that the smallest real pole of the map 
\begin{equation}\label{zetPF}
s\mapsto \left\langle \zeta_{P_{\bm{F}}}(s), \psi\right\rangle
\end{equation} 
and the corresponding order do not depend on the choice of $\psi$, and also that they equal the values $p/2$ and 1, respectively. \\

In the case that  the condition $(\mathcal{H}_3)$ is assumed, it may furthermore be imposed without loss of generality, even if it means making it slightly bigger, that the support of $\psi\ge 0$   intersects non trivially the variety $\mathcal{Z}_{\R}(\bm{F})$ (this can be done since the condition $(\mathcal{H}_2)$ is assumed). In view of the point (X2) in Proposition~\ref{proprlctgene}, the goal is then to prove that the real log-canonical threshold is again the pair $(p/2, 1)$ by establishing this claim for the above--defined  map~\eqref{zetPF}.\\

To prove these claims, note that the smooth complete intersection assumption implies that the function $$\bm{x}=\left(x_1, \dots, x_n\right)\mapsto \left(F_1(\bm{x}), \dots, F_p(\bm{x}), x_{p+1}, \dots , x_n\right)$$realises a $\mathcal{C}^1$-diffeomorphism between $K$ and its image (its Jacobian is $\pm\left\|\bigwedge_{i=1}^p \nabla F_i(\bm{x})\right\|^2$). As a consequence, the corresponding change of variables yields that for any real $\sigma$ strictly smaller than the smallest real pole of the map~\eqref{zetPF}, 
\begin{align*}
\left\langle \zeta_{P_{\bm{F}}}(\sigma), \psi\right\rangle\;=\; \int_{\R^n}\frac{\psi(\bm{x})}{\left\|P_{\bm{F}}(\bm{x})\right\|^{2\sigma}}\cdot\textrm{d}\bm{x}\;\underset{\eqref{normpolynsq}}{=}\;\int_{\R^p} \frac{\widetilde{\psi}_{\bm{F}}(\bm{y})}{\left\|\bm{y}\right\|^{2\sigma}}\cdot\textrm{d}\bm{y},
\end{align*}
where $\widetilde{\psi}_{\bm{F}}$ is a smooth nonnegative test function taking a strictly positive value at the origin (this follows from the assumption that the support of $\psi$ intersects non trivially the variety $\mathcal{Z}_{\R}(\bm{F})$ from the above discussion). Thus, working in polar coordinates, one obtains the existence of a nonnegative map $\xi_{\psi}^{\bm{F}}$ in $\mathcal{C}_c^\infty(\R)$ taking a strictly positive value at the origin such that 
\begin{align*}
\left\langle \zeta_{P_{\bm{F}}}(\sigma), \psi\right\rangle\;&=\; \int_{\R}\xi_{\psi}^{\bm{F}}(r)\cdot r^{p-1-2\sigma}\cdot \textrm{d}r.
\end{align*}

It is clear that the smallest real pole of the map defined by the right-hand side is $p/2$ and that it has order 1, regardless of the choice of the smooth compactly supported map $\xi_{\psi}^{\bm{F}}\ge 0$ such that $\xi_{\psi}^{\bm{F}}(0)>0$. This completes the proof of the proposition.
\end{proof}

\begin{proof}[Completion of the proof of Theorem~\ref{volestim}] The final goal in this chapter is to establish the few remaining claims in the theorem which have not already been proved. To this end, let  $P_{\bm{F}}(\bm{x})$ be as in~\eqref{normpolynsq}. Denote by $\zeta_{P_{\bm{F}}}$ the corresponding distribution as defined in~\eqref{zetaP} and recall that $\zeta_{\bm{F}}$ denotes the distribution introduced in~\eqref{distrizetaR}. It is clear from these definitions that when tested against any smooth compactly supported function, the smallest pole of the former distribution is obtained by dividing by 2  the smallest pole of the latter, and that the orders coincide. \\

In view of this observation, it is immediate from Proposition~\ref{smoothcase} that in the case of a smooth complete intersection, the real log canonical threshold of the distribution $\zeta_{\bm{F}}$ is, in the notations of the theorem, the pair $(p,1)$.\\

Similarly, to prove that the pairs $(r_{\bm{F}}(K),m_{\bm{F}}(K))$ and $(r_{\bm{F}},m_{\bm{F}})$ can be determined by resolving singularities in a finite number of steps and that their first components are strictly positive rationals, it is enough to prove this statement for the real log-canonical thresholds $\R\textrm{-}LCT_{P_{\bm{F}}}(K)$ and $\R\textrm{-}LCT_{P_{\bm{F}}}$ (i.e.~for the polynomial $P_{\bm{F}}(\bm{x})$). The claim then immediately follows from the points (X1) and (X2) in Proposition~\ref{proprlctgene} as far as $\R\textrm{-}LCT_{P_{\bm{F}}}(K)$ is concerned. As for  $\R\textrm{-}LCT_{P_{\bm{F}}}$, one needs to additionally consider the equation~\eqref{rlctpolezetaprim} to reach the same conclusion. \\

To prove that the rationals $-r_{\bm{F}}(K)/2=-r_{P_{\bm{F}}}(K)$ and $-r_{\bm{F}}/2=-r_{P_{\bm{F}}}$ are roots of the  Sato-Bernstein polynomial $B_{\bm{F}}(s)$ with respective multiplicities at most $m_{\bm{F}}(K)=m_{P_{\bm{F}}}(K)$ and $m_{\bm{F}}=m_{P_{\bm{F}}}$, apply  the point (X3) in Proposition~\ref{proprlctgene} in the case of the homogeneous form $P_{\bm{F}}(\bm{x})$. Then,
\begin{itemize}
\item the claim for the pair $\left(-r_{P_{\bm{F}}}(K), m_{P_{\bm{F}}}(K)\right)$ follows from the observation that the global Sato-Bernstein polynomial $B_{P_{\bm{F}}}\left(K, s\right)$ with respect to the domain $K$ divides the Sato-Bernstein polynomial $B_{P_{\bm{F}}}\left( s\right)$ of the homogeneous form $P_{\bm{F}}(\bm{x})$;

\item in the case of the pair $\left(-r_{P_{\bm{F}}}, m_{P_{\bm{F}}}\right)$, it is implied by the same division property satisfied by the polynomial $B_{P}(\Omega, s)$, where $\Omega\in\mathfrak{O}$, with the additional need to include the conclusion of Lemma~\ref{lempolemulti} .
\end{itemize}

To complete the proof of Theorem~\ref{volestim}, it remains to notice that Proposition~\ref{locsmalpole} applied to the homogeneous form $P_{\bm{F}}(\bm{x})$ of degree $q=2d$ and the relation $r_{\bm{F}}/2=r_{P_{\bm{F}}}$ yield that the rational $r_{\bm{F}}$ lies in the interval $(0, n/d]$. 
\end{proof}

\chapter[Log--Canonical Threshold and the Sato--Bernstein Polynomial]{Log--Canonical Threshold, Roots of the Sato--Bernstein Polynomial and Related Multiplicities}\label{lctrootpole} 

\vspace{20mm}

The goal of this section is to establish Theorem~\ref{lctordermultplicity}. Let $P(\bm{z})\in\C[\bm{z}]$ be a (non--necessarily homogeneous) complex polynomial in $n$--variables $\bm{z}=\left(z_1, \dots, z_n\right)$. Define its \emph{(complex) log--canonical threshold} as the pair 
\begin{equation*}\label{clct}
\C\textrm{-}LCT_P\;=\;\left(\hat{r}_P, \hat{m}_P\right),
\end{equation*}
where $\hat{r}_P$ is the smallest real pole of the complex zeta distribution $\widehat{\zeta}_{P}$ defined in~\eqref{zetadistribcompl} and where $ \hat{m}_P$ is its order (in accordance with the previous section, this slightly extends the definition of the complex log--canonical threshold given in \S 1.3 of the Introduction, where only the quantity $\hat{r}_P$ was considered). The existence of the pair $\C\textrm{-}LCT_P$ follows from the fact that $\widehat{\zeta}_{P}$ admits a meromorphic extension to the entire complex plane. To see this, let $B_p(s)\in\Q[s]$ be the Sato--Bernstein polynomial of $P(\bm{z})$ and let $\mathcal{D}(\bm{z}, s, \bm{\partial})$ be a differential operator such that the identity~\eqref{sbrelation} holds with $B_p(s)$. If $\mathcal{D}(\bm{z}, s, \bm{\partial})$ is expanded as in~\eqref{expdiffop}, define its conjugate by setting, with obvious notations, $$\overline{\mathcal{D}}(\bm{z}, s, \bm{\partial})\;=\; \sum_{\bm{\alpha}\in\N_0^n} \overline{a_{\bm{\alpha}}(\bm{z},s)}\cdot \overline{\bm{\partial}}^{\bm{\alpha}}.$$ It is then easily seen that the operator $\mathcal{D}(\bm{z}, s, \bm{\partial})$ commutes with its conjugate $\overline{\mathcal{D}}(\bm{z}, s, \bm{\partial})$. Moreover, $$\left(\mathcal{D}(\bm{z}, s, \bm{\partial}) P(\bm{z})^{s+1}\right)\cdot\left(\overline{\mathcal{D}}(\bm{z}, s, \bm{\partial})  \overline{P(\bm{z})}^{s+1}\right)\;=\; \left(\mathcal{D}(\bm{z}, s, \bm{\partial})\overline{\mathcal{D}}(\bm{z}, s, \bm{\partial})\right) \left|P(\bm{z})\right|^{2(s+1)},$$ which implies that for any $\psi$ in $\mathcal{C}_{c}^{\infty}(\C^n)$,
\begin{equation}\label{sbzetameorextcompl}
\left|B_{P}(s)\right|^2\cdot \left\langle \widehat{\zeta}_{P}\left(s\right), \psi\right\rangle\; =\;  \left\langle \zeta_{P}\left(s-1\right), \left(\mathcal{D}(\bm{z}, s, \bm{\partial})\overline{\mathcal{D}}^*(\bm{z}, s, \bm{\partial})\right) \psi\right\rangle.
\end{equation}
The meromorphic continuation of the distribution $\widehat{\zeta}_{P}$ then follows in the same way as in the proof of Lemma~\ref{propzetamerogro}.\\

It is known, see, e.g., \cite[\S\S 10.6 \&  10.7]{kollar} and~\cite[\S 3]{saito}, that the smallest pole $\hat{r}_P$ is a well--defined rational number such that 
\begin{equation}\label{rhatineeq}
0\;<\; \hat{r}_P\;\le\; 1. 
\end{equation}
In the case that $P(\bm{z})$ is a polyomial \emph{with real coefficients}, it always holds that $$r_P\;\ge\; \hat{r}_P,$$where $r_P$ is the first component of the \emph{real} log--canonical threshold defined in~\eqref{rlctpolezeta}. See~\cite{saito} for a justification of this claim, where it is also shown that the inequality may be strict, and that it may also happen  that $r_P>1$ in the case that $\mathcal{Z}_\R(P)$ is an algebraic variety of dimension strictly less than $n-1$. It is nevertheless the case that the (global) complex and real Sato--Bernstein polynomials of $P(\bm{z})$ which, in the notation of~\eqref{locglobsb} and~\eqref{globalsbomegabis}, are respectively $B_P(s)$ and $B_P(\R^n, s)$, coincide. To see this, note first that the complexification of a relation such as~\eqref{sbrelation} valid  for real variables shows that $B_P(s)$ divides $B_P(\R^n, s)$. Conversely, if such a relation is true for complex variables, then restricting to real variables and taking the real part of each side of the equation shows that $B_P(\R^n, s)$ divides $B_P(s)$ when the polynomial $P(\bm{z})$ is assumed to have real coefficients, whence the claim.\\

From the just established relation $B_P(s) = B_P(\R^n, s)$ and from the fact that $r_P$  can be strictly larger than $\hat{r}_P$, one deduces there can be no real analogue of Theorem~\ref{lctordermultplicity} when considering the real zeta function of a polynomial $P(\bm{x})\in\R[\bm{x}]$.\\

Revert to the general case where $P(\bm{z})\in\C[\bm{z}]$. The first step in the proof of Theorem~\ref{lctordermultplicity} is the following complex analogue of Lemma~\ref{lemmonorlct}. Given $Z\subset\C^n$ and an element $\psi$ in  $\mathcal{C}_{c}^{\infty}(\C^n)$, define, whenever it exists, $\C\textrm{-}LCT_P(Z| \psi)$ as the pair $\left(\hat{r}(\psi), \hat{m}(\psi)\right)$, where $\hat{r}(\psi)$ is the smallest real pole of the map defined for $\textrm{Re}(s)\le 0$ by $s\mapsto \int_Z\left|P(\bm{z})\right|^{-2s}\cdot\psi(\bm{z}, \overline{\bm{z}})\cdot\textrm{d}\bm{z}\wedge\textrm{d}\overline{\bm{z}}$ and where $\hat{m}(\psi)$ is its order. Let $\C\textrm{-}LCT_P(Z| \psi)=(\infty, -)$ if this is undefined. When $Z=\C^n$, set for the sake of simplicity of notation $\C\textrm{-}LCT_{P}\left(\C^n, \psi\right)= \C\textrm{-}LCT_P(\psi)$.

\begin{lem}\label{lemclctloc}
Assume that $\psi$ is an element in $\mathcal{C}_{c}^{\infty}(\C^n)$ not vanishing in some neighbourhood of the origin. Let $\bm{k}=(k_1, \dots, k_n)\in\N_0^n$ and  $\bm{h}=(h_1, \dots, h_n)\in\N_0^n$ be $n$--tuples of non--negative integers. Then, with the same convention as in  Lemma~\ref{lemmonorlct} (and with a slight abuse of notations), $$\C\textrm{-}LCT_{|\bm{z}^{\bm{k}}|}\left(\, |\bm{z}^{2\bm{h}}|\cdot\psi\right)\;=\; \left(\hat{r}, \hat{m}\right),$$ where $$ \hat{r}\;=\; \min_{1\le i\le n}\frac{h_i+1}{k_i} \qquad \qquad  \textrm{and}\qquad \qquad  \hat{m}\;=\; \#\left\{1\le i\le n\; :\; \frac{h_i+1}{k_i}\;=\; r\right\}.$$ 
\end{lem}

\begin{proof}
Decompose each coordinate $z_i$ ($1\le i \le n$) of the point $\bm{z}\in\C^n\backslash \{\bm{0}\}$ in polar form as $z_i=\rho_i e^{i\theta_i}$, where $\rho_i>0$ and $\theta_i\in [0, 2\pi)$. Set $Y_n^+=(0, \infty)^n$, $\bm{\rho}=\left(\rho_1, \dots, \rho_n\right)\in Y_n^+$ and $\bm{\theta}=\left(\theta_1, \dots, \theta_n\right)\in [0, 2\pi)^n$ and note that $\textrm{d}\bm{z}\wedge\textrm{d}\overline{\bm{z}}=\bm{\rho}\cdot \textrm{d}\bm{\rho}\wedge \textrm{d}\bm{\theta}$. Then, for $s\in\C$ such that the integrals exist, 
\begin{equation}\label{polarlct}
\int_{\C^n}\left|\bm{z}^{\bm{k}}\right|^{-2s}\cdot \left|\bm{z}^{2\bm{h}}\right|\cdot\psi(\bm{z}, \overline{\bm{z}})\cdot \textrm{d}\bm{z}\wedge\textrm{d}\overline{\bm{z}}\;=\; \int_{Y_n^+}\bm{\rho}^{2\bm{h}+\bm{1}-2s\bm{k}}\cdot \tilde{\psi}(\bm{\rho})\cdot\textrm{d}\bm{\rho}.
\end{equation} 
Here, $\bm{1}=\left(1, \dots, 1\right)\in\N^n$
and $\tilde{\psi}$ is a complex valued bounded, smooth  function supported in the intersection of a compact neighbourhood of the origin with  the closure of $Y_n^+$ . Furthermore, the map $\tilde{\psi}$ does not vanish when $\bm{\rho}\in Y_n^+$ has small enough norm. Thus, taking the minimum with respect to the order introduced in~\eqref{ordering},
\begin{align*}
&\C\textrm{-}LCT_{|\bm{z}^{\bm{k}}|}\left(|\bm{z}^{\bm{h}}|\cdot\psi\right)\\ & \qquad\qquad = \; \min\left\{\R\textrm{-}LCT_{\bm{\rho}^{\bm{2k}}}\left(Y_n^+ | \bm{\rho}^{2\bm{h}+\bm{1}}\cdot\textrm{Re}\left(\tilde{\psi}\right)\right),\, \R\textrm{-}LCT_{\bm{\rho}^{\bm{2k}}}\left(Y_n^+ | \bm{\rho}^{2\bm{h}+\bm{1}}\cdot\textrm{Im}\left(\tilde{\psi}\right)\right)\right\},
\end{align*} 
where either $\textrm{Re}\left(\tilde{\psi}\right)(\bm{0})\neq 0$ or $\textrm{Im}\left(\tilde{\psi}\right)(\bm{0})\neq 0$ . The claim is therefore reduced to Lemma~\ref{lemmonorlct}.
\end{proof}

By analogy with the real case, define the \emph{local (complex) zeta distribution $\widehat{\zeta}_{P, \bm{z}_0}$ at $\bm{z}_0\in\C^n$} by setting $\left\langle \widehat{\zeta}_{P, \bm{z}_0}, \varphi_{\bm{z}_0}\right\rangle= \left\langle \widehat{\zeta}_{P}, \varphi_{\bm{z}_0}\right\rangle$ for any smooth  function $\varphi_{\bm{z}_0}$ compactly supported in a small enough open neigbourhood of $\bm{z}_0$. This neighbourhood is chosen such that a resolution of singularities map is well-defined over it when $P(\bm{z}_0)=0$ (cf.~Theorem~\ref{hironaka}) and such that the polynomial $P(\bm{z})$ does not vanish over it when $P(\bm{z}_0)\neq 0$.\\

Let furthermore $$\C\textrm{-}LCT_{P, \bm{z}_0}\;=\;\left(\hat{r}_{P, \bm{z}_0},\, \hat{m}_{P, \bm{z}_0}\right)$$denote the \emph{(complex) log--canonical threshold of the polynomial $P(\bm{z})$ at $\bm{z}_0$}~: it is here again the data of the smallest real pole $\hat{r}_{P, \bm{z}_0}$ of $\widehat{\zeta}_{P, \bm{z}_0}$ and of its order $\hat{m}_{P, \bm{z}_0}$, provided it is well--defined.

\begin{lem}\label{lemcomloct}
Let $\bm{z}_0\in\C^n$. If $P(\bm{z}_0)=0$, let $\hat{g}~:  \widehat{\mathcal{M}}\rightarrow \widehat{\mathcal{W}}$ be an analytic morphism obtained by applying a resolution of singularities to $P(\bm{z})$ at $\bm{z}_0$, where $\widehat{\mathcal{M}}$ is a complex manifold and where  $\widehat{\mathcal{W}}\subset\C^n$ is a neighbourhood of $\bm{z}_0$. Then,
 \begin{itemize}
\item[(Y1)] when $P(\bm{z}_0)\neq0$, the pair $\C\textrm{-}LCT_{P, \bm{z}_0}= \left(\hat{r}_{P, \bm{z}_0},\, \hat{m}_{P, \bm{z}_0}\right)$ is not defined, and thus set to $\left(\infty, -\right)$. If, however, $P(\bm{z}_0)=0$, it is well--defined, in which case  the preimage $\hat{g}^{-1}\left(\left\{\bm{z}_0\right\}\right)$  is covered by finitely many charts $\left\{\widehat{\mathcal{M}}_{\bm{y}_0}\right\}_{\bm{y}_0\in \hat{I}(\bm{z}_0)}$. \\
\end{itemize} 

Assuming from now on that $P(\bm{z}_0)=0$, let $\hat{J}(\bm{z}_0)\subset \hat{I}(\bm{z}_0)$ be a maximal index set such that for each $\bm{y}_0\in \hat{J}(\bm{z}_0)$, the set $\widehat{\mathcal{M}}_{\bm{y}_0}$ has positive measure. Then,

\begin{itemize}
\item[(Y2)] for each index $\bm{y}_0\in \hat{J}(\bm{z}_0)$,  there exists an integer $\hat{K}(\bm{y}_0)\ge 1$ such that for any $\varphi_{\bm{z}_0}$ smooth and compactly supported in a small enough neighbourhood of $\bm{z}_0$, the meromorphic map $\left\langle \widehat{\zeta}_{P, \bm{z}_0}, \varphi_{\bm{z}_0}\right\rangle$ can be decomposed as the sum of $\hat{K}(\bm{y}_0)$ integrals of the form $$\int_{\C^n}\left|\bm{y}^{\bm{k}_j(\bm{y}_0)}\right|^{-2s}\cdot \left|\bm{y}^{2\bm{h}_j(\bm{y}_0)}\right|\cdot\psi_j^{(\bm{y}_0)}(\bm{y}, \overline{\bm{y}})\cdot \textrm{d}\bm{y}\wedge\textrm{d}\overline{\bm{y}}.$$ Here, the smooth maps $\psi_j^{(\bm{y}_0)}$ are supported in a neighbourhood of $\bm{y}_0$ for $1\le j\le \hat{K}(\bm{y}_0)$. Furthermore, the integer $n$--tuples  $$\bm{k}_j(\bm{y}_0)=\left(k_j^{(1)}(\bm{y}_0), \dots, k_j^{(n)}(\bm{y}_0)\right) \qquad \textrm{and}\qquad \bm{h}_j(\bm{y}_0)=\left(h_j^{(1)}(\bm{y}_0), \dots, h_j^{(n)}(\bm{y}_0)\right)$$are obtained from the monomial forms~\eqref{monoform} and~\eqref{monoformbis} of the resolution of singularities $\hat{g}$ restricted to $\widehat{\mathcal{M}}_{\bm{y}_0}$. 

\item[(Y3)] moreover, it holds that
\begin{equation}\label{formullct} 
\hat{r}_{P, \bm{z}_0}\;=\; \min_{\bm{y}_0\in \hat{J}(\bm{x}_0)}\;\min_{1\le j\le \hat{K}(\bm{y}_0)}\;\min_{1\le l\le n}\; \frac{h_j^{(l)}(\bm{y_0})+1}{k_j^{(l)}}
\end{equation} 
and 
\begin{equation}\label{formullctm}
\hat{m}_{P, \bm{z}_0}\;=\;\max_{\bm{y}_0\in \hat{J}(\bm{x}_0)}\;\max_{1\le j\le \hat{K}(\bm{y}_0)}\;\#\left\{1\le l\le n\; :\; \frac{h_j^{(l)}(\bm{y_0})+1}{k_j^{(l)}}\;=\; \hat{r}_{\bm{z}_0}(P)\right\}
\end{equation}
with the same conventions as in Lemma~\ref{lemmonorlct} to deal with the case of vanishing denominators. 
\end{itemize}
\end{lem}

The formula for $\hat{r}_{P, \bm{z}_0}$ above is known --- see, e.g., \cite[\S 10.7]{kollar}. The main novelty is the expression for the multiplicity $\hat{m}_{P, \bm{z}_0}$. 

\begin{proof}
Note first that the local distribution $\widehat{\zeta}_{P, \bm{z}_0}$ can be meromorphically continued  to the entire complex plane since the global distribution $\widehat{\zeta}_{P}$ can be so. The lemma then turns out to be a rephrasing of the results established as part of the proof of Lemma~\ref{locrlct} upon considering a complex resolution of singularities. The only change leading to the specific forms of the above statements (Y2) and (Y3) is that the equation~\eqref{eqtdeczeat} now becomes 
\begin{align*}
&\left\langle \widehat{\zeta}_{P, \bm{z}_0}(\sigma), \varphi_{\bm{z}_0}\right\rangle \\ &\qquad =\;  \sum_{\bm{y}_0\in \hat{J}(\bm{x}_0)} \int_{\widehat{\mathcal{M}}_{\bm{y}_0}} \left|P(g\left(\bm{y}\right))\right|^{-2\sigma}\cdot\varphi_{\bm{z}_0}\left(g(\bm{y}), \overline{g(\bm{y})}\right)\cdot \left|\textrm{Jac}_g(\bm{y})\right|^2\cdot \pi_{\bm{y}_0}\left(\bm{y}\right)\cdot\textrm{d}\bm{y}\wedge\textrm{d}\overline{\bm{y}},
\end{align*}
where $\left|\textrm{Jac}_g(\bm{y})\right|^2$ is the (real) Jacobian of the map of resolution of singularities con\-si\-de\-red as a smooth map of real $2n$--dimensional manifolds. This is enough to justify (Y2). The calculation of the complex log--canonical threshold in (Y3) is then easily seen to be reduced to the determination of the minimum over all indices $1\le j \le \hat{K}(\bm{y}_0)$ and over all $\bm{y}_0\in \hat{J}(\bm{z}_0)$ of the pair $\C\textrm{-}LCT_{|\bm{z}^{\bm{k}_j(\bm{y}_0)}|}\left(|\bm{z}^{2\bm{h}_j(\bm{y}_0)}|\cdot\psi\right)$ for some smooth function $\psi$ satisfying the assumptions of  Lemma~\ref{lemclctloc} (here again with respect to the order introduced in~\eqref{ordering}). \\

This concludes the proof of Lemma~\ref{lemcomloct}. 
\end{proof}

The following is a local version of  Theorem~\ref{lctordermultplicity} which is needed to prove it in full generality.

\begin{prop}\label{proplocthmlctsb}
Let $\bm{z}_0\in\C^n$ be such that $P(\bm{z}_0)=0$. Then, the log--canonical threshold $\C\textrm{-}LCT_{P, \bm{z}_0}=\left(\hat{r}_{P, \bm{z}_0}, \hat{m}_{P, \bm{z}_0}\right)$ is such that $-\hat{r}_{P, \bm{z}_0}$ is the largest root of the local polynomial $B_{P, \bm{z}_0}(s)$ with multiplicity \emph{exactly} $\hat{m}_{P, \bm{z}_0}$.
\end{prop}

\begin{proof}
The  statement is established with the help of two lemmata, the first one showing that the multiplicity under consideration is \emph{at least}  $\hat{m}_{P, \bm{z}_0}$ and the second one that it is \emph{at most} this same quantity.

\begin{lem}\label{lem1167}
Under the assumptions of the proposition, \sloppy the local log--canonical threshold $\C\textrm{-}LCT_{P, \bm{z}_0}=\left(\hat{r}_{P, \bm{z}_0},\, \hat{m}_{P, \bm{z}_0}\right)$ is such that $-\hat{r}_{P, \bm{z}_0}$ is the largest root of the local Sato--Bernstein polynomial $B_{P, \bm{z}_0}(s)$ with multiplicity \emph{at least} $\hat{m}_{P, \bm{z}_0}$.
\end{lem}

\begin{proof}
Let $\mathcal{D}_{\bm{z}_0}(\bm{z}, s, \bm{\partial})$ be a differential operator defined in a neighbourhood $Z_{\bm{z}_0}\subset\C^n$ of $\bm{z}_0$ such that for all $\bm{z}\in Z_{\bm{z}_0}$, 
\begin{equation}\label{eqt1}
B_{P, \bm{z}_0}(s) P(\bm{z})^{s}= \mathcal{D}_{\bm{z}_0}(\bm{z}, s, \bm{\partial})\cdot P(\bm{z})^{s+1}. 
\end{equation}

Let $\psi$ be in $\mathcal{C}_{c}^{\infty}(\C^n)$ such that the meromorphic map $s\mapsto  \left\langle\widehat{\zeta}_{P, \bm{z}_0}(s), \psi\right\rangle$ achieves the pair $\C\textrm{-}LCT_{P, \bm{z}_0}$. Since the dual operator $\mathcal{D}^*_{\bm{z}_0}(\bm{z}, s, \bm{\partial})$ is polynomial in $s$ and analytic in $\bm{z}$, the real $$M\left(\psi, \hat{r}_{P, \bm{z}_0}\right)\;=\; \max_{\underset{\bm{z}\in\textrm{Supp }\psi}{\left|\sigma-\hat{r}_{P, \bm{z}_0}\right|\le 1}} \left|\mathcal{D}^*_{\bm{z}_0}(\bm{z}, -\sigma, \bm{\partial})\psi\left(\bm{z}, \overline{\bm{z}}\right)\right|$$is well--defined. Fixing $\varphi\ge 0$ in $\mathcal{C}_{c}^{\infty}(\C^n)$ such that $\varphi\left(\bm{z}, \overline{\bm{z}}\right)\ge M\left(\psi, \hat{r}_{P, \bm{z}_0}\right)$ for all $\bm{z}\in\textrm{Supp } \psi$, one obtains that for all reals $\sigma\in\left(\hat{r}_{P, \bm{z}_0}-1,\,\hat{r}_{P\bm{z}_0}\right)$, 
\begin{align*}
\left|B_{P, \bm{z}_0}(-\sigma)\cdot \left\langle\widehat{\zeta}_{P, \bm{z}_0}(\sigma), \psi\right\rangle \right| \;&=\; \left|\int_{\C^n}B_{P, \bm{z}_0}(-\sigma)\cdot  P(\bm{z})^{-\sigma}\cdot\left(\overline{P(\bm{z})}\right)^{-\sigma}\cdot \psi\left(\bm{z}, \overline{\bm{z}}\right) \cdot\textrm{d}\bm{z}\wedge\textrm{d}\overline{\bm{z}}  \right|\\
&\underset{\eqref{eqt1}}{=}\; \left|\int_{\C^n} P(\bm{z})\cdot\left|P(\bm{z})\right|^{-2\sigma}\cdot \mathcal{D}^*_{\bm{z}_0}(\bm{z}, -\sigma, \bm{\partial})\psi\left(\bm{z}, \overline{\bm{z}}\right) \cdot\textrm{d}\bm{z}\wedge\textrm{d}\overline{\bm{z}}  \right|\\
&\le\; \int_{\C^n} \left|P(\bm{z})\right|^{-2(\sigma-1/2)}\varphi\left(\bm{z}, \overline{\bm{z}}\right) \cdot\textrm{d}\bm{z}\wedge\textrm{d}\overline{\bm{z}} \\
&=\; \left\langle\widehat{\zeta}_{P, \bm{z}_0}\left(\sigma-\frac{1}{2}\right), \varphi\right\rangle.
\end{align*}
From the definition of  $\hat{r}_{P, \bm{z}_0}$ as the smallest real pole of the local zeta distribution, this last quantity remains bounded as $\sigma$ tends to $\hat{r}_{P, \bm{z}_0}$, which is easily seen to imply the claim.
\end{proof}

\begin{lem}\label{lem1052}
Under the assumptions of the proposition, \sloppy the local log--canonical threshold $\C\textrm{-}LCT_{P, \bm{z}_0}=\left(\hat{r}_{P, \bm{z}_0},\, \hat{m}_{P, \bm{z}_0}\right)$ is such that $-\hat{r}_{P, \bm{z}_0}$ is the largest root of the local Sato--Bernstein polynomial $B_{P, \bm{z}_0}(s)$ with multiplicity \emph{at most} $\hat{m}_{P, \bm{z}_0}$.
\end{lem}

\begin{proof}
The proof relies on Lichtin's refinement of Kashiwara's Theorem~\ref{kashiroot} as found in~\cite{lichtin}.\\

Consider a resolution of singularity $\hat{g}$ of the polynomial $P(\bm{z})$ at $\bm{z}=\bm{z}_0$ as in Lemma~\ref{lemcomloct}, and keep the notations of this statement. Set 
\begin{equation}\label{defqhat}
\hat{\omega}~: \bm{y}\in\widehat{\mathcal{M}}\mapsto \left(P\circ \hat{g}\right)(\bm{y})
\end{equation} 
and let $B_{\hat{\omega}}(s)$ denote the Sato--Bernstein polynomial of the analytic map $\hat{\omega}$ (it is defined as in the polynomial case by allowing the differential operators $\mathcal{D}(\bm{y}, s, \bm{\partial})$ to be analytic in the variable $\bm{y}$; furthermore, it still satisfies the conclusions of Kashiwara's Theorem~\ref{kashiroot}). It then follows from~\cite[Equation~(4.6)]{lichtin} that there exists an integer $a\ge 0$ such that $B_{P, \bm{z}_0}(s)$ divides the product $\prod_{m=0}^a B_{\hat{\omega}}(s+m)$.\\

It is known that a local version of the inequality~\eqref{rhatineeq} holds; namely, that $-\hat{r}_{P, \bm{z}_0}$ is a rational number in $[-1, 0)$ which is, furthermore, the largest root of the polynomial $B_{P, \bm{z}_0}(s)$  --- see~\cite[\S\S 10.6 \& 10.7]{kollar} and~\cite[\S 3]{saito}. Since the roots of the polynomial $B_{\hat{\omega}}(s)$ are negative, this implies that the multiplicity of $-\hat{r}_{P, \bm{z}_0}$ as a  root of $B_{P, \bm{z}_0}(s)$ is bounded above by its muliplicity as a root of the factor $B_{\hat{\omega}}(s)$  in the above product (which factor is obtained when $m=0$).\\

Assume that $\left\{\widehat{\mathcal{N}}_{\bm{y}_0}\right\}_{\bm{y}_0}\subset \widehat{\mathcal{M}}$ is a finite collection of charts covering an open neighbourhood of the preimage $\hat{g}^{-1}(\{\bm{z}_0\})$. Assume also that the restriction of the map $\hat{\omega}$ to each $\widehat{\mathcal{N}}_{\bm{y}_0}$ admits a \sloppy monomial form such as in~\eqref{monoform} and~\eqref{monoformbis} with exponents, say, $\bm{k}(\bm{y}_0)=\left(k^{(1)}(\bm{y}_0), \dots, k^{(n)}(\bm{y}_0)\right)$ and $\bm{h}(\bm{y}_0)=\left(h^{(1)}(\bm{y}_0), \dots, h^{(n)}(\bm{y}_0)\right)$ in $\N_0^n$. Then, as stated in~\cite[Equation~(4.13)]{lichtin}, the Sato--Bernstein polynomial of such a restriction is 
\begin{equation}\label{berstereslichtin}
\prod_{i=1}^n\prod_{b=1}^{k^{(i)}(y_0)}\left(s+\frac{h^{(i)}(\bm{y}_0)+b}{k^{(i)}(\bm{y}_0)}\right),
\end{equation} 
where a factor is conventionally taken as 1 when  $k^{(i)}(\bm{y}_0)=0$. Furthermore --- this follows from~\cite[Point (3) p.301]{lichtin}, the polynomial $B_{\hat{\omega}}(s)$ (of the "unrestricted" map $\hat{\omega}$ defined in~\eqref{defqhat}) divides any nonzero element lying in the intersection of the ideals generated by the polynomials~\eqref{berstereslichtin} as $\bm{y}_0$ varies. (It should be noted here that for the purpose of his work, it is enough for Lichtin to choose in the aforementioned reference the product of all such polynomials as a nonzero element in this intersection. The present proof, however, requires to take below their lowest common multiple.)\\

Specialise these claims to the case where the collection of charts $\left\{\widehat{\mathcal{N}}_{\bm{y}_0}\right\}_{\bm{y}_0}$ is the neighbourhood of the preimage $\hat{g}^{-1}(\{\bm{z}_0\})$ induced by the charts $\left\{\widehat{\mathcal{M}}_{\bm{y}_0}\right\}_{\bm{y}_0\in \hat{I}(\bm{z}_0)}$ in Lemma~\ref{lemcomloct} (where each chart $\widehat{\mathcal{M}}_{\bm{y}_0}$ is itself a union of at most $\hat{K}(\bm{y}_0)\ge 1$ subsets). One thus obtains that the polynomial $B_{\hat{\omega}}(s)$ divides 
\begin{equation}\label{lcmpolyn}
\textrm{lcm}\left\{ \prod_{i=1}^n\prod_{b=1}^{k_j^{(i)}(y_0)}\left(s+\frac{h_j^{(i)}(\bm{y}_0)+b}{k_j^{(i)}(\bm{y}_0)}\right)\; :\; \bm{y}_0\in \hat{I}(\bm{z}_0), \; 1\le j\le \hat{K}(\bm{y}_0)\right\}.
\end{equation}
The equation~\eqref{formullct} shows that the pole $\hat{r}_{P, \bm{z}_0}$ is the opposite of the largest root of the polynomial~\eqref{lcmpolyn} (which largest root is in particular attained when $b=1$). The multiplicity of this root is clearly given by the equation~\eqref{formullctm} expressing the order of the quantity $\hat{r}_{P, \bm{z}_0}$ considered as a pole   of the local zeta distribution. As a consequence, the equation~\eqref{formullctm} also provides an upper bound for the multiplicity of $-\hat{r}_{\bm{z}_0}$ as the root of $B_{\hat{\omega}}(s)$, and therefore as a root of the local Sato--Bernstein polynomial $B_{P, \bm{z}_0}(s)$.  \\

This establishes Lemma~\ref{lem1052}.
\end{proof}
Proposition~\ref{proplocthmlctsb} follows from Lemmata~\ref{lem1167} and~\ref{lem1052}.
\end{proof}

\begin{proof}[Proof of Theorem~\ref{lctordermultplicity}] 
Let $\bm{z}_0\in\C^n$. Consider an open neighbourhood $Z_{\bm{z}_0}\subset\C^n$ of $\bm{z}_0$ obtained as the intersection of a neighbourhood where the local Sato--Bernstein polynomial $\widehat{B}_{P, \bm{z}_0}(s)$ exists and of a neighbourhood supporting the test functions against which the  above defined local zeta distribution $\widehat{\zeta}_{P, \bm{z}_0}$ is evaluated. \\ 

Then, on the one hand, it follows from the equation~\eqref{locglobsb} that there exists $\bm{z}_0\in\C^n$ such that the largest root of the local polynomial $B_{P, \bm{z}_0}(s)$ and its multiplicity coincide with those of the global polynomial $B_P(s)$ (recall here that the relation~\eqref{locglobsb} is only valid when the ground field is $\K=\C$).\\ 

On the other, the global zeta distribution $\widehat{\zeta}_{P}$ can be decomposed with the help of a locally finite partition of unity into a sum of local zeta distributions $\widehat{\zeta}_{P, \bm{z}_0}$, each evaluated against test functions supported in the above defined open sets $Z_{\bm{z}_0}$. From Proposition~\ref{proplocthmlctsb} and from the point (Y1) in Lemma~\ref{lemcomloct}, the  log--canonical thresholds of all these local distributions at points $\bm{z}_0$ of the algebraic variety $\mathcal{Z}_{\C}(P)$ determine the largest roots of the corresponding local Sato--Bernstein polynomials and their multiplicities.\\

The theorem then immediately follows upon combining the above two points.
\end{proof}

The local version of Theorem~\ref{lctordermultplicity} obtained in Proposition~\ref{proplocthmlctsb} implies the following neat statement relating, on the one hand the largest root of the local Sato--Bernstein polynomial and its multiplicity and, on the other, the growth rate of the weighted volume of the \emph{complex} domain around a point in the variety $\mathcal{Z}_\C(P)$. In order to state the result, given a vector $\bm{z}=\left(z_1, \dots, z_n\right)\in\C^n$, set (with a slight abuse of notation) $$\left\|\bm{z}\right\|\;=\; \sqrt{\sum_{i=1}^n\left|z_i\right|^2}.$$

\begin{thm}\label{lctvolcomplex}
Let $P(\bm{z})\in\C[\bm{z}]$ be polynomial in $n\ge 2$ variables and let $\bm{z}_0\in\C^n$. Let $\hat{r}_{P, \bm{z}_0}\in (0,1]$ be the opposite of the largest root of the local Sato--Bernstein polynomial $B_{P, \bm{z}_0}(s)$. Denote by $\hat{m}_{P, \bm{z}_0}\ge 1$ the multiplicity of this root.  Then, there exists $\delta_P(\bm{z}_0)>0$ such that for any map $\psi_{\bm{z}_0}\ge 0$ in $\mathcal{C}_{c}^{\infty}(\C^n)$ taking a positive value at $\bm{z}_0$ and supported in the ball determined by the inequality $\left\|\bm{z}-\bm{z}_0\right\|<\delta_P(\bm{z}_0)$, it holds that 
\begin{equation}\label{volsblocvol}
\int_{\C^n}\chi_{\left\{\left|P(\bm{z})\right|\le c\right\}}\cdot\psi_{\bm{z}_0}(\bm{z}, \overline{\bm{z}})\cdot \textrm{d}\bm{z}\wedge\textrm{d}\overline{\bm{z}}\;=\; \gamma_{P, \bm{z}_0}(\psi_{\bm{z}_0})\cdot c^{\hat{r}_{P, \bm{z}_0}}\cdot R_{\psi_{\bm{z}_0}}\left(\left|\log c\right|\right)+O\left(c^{\hat{r}_{P, \bm{z}_0}+\varepsilon}\right)
\end{equation} 
for some $\varepsilon>0$ as the parameter $c$ tends to zero. In this relation, $\gamma_{P, \bm{z}_0}(\psi_{\bm{z}_0})>0$ is a positive constant and $R_{\psi_{\bm{z}_0}}(x)\in\R[x]$ is a univariate monic polynomial  with degree $\hat{m}_{P, \bm{z}_0}-1$.

\end{thm}

Upon considering an appropriate partition of unity, this statement and its proof can be extended to the case where the volume estimates are concerned with a union of finitely many balls with radii $\delta_P(\bm{z}_0)>0$ centered at points $\bm{z}_0$ lying in the complex variety $\mathcal{Z}_\C(P)$. The growth rate of the volume estimates is then a function of the largest root of the lowest common multiple of the corresponding local Sato--Bernstein polynomials, and of its multiplicity.

\begin{proof} The proof of the statement is sketched as it follows upon repeating various calculations done previously. \\

Let $\delta_P(\bm{z}_0)>0$ be chosen so that the local Sato--Bernstein polynomial $B_{P, \bm{z}_0}(s)$ exists over the complex ball of radius $\delta_P(\bm{z}_0)$ centered at $\bm{z}_0$ and such that the local zeta distribution $\widehat{\zeta}_{P, \bm{z}_0}$ is well-defined when tested against smooth functions supported in this same ball.\\

As in Lemma~\ref{lemloir}, the relation~\eqref{volsblocvol} is a consequence of Chambert--Loir's and Tschinkel's Tauberian Theorem~\cite[Appendix A, Theorem A.1]{chambloitshcik} applied in the fol\-lo\-wing form~: define a measure $\mu_{\psi_{\bm{z}_0}}$ in $\C^n$ by setting $\mu_{\psi_{\bm{z}_0}}(A)= \int_A \psi_{\bm{z}_0}(\bm{z}, \overline{\bm{z}})\cdot\textrm{d}\bm{z}\wedge\textrm{d}\overline{\bm{z}}$ for any Borel set $A\subset\C^n$. The assumption that $\psi_{\bm{z}_0}(\bm{z}_0)>0$ implies that the meromorphic map $s\mapsto  \left\langle\widehat{\zeta}_{P, \bm{z}_0}(s), \psi\right\rangle$ achieves the pair $\C\textrm{-}LCT_{P, \bm{z}_0}$ (the argument follows from the complex analogue of the  point (W1) in Lemma~\ref{locrlct}, which analogue is a direct consequence of Lemma~\ref{lemcomloct}). The aforementioned Tauberian Theorem then yields that the relation~\eqref{volsblocvol} holds provided that the assumptions (a) and (b) stated in the proof of Lemma~\ref{lemloir} are satisfied for the local complex zeta function $\widehat{\zeta}_{P, \bm{z}_0}$. This is indeed easily verified in the same way as in this proof upon working with the complex analogue~\eqref{sbzetameorextcompl} of the Sato--Bernstein relation.\\

This is enough to establish Theorem~\ref{lctvolcomplex} in view of Proposition~\ref{proplocthmlctsb}, which states that the largest root of the Sato--Bernstein polynomial $B_{P, \bm{z}_0}(s)$ is the opposite of the smallest real pole of $\widehat{\zeta}_{P, \bm{z}_0}$,  the order and the multiplicity coinciding.
\end{proof}

\newpage
$\quad$
\newpage
\chapter{Counting Solutions to Generic Homogeneous Forms Inequalities}\label{ranunidisfamhomfor} 

\vspace{10mm}

The goal in this chapter is to establish Theorems~\ref{athmarg} and~\ref{athmarg2}. To this end, given an integer $n\ge 2$, denote by $\mu_n$ the Haar probability measure on the space $X_n$ of uni\-mo\-vdular lattices in $\R^n$ identified with the homogeneous space $ \textrm{SL}_n(\R)/\textrm{SL}_n(\Z)$ through the map $\mathfrak{g}\in SL_n(\R)/SL_n(\Z)\mapsto \mathfrak{g}\cdot\Z^n\in X_n$.  Recall that $\V_n$ denotes the $n$--dimensional Lebesgue measure. The proofs 
rely on the sharp volume estimates obtained in Chapter~\ref{secvolestim} and on the following lemma.

\begin{lem}[Athreya \& Margulis, 2009]\label{latticeintersec}
Let $A\subset\R^n$ be a Lebesgue measurable set such that $\V_n(A)>0$. Then, $$\mu_n\left(\left\{\Lambda\in X_n : \Lambda\cap A=\emptyset\right\}\right)\;\le\;\frac{c_n}{\V_n(A)}$$ for some constant $c_n>0$ depending only on the dimension $n\ge 2$.
\end{lem}

\begin{proof}
See~\cite[Theorem~2.2]{athreyamargloglaw}.
\end{proof}

For the sake of simplicity of notation, given real numbers $a, b$, set
\begin{equation*}
\overline{\mathcal{S}}_{\bm{F}}(a, b)\;=\; \mathcal{S}_{\bm{F}}\left(B_n, I_n, a, b\right),
\end{equation*}
where $B_n$ denotes the unit Eulidean ball and where the set on the right-hand side is defined in~\eqref{ensalg0}. When $a=T$ is a (large) parameter and when $b=a(T)$ is a function of this parameter, this is also the set $\underline{\mathcal{S}}_{\bm{F}}\left(B_n, b(T)\right)$ defined in~\eqref{ensalg0fafaddebis}. 

\begin{proof}[Proof of Theorem~\ref{athmarg}] 
Given reals $T, b>0$, from the homogeneity of degree $d$ of the forms composing the $p$--tuple $\bm{F}(\bm{x})$, it holds that 
\begin{align}
\V_n\left(\overline{\mathcal{S}}_{\bm{F}}(T, b)\right)\; &=\; \V_n\left(\left\{\bm{x}\in B_n(T)\; :\; 0<\left\|\bm{F}\left(\bm{x}\right)\right\|<b\right\}\right)\nonumber\\ 
&=\; \int_{\R^n}\chi_{\left\{\left\|\bm{F}(\bm{x})\right\|<b\right\}}\cdot\chi_{\left\{\left\|\bm{x}\right\|\le T\right\}}\cdot\textrm{d}\bm{x} \nonumber\\
&=\; b^{n/d}\cdot \int_{\R^n} \chi_{\left\{\left\|\bm{F}(\bm{x})\right\|<1\right\}}\cdot\chi_{\left\{\left\|\bm{x}\right\|\le b^{-1/d}\cdot T\right\}}\cdot\textrm{d}\bm{x}\nonumber\\
&=\; b^{n/d}\cdot \V_n\left(\overline{\mathcal{S}}_{\bm{F}}\left(b^{-1/d}T, 1\right)\right).\label{changevarvol}
\end{align}
Thus, defining for $j\ge 0$ the sets $$A_{\bm{F}}(j)\;=\; \left\{\bm{x}\in B_n\left(f(2^j)\right)\; :\; 0<\left\|\bm{F}(\bm{x})\right\|<2^{-j}\right\},$$ Lemma~\ref{latticeintersec} implies that $$\mu_n\left(\left\{\Lambda\in X_n : \Lambda\cap A_{\bm{F}}(j)=\emptyset\right\}\right)\;\underset{\eqref{formulevol}\, \& \,\eqref{changevarvol}}{\ll}\;\; \frac{2^{jr}}{f\left(2^j\right)^{n-rd}\cdot\left|\log \left(f\left(2^j\right)^d\cdot 2^{-j}\right)\right|^{m-1}}\cdotp$$By assumption, the right--hand side is the general term of a convergent series. As a consequence,  the Borel--Cantelli Lemma yields that for $\mu_n$--almost every $\Lambda\in X_n$, there exists an integer $j(\Lambda)\ge 0$ such that for all $j\ge j(\Lambda)$, one has  that $\Lambda\cap A_{\bm{F}}(j)\neq \emptyset$. Since the probability measure $\mu_n$ descends from the Haar measure on $\textrm{SL}_n\left(\R\right)$, this is easily seen to imply the existence, for almost all $\mathfrak{g}\in SL_n\left(\R\right)$,  of an integer $j_{\mathfrak{g}}\ge 1$ and of a vector $\bm{m}\in\Z^n\backslash\left\{\bm{0}\right\}$ such that for all $j\ge j_{\mathfrak{g}}$, 
\begin{equation}\label{rel2j}
\left\|\left(\bm{F}\circ \mathfrak{g}\right)\left(\bm{m}\right)\right\|< 2^{-j} \qquad \textrm{and}\qquad \left\|\mathfrak{g}\cdot\bm{m}\right\|\le f\left(2^j\right).
\end{equation}

Fix $\varepsilon\in \left(0, 2^{-j_{\mathfrak{g}}+1}\right)$ and let $j\ge j_{\mathfrak{g}}$ be the integer such that $2^{-j}\le\varepsilon<2^{-j+1}$. From the growth condition~\eqref{growthfcth} and the assumption that the function $f$ is non--decreasing, one infers the existence of an integer $j_f\ge 0$ and of a constant $\kappa_f\ge 1$ such that for all $j\ge j\left(\mathfrak{g}, f\right)=\max\left\{j_f, j_{\mathfrak{g}}\right\}$, the inequalities~\eqref{rel2j} yield that $$\left\|\left(\bm{F}\circ \mathfrak{g}\right)\left(\bm{m}\right)\right\|< \varepsilon \qquad \textrm{and}\qquad \left\|\mathfrak{g}\cdot\bm{m}\right\|\le \kappa_f\cdot f\left(\varepsilon^{-1}\right).$$Define $\kappa_{\mathfrak{g}^{-1}}>0$ as the operator norm of the linear map determined by $\mathfrak{g}^{-1}$ in such a way that $\left\|\bm{m}\right\|\le \kappa_{\mathfrak{g}^{-1}}\cdot \left\|\mathfrak{g}\cdot\bm{m}\right\|$. Then, upon setting $\kappa_{\bm{F}}(\mathfrak{g}, f)=\kappa_f\cdot \left(\kappa_{\mathfrak{g}^ {-1}}\right)^{-1}$, one obtains $$\left\|\left(\bm{F}\circ \mathfrak{g}\right)\left(\bm{m}\right)\right\|< \varepsilon \qquad \textrm{and}\qquad \left\|\bm{m}\right\|\le \kappa_{\bm{F}}(\mathfrak{g}, f)\cdot f\left(\varepsilon^{-1}\right).$$The proof is then complete upon defining the quantity $\varepsilon_{\bm{F}}(\mathfrak{g}, f)$ in the statement of Theorem~\ref{athmarg} as $\varepsilon_{\bm{F}}(\mathfrak{g}, f)= 2^{- j\left(\mathfrak{g}, f\right)+1}$.
\end{proof}

As for Theorem~\ref{athmarg2}, it is a particular case of the following more general result. To state it, given a real number $x$,  set $x_{+}=\max\left\{x, 0\right\}\ge 0$ and $x_-=\max\left\{-x, 0\right\}\ge 0$ so that $x=x_{+}-x_{-}$ and $|x|=x_{+}+x_{-}$.

\begin{thm}\label{thmprobacountvol}
Let $P(\bm{x})\in\R[\bm{x}]$ be a homogeneous polynomial of degree $q\ge 1$ in $n\ge 3$ variables. Given real numbers $a<b$, a parameter $T\ge 1$ and a matrix $\mathfrak{g}\in SL_n(\R)$, set $$\widetilde{\mathcal{N}}_{P}(\mathfrak{g},  a, b, T)\; =\; \#\left\{\bm{m}\in\Z^n \cap B_n(T) : a<\left(P\circ \mathfrak{g}\right)\left(\bm{m}\right)\le b\right\}.$$ 

Assume that whenever $b_{\pm}>a_{\pm}$, it holds that $n>q\cdot r_P^{\pm}$, where the quantities $r_P^{\pm}$ are defined in the equation~\eqref{latedef}. Then, for  almost all $\mathfrak{g}\in SL_n(\R)$, there exists exponents $\eta_P^{\pm}>0$ depending only on the polynomial $P(\bm{x})$ such that, as $T$ tends to infinity,
\begin{align}
\frac{\widetilde{\mathcal{N}}_{P}(\mathfrak{g}, a, b, T)}{T^n}\; =\; & \left(\gamma_{P, \mathfrak{g}}^{+}+O\left(f_P^{+}\left(\frac{b_+}{T^q}\right)^{-\eta_P^{+}}\right)\right)\cdot \left(\frac{b_+}{T^q}\right)^{r_P^{+}}\cdot\left|\log \left(\frac{b_+}{T^q}\right)\right|^{m_P^{+}-1}\nonumber \\
&\qquad -\left(\gamma_{P, \mathfrak{g}}^{+}+O\left(f_P^{+}\left( \frac{a_+}{T^q}\right)^{-\eta_P^{+}}\right)\right)\cdot  \left(\frac{a_+}{T^q}\right)^{r_P^{+}}\cdot\left|\log \left(\frac{a_+}{T^q}\right)\right|^{m_P^{+}-1}\nonumber \\
&-\left(\gamma_{P, \mathfrak{g}}^{-}+O\left(f_P^{-}\left(\frac{b_-}{T^q} \right)^{-\eta_P^{-}}\right)\right)\cdot  \left(\frac{b_-}{T^q} \right)^{r_P^{-}}\cdot\left|\log \left(\frac{b_-}{T^q}\right)\right|^{m_P^{-}-1}\nonumber \\
&\qquad +\left(\gamma_{P, \mathfrak{g}}^{-}+O\left(f_P^{-}\left(\frac{a_-}{T^q} \right)^{-\eta_P^{-}}\right)\right)\cdot  \left(\frac{a_-}{T^q} \right)^{r_P^{+}}\cdot\left|\log \left(\frac{a_-}{T^q}\right)\right|^{m_P^{-}-1}.\label{soughtestim}
\end{align}
In this relation, the functions $f_P^{\pm}$ and the constants $\gamma_{P, \mathfrak{g}}^{\pm}\ge 0$ are those defined in Corollary~\ref{zetapm}. Furthermore, the implicit constants are allowed to depend on $\mathfrak{g}$.
\end{thm}

Admissible ranges for the exponents $\eta_P^{\pm}>0$ can easily be derived from the proof of the theorem. They have not been recorded here as they depend on a rather tedious distinction of cases depending on the form of the maps $f_P^\pm$ and on the signs of the real numbers $a$ and $b$.\\

\begin{proof}[Deduction of Theorem~\ref{athmarg2} from Theorem~\ref{thmprobacountvol}] To see why the above statement implies The\-orem~\ref{athmarg2}, take $ P(\bm{x})=P_{\bm{F}}(\bm{x})=\left\|\bm{F}(\bm{x})\right\|^2$. Then,  the volume estimate~\eqref{volestimformsign} in Corollary~\ref{zetapm} reduces to an estimate for the quantity $\V_n\left(\underline{\mathcal{S}}_{\bm{F}}(B_n, b(T))\right)$ when $q=2d$ and is only relevant when all signs are positive. Furthermore, the pair $\left(r^{+}_P, m^{+}_P\right)=\left(r^{+}_{P_{\bm{F}}}, m^{+}_{P_{\bm{F}}}\right)$ therein becomes the real log-canonical threshold, defined in~\eqref{rlctpolezeta}, of the polynomial $P_{\bm{F}}(\bm{x})$. This pair also equals $(r_{\bm{F}}/2, m_{\bm{F}})$, where $r_{\bm{F}}>0$ is the smallest real pole of the distribution $\zeta_{\bm{F}}$ defined in~\eqref{distrizetaR} and where $m_{\bm{F}}\ge 1$ is its order. From Theorem~\ref{volestim}, the pair $(r_{\bm{F}},m_{\bm{F}})$ is therefore the one appearing in the volume estimate~\eqref{formulevol}. To obtain Theorem~\ref{athmarg2}, it then remains to notice that, for this choice of the homogeneous form $P_{\bm{F}}(\bm{x})$,  all quantities appearing with a negative sign in the above  Theorem~\ref{thmprobacountvol} can be taken as zero. 
\end{proof}

\begin{proof}[Proof of Theorem~\ref{thmprobacountvol}]
The proof relies on the refinement by Kelmer and Yu~\cite[Theo\-rem~1]{dubyu} of Theorem~1.2 by Athreya and Margulis in~\cite{athrmargbib}. The latter is concerned with the case of  unimodular twists of indefinite quadratic forms satisfying the volume estimate~\eqref{volsqgtab} and the former with unimodular twists of the diagonal homogeneous forms defined in~\eqref{diagkelyu}. The additional --- and essential -- ingredient introduced here to deal with the most general case of any homogeneous form is the sharp volume estimate obtained in Corollary~\ref{zetapm}, the form of which requires a much finer analysis of the error terms.\\ 

Given $\mathfrak{g}\in SL_n(\R)$, set as in~\eqref{defslast} $$\widetilde{\mathcal{S}}_{P}(\mathfrak{g}, a, b, T)\; =\; \left(\left(P\circ\mathfrak{g}\right)^{-1}\left(\left[a, b\right]\right)\right)\;\cap\; B_n(T)\; =\; \left\{\bm{x}\in B_n(T)\; :\; a< \left(P\circ\mathfrak{g}\right)(\bm{x}) \le b\right\}$$  and  note that 
\begin{align*}
\widetilde{\mathcal{N}}_{P}(\mathfrak{g},  a, b, T)\; &=\;\#\left(\widetilde{\mathcal{S}}_{P}(\mathfrak{g}, a, b, T)\;\cap\;\Z^n\right)\\
& =\; \#\left(\left(P^{-1}\left(\left(a, b\right]\right)\right)\;\cap\; \left(\mathfrak{g}\cdot B_n(T)\right)\;\cap\;\left(\mathfrak{g}\cdot\Z^n\right)\right)\\
&=\; \#\left\{\bm{m}\in \mathfrak{g}\cdot\Z^n\cap \mathfrak{g}\cdot B_n(T)\; :\; a<P(\bm{x})\le b\right\}.
\end{align*}  
It is enough to establish Theorem~\ref{thmprobacountvol} when the transformation $\mathfrak{g}$ is assumed to lie in any given compact subset $B$ of $\textrm{SL}_n(\R)$. Let $$\mathfrak{B}_{P, B}(a, b, T)\;=\; \left\{\mathfrak{g}\in B\; :\; \left|\widetilde{\mathcal{N}}_{P}(\mathfrak{g}, a, b, T)- \textrm{Vol}_n\left(\widetilde{\mathcal{S}}_{P}(\mathfrak{g}, a, b, T)\right)\right|\ge k_P(\mathfrak{g}, a, b, T)\right\},$$where $T\mapsto k_P(\mathfrak{g}, a, b, T)$  is some increasing function to be defined. Throughout this proof, all implicit multiplicative constants are allowed to depend on $B$.\\

Given an integer $j\ge 1$, let $\mathcal{O}_j$ be the $(1/j)$--neighbourhood of the identity in $\textrm{SL}_n(\R)$ induced by the operator norm (with respect to the Euclidean norm). It is proved in~\cite[Lemma~2.1]{dubyu} that there exists a finite set $I_j\subset K$ with cardinality 
\begin{equation}\label{cardinality}
\#I_j\;\ll\; j^{(n+2)(n-1)/2}
\end{equation} such that $B\subset \bigcup_{\mathfrak{h}\in I_j}\left(\mathfrak{h}\cdot \mathcal{O}_j\right)$. Fixing $j\ge 1$ and $\mathfrak{h}\in I_j$, let $$\widetilde{\mathcal{S}}_{P}^{(1)}(\mathfrak{h}, a, b, j)\; =\; \left(\left(P\circ\mathfrak{h}\right)^{-1}\left(\left(a, b\right]\right)\right)\;\bigcap\; B_n\!\left(x_1(j, \sigma)\right)$$ and $$\widetilde{\mathcal{S}}_{P}^{(2)}(\mathfrak{h}, a, b, j)\;=\;\left(\left(P\circ\mathfrak{h}\right)^{-1}\left(\left(a, b\right]\right)\right)\;\bigcap\; B_n\!\left(x_2(j, \sigma)\right),$$ where $\sigma>0$ is a parameter to be adjusted later and where $$x_1(j, \sigma)\;=\; \left(1-\frac{1}{j}\right)\cdot j^{\sigma} \qquad \textrm{and}\qquad x_2(j, \sigma)\;=\; \left(1+\frac{1}{j}\right)\cdot\left(j+1\right)^{\sigma}.$$ Thus, whenever $ j^{\sigma}\le T< (j+1)^{\sigma}$ and $\mathfrak{g}\in \mathfrak{h}\cdot \mathcal{O}_j$, it holds that 
\begin{equation}\label{inclu0}
\widetilde{\mathcal{S}}_{P}^{(1)}(\mathfrak{h}, a, b, j)\;\subset\; \left(P^{-1}\left(\left(a, b\right]\right)\right)\;\cap\; \left(\mathfrak{g}\cdot B_n(T)\right) \;\subset\; \widetilde{\mathcal{S}}_{P}^{(2)}(\mathfrak{h}, a, b, j).
\end{equation}

Denote by $ \widetilde{\mu}_n$ the Haar measure on $\textrm{SL}_n(\R)$. As detailed in~\cite[p.12]{athrmargbib}, it follows from Roger's second moment formula for the Siegel transform combined with Chebyshev's inequality that, whenever $n\ge 3$, given any measurable set $ \mathcal{A}\subset\R^n$ and any real parameter $M\ge 1$, it holds that 
\begin{equation}\label{roger}
\widetilde{\mu}_n\left(\mathfrak{M}_B\left(\mathcal{A}, M\right)\right)\ll\frac{\V_n\left(\mathcal{A}\right)}{M^2},
\end{equation}
where
$$\mathfrak{M}_B\left(\mathcal{A}, M\right) = \left\{\mathfrak{g}\in B : \left|\#\left(\left(\mathfrak{g}\cdot\Z^n\right)\cap \mathcal{A}\right)-\V_n\left(\mathcal{A}\right)\right|> M\right\}$$and where the implicit constant depends only on the dimension $n$. Moreover, \cite[Theorem~6]{dubyu} establishes that for any $j\ge 1$, 
\begin{equation}\label{inclu1}
\bigcup_{j^{\sigma}\le T<(j+1)^{\sigma}}\mathfrak{B}_{P, B}(a, b, T)\;\subset\;\bigcup_{\mathfrak{h}\in I_j} \bigcup_{l=1}^2 \mathfrak{M}\left(\widetilde{\mathcal{S}}_{P}^{(l)}(\mathfrak{h}, a, b, j), \widetilde{k}_P(\mathfrak{h}, a, b, j)\right),
\end{equation} 
where 
\begin{equation}\label{inclu2}
\widetilde{k}_P(\mathfrak{h}, a, b, j)\;=\; k_P\left(\mathfrak{h}, a, b, j^{\sigma}\right) - \V_n\left(\left(\widetilde{\mathcal{S}}_{P}^{(2)}(\mathfrak{h}, a, b, j)\right)\Big\backslash \left(\widetilde{\mathcal{S}}_{P}^{(1)}(\mathfrak{h}, a, b, j)\right)\right).
\end{equation}
Keeping the notations of Corollary~\ref{zetapm}, an elementary change of variables as in the proof of Theorem~\ref{athmarg} (see the equation~\eqref{changevarvol}) yields that 
\begin{align*}
&\frac{\V_n\left(\widetilde{\mathcal{S}}_{P}^{(1)}(\mathfrak{h}, a, b, j)\right)}{x_1(j, \sigma)^n} \; =\; \\ 
&\;\;\; \left(\gamma_{P, \mathfrak{h}}^{+}+O\left(\left(f_P^{+}\left(\frac{b_+}{x_1(j, \sigma)}\right)\right)^{-\delta_P^{+}}\right)\right)\cdot \left(\left(\frac{b_{+}}{x_1(j, \sigma)^q}\right)^{r^{+}_P}\cdot\left|\log\left(\frac{b_+}{x_1(j, \sigma)^q}\right)\right|^{m^{+}_P-1}\right) \\
& \;\;\;\quad  - \left(\gamma_{P, \mathfrak{h}}^{+}+O\left(\left(f_P^{+}\left(\frac{a_+}{x_1(j, \sigma)}\right)\right)^{-\delta_P^{+}}\right)\right)\cdot \left(\left(\frac{a_{+}}{x_1(j, \sigma)^q}\right)^{r^{+}_P}\cdot\left|\log\left(\frac{a_+}{x_1(j, \sigma)^q}\right)\right|^{m^{+}_P-1}\right)\\
&\;\;\; -\left(\gamma_{P, \mathfrak{h}}^{-}+O\left(\left(f_P^{-}\left(\frac{b_-}{x_1(j, \sigma)}\right)\right)^{-\delta_P^{-}}\right)\right)\cdot \left(\left(\frac{b_{-}}{x_1(j, \sigma)^q}\right)^{r^{-}_P}\cdot\left|\log\left(\frac{b_-}{x_1(j, \sigma)^q}\right)\right|^{m^{-}_P-1}\right) \\
& \;\;\;\quad  + \left(\gamma_{P, \mathfrak{h}}^{-}+O\left(\left(f_P^{-}\left(\frac{a_-}{x_1(j, \sigma)}\right)\right)^{-\delta_P^{-}}\right)\right)\cdot \left(\left(\frac{a_{-}}{x_1(j, \sigma)^q}\right)^{r^{-}_P}\cdot\left|\log\left(\frac{a_-}{x_1(j, \sigma)^q}\right)\right|^{m^{-}_P-1}\right),
\end{align*}
where the implicit constants are uniform over the compact set $B$. The same formula holds for $\V_n\left(\widetilde{\mathcal{S}}_{P}^{(2)}(\mathfrak{h}, a, b, j)\right)$ upon replacing $x_1(j, \sigma)$ with $x_2(j, \sigma)$. As a consequence, a tedious but not difficult asymptotic expansion calculation leads one to
\begin{align}
\V_n&\left(\left(\widetilde{\mathcal{S}}_{P}^{(2)}(\mathfrak{h}, a, b, j)\right)\Big\backslash \left(\widetilde{\mathcal{S}}_{P}^{(1)}(\mathfrak{h}, a, b, j)\right)\right)\nonumber\\
&\qquad\ll\;\;\;  \left(b_{+}^{r^{+}_P}-a_{+}^{r^{+}_P}\right)\cdot\max\left\{j^{-1}, \, \left(f_P^{+}\left(j^{\sigma}\right)\right)^{-\delta_P^{+}} \right\}\cdot j^{\sigma\left(n-q\cdot r^+_P\right)}\cdot\left(\log j\right)^{m^+_P-1}\nonumber\\
& \qquad\quad\quad+ \left(b_{-}^{r^{-}_P}-a_{-}^{r^{-}_P}\right)\cdot\max\left\{j^{-1}, \, \left(f_P^{-}\left(j^{\sigma}\right)\right)^{-\delta_P^{-}} \right\}\cdot j^{\sigma\left(n-q\cdot r^-_P\right)}\cdot\left(\log j\right)^{m^-_P-1}.\label{volerror}
\end{align}
In this bound, the implicit constant may depend on $a$ and on $b$.\\

The proof is now split into several cases~:
\begin{itemize}
\item Case A~:  $b_{+}^{r^{+}_P}>a_{+}^{r^{+}_P}$ and $b_{-}^{r^{-}_P}=a_{-}^{r^{-}_P}$
, and either (i)~: $f_P^{+}(T)=T$ or (ii)~: $f_P^{+}(T)=\log T$;
\item Case B~:  $b_{+}^{r^{+}_P}=a_{+}^{r^{+}_P}$ and $b_{-}^{r^{-}_P}>a_{-}^{r^{-}_P}$
, and either (i)~: $f_P^{-}(T)=T$ or (ii)~: $f_P^{-}(T)=\log T$;
\item Case C~:  $b_{+}^{r^{+}_P}>a_{+}^{r^{+}_P}$ and $b_{-}^{r^{-}_P}>a_{-}^{r^{-}_P}$.
\end{itemize}

It is enough to establish the sought estimate~\eqref{soughtestim} in Case (A). Indeed, the proof and the result in Case (B) then follow upon switching signs and Case (C) is obtained by adding the appropriate estimates obtained in the first two cases.\\

Assume therefore that (A)~: $b_{+}^{r^{+}_P}>a_{+}^{r^{+}_P}$ and $b_{-}^{r^{-}_P}=a_{-}^{r^{-}_P}$, and also that (i)~: $f_P^{+}(T)=T$ hold. Letting 
\begin{equation}\label{gamma1}
\sigma\;\ge\; 1/\delta_P^{+},
\end{equation} 
one obtains 
\begin{equation*}
\V_n\left(\left(\widetilde{\mathcal{S}}_{P}^{(2)}(\mathfrak{h}, a, b, j)\right)\Big\backslash \left(\widetilde{\mathcal{S}}_{P}^{(1)}(\mathfrak{h}, a, b, j)\right)\right)\;\ll  \; j^{\sigma\left(n-q\cdot r^+_P\right)-1}\cdot \left(\log j \right)^{m_P^+-1}. 
\end{equation*}
Then, the inequality~\eqref{roger} combined with the relations~\eqref{cardinality}, \eqref{inclu0} and~\eqref{inclu1} 
yields that 
\begin{equation*}
\widetilde{\mu}_n\left(\bigcup_{j^{\sigma}\le T<(j+1)^{\sigma}}\mathfrak{B}_{P, B}(a, b, T)\right)\;\ll\; j^{(n+2)(n-1)/2}\cdot\frac{\V_n\left(\widetilde{\mathcal{S}}_{P}^{(2)}(\mathfrak{h}, a, b, j)\right)}{\left(\widetilde{k}_P(\mathfrak{h}, a, b, j)\right)^2},
\end{equation*}
where
\begin{equation*}
\V_n\left(\widetilde{\mathcal{S}}_{P}^{(2)}(\mathfrak{h}, a, b, j)\right)\; \ll\; j^{\sigma\left(n-q\cdot r^+_P\right)}\cdot \left(\log j \right)^{m_P^+-1} .
\end{equation*}
In these inequalities, the implicit constants are, from Corollary~\ref{zetapm}, uniform over the compact set $B$. Define $k_P\left(\mathfrak{h}, a, b, j\right)=j^{\gamma}$ for some exponent $\gamma>0$ satisfying the inequa\-li\-ties
\begin{equation}\label{gamma2}
n-q\cdot r^+_P-\frac{1}{\sigma}\;<\;\gamma\;<\; n-q\cdot r^+_P.
\end{equation} 
Then, the equation~\eqref{inclu2} implies that 
\begin{align*}
\widetilde{\mu}_n\left(\bigcup_{j^{\sigma}\le T<(j+1)^{\sigma}}\mathfrak{B}_{P, B}(a, b, T)\right)\;&\ll\;\frac{ j^{(n+2)(n-1)/2}\cdot j^{\sigma (n-q\cdot r^+_P)}\cdot \left(\log j\right)^{m_P^+-1}}{\left(j^{\gamma\sigma}+O\left(j^{\sigma (n-q\cdot r^+_P)-1}\cdot \left(\log j\right)^{m_P^+-1}\right)\right)^2},\\
&\underset{\eqref{gamma2}}{\ll}\; j^{(n+2)(n-1)/2+\sigma(n-q\cdot r^+_P-2\gamma)}\cdot \left(\log j\right)^{m_P^+-1}.
\end{align*}
This is the general term of a convergent series whenever 
\begin{equation}\label{gamma3}
\gamma\;>\; \frac{n-q\cdot r^+_P}{2}+\frac{n(n+1)}{4\sigma}\cdotp
\end{equation} 
Under the assumption that $n-q\cdot r^+_P>0$, taking into account the constraints~\eqref{gamma1}, \eqref{gamma2} and~\eqref{gamma3}, this occurs for any 
\begin{equation}\label{gammacontr}
\gamma\;\in\;\left(n-q\cdot r^+_P-\frac{1}{\sigma_0}, n-q\cdot r^+_P\right), \quad\textrm{where} \quad \sigma_0\;=\;\max\left\{\frac{1}{\delta_{P}^+}, \frac{n(n+1)+4}{2\cdot(n-q\cdot r^+_P)}\right\}.
\end{equation}

From the Borel--Cantelli Lemma, the upper limit set formed by the family $\left(\mathfrak{B}_{P, B}(a, b, T)\right)_{T\ge 1}$ as $T$ tends to infinity has then  zero $\widetilde{\mu}_n$--measure. In other words, for almost all $\mathfrak{g}\in B$, 
\begin{align*}
\left|\widetilde{\mathcal{N}}_{P}(\mathfrak{g}, a, b, T)-\textrm{Vol}_n\left(\widetilde{\mathcal{S}}_{P}(\mathfrak{g}, a, b, T)\right)\right|\;\ll\;  T^\gamma
\end{align*}
with an implicit constant depending on the choice of $\mathfrak{g}$ and on the reals $a$ and $b$. To complete the proof of Theorem~\ref{thmprobacountvol} in Case A.(i), it is then enough to notice that the above estimate and  Corollary~\ref{zetapm} imply that
\begin{align}
&\left|\widetilde{\mathcal{N}}_{P}(\mathfrak{g}, a, b, T) - \gamma_{P, \mathfrak{g}}^+\cdot T^n\cdot  \left( \left(\frac{b_{+}}{T^q}\right)^{r^{+}_P}\cdot\left|\log\left(\frac{b_+}{T^q}\right)\right|^{m_P^{+}-1}-\left(\frac{a_{+}}{T^q}\right)^{r^{+}_P}\cdot \left|\log\left(\frac{a_+}{T^q}\right)\right|^{m_P^{+}-1}\right)\right|\; \nonumber \\ 
&\le \; \left|\widetilde{\mathcal{N}}_{P}(\mathfrak{g}, a, b, T) -\textrm{Vol}_n\left(\widetilde{\mathcal{S}}_{P}(\mathfrak{g}, a, b, T)\right)\right| + \nonumber \\
&\left| \gamma_{P, \mathfrak{g}}^+\cdot T^n\cdot  \left( \left(\frac{b_{+}}{T^q}\right)^{r^{+}_P}\cdot\left|\log\left(\frac{b_+}{T^q}\right)\right|^{m_P^{+}-1}-\left(\frac{a_{+}}{T^q}\right)^{r^{+}_P}\cdot \left|\log\left(\frac{a_+}{T^q}\right)\right|^{m_P^{+}-1}\right)-\textrm{Vol}_n\left(\widetilde{\mathcal{S}}_{P}(\mathfrak{g}, a, b, T)\right)\right|\nonumber\\
& \ll\; T^\gamma + T^{n- q\cdot r^+_P-\delta_P^+}\cdot\left(\log T\right)^{-1+m_P^+}. \label{inegtriconcl}
\end{align}
Theorem~\ref{thmprobacountvol} thus follows in this case upon choosing any $\eta_P^+< \min\left\{n- q\cdot r^+_P-\gamma, \delta_P^+\right\}$ when $\gamma$ satisfies the conditions~\eqref{gammacontr}; that is, upon choosing any $\eta_P^+\in \left(0, 1/\sigma_0\right)$.\\

Assume now that the remaining case holds, namely that (A)~: $b_{+}^{r^{+}_P}>a_{+}^{r^{+}_P}$ and \mbox{$b_{-}^{r^{-}_P}=a_{-}^{r^{-}_P}$}, and also that (ii)~: $f_P^{+}(T)=\log T$. Then, the error term~\eqref{volerror} becomes 
\begin{align*}
\V_n&\left(\left(\widetilde{\mathcal{S}}_{P}^{(2)}(\mathfrak{h}, a, b, j)\right)\Big\backslash \left(\widetilde{\mathcal{S}}_{P}^{(1)}(\mathfrak{h}, a, b, j)\right)\right)\;\ll\; j^{\sigma\left(n-q\cdot r^{+}_P\right)}\cdot\left(\log j\right)^{m_P^+-1-\delta_P^+}
\end{align*} 
for any choice of the parameter $\sigma>0$. Set $k_P\left(\mathfrak{h}, a, b, j\right)=j^{n-q\cdot  r^+_P}\cdot\left(\log j\right)^{m^+_P-1-\eta_P^+}$ for some 
\begin{equation}\label{contretaP+}
0\;<\; \eta_P^+\;<\; \delta_P^+. 
\end{equation}
Then, as above, one obtains
\begin{align*}
&\widetilde{\mu}_n\left(\bigcup_{j^{\sigma}\le T<(j+1)^{\sigma}}\mathfrak{B}_{P, B}(a, b, T)\right)\;\\
&\qquad\qquad\qquad\qquad\qquad\quad\ll\;\frac{ j^{(n+2)(n-1)/2}\cdot j^{\sigma (n-q\cdot r^+_P)}\cdot \left(\log j\right)^{m_P^+-1}}{\left(j^{n-q\cdot r_P^+}\cdot\left(\log j\right)^{m_P^+-1-\eta_P^+}+O\left(j^{n-q\cdot r^+_P}\cdot\left(\log j\right)^{m_P^+-1-\delta_P^+}\right)\right)^2},\\
&\qquad\qquad\qquad\qquad\qquad\quad\ll\; j^{(n+2)(n-1)/2-\sigma(n-q\cdot r^+_P)}\cdot \left(\log j\right)^{-m_P^++1-2\eta_P^+}.
\end{align*}
Under the assumption that $n-q\cdot r^+_P>0$, this is the general term of a convergent series for $\sigma$ large enough. From the Borel--Cantelli Lemma, this yields that for almost all $\mathfrak{g}\in B$,
\begin{align*}
\left|\widetilde{\mathcal{N}}_{P}(\mathfrak{g}, a, b, T)-\textrm{Vol}_n\left(\widetilde{\mathcal{S}}_{P}(\mathfrak{g}, a, b T)\right)\right|\;\ll\;  T^{n-q\cdot r^+_P}\cdot\left(\log T\right)^{m_P^+-1-\eta_P^+}
\end{align*}
with an implicit constant depending on the choice of $\mathfrak{g}$ and on the reals $a$ and $b$. The proof is then concluded as in~\eqref{inegtriconcl}. The upshot is that Theorem~\ref{thmprobacountvol} is obtained in this case upon choosing any exponent $\eta_P^+$ in the range determined by~\eqref{contretaP+}.\\

This completes the proof of Theorem~\ref{thmprobacountvol}.
\end{proof}

\chapter{Deterministic Counting in Large Domains}\label{secdeterlarge} 

Let $K$ be a globally semianalytic subset of $\R^n$ meeting the assumptions $(\mathcal{H}_1)$ and $(\mathcal{H}_2)$ stated in the Introduction. The goal of this section is to establish a more precise version of Theorem~\ref{sarnak} in the case that the domain $\mathcal{S}^{\dag}_{\bm{F}}(K, T, \alpha)$ defined in~\eqref{colsfsar} is "large enough" in the sense that $0<\alpha\le 1$. Indeed, when, moreover, the assumption $(\mathcal{H}_3)$ holds (one then deals with what is referred to as Case (2) in Section~\ref{sec0.1intro} of the Introduction), it is then possible to determine the precise asymptotic order for the counting function $\mathcal{N}_{\bm{F}}^{\dag}(K, T, \alpha)$ defined in~\eqref{nbsfsar} when $0<\alpha<1$, and to also obtain an upper bound for this function when $\alpha=1$. When one assumes instead the stronger conditon that the set $K$ is star-shaped with respect to the origin and that it contains the origin in its interior (one then deals with what is referred to as Case (3) in  Section~\ref{sec0.1intro} of the Introduction), it becomes possible to get an  asymptotic expansion for $\mathcal{N}_{\bm{F}}^{\dag}(K, T, \alpha)$ when $0<\alpha<1$.\\

To state the main result proved in this chapter, let $P(\bm{x})\in\R\left[\bm{x}\right]$ denote a homogeneous form of degree $q\ge 2$. By analogy with the above notations, let 
\begin{equation}\label{defspalpT} 
\mathcal{S}_{P}^{\dag}(K, T, \alpha)\;=\; \left\{\bm{x}\in T\cdot K\; :\; \left|P\left(\bm{x}\right)\right|\le T^{q-\alpha}\right\}
\end{equation}
and 
\begin{equation}\label{defspalpTbis}
\mathcal{N}_{P}^{\dag}(K, T, \alpha)\;=\; \#\left(\Z^n\cap \mathcal{S}_{P}^{\dag}(K, T, \alpha)\right).
\end{equation}

Theorem~\ref{sarnak} in the case where $0<\alpha\le 1$ is then a consequence of Case (1') in the statement below~: it indeed suffices to apply it to the polynomial $P(\bm{x})=\left\|\bm{F}\left(\bm{x}\right)\right\|^2$  after substituting the parameters $\left(2d, 2\alpha\right)$ for $\left(q, \alpha\right)$ in the above definition~\eqref{defspalpT}. The resulting conclusion generalises that of Theorem~\ref{sarnak}  inasmuch as it is valid \emph{without assuming the flatness condition~\eqref{cndiestimvolsa1reutilise}}. 

\begin{thm}[Asymptotic behaviour of the counting function of homogeneous form inequalities in large domains]\label{countinglargedomains} 
Assume  that the set $K\subset\R^n$ is globally semianalytic and that it satisfies the assumptions $(\mathcal{H}_1)$ and $(\mathcal{H}_2)$. \\

Then, the counting function $\mathcal{N}_{P}^{\dag}(K, T, \alpha)$ can be estimated as follows~:

\begin{itemize}

\item[(1')\,--] assume that the condition $(\mathcal{H}_{3})$ holds. Then,  as the parameter $T$ tends to infinity, one has that
\begin{equation*}\label{countchap5formulevolbis}
\mathcal{N}_{P}^{\dag}(K, T, \alpha)\;\asymp\; \V_n\left(\mathcal{S}_{P}^{\dag}(K, T, \alpha)\right)
\end{equation*}
provided that $\alpha\in (0,1)$. When $\alpha=1$, it holds that 
\begin{equation}\label{vol-0-alpah-1}
\mathcal{N}_{P}^{\dag}(K, T, 1)\;\ll\;\V_n\left(\mathcal{S}_{P}^{\dag}(K, T, 1)\right).
\end{equation} 

\item[\emph{(2')\,--}] assume the stronger condition that the set $K$ is star-shaped with respect to the origin and that it contains the origin in its interior.  Then, provided that $\alpha\in (0,1)$, one has that as the parameter $T$ tends to infinity, 
\begin{equation*}\label{vol0-alpah-1}
\mathcal{N}_{P}^{\dag}(K, T, \alpha)\;=\;\left(1+O\left(\frac{1}{T^{\delta}}\right)\right)\cdot \V_n\left(\mathcal{S}_{P}^{\dag}(K, T, \alpha)\right)
\end{equation*}  
for some  $\delta>0$. 
\end{itemize}
\end{thm}

In this statement, the volumes of the sets $\mathcal{S}_{P}^{\dag}(K, T, \alpha)$ (where $0<\alpha\le  1$) are respectively given by Cases (2) and (3) of Theorem~\ref{volestim} stated in the Introduction.\\

The rest of this chapter is devoted to the proof of Theorem~\ref{countinglargedomains}, starting with the inequality~\eqref{vol-0-alpah-1}.

\begin{proof}[Proof of the volume estimate~\eqref{vol-0-alpah-1}] Assume that the conditions $(\mathcal{H}_1)-(\mathcal{H}_3)$ are met. An elementary packing argument, a variant of which appears in~\cite[Appendix~1]{sarnakprinci}, shows that the sought estimate holds for any value of 
\begin{equation}\label{alphavalue}
\alpha\in [0,1].
\end{equation} 
Indeed, fix $\eta\in (0, 1/4)$ and assume that $\bm{m}\in \Z^n\cap \mathcal{S}_{P}^{\dag}(K, T, \alpha)$. Let $\bm{u}\in\R^n$ be a vector with Euclidean norm at most $\eta$. Consider then the Taylor expansion $$P\left(\bm{m}+\bm{u}\right)\;=\; P(\bm{m})+\sum_{k=1}^q\frac{1}{k!}\cdot\left(\left(\bm{u\cdot}\nabla\right)^kP\right)(\bm{m}),$$where $\bm{u\cdot}\nabla$ is short-hand notation for the differential operator $\sum_{i=1}^{n}u_i\cdot\partial_i$ when $\bm{u}=\left(u_1, \dots, u_n\right)$ and $\nabla=\left(\partial_1, \dots, \partial_n\right)$, and where the power applied to it must be understood as a formal multiplication process. \\

For any given integer $1\le k\le q$, the map $\bm{x}\mapsto \left(\left(\bm{u\cdot}\nabla\right)^kP\right)(\bm{x})$ is homogeneous of degree $q-k$. As a consequence, upon decomposing the integer vector $\bm{m}$ as $\bm{m}=T\cdot \bm{y}$, where $\bm{y}\in K$, one obtains that $$\left|P\left(\bm{m}+\bm{u}\right)\right|\;\le\; T^{q-\alpha}+ \eta\cdot c_P(K)\cdot T^{q-1}\;\underset{\eqref{alphavalue}}{\le}\; \left(1+ \eta\cdot c_P(K)\right)\cdot T^{q-\alpha}$$with a constant $c_P(K)>0$  determined by the sup norms of the polynomial $P(\bm{x})$ and of its partial derivatives over the set $K$. This shows that the pairwise disjoint Euclidean balls with radius $\eta$ centered at the integer points in $\mathcal{S}_{P}^{\dag}(K, T, \alpha)$  are contained in, say, the set $\underline{\mathcal{S}}_P\left(K^{(\eta)},\; (1+\eta\cdot c_P(K))\cdot T^{q-\alpha}\right)$ defined explicitly in~\eqref{ensalgbis}. Here, $K^{(\eta)}$ denotes the $\eta$-neighbourhood of the set $K$. Thus, 
\begin{equation}\label{nbvollardom}
\mathcal{N}_{P}^{\dag}(K, T, \alpha)\;\le\; \V_n\left(\underline{\mathcal{S}}_P\!\left(K^{(\eta)},\; (1+\eta\cdot c_P(K))\cdot T^{q-\alpha}\right)\right).
\end{equation} 
If  $\eta>0$ is chosen small enough, it is immediate that the set $K^{(\eta)}$ also satisfies the assumption $(\mathcal{H}_3)$ with the same compact set $C$ and the same open set $U$ as those introduced in the support restriction condition~\eqref{suppproper} for $K$. Consequently, with the notations of $(\mathcal{H}_3)$, one deduces that the pair $\left(r_P\left(K^{(\eta)}\right), m_P\left(K^{(\eta)}\right)\right)$ equals $\left(r_P(K), m_P(K)\right)$. It then follows from Case (2) in Theorem~\ref{volestim} that $$\V_n\left(\underline{\mathcal{S}}_P\left(K^{(\eta)},\; (1+\eta\cdot c_P(K))\cdot T^{q-\alpha}\right)\right)\;\ll\; \V_n\left(\mathcal{S}_P\!\left(K,\; T^{q-\alpha}\right)\right)\; = \;\V_n\left(\underline{\mathcal{S}}_{P}^{\dag}(K, T, \alpha)\right).$$ In view of the inequality~\eqref{nbvollardom}, this completes the proof.
\end{proof}

Establishing the remaining statements part of Theorem~\ref{countinglargedomains} requires more work. To do so, assume from now on and until the end of the chapter that $\alpha\in (0,1)$.

\section{Smooth Counting}

It is convenient to work in the first place with a smooth version of the counting function $\mathcal{N}_{P}^{\dag}(K, T, \alpha)$. To this end, fix a smooth map $\psi\ge 0$ in $\mathcal{C}_{c}^{\infty}(\R^n)$ whose support will be specified later to approximate suitably the set $K$. The proof of the remaining cases in Theorem~\ref{countinglargedomains} requires to also approximate the characteristic function of an interval by a function whose Fourier transform behaves "nicely". This is a well--studied topic where the extremal functions introduced by Beurling and Selberg are often used (see the discussion and the references in~\cite{CGT} for more details). In the present situation however, it is more relevant to work with the map obtained from the following (particular case of a) result due to Colzani, Gigante and Travaglini in~\cite[Theorem~1.1]{CGT}. \\

Before stating it, recall that a function $\theta~: \R_+\rightarrow\R_+$ has fast decay at infinity if for any $\beta>0$, there exists a constant $c_{\beta}>0$ such that for all $t\ge 0$, it holds that 
\begin{equation}\label{fast}
0\;\le\; \theta(t)\;\le\; \frac{c_{\beta}}{\left(1+t\right)^{\beta}}\cdotp
\end{equation} 
Also, a holomorphic function $F$ is said to be of exponential type $R>0$ if $$R\;=\; \limsup_{r\rightarrow\infty}\;\max_{\left|z\right|=r}\;\frac{\log\left|F(z)\right|}{\left|z\right|}\cdotp$$Equivalently, this is saying that for any $\varepsilon>0$, there exists a constant $\kappa(\varepsilon)>0$ such that $\left|F(z)\right|\le \kappa(\varepsilon)\cdot e^{(R+\varepsilon)\cdot \left|z\right|}$.

\begin{thm}[Colzani, Gigante \& Travaglini, 2011]\label{thmcolzani}
Let $c\in \left(0,1\right)$. Set $I\left(c\right)=\left[-c, c\right]$. Then, there exists a positive function $\theta$ with fast decay at infinity such that for any $R>0$, there exist integrable entire functions $A_{(R, c)}$ and $B_{(R, c)}$ satisfying the following properties~: these two functions are of exponential type $R$ and have their Fourier transforms supported in the interval $\left[-R, R\right]$; furthermore, for all $x\in\R$, 
\begin{equation}\label{ineqartchiR}
A_{(R, c)}(x)\;\le\; \chi_{I\left(c\right)}(x)\;\le\; B_{(R, c)}(x)
\end{equation} 
and 
\begin{equation}\label{ineqartchiRbis}
0\;\le\;  B_{(R, c)}(x)-A_{(R, c)}(x)\;\le\;\theta\left(R\cdot \left|\left|x\right|-c\right|\right).
\end{equation}   
\end{thm}

The main point of this statement is that it allows for an easy control of the gap $\left|B_{(R, c)}(x)-A_{(R, c)}(x)\right|$ by the distance from $x\in\R$ to the \emph{boundary} of the interval $I\left(c\right)$~: roughly speaking, the quantity $\left|B_{(R, c)}(x)-A_{(R, c)}(x)\right|$ is approximately equal to 1 at points $x$ within distance $1/R$ of the interval $I(c)$ and essentially zero at larger distances.\\

Let $T\ge 1$. The proof of Theorem~\ref{countinglargedomains} mainly  reduces to the estimate for $\widetilde{\mathcal{N}}_{P}\left(\psi, T, \alpha\right)$ obtained in Proposition~\ref{propsmoothcountinglargdom} below, which takes into account the support restriction condition~\eqref{suppproper} induced by the assumption $\left(\mathcal{H}_3\right)$. To state it, fix a compact set $\mathcal{K}$ contained in the open set $U$ and containing the compact set $K$ in its interior, where the set $U$ is defined as part of  the condition~\eqref{suppproper}. Assume then that the map $\psi$ is such that 
\begin{equation}\label{2emesupprestrispi}
C\;\subset\; \textrm{Supp}\;\psi\;\subset\; \mathcal{K},\qquad 0\;\le\; \psi\;\le\; 1 \qquad \textrm{and}\qquad \max_{\bm{x}\in \mathcal{K}}\;\psi(\bm{x})\;=\; 1,
\end{equation} 
where the set $C$ is also defined in~\eqref{suppproper}. Take 
\begin{equation}\label{defN_P}
c\;=\; T^{-\alpha} \qquad \textrm{and}\qquad R\;=\;\frac{T}{2\cdot \Delta_P(\mathcal{K})}, \qquad \textrm{where}\qquad \Delta_P(\mathcal{K})\;=\;\max_{\bm{x}\in \mathcal{K}}\left\|\left(\nabla P\right)(\bm{x})\right\|,
\end{equation}
and define $E_{(T, \alpha)}$ as being either of the maps $A_{(T/(2\Delta_P(\mathcal{K})),\; T^{-\alpha})}$ or $B_{(T/(2\Delta_P(\mathcal{K})),\; T^{-\alpha})}$ in  Theorem~\ref{thmcolzani} (i.e.~the maps obtained with the above choice for the parameters $R$ and $c$). \\

Let then 
\begin{equation}\label{sumcountingsmooth}
\widetilde{\mathcal{N}}_{P}\left(\psi, T, \alpha\right)\;=\; \sum_{\bm{k}\in\Z^n} G^{\left(P, \psi\right)}_{T, \alpha}\left(\bm{k}\right),  
\end{equation}
where, given $\bm{x}\in\R^n$, 
\begin{equation}\label{defgppsitalp}
G^{\left(P, \psi\right)}_{T, \alpha}\left(\bm{x}\right)\;=\; \psi\left(\frac{\bm{x}}{T}\right)\cdot E_{(T, \alpha)}\left(P\left(\frac{\bm{x}}{T}\right)\right).
\end{equation}
Given $c\in (0,1)$, set furthermore, as in~\eqref{defmuvolweighted},
\begin{equation}\label{defmuppsit}
\mu_P\left(\psi, c\right)\;=\;\int_{\R^n}\chi_{\left\{\left|P(\bm{x})\right|<c\right\}}\cdot \psi(\bm{x})\cdot \textrm{d}\bm{x}.
\end{equation}
It is useful to record right here the asymptotic behaviour at the origin of the map $c>0\mapsto \mu_P\left(\psi, c\right)$. This can be  inferred from Lemma~\ref{lemloir} upon choosing, with  the notation therein, the  compact set $\mathcal{K}$ as the one defined above and the collection of maps $\mathcal{F}$ as the singleton $\left\{\psi\right\}$. Indeed, note first that the assumptions of Lemma~\ref{lemloir} are met from the last claim in the point (X1) of Proposition~\ref{proprlctgene}. From this point (X1) again and from the conclusions of the lemma, one then obtains the existence of a rational $\rho> 0$ and of an integer $m\ge 0$,  both independent of $\psi$ as chosen above, satisfying the estimate
\begin{equation}\label{asympmupsansome}
\mu_P\!\left(\psi, c\right)= \left(\Theta_P\!\left(\psi\right)+ h_\psi(c)\right)\cdot c^{\rho}\cdot \left|\log c\right|^m,\; \textrm{where} \; \left|h_\psi(c)\right|\le (m+1)\cdot\left(C_N\!\left(\psi, \mathcal{K}\right)\right)^M\cdot\eta\!\left(c\right).
\end{equation}
In this relation,
\begin{itemize}
\item the integers $M, N\ge  1$ depend only on $\mathcal{K}$ and on the polynomial $P(\bm{x})$;
\item $C_N\left(\psi, \mathcal{K}\right)$ is 
the constant defined in~\eqref{defcnfk}, namely
\begin{equation}\label{defcnfkbis}
C_N\left(\psi, \mathcal{K}\right)\;=\;  \max_{0\le \left|\bm{k}\right|\le N} \; \max_{\bm{x}\in\mathcal{K}}\; \left|\frac{\partial^{\left|\bm{k}\right|}\psi}{\partial\bm{x}^{\bm{k}}}(\bm{x})\right|.
\end{equation}
The last relation in~\eqref{2emesupprestrispi} implies that $C_N\left(\psi, \mathcal{K}\right)\ge 1$;
\item from the uniformity claim (V3) in Lemma~\ref{lemloir}, 
\begin{equation}\label{asympmupsansomebis}
0\;<\; \Theta_P\!\left(\psi\right)\;\le\; \left(C_N\!\left(\psi, \mathcal{K}\right)\right)^M;
\end{equation}
\item the quantity $\eta\left(c\right)$ tends to 0 as $c\rightarrow 0^+$ uniformly in the choice of $\psi$ meeting the relations~\eqref{2emesupprestrispi} under the assumption $(\mathcal{H}_3)$. Indeed, the factor $(m+1)\cdot\left(C_N\!\left(\psi, \mathcal{K}\right)\right)^M\cdot\eta\!\left(c\right)$ in~\eqref{asympmupsansome} is obtained upon adding upper bounds for two maps depending on $c$~: 
\begin{itemize}
\item on the one hand the map 
\begin{equation}\label{lastdernier}
c\mapsto \Theta_P\!\left(\psi\right)\cdot \left(R_{\psi}\left(\left|\log c\right|\right)-\left|\log c\right|^m\right)/\left|\log c\right|^m,
\end{equation} 
where $R_{\psi}$ is, with the present notations, the degree $m$ monic polynomial appearing in~\eqref{asympmu0}; furthermore, from the uniformity claim (V3) in Lemma~\ref{lemloir}, each of the $m$ coefficients obtained when expanding this expression is, up to a multiplicative constant depending on $\mathcal{K}$, bounded above by $\left(C_N\!\left(\psi, \mathcal{K}\right)\right)^M$;
\item on the other hand, the map determined by the error term present in~\eqref{asympmu0}, which is of the form $O\left(c^{\varepsilon}\right)$ for some fixed $\varepsilon>0$. Here, from the uniformity claim (V2) in Lemma~\ref{lemloir} and under $(\mathcal{H}_3)$, the implicit constant can be taken up to a multiplicative constant depending on $\mathcal{K}$ and on $\varepsilon$, and uniformly in any $\psi$ meeting the conditions stated  in~\eqref{2emesupprestrispi}, as $C_N\left(\psi, \mathcal{K}\right)$.
\end{itemize} 
Upon adding the resulting form of the error term in the latter case with the upper bound obtained for the map~\eqref{lastdernier}, this yields the required shape $(m+1)\cdot\left(C_N\!\left(\psi, \mathcal{K}\right)\right)^M\cdot\eta\!\left(c\right)$ with the sought uniformity claim on $\eta$.
\end{itemize}
 

\begin{prop}\label{propsmoothcountinglargdom} Assume that $\psi$ is a map in $\mathcal{C}_{c}^{\infty}(\R^n)$  satisfying the relations~\eqref{2emesupprestrispi} and ta\-king  strictly positive values over the set $C$ defined therein. Let $\alpha\in (0,1)$. Then, under the assumption of the homogeneity of the polynomial $P(\bm{x})$, there exist integers $M, N\ge 1$  such that 
\begin{equation}\label{soughtestimate}
\left|\widetilde{\mathcal{N}}_{P}\left(\psi, T, \alpha\right)- T^n\cdot \mu_P\left(\psi, T^{-\alpha}\right)\right|\;\ll\; C_{N}(\psi, \mathcal{K})^{M}\cdot T^{n-\delta}\cdot \mu_P\left(\psi, T^{-\alpha}\right)
\end{equation} 
for any $\delta\in (0, 1-\alpha)$. In this inequality, the constant $C_{N}(\psi, \mathcal{K})$ is the one defined in~\eqref{defcnfkbis}. Furthermore, the implicit constant, the integers $M$ and $N$ and the choice of the exponent $\delta$ do not depend on $\psi$. 
\end{prop}

\begin{proof}[Deduction of Theorem~\ref{countinglargedomains} from Proposition~\ref{propsmoothcountinglargdom} when $0<\alpha<1$] Several steps  in this deduction already appeared in Section~\ref{tauberian} so that the argument is here sketched highlighting the parts which are specific to this setup. \\

Consider first Case (1') of Theorem~\ref{countinglargedomains}. Fix $0<\alpha<1$ and let $\psi_1, \psi_2$ be two maps  in $\mathcal{C}_c^{\infty}(\R^n)$ satisfying the conditions~\eqref{2emesupprestrispi}  and the inequalities 
\begin{equation*}
0\;\le\; \psi_1\;\le\; \chi_K\;\le\; \psi_2.
\end{equation*}
Under the assumption $(\mathcal{H}_3)$, the relation~\eqref{asympmu0} part of Lemma~\ref{lemloir} and the expansion~\eqref{decompvol0} yield that 
\begin{equation*}
T^n\cdot \mu_P\left(\psi_1, T^{-\alpha}\right) \;\asymp\;  \V_n\left(\mathcal{S}_{P}^{\dag}(K, T, \alpha)\right) \;\asymp\; T^n\cdot \mu_P\left(\psi_2, T^{-\alpha}\right).
\end{equation*} 
As a consequence, one infers from Proposition~\ref{propsmoothcountinglargdom} that 
\begin{align*}
\V_n\left(\mathcal{S}_{P}^{\dag}(K, T, \alpha)\right)\;\ll\; \widetilde{\mathcal{N}}_{P}\left(\psi_1, T, \alpha\right)\; &\le\;\mathcal{N}_{P}^{\dag}(K, T, \alpha)\\
&\le\;  \widetilde{\mathcal{N}}_{P}\left(\psi_2, T, \alpha\right)\;\ll\; \V_n\left(\mathcal{S}_{P}^{\dag}(K, T, \alpha)\right).
\end{align*}
Here, the left-most inequality relies on the fact that the exponent $\delta$ in~\eqref{soughtestimate} is, by assumption, strictly positive. This establishes the conclusion of Theorem~\ref{countinglargedomains} in Case~(1') when $0<\alpha<1$. \\

As for Case (2'), it relies on the arguments developed in the proof of Theorem~\ref{mainthmvolomega} (in particular in  those present in the proofs of Lemmata~\ref{estimcoefdomomeggT} and~\ref{estimerreurunif}). To see this, fix an integer $k\ge 1$. Let $\psi_{k}^-$ and $\psi_{k}^+$ be smooth maps in $\mathcal{C}_{c}^{\infty}(\R^n)$ pointwise monotonic as $k$ tends to infinity such that 
\begin{equation}\label{ineqchipsiepsi} 
\chi_{K\left(1-1/k\right)}\;\le\; \psi_{k}^-\;\le\; \chi_{K}\;\le\; \psi_{k}^+\;\le\;  \chi_{K\left(1+1/k\right)}
\end{equation}
(recall here that, given $\lambda>0$, one sets $K(\lambda)=\lambda\cdot K$ and that the inequalities between the characteristic functions are indeed valid under the assumptions that $K$ is star-shaped with respect to the origin). Upon choosing $k$ larger than some integer $k_0$, assume furthermore that the maps $\psi_{k}^{\pm}$ also meet the conditions stated in~\eqref{2emesupprestrispi}. From~\cite[Theorem~1.4.1, p.25]{horm}, they can be chosen so that the constant $C_{N}\left(\psi_{k}^\pm, \mathcal{K}\right)$ defined in~\eqref{defcnfkbis} 
is bounded above as 
\begin{equation}\label{upoucpsikpm}
C_{N}\left(\psi_{k}^\pm, \mathcal{K}\right)\;\ll\; k^{N}
\end{equation} 
for some implicit constant depending only on $N\ge 1$. Given $k\ge k_0$, set then $$\widetilde{\mathcal{N}}_{P, k}^{-}\left(T, \alpha\right)\;=\; \widetilde{\mathcal{N}}_{P}\left(\psi_{k}^-, T, \alpha\right)\qquad \textrm{when} \qquad E_{(T, \alpha)}=A_{(T/(2\Delta_P(\mathcal{K})),\; T^{-\alpha})}$$ and $$\widetilde{\mathcal{N}}_{P, k}^+\left(T, \alpha\right)\;=\; \widetilde{\mathcal{N}}_{P}\left(\psi_{k}^+, T, \alpha\right)\qquad \textrm{when} \qquad E_{(T, \alpha)}=B_{(T/(2\Delta_P(\mathcal{K})), \; T^{-\alpha})}$$ so that 
\begin{equation}\label{sandwicoundag}
\widetilde{\mathcal{N}}_{P, k}^-\left(T, \alpha\right)\;\le \; \mathcal{N}_{P}^{\dag}\left(T, \alpha\right) \;\le\; \widetilde{\mathcal{N}}_{P, k}^+\left(T, \alpha\right).
\end{equation}
To emphasize the analogy with the proof of Theorem~\ref{mainthmvolomega}, given $c\in (0,1)$, let furthermore 
\begin{equation}\label{defnup}
\nu_P\left(K, c\right)\; =\; \int_{K} \chi_{\left\{\left|P(\bm{x})\right|\le c\right\}}\cdot\textrm{d}\bm{x}
\end{equation} 
and $$\nu^\pm_{P}\left( k, c\right)\; =\;\mu_{P}\left(\psi_{k}^\pm , c\right)\;=\; \int_{\R^n} \chi_{\left\{\left|P(\bm{x}\right|\le c\right\}}\cdot\psi_{k}^\pm(\bm{x})\cdot\textrm{d}\bm{x}.$$ Then, for any $k\ge 1$,
\begin{align*} 
\left|\widetilde{\mathcal{N}}_{P, k}^{\pm}\left(T, \alpha\right)-\V_n\left(\mathcal{S}_{P}^{\dag}(K, T, \alpha)\right)\right|\;&=\; \left|\widetilde{\mathcal{N}}_{P, k}^{\pm}\left(T, \alpha\right)-T^n\cdot \nu_P\left(K, T^{-\alpha}\right)\right|\\
&\underset{\eqref{ineqchipsiepsi}}{\le}\; \left|\widetilde{\mathcal{N}}_{P, k}^{\pm}\left(T, \alpha\right)-T^n \cdot \nu^\pm_{P}\left( k, T^{-\alpha}\right) \right| \\
&\qquad \qquad \qquad \qquad \quad   + T^n\cdot\left|\nu^+_{P}\left( k, T^{-\alpha}\right)-\nu^-_{P}\left( k, T^{-\alpha}\right)\right|.
\end{align*}

Since $\chi_K\le \psi^+_{k_0}$, repeating the calculations done in the proof of Lemma~\ref{estimcoefdomomeggT}, one obtains that for any $c>0$ and any $k\ge k_0$, $$\left|\nu^+_{P}\left( k, c\right)-\nu^-_{P}\left( k, c\right)\right|\;\ll\; \frac{\nu^+_{P}\left( k_0, c\right)}{k}$$with an implicit constant independent of $k$ and $c$. From Proposition~\ref{propsmoothcountinglargdom}, this implies that 
\begin{align}\label{ineqnpmpktavolspat} 
\left|\widetilde{\mathcal{N}}_{P, k}^{\pm}\left(T, \alpha\right)-\V_n\left(\mathcal{S}_{P}^{\dag}(K, T, \alpha)\right)\right|\;&\ll\; T^n\cdot \left(C_{N}(\psi_k^{\pm}, \mathcal{K})^{M}\cdot \frac{\nu_P^+\left(k_0, T^{-\alpha}\right)}{T^\delta}+ \frac{\nu^+_{P}\left( k_0, T^{-\alpha}\right)}{k}\right)
\end{align}
for some $\delta>0$ and some integers $M, N\ge 1$. Here, the inclusions \mbox{$\textrm{Supp }\psi_{k_0}^+\subset K(1+1/k_0) \subset K(2) $} yield that
\begin{align*}
\nu_P^+\left(k_0, T^{-\alpha}\right)\; &\underset{\eqref{defnup}}{\le}\; \frac{\nu_P\left(K, 2^{-q}\cdot T^{-\alpha}\right)}{2^n}\; =\;\frac{\V_n\left(\mathcal{S}_{P}^{\dag}(K, 2^{q/\alpha} T, \alpha)\right)}{2^n\cdot \left(2^{q/\alpha} T\right)^n}\\
&\ll\; \frac{\V_n\left(\mathcal{S}_{P}^{\dag}(K, T, \alpha)\right)}{T^n}\;= \nu_P\left(K, T^{-\alpha}\right),
\end{align*} 
where the second last relation is easily deduced from the explicit volume estimate obtained in Case (3) of Theorem~\ref{volestim}. It is then a consequence of the inequalities~\eqref{upoucpsikpm} and~\eqref{ineqnpmpktavolspat} that $$\left|\widetilde{\mathcal{N}}_{P, k}^{\pm}\left(T, \alpha\right)-\V_n\left(\mathcal{S}_{P}^{\dag}(K, T, \alpha)\right)\right|\;\ll\; \V_n\left(\mathcal{S}_{P}^{\dag}(K, T, \alpha)\right)\cdot \left(\frac{k^{MN}}{T^\delta}+\frac{1}{k}\right).$$In view of the inequalities~\eqref{sandwicoundag}, specialising $k$ to be the integer part of $T^{\delta/(1+MN)}$ provides the sought power saving. This completes the deduction of Theorem~\ref{countinglargedomains} from Proposition~\ref{propsmoothcountinglargdom} when $\alpha\in (0,1)$.  
\end{proof}

The proof of Proposition~\ref{propsmoothcountinglargdom} relies on the Poisson summation formula, which implies that the counting function~\eqref{sumcountingsmooth} can be expanded as 
\begin{equation}\label{poissoncountpsiat}
\widetilde{\mathcal{N}}_{P}\left(\psi, T, \alpha\right)\;=\; \sum_{\bm{k}\in\Z^n} \widehat{G}^{\left(P, \psi\right)}_{T, \alpha}\left(\bm{k}\right). 
\end{equation}
Here, $\widehat{G}^{\left(P, \psi\right)}_{T, \alpha}$ denotes the Fourier transform of the map $G^{\left(P, \psi\right)}_{T, \alpha}$ defined in~\eqref{defgppsitalp}. Expli\-ci\-tly,  an elementary change of variables shows that, given $\bm{\xi}\in\R^n$, 
\begin{align}\label{fouriercoeff1}
\widehat{G}^{\left(P, \psi\right)}_{T, \alpha}\left(\bm{\xi}\right)= T^n\cdot \widehat{H}^{\left(P, \psi\right)}_{T, \alpha}\left(\bm{\xi}\right),
\end{align}
where
\begin{align}\label{fouriercoeff1b}
\widehat{H}^{\left(P, \psi\right)}_{T, \alpha}\left(\bm{\xi}\right)\;=\; \int_{\R^n}\e\left(-T\bm{\xi\cdot x}\right)\cdot E_{(T, \alpha)}\left(P\left(\bm{x}\right)\right)\cdot \psi\left(\bm{x}\right)\cdot \textrm{d}\bm{x}.
\end{align}

The term $\widehat{G}^{\left(P, \psi\right)}_{T, \alpha}(\bm{0})$ obtained when $\bm{k}=\bm{0}$ is referred to as the leading term in the sum~\eqref{poissoncountpsiat}, and the remaining sum as the error term. These two quantities are analysed separately, starting with the latter.

\section{Non--stationary Phase Analysis and Estimate of the Error Term}

The following claim essentially establishes that the error term has fast decay in the variable $T$.

\begin{lem}[Asymptotic behavior of the error term]\label{faserrorfourier}
Let $\alpha\in (0,1)$ and $\eta\in (0,1)$. Fix an integer $j\ge 1$ such that $j\eta>n$. Assume that $\psi$ is a smooth map supported in the compact set $\mathcal{K}$ introduced in~\eqref{2emesupprestrispi} and that $T$ is a parameter larger than $2^{1/(1-\eta)}$.  Then, $$\left|\sum_{\underset{\bm{k}\in\Z^n}{\bm{k}\neq\bm{0}}} \widehat{G}^{\left(P, \psi\right)}_{T, \alpha}\left(\bm{k}\right)\right|\;\ll\; C_j\!\left(\psi, \mathcal{K}\right)\cdot T^{n-j\eta}\cdot \log(T), $$where the constant $C_j\!\left(\psi, \mathcal{K}\right)$ is defined in~\eqref{defcnfkbis} and where the implicit constant does not depend   on the smooth map $\psi$.
\end{lem}

The proof of this statement relies on the following effective non-stationary phase estimate which can be found in~\cite[Theorem~7.7.1]{horm}.

\begin{prop}[Effective non--stationary phase estimate]\label{effnonstaphase}
Let $j\ge 1$ be an integer and let $\psi$ be in $\mathcal{C}_{c}^{j}(\R^n)$ with support contained in a compact set $\mathcal{K}$. Assume that the real--valued function $f$ is $j+1$ times continuously differentiable in an open neighbourhood $\mathcal{V}$ of $\mathcal{K}$ . Then, for any real parameter $\lambda>0$, it holds that $$\left|\int_{\R^n}e\left(\lambda f(\bm{x})\right)\cdot\psi(\bm{x})\cdot \textrm{d}\bm{x}\right|\;\ll\; \frac{C_j\left(\psi, \mathcal{K}\right)}{\lambda^j \cdot \delta_\psi( f)^j\cdot \min\left\{1, \delta_\psi( f)^j\right\}}$$ assuming that the quantity $$\delta_\psi(f)\;=\;\inf_{\bm{x}\in\mathcal{K}}\;\left\|\nabla f\left(\bm{x}\right)\right\|$$should be  strictly positive. In the above inequality, the implicit constant is independent from $\psi$ and $\lambda$; furthermore, it remains bounded as long as the function $f$ remains in a bounded subset of the set of $j+1$ times continuously differentiable functions on $\mathcal{V}$.
\end{prop}

\begin{proof}[Proof of Lemma~\ref{faserrorfourier}]
Fix $\bm{k}\neq\bm{0}$ in $\Z^n$. Then, by Fourier inversion, 
\begin{equation}\label{fourierinv} 
\widehat{H}^{\left(P, \psi\right)}_{T, \alpha}\left(\bm{k}\right)\;=\; \int_{-T/(2\Delta_P(\mathcal{K}))}^{T/(2\Delta_P(\mathcal{K}))}\widehat{E}_{(T, \alpha)}(t)\cdot\left(\int_{\R^n}e\left(tP(\bm{x})-T\bm{k\cdot x}\right)\cdot\psi(\bm{x})\cdot \textrm{d}\bm{x}\right)\cdot\textrm{d}t,
\end{equation} 
where the real $\Delta_P(\mathcal{K})$ is defined in~\eqref{defN_P} and $\widehat{H}^{\left(P, \psi\right)}_{T, \alpha}\left(\bm{k}\right)$ in~\eqref{fouriercoeff1b}. Let $\eta\in (0,1)$ and $T\ge 2^{1/(1-\eta)}$ be as in the assumption. Fix $t\in [-T/(2\Delta_P(\mathcal{K})), T/(2\Delta_P(\mathcal{K}))]$. Then, given $\bm{x}\in\textrm{Supp }\psi \subset \mathcal{K}$, 
\begin{align}\label{ineqTktNp}
\left\|t\cdot\nabla P (\bm{x}) -T\bm{k}\right\| \;\ge\; T\cdot \left\|\bm{k}\right\|- \left|t\right|\cdot  \left\|\nabla P (\bm{x})\right\|   \;\underset{\eqref{defN_P}}{\ge}\;T\cdot \left\|\bm{k}\right\|- \left|t\right|\cdot \Delta_P(\mathcal{K}). 
\end{align}
As a consequence, the inequality $\left\|t\cdot\nabla P (\bm{x}) -T\bm{k}\right\| \ge \left(T\cdot \left\|\bm{k}\right\|\right)^{\eta}$ is met as soon as  $$\left\|\bm{k}\right\|\;\ge\; \frac{\left\|\bm{k}\right\|^{\eta}}{T^{1-\eta}}+\frac{1}{2}\cdotp$$In turn, under the assumption that $T\ge 2^{1/(1-\eta)}$, this condition is fulfilled for all $\bm{k}\in\Z^n\backslash\left\{\bm{0}\right\}$. Proposition~\ref{effnonstaphase} (applied with $\lambda=1$) then implies that for any integer $j\ge 1$, 
\begin{equation}\label{fastdecay}
\left|\int_{\R^n}e\left(tP(\bm{x})-T\bm{k\cdot x}\right)\cdot\psi(\bm{x})\cdot \textrm{d}\bm{x}\right|\;\ll\; \frac{C_j\!\left(\psi, \mathcal{K}\right)}{\left(T\cdot \left\|\bm{k}\right\|\right)^{j\eta}}\cdotp
\end{equation} 
Under the assumption $j\eta>n$, one thus obtains that 
\begin{align}
\left|\sum_{\underset{\bm{k}\in\Z^n}{\bm{k}\neq\bm{0}}} \widehat{G}^{\left(P, \psi\right)}_{T, \alpha}\left(\bm{k}\right)\right|\;&\underset{\eqref{fouriercoeff1}}{\le}\; T^n\cdot \sum_{\underset{\bm{k}\in\Z^n}{\bm{k}\neq\bm{0}}} \left|\widehat{H}^{\left(P, \psi\right)}_{T, \alpha}\left(\bm{k}\right)\right|\nonumber \\
&\underset{\eqref{fourierinv}\, \&\, \eqref{fastdecay}}{\ll}\; C_j\!\left(\psi, \mathcal{K}\right)\cdot T^{n-j\eta}\cdot\left( \sum_{\underset{\bm{k}\in\Z^n}{\bm{k}\neq\bm{0}}} \left\|\bm{k}\right\|^{-j\eta}\right)\nonumber\\
&\qquad \qquad \qquad \qquad \qquad \qquad \times \left( \int_{ \left|t\right|\le T/(2\Delta_P(\mathcal{K}))}\left|\widehat{E}_{(T, \alpha)}(t)\right|\cdot\textrm{d}t \right)\nonumber \\
&\ll \; C_j\!\left(\psi, \mathcal{K}\right)\cdot T^{n-j\eta}\cdot\left( \int_{\left|t\right|\le T/(2\Delta_P(\mathcal{K}))}\left|\widehat{E}_{(T, \alpha)}(t)\right|\cdot\textrm{d}t \right). \label{estimehattalph}
\end{align}
To evaluate this last factor, let $t$ be any real lying in the interval $\left[-T/(2\Delta_P(\mathcal{K})), T/(2\Delta_P(\mathcal{K}))\right]$. From the elementary estimate  $$\left|\widehat{\chi}_{I\left(T^{-\alpha}\right)}(t)\right|\; =\; \left|\frac{\sin\left(2\pi T^{-\alpha}t\right)}{\pi t}\right|\;\le\; \min\left\{T^{-\alpha}, \; \frac{1}{\left|\pi t\right|} \right\},$$ one deduces that 
\begin{align*}
\left|\widehat{E}_{(T, \alpha)}(t)\right|\;&\le\; \left|\widehat{\chi}_{I\left(T^{-\alpha}\right)}(t)\right|+\left|\left(\widehat{E}_{(T, \alpha)}-\widehat{\chi}_{I\left(T^{-\alpha}\right)}\right)(t)\right|\\
& \le\; \left|\widehat{\chi}_{I\left(T^{-\alpha}\right)}(t)\right|+\int_{\R} \left|E_{(T, \alpha)}(x)-\chi_{I\left(T^{-\alpha}\right)}(x)\right|\cdot \textrm{d}x\\
&\underset{\eqref{ineqartchiR}\;\&\;\eqref{ineqartchiRbis}}{\le}\; \min\left\{2, \;\frac{1}{\left|t\right|}\right\}+\int_{\R} \theta\left(\frac{T}{2\Delta_P(\mathcal{K})}\cdot\left|\left|x\right|-T^{-\alpha}\right|\right)\cdot \textrm{d}x\\
&\ll \min\left\{2, \;\frac{1}{\left|t\right|}\right\}+\frac{1}{T}
\end{align*}
with an implicit constant depending on the integral of $\theta$ over its domain of definition. Inserting this estimate in~\eqref{estimehattalph} then concludes the proof.
\end{proof}

\section{Analysis of the Leading Term}

The estimation of the leading term $\widehat{G}^{\left(P, \psi\right)}_{T, \alpha}\left(\bm{0}\right)$ in the sum~\eqref{poissoncountpsiat} reduces to the following claim~:

\begin{lem}[Asymptotic behavior of the leading term]\label{asympleading}
Keep the assumptions and the notations of Proposition~\ref{propsmoothcountinglargdom}, assuming in particular that $\alpha\in (0,1)$. Then, there exist integers $M, N\ge 1$  such that for any $\delta\in (0, 1-\alpha)$, it holds that $$\left|\widehat{H}^{\left(P, \psi\right)}_{T, \alpha}\left(\bm{0}\right)-\mu_P\left(\psi, T^{-\alpha}\right)\right|\;\ll\; \left(C_{N}\left(\psi, \mathcal{K}\right)\right)^{M}\cdot \frac{  \mu_P\left(\psi, T^{-\alpha}\right)}{T^\delta}\cdotp$$In this inequality, the quantity $\widehat{H}^{\left(P, \psi\right)}_{T, \alpha}\left(\bm{0}\right)$ is defined in~\eqref{fouriercoeff1b}, the map $\mu_P\left(\psi, \,\cdot\,\right)$ in~\eqref{defmuppsit} and the constant $C_{N}\left(\psi, \mathcal{K}\right)$ in~\eqref{defcnfkbis}. Also, the implicit constant is independent of the smooth map $\psi$.
\end{lem}

\begin{proof}It follows from Theorem~\ref{thmcolzani} that
\begin{align}
&\left|\widehat{H}^{\left(P, \psi\right)}_{T, \alpha}\left(\bm{0}\right) - \mu_P\left(\psi, T^{-\alpha}\right)\right|\nonumber \\
& \qquad \qquad \qquad  \le\; \int_{\R^n}\left|B_{(T/(2\Delta_P(\mathcal{K})), \;T^{-\alpha})}\left(P(\bm{x})\right)- A_{(T/(2\Delta_P(\mathcal{K})), \;T^{-\alpha})}\left(P(\bm{x})\right)\right|\cdot \psi\left(\bm{x}\right)\cdot\textrm{d}\bm{x}\nonumber\\
& \qquad \qquad \qquad  \le\; \int_{\R^n} \theta\left(\frac{T}{2\Delta_P(\mathcal{K})}\cdot\left|\left|P(\bm{x})\right|-T^{-\alpha}\right|\right)\cdot \psi\left(\bm{x}\right)\cdot\textrm{d}\bm{x},\label{inttheors}
\end{align}
where the constant $\Delta_P(\mathcal{K})$ is defined in~\eqref{defN_P}. \\

Let $\delta>0$ stand for some small exponent, the range of which is to be determined. Define the set  $$\mathfrak{S}_P\left(T, \alpha\right)\;=\;\left\{\bm{x}\in\R^n\; :\; \left|\left|P(\bm{x})\right|-T^{-\alpha}\right|\le\frac{2\cdot \Delta_P(\mathcal{K})}{T^{1-\delta}} \right\}$$and decompose the integral in~\eqref{inttheors} as 
\begin{align} 
&\int_{\R^n} \theta\left(\frac{T}{2\Delta_P(\mathcal{K})}\cdot\left|\left|P(\bm{x})\right|-T^{-\alpha}\right|\right)\cdot \psi\left(\bm{x}\right)\cdot\textrm{d}\bm{x}\nonumber \\
&\qquad\qquad =  \int_{\R^n} \chi_{\mathfrak{S}_P\left(T, \alpha\right)}\left(\bm{x}\right)\cdot \theta\left(\frac{T}{2\Delta_P(\mathcal{K})}\cdot\left|\left|P(\bm{x})\right|-T^{-\alpha}\right|\right)\cdot \psi\left(\bm{x}\right)\cdot\textrm{d}\bm{x}\; \nonumber \\
&\qquad\qquad\qquad +\int_{\R^n} \left(1-\chi_{\mathfrak{S}_P\left(T, \alpha\right)}\left(\bm{x}\right)\right)\cdot\theta\left(\frac{T}{2\Delta_P(\mathcal{K})}\cdot\left|\left|P(\bm{x})\right|-T^{-\alpha}\right|\right)\cdot \psi\left(\bm{x}\right)\cdot\textrm{d}\bm{x}.\label{aestimer}
\end{align}
The first integral on the right--hand side of this equation can be bounded as follows~: 
\begin{align}
&\int_{\R^n} \chi_{\mathfrak{S}_P\left(T, \alpha\right)}\left(\bm{x}\right)\cdot \theta\left(\frac{T}{2\Delta_P(\mathcal{K})}\cdot\left|\left|P(\bm{x})\right|-T^{-\alpha}\right|\right)\cdot \psi\left(\bm{x}\right)\cdot\textrm{d}\bm{x}\nonumber\\
&\qquad\qquad\qquad\qquad\qquad \qquad \le \; \left(\max_{x\ge 0}\theta(x)\right)\cdot \int_{\R^n} \chi_{\mathfrak{S}_P\left(T, \alpha\right)}(\bm{x})\cdot \psi\left(\bm{x}\right)\cdot\textrm{d}\bm{x}.\label{firstupbound}
\end{align}
Introducing the map $\mu_P\left(\psi, \,\cdot\,\right) $ defined in~\eqref{defmuppsit}, this last integral becomes
\begin{align}
&\int_{\R^n} \chi_{\mathfrak{S}_P\left(T, \alpha\right)}(\bm{x})\cdot \psi\left(\bm{x}\right)\cdot\textrm{d}\bm{x}\nonumber\\
&\qquad\qquad=\; \mu_P\left(\psi, \frac{1}{T^{\alpha}}+\frac{2\Delta_P(\mathcal{K})}{T^{1-\delta}}\right) - \mu_P\left(\psi, \frac{1}{T^{\alpha}}-\frac{2\Delta_P(\mathcal{K})}{T^{1-\delta}}\right)\nonumber\\
&\qquad\qquad\underset{\eqref{asympmupsansome}\;\&\;\eqref{asympmupsansomebis}}{\ll}\; \left(C_N\left(\psi, \mathcal{K}\right)\right)^M\cdot \frac{1}{T^{1-\delta}}\cdot\frac{\left(\log T\right)^m}{T^{\alpha (\rho-1)}}\label{ineqdelatalp1-}
\end{align}
with an implicit constant depending on the parameters $\alpha$, $\rho, m$ and on the polynomial $P(\bm{x})$. Here, the inequality~\eqref{ineqdelatalp1-} holds as long as 
\begin{equation}\label{rangedelta}
0\;<\;\delta\;<\; 1-\alpha
\end{equation}  
since this condition guarantees that $T^{-1+\delta}=o\left(T^{-\alpha}\right)$ as $T$ tends to infinity. Under this assumption, the upper bound~\eqref{firstupbound} implies that 
\begin{align}\label{firstfinupbbis}
&\int_{\R^n} \chi_{\mathfrak{S}_P\left(T, \alpha\right)}\left(\bm{x}\right)\cdot \theta\left(\frac{T}{2\Delta_P(\mathcal{K})}\cdot\left|\left|P(\bm{x})\right|-T^{-\alpha}\right|\right)\cdot \psi\left(\bm{x}\right)\cdot\textrm{d}\bm{x}\nonumber \\
&\qquad\qquad\qquad\qquad \qquad \qquad \qquad\qquad \ll\; \left(C_N\left(\psi, \mathcal{K}\right)\right)^M\cdot \frac{\left(\log T\right)^m}{T^{1-\delta+\alpha (\rho-1)}}\cdotp
\end{align}
As for the second term on the right--hand side of the equation~\eqref{aestimer}, one infers from the definition of the set $\mathfrak{S}_P\left(T, \alpha\right)$ and from the fast decay of the function $\theta$ that for any large $\beta>0$,
\begin{align}
\int_{\R^n} \left(1-\chi_{\mathfrak{S}_P\left(T, \alpha\right)}\left(\bm{x}\right)\right)\cdot\theta\left(\frac{T}{2\Delta_P(\mathcal{K})}\cdot\left|\left|P(\bm{x})\right|-T^{-\alpha}\right|\right)\cdot \psi\left(\bm{x}\right)\cdot\textrm{d}\bm{x}  \; \ll\; T^{-\beta\delta}\label{dinalfhju}
\end{align}
with an implicit constant depending only on the choice of $\beta$ and on the volume of the compact set $\mathcal{K}$ containing the support of the map $\psi$ (recall that this map is assumed to meet the conditions stated in~\eqref{2emesupprestrispi}).\\

Under the restriction~\eqref{rangedelta}, one can choose $\beta>0$ large enough so that  inequalities~\eqref{firstfinupbbis} and~\eqref{dinalfhju}  imply that $$\left|\widehat{H}^{\left(P, \psi\right)}_{T, \alpha}\left(\bm{0}\right) - \mu_P\left(\psi, T^{\alpha}\right)\right|\;\ll\;\frac{ \left(C_{N}\left(\psi, \mathcal{K}\right)\right)^M}{T^{1-\delta+\alpha(\rho-1)-\eta}}\cdotp$$Since the restriction~\eqref{rangedelta} also guarantees that $1+\alpha(\rho-1)-\delta>\alpha\rho$ when $\alpha<1$, in view of the asymptotic expansion~\eqref{asympmupsansome}, this provides the required power saving  in the error term upon choosing $\eta>0$ small enough. This completes the proof of the lemma.
\end{proof}

\begin{proof}[Completion of the proof of Proposition~\ref{propsmoothcountinglargdom}] In view of the Poisson summation formula~\eqref{poissoncountpsiat}, the sought estimate~\eqref{soughtestimate} is an immediate consequence of  Lemmata~\ref{faserrorfourier} and~\ref{asympleading} when $\alpha\in (0,1)$.
\end{proof}

%

\chapter{Deterministic Counting in Thin Domains}\label{secdeternarrow}

\vspace{10mm}

Throughout this final chapter, let $\alpha>1$ be a real and let $\bm{F}(\bm{x})=\left(F_1(\bm{x}), \dots, F_p(\bm{x})\right)$ be a $p$-tuple  of homogeneous forms in $n\ge 2$ variables, each of degree $d\ge 2$. Assume that  the set $K\subset\R^n$ meets the conditions $(\mathcal{H}_1)$, $(\mathcal{H}_2)$ and $(\mathcal{H}_3)$ stated in the Introduction (Chapter 1) and that it is furthermore semialgebraic (recall that this means that it is given as a finite union of sets that can be defined as finitely many equalities and inequalities involving polynomial maps). Given a large parameter $T\ge 1$, recall also the definitions of the set $\mathcal{S}^{\dag}_{\bm{F}}\left(K,T, \alpha\right)$ and of the corresponding counting function $\mathcal{N}^{\dag}_{\bm{F}}\left(K,T, \alpha\right)$ given in Chapter 1, namely
\begin{align}\label{volFchap6}
\mathcal{S}^{\dag}_{\bm{F}}\left(K,T, \alpha\right) \;=\; \left\{\bm{x}\in T\cdot K \; :\; \left\|\bm{F}(\bm{x})\right\|\le T^{d-\alpha}\right\}
\end{align}
and 
\begin{equation}\label{bnFchap6} 
\mathcal{N}^{\dag}_{\bm{F}}\left(K,T, \alpha\right)\;=\; \#\left(\mathcal{S}^{\dag}_{\bm{F}}\left(K,T, \alpha\right)\cap \Z^n\right).
\end{equation}
The main goal of this final chapter is to establish Theorem~\ref{thmgeneralsarnak} below, which comprises Theorem~\ref{sarnak} stated in the Introduction as a particular case when $\alpha>1$ (the case $\alpha\le 1$   being settled in the previous chapter). \\

To this end, first are recalled further notations and definitions introduced in Chapter~1. When $\bm{v}\in\Sph^{n-1}$ and $\sigma\in\R$, the set $\bm{v}^{\perp}(\sigma)$ denotes the affine hyperplane
\begin{equation}\label{defvperptau}
\bm{v}^{\perp}(\sigma)\; :=\; \left\{\bm{x}\in\R^n \; :\; \bm{v\cdot x}=\sigma\right\}.
\end{equation}
Fixing a compact set $ \mathcal{K}$ containing $K$ in its interior and contained in the open set $U$ defined as part of the support restriction condition~\eqref{suppproper}, let for any given  $\varepsilon>0$
\begin{equation*}\label{defmupcvte}
\mu_{\bm{F}, \mathcal{K}}\left(\bm{v}, \sigma, \varepsilon\right)\;=\; \V_{n-1}\left(\left\{\bm{x}\in \mathcal{K}\cap \bm{v}^{\perp}(\sigma)\; :\; \left\|\bm{F}(\bm{x})\right\|\le\varepsilon\right\}\right)
\end{equation*}
and 
\begin{equation}\label{defmpueps} 
M_{\bm{F}}\left(\mathcal{K}, \varepsilon\right)\;=\; \sup_{\bm{v}\in\Sph^{n-1}}\;\sup_{\sigma\in\R}\;\mu_{\bm{F},\mathcal{K}}\left(\bm{v}, \sigma, \varepsilon\right).
\end{equation}
The "biggest level of flatness" that can be achieved when intersecting the sublevel set $\left\{\bm{x}\in \mathcal{K}\; :\; \left|\bm{F}(\bm{x})\right|\le\varepsilon\right\}$  with affine hyperplanes is quantified by the value of the real
\begin{equation}\label{conditionforall} 
q_{\bm{F}}(\mathcal{K})\;=\; \liminf_{\varepsilon\rightarrow 0^+} \left(\frac{\log M_{\bm{F}}\left(\mathcal{K}, \varepsilon\right)}{\log \varepsilon}\right).
\end{equation}
This is a well--defined real number under the assumption $(\mathcal{H}_2)$.\\

The dimension $\dim\left(\mathcal{Z}_{\mathcal{K}}(\bm{F})\right)$ of  the algebraic variety $\mathcal{Z}_\R(\bm{F})$ over the set $\mathcal{K}$  is denoted by $\tau_{\bm{F}}(\mathcal{K})$ and its codimension by $\widehat{\tau}_{\bm{F}}(\mathcal{K})$, namely by $\widehat{\tau}_{\bm{F}}(\mathcal{K})= n- \dim\left(\mathcal{Z}_{\mathcal{K}}(\bm{F})\right)$. Recall that in the case that the variety $\mathcal{Z}_{\mathcal{K}}(\bm{F})$ has smooth complete intersection, it holds that $\widehat{\tau}_{\bm{F}}(\mathcal{K})=p$.\\

The main result in this chapter can now be stated as follows~:
\begin{thm}\label{thmgeneralsarnak}
Let $\alpha>1$ be a real. Assume that $K$ is a semialgebraic set meeting the assumptions $(\mathcal{H}_1)$, $(\mathcal{H}_2)$ and $(\mathcal{H}_3)$ stated in the Introduction. Let it be contained in the interior of a compact set $\mathcal{K}\subset U$, where  the open set $U$ is defined in the support restriction condition~\eqref{suppproper}. Also, let $r_{\bm{F}}(K)$ be the nonnegative real introduced in the assumption $(\mathcal{H}_3)$. Then, with the above  notations, the following claims are verified~:
 \begin{itemize}
\item[(Z1)] The upper bound
\begin{equation}\label{estimvolsa}
\mathcal{N}^{\dag}_{\bm{F}}\left(K,T, \alpha\right)\;\ll\; \V_n\left(\mathcal{S}^{\dag}_{\bm{F}}\left(K,T, \alpha\right)\right)+T^{\tau_{\bm{F}}(\mathcal{K})-\delta_{\bm{F}}(\mathcal{K}, K)}
\end{equation}
holds for some exponent $\delta_{\bm{F}}(\mathcal{K}, K)>0$ depending on the polynomial $P(\bm{x})$ and on the compact sets $K\subset \mathcal{K}$ if either of the following two mutually exclusive assumptions are satisfied~: 
\begin{itemize}
\item[(i)] in the generic case when 
\begin{equation}\label{cndiestimvolsa1}
q_{\bm{F}}(\mathcal{K})\;>\; r_{\bm{F}}(K)-1,
\end{equation}
one requires that the measure of flatness $q_{\bm{F}}(\mathcal{K})$ should be bounded below  in this sense~:
\begin{equation}\label{condiquP}
q_{\bm{F}}(\mathcal{K})\;>\; n-1-r_{\bm{F}}(K)\cdot\frac{\tau_{\bm{F}}(\mathcal{K})}{\widehat{\tau}_{\bm{F}}(\mathcal{K})}\cdotp
\end{equation}
One can then choose any exponent $\delta_{\bm{F}}(K, \mathcal{K})$ such that
\begin{equation}\label{boundelat}
0\;<\;\delta_{\bm{F}}(K, \mathcal{K})\;<\;  \tau_{\bm{F}}(\mathcal{K})-\frac{n\cdot\left(n-1-q_{\bm{F}}(\mathcal{K})\right)}{n-1-q_{\bm{F}}(\mathcal{K})+r_{\bm{F}}(K)}\cdotp
\end{equation}

\item[(ii)]  in the degenerate case when 
\begin{equation}\label{cndiestimvolsa2}
q_{\bm{F}}(\mathcal{K})\;\le\; r_{\bm{F}}(K)-1,
\end{equation}
one requires that  the measure of flatness $q_{\bm{F}}(\mathcal{K})$ should be bounded below in this sense~: 
\begin{equation}\label{cndiestimvolsa3}
q_{\bm{F}}(\mathcal{K})\;>\; \widehat{\tau}_{\bm{F}}(\mathcal{K})-1.
\end{equation}
One can then choose any exponent $\delta_{\bm{F}}(K, \mathcal{K})$ such that
\begin{equation}\label{boundelatbis}
0\;<\;\delta_{\bm{F}}(\mathcal{K}, K)\;<\;  q_{\bm{F}}(\mathcal{K})-\widehat{\tau}_{\bm{F}}(\mathcal{K})+1.\\
\end{equation}
\end{itemize}

In the inequality~\eqref{estimvolsa},  the volume term on the right--hand side is explicitly determined by Case (2) in Theorem~\ref{volestim}.

\item[(Z2)] In the case that that the algebraic variety $\mathcal{Z}_{\R}(\bm{F})$ is  not contained in a linear subspace of dimension $n-p$ and that, furthermore, the set of homogeneous forms $\bm{F}(\bm{x})$ has smooth complete intersection over $\mathcal{K}$ (recall that this is understood in the sense that the map~\eqref{mapsmoothcominter} does not vanish over $\mathcal{K}$), then the conditions~\eqref{cndiestimvolsa1}  and~\eqref{condiquP} are both met. As a consequence, an estimate of the form~\eqref{estimvolsa} indeed holds.
\end{itemize}
\end{thm}

Claim (Z2) is established in Section~\ref{smoothcasevariety} below. Most of the work in proving the above statement deals with establishing Claim (Z1). To this end, it is enough to consider the case of $p=1$ homogeneous form, say $P(\bm{x})$, of degree denoted by $q\ge 1$. The general case stated in (Z1) then follows upon applying the results obtained in this particular situation to the polynomial $P_{\bm{F}}(\bm{x})=\left\|\bm{F}(\bm{x})\right\|^2$ and upon squaring all relations involving this polynomial. Explicitly, this is saying that it suffices to replace in  what follows the polynomial $P(\bm{x})$ with $P_{\bm{F}}(\bm{x})$, the real $\varepsilon>0$ with $\varepsilon^2$, the integer $q$ with $2d$ and the real $\alpha>1$ with $2\alpha$. The details of this elementary verification are left to the reader.\\

In view of this reduction, it is natural to establish (Z1)  employing the notations corresponding to the analogues in the case $p=1$ of the set $\mathcal{S}^{\dag}_{\bm{F}}\left(K,T, \alpha\right)$ 
and of the counting function $\mathcal{N}^{\dag}_{\bm{F}}\left(K,T, \alpha\right)$
. These are respectively the set $\mathcal{S}^{\dag}_{P}\left(K,T, \alpha\right)$ and the counting function $\mathcal{N}^{\dag}_{P}\left(K,T, \alpha\right)$ introduced in Chapter~\ref{secdeterlarge} --- see the definitions~\eqref{defspalpT} and~\eqref{defspalpTbis}.\\

In order to establish the claim (Z1), let $j\ge 1$ be an integer and let $\varphi$ be in $\mathcal{C}_c^j(\R)$ and $\psi$ in $\mathcal{C}_c^j(\R^n)$ such that 
\begin{equation}\label{bornesphi}
\chi_{[-1, 1]}\;\le\; \varphi \;\le\; \chi_{[-2, 2]}\qquad \textrm{and}\qquad\chi_{K}\;\le\; \psi\;\le\; \chi_{\mathcal{K}}.
\end{equation} 

As the proof relies on semialgebraic geometry, the maps $\varphi$ and $\psi$ are also assumed to be semialgebraic (meaning that their graphs are semialgebraic sets). The existence of such maps is guaranteed, for instance, by~\cite[Proposition~4.8]{kaw}.\\

Fix
\begin{equation}\label{bornevarespi} 
\varepsilon\;\ge\; T^{-\alpha}.
\end{equation} 
It then follows from the Poisson Summation formula and from the homogeneity of the polynomial $P(\bm{x})$ that 
\begin{align}
\mathcal{N}^{\dag}_{P}\left(K,T, \alpha\right)\;& \le\; \sum_{\bm{k}\in\Z^n}\psi\left(\frac{\bm{k}}{T}\right)\cdot \varphi\left(\frac{P(\bm{k})}{\varepsilon\cdot  T^{q}}\right)\nonumber\\
&=\; T^n\cdot \sum_{\bm{k}\in\Z^n}\int_{\R^n}\psi(\bm{x})\cdot\varphi\left(\varepsilon^{-1}\cdot P(\bm{x})\right)\cdot e\left(-T\bm{k\cdot x}\right)\cdot\textrm{d}\bm{x}\nonumber\\
&=\;T^n\cdot \int_{\R^n}\psi(\bm{x})\cdot\varphi\left(\varepsilon^{-1}\cdot P(\bm{x})\right)\cdot\textrm{d}\bm{x}\nonumber\\
&\qquad \qquad + T^n\cdot \sum_{\bm{k}\in\Z^n\backslash\{\bm{0}\}}\int_{\R^n}\psi(\bm{x})\cdot\varphi\left(\varepsilon^{-1}\cdot P(\bm{x})\right)\cdot e\left(-T\bm{k\cdot x}\right)\cdot\textrm{d}\bm{x}.\label{resultingsum}
\end{align}
Upon specialising the real $\varepsilon>0$ and the maps $\varphi$ and $\psi$, the integral in the former term in this last expression is related in Section~\ref{compleproof} to the volume of the set of points $\bm{x}\in K$ such that $\left|P(\bm{x})\right|<T^{-\alpha}$. The analysis of the oscillatory integrals and of the resulting sum in the latter term in~\eqref{resultingsum} constitutes the main substance of the proof and will be carried out  along these lines~: by a non--stationary phase analysis developed in Section~\ref{nonstationarylastsec} below, the sum is first reduced to a finite one up to an error term with rapid decay in the quantity $\left\|T\bm{k}\right\|$, provided that the integer $j\ge 1$ defining the regularity of the maps $\varphi$ and $\psi$ is large enough.  To deal with this finite sum, upon decomposing in polar coordinates the nonzero vector $T\bm{k}$ as 
\begin{equation}\label{decompopolar} 
-T\bm{k}=\lambda \bm{v},\qquad\textrm{where} \qquad \lambda>0\qquad \textrm{and}\qquad \bm{v}\in\Sph^{n-1},
\end{equation}
the goal is to obtain, under suitable assumptions, a uniform decay estimate for the Fourier coefficient appearing in~\eqref{resultingsum}; namely an estimate of the form
\begin{equation}\label{soughtestimatefourier} 
\int_{\R^n}\psi(\bm{x})\cdot\varphi\left(\varepsilon^{-1}\cdot P(\bm{x})\right)\cdot e\left(\lambda\bm{v\cdot x}\right)\cdot\textrm{d}\bm{x}\;\ll\; \frac{\varepsilon^{\nu}}{\lambda^\gamma}
\end{equation}
for some positive exponents $\gamma$ and $\nu$. The crucial requirement here, which prevents the reduction of the problem to a (more or less standard) local harmonic analysis estimate of oscillatory integrals, is that  the implicit constant in the inequality above must be independant of the choice of the unit vector $\bm{v}\in\Sph^{n-1}$. This issue is overcome in Section~\ref{doublegfl} with the help of  o--minimal theory. The sum appearing in~\eqref{resultingsum}, when taken over finitely many nonzero integer vectors $\bm{k}$, can then be estimated upon specialising $\lambda$ to the value $\left\|T\bm{k}\right\|$ in~\eqref{soughtestimatefourier} and upon optimising the value of $\varepsilon$ under the constraint~\eqref{bornevarespi}. This will be seen to imply Theorem~\ref{thmgeneralsarnak}.

\section{Non--Stationary Phase Analysis and Reduction to a Finite Sum}\label{nonstationarylastsec}

This section deals with the problem of reducing the infinite sum appearing in~\eqref{resultingsum} to a finite one through a non--stationary phase analysis.

\begin{lem}\label{lemnonstatchap6}
Let $\beta>0$ be a large parameter and let $\kappa\in (0,1)$ be a small one. Fix $T\ge 2$ and $\varepsilon\in (0, 1/T)$ and assume that  $\varphi\ge 0$ is in $\mathcal{C}_c^j(\R)$ and $\psi\ge 0$  in $\mathcal{C}_c^j(\R^n)$, where 
\begin{equation}\label{boundjreg} 
j\;\ge\; 1+\frac{n+\beta\cdot (1-\kappa)}{\kappa}\cdotp
\end{equation} 
Then, there exists a constant $a_P(\psi, \kappa)>0$ depending on the parameter $\kappa>0$  and  on the sup norm of the polynomial $P(\bm{x})$ over the support of $\psi$ such that 
\begin{align*}
& \sum_{\bm{k}\in\Z^n\backslash\{\bm{0}\}}\int_{\R^n}\psi(\bm{x})\cdot\varphi\left(\varepsilon^{-1}\cdot P(\bm{x})\right)\cdot e\left(-T\bm{k\cdot x}\right)\cdot\textrm{d}\bm{x}\\
&\quad  \quad  \ll\;  \left(\varepsilon T\right)^\beta + \sum_{1\le \left\|\bm{k}\right\|\le  a_P(\psi, \kappa)\cdot \left(T\varepsilon\right)^{-1/(1-\kappa)}}\left|\int_{\R^n}\psi(\bm{x})\cdot\varphi\left(\varepsilon^{-1}\cdot P(\bm{x})\right)\cdot e\left(-T\bm{k\cdot x}\right)\cdot\textrm{d}\bm{x}\right|.
\end{align*}
Here, the implicit constant depends on the parameters $j, \beta, \kappa$ and on the weight functions $\psi$ and $\varphi$.
\end{lem}

\begin{proof}
Let $\bm{k}\in\Z^n\backslash\{\bm{0}\}$ be an integer vector and let $L\ge 1$ be a real. By Fourier inversion, 
\begin{align*}
&\int_{\R^n}\psi(\bm{x})\cdot\varphi\left(\varepsilon^{-1}\cdot P(\bm{x})\right)\cdot e\left(-T\bm{k\cdot x}\right)\cdot\textrm{d}\bm{x} \\
&\qquad\qquad=\; \underbrace{\int_{\left|t\right|\le L}\hat{\varphi}(t)\cdot\left(\int_{\R^n}\psi(\bm{x})\cdot e\left(t\varepsilon^{-1}\cdot P(\bm{x})-T\bm{k\cdot x}\right)\cdot\textrm{d}\bm{x}\right)\cdot\textrm{d}t}_{=I_1(L)} \\
&\qquad\qquad\qquad + \underbrace{\int_{\left|t\right|> L}\hat{\varphi}(t)\cdot\left(\int_{\R^n}\psi(\bm{x})\cdot e\left(t\varepsilon^{-1}\cdot P(\bm{x})-T\bm{k\cdot x}\right)\cdot\textrm{d}\bm{x}\right)\cdot\textrm{d}t}_{=I_2(L)}.
\end{align*}
Since $\varphi$ is in $\mathcal{C}_c^j(\R)$, repeated integration by parts show that its Fourier transform $t\mapsto\hat{\varphi}(t)$ decays as $O\left(t^{-j}\right)$ when $t\rightarrow\infty$. As a consequence, 
\begin{equation}\label{boundi2}
\left|I_2(L)\right|\;\le\; \left(\int_{\R^n}\psi(\bm{x})\cdot\textrm{d}\bm{x}\right)\cdot \left(\int_{\left|t\right|> L}\left|\hat{\varphi}(t)\right|\cdot\textrm{d}t\right)\;\ll\; L^{-j+1}.
\end{equation}

As for the integral $I_1(L)$, it can be estimated with the help of the effective non--stationary phase result stated in Proposition~\ref{effnonstaphase} (upon letting $\lambda=1$ and \mbox{$f(\bm{x})=t\varepsilon^{-1}\cdot P(\bm{x})-T\bm{k\cdot x}$} therein). To this end, given $t\in \left[-L, L\right]$ and $\bm{x}\in\textrm{Supp}\, \psi$, note first that, upon setting $\Delta_P(\mathcal{K})=\max_{\bm{x}\in\mathcal{K}}\left\|\nabla P(\bm{x})\right\|$, it holds that
\begin{equation}\label{boundi1}
\left\|t\varepsilon^{-1}\cdot\nabla P(\bm{x})-T\bm{k}\right\|\;\underset{\eqref{bornesphi}}{\ge}\; T\left\|\bm{k}\right\|-L\varepsilon^{-1}\cdot \Delta_P(\mathcal{K}).
\end{equation}  
Under the assumption that $T\ge 2$, this quantity is larger than $(T\left\|\bm{k}\right\|)^{\kappa}$ whenever 
\begin{equation}\label{boundnormk}
\left\|\bm{k}\right\|\;\ge\; \frac{L\cdot \Delta_P(\mathcal{K})}{T\varepsilon\cdot\left(1-2^{-(1-\kappa)}\right)}\cdotp
\end{equation} 
Proposition~\ref{effnonstaphase} then yields that 
\begin{equation}\label{lastbound98}
\left|I_1(L)\right|\;\ll\; \left(T\left\|\bm{k}\right\|\right)^{-j\kappa}
\end{equation} 
with an implicit constant depending on the map $\psi$ and on the integer $j$. Specialise the inequalities~\eqref{boundi2} and~\eqref{boundi1} to the case that $L=\left\|\bm{k}\right\|^{\kappa}$ and note that the condition~\eqref{boundnormk} then amounts to 
\begin{equation}
\left\|\bm{k}\right\|\;\ge\;  a_P(\psi, \kappa)\cdot\left(T\varepsilon\right)^{-1/(1-\kappa)}, \quad \textrm{where} \quad a_P(\psi, \kappa) \;=\; \left(\frac{\Delta_P(\mathcal{K})}{1-2^{-(1-\kappa)}}\right)^{1/(1-\kappa)}.
\end{equation}
One then obtains that
\begin{align*}
& \sum_{\bm{k}\in\Z^n\backslash\{\bm{0}\}}\int_{\R^n}\psi(\bm{x})\cdot\varphi\left(\varepsilon^{-1}\cdot P(\bm{x})\right)\cdot e\left(-T\bm{k\cdot x}\right)\cdot\textrm{d}\bm{x}\\
&\qquad \qquad \le\;  \sum_{1\le \left\|\bm{k}\right\|\le  a_P(\psi, \kappa)\cdot \left(T\varepsilon\right)^{-1/(1-\kappa)}}\left|\int_{\R^n}\psi(\bm{x})\cdot\varphi\left(\varepsilon^{-1}\cdot P(\bm{x})\right)\cdot e\left(-T\bm{k\cdot x}\right)\cdot\textrm{d}\bm{x}\right|\\
&\qquad \qquad \qquad\qquad\qquad\;+\;  \sum_{\left\|\bm{k}\right\|>  a_P(\psi, \kappa)\cdot \left(T\varepsilon\right)^{-1/(1-\kappa)}}\left(I_1\left(\left\|\bm{k}\right\|^{\kappa}\right)+I_2\left(\left\|\bm{k}\right\|^{\kappa}\right)\right).
\end{align*}
Here, given any integer $j$ satisfying the bound~\eqref{boundjreg}, 
\begin{align*}
\sum_{\left\|\bm{k}\right\|>  a_P(\psi, \kappa)\cdot \left(T\varepsilon\right)^{-1/(1-\kappa)}}\left(I_1\left(\left\|\bm{k}\right\|^{\kappa}\right)+I_2\left(\left\|\bm{k}\right\|^{\kappa}\right)\right)&\underset{\eqref{boundi2}\&\eqref{lastbound98}}{\ll} \sum_{\left\|\bm{k}\right\|>  a_P(\psi, \kappa)\cdot \left(T\varepsilon\right)^{-1/(1-\kappa)}} \left\|\bm{k}\right\|^{-\kappa (j-1)}\\
&\ll \; \left(T\varepsilon\right)^{(\kappa\cdot(j-1)-n)/(1-\kappa)}\\
& \underset{\eqref{boundjreg}}{\ll}\; \left(T\varepsilon\right)^{\beta}
\end{align*}
since $T\varepsilon<1$ by assumption. This  concludes the proof of the lemma.
\end{proof}

\section{Geometric Tomography on Semi--Algebraic Sets}\label{doublegfl}

In view of Lemma~\ref{lemnonstatchap6} and of the discussion held in the introduction to this chapter, the problem is reduced to obtaining uniform decay estimates for the oscillatory integral appearing in the relation~\eqref{resultingsum} when the integer vector $\bm{k}$ varies in a bounded domain. To this end, the integral is first expressed as the Fourier transform of the so--called Gel'fand--Leray function in Section~\ref{subsecdecfour} below. The properties of this function relevant to the proof of Theorem~\ref{thmgeneralsarnak} are then established in the following Section~\ref{domdefgl}.

\subsection{Uniform Fourier Decay of the Gel'fand--Leray Function}\label{subsecdecfour}

Define the \emph{Gel'fand--Leray function} $g_P$ depending on the variables $\left(\bm{v}, \varepsilon, \sigma\right)\in\Sph^{n-1}\times\R^2$ and on the weight functions $\psi$ and $\varphi$ as the map 
\begin{equation}\label{defgammamap0}
\left(\bm{v}, \varepsilon, \sigma\right)\in\Sph^{n-1}\times\R^2\;\mapsto\; \left\langle g_P \left(\bm{v}, \varepsilon, \sigma\right), \left(\varphi, \psi\right)\right\rangle\;=\; \int_{\bm{v}^{\perp}(\sigma)}\psi\left(\bm{x}\right)\cdot \varphi\left(\varepsilon^{-1}\cdot P(\bm{x})\right)\cdot\textrm{d}\bm{x}, 
\end{equation}
where the affine subspace $\bm{v}^{\perp}(\sigma)$ is defined in~\eqref{defvperptau} (this slightly unusual notation for the function $g_P$ takes into account the fact that it depends both on $\varphi$ and $\psi$). It satifies the property that, given $\varepsilon>0$, $\lambda\ge 1$ and $\bm{v}\in\Sph^{n-1}$, the oscillatory integral in~\eqref{resultingsum} can be decomposed after an elementary change of variables as the Fourier transform of the Gel'fand--Leray function as follows~:
\begin{equation}\label{linkoscigelfler}
\int_{\R^n}\psi\left(\bm{x}\right)\cdot\varphi\left(\varepsilon^{-1}\cdot P(\bm{x})\right)\cdot e\left(\lambda\bm{v\cdot x}\right)\cdot\textrm{d}\bm{x}\;=\; \int_{\R}\left\langle g_P \left(\bm{v}, \varepsilon, \sigma\right), \left(\varphi, \psi\right)\right\rangle\cdot e\left(\lambda\sigma\right)\cdot\textrm{d}\sigma.
\end{equation}
This transformation can be seen as a formalisation of the principle of geometric to\-mo\-gra\-phy according to which the properties of a sufficiently nice set can be understood from the knowledge of the volume of their slices with lower dimensional families of sets. \\

The main property of interest satisfied by the Gel'fand--Leray function is contained in the following lemma, which is proved in the next section.

\begin{lem}\label{boundedvargelfler}
Recall that the weight functions $\psi$ and $\varphi$ introduced in~\eqref{bornesphi} are assumed to be semialgebraic. Then, there exists an integer $b_P\left(\varphi, \psi\right)\ge 1$ depending on these weights and on the polynomial $P(\bm{x})$ only such that for any given real number $\varepsilon>0$ and any vector $\bm{v}\in\Sph^{n-1}$, the real line can be  partitioned into at most $b_P\left(\varphi, \psi\right)$ intervals on the interior of which, whenever nonempty, the map 
\begin{equation}\label{defgammamap}
\sigma\in\R\;\mapsto\; \left\langle g_P \left(\bm{v}, \varepsilon, \sigma\right), \left(\varphi, \psi\right)\right\rangle 
\end{equation}
is continuously differentiable and monotonic.
\end{lem}

This immediately implies the main ingredient in the proof of the counting estimate in  Theorem~\ref{thmgeneralsarnak}~:

\begin{coro}\label{coromainestiunif}
Under the assumptions of the previous lemma, given $\varepsilon>0$, set $$M_P\left(\psi, \varphi, \varepsilon\right)\;=\; \sup_{\bm{v}\in\Sph^{n-1}}\;\sup_{ \sigma\in\R}\;  \left\langle g_P \left(\bm{v}, \varepsilon, \sigma\right), \left(\varphi, \psi\right)\right\rangle.$$ Then, it holds that 
\begin{equation}\label{ineqmaingl}
\int_{\R^n}\psi\left(\bm{x}\right)\cdot\varphi\left(\varepsilon^{-1}\cdot P(\bm{x})\right)\cdot e\left(\lambda\bm{v\cdot x}\right)\cdot\textrm{d}\bm{x}\;\ll\;\frac{M_P\left(\psi, \varphi, \varepsilon\right)}{\lambda}
\end{equation}
for some implicit constant independent of the choice of $\bm{v}\in\Sph^{n-1}$ and of $\lambda>0$.
\end{coro}

\begin{proof}
Express the oscillatory integral on the left--hand side of the inequality~\eqref{ineqmaingl} as the Fourier transform of the Gel'fand--Leray function as in the equation~\eqref{linkoscigelfler}. Then, decompose this Fourier transform as 
\begin{align}\label{decointergelfler}
\int_{\R}\left\langle g_P \left(\bm{v}, \varepsilon, \sigma\right), \left(\varphi, \psi\right)\right\rangle\cdot e\left(\lambda\sigma\right)\cdot\textrm{d}\sigma\;=\; \sum_{j=1}^{c}\left(\int_{I_{j}}\left\langle g_P \left(\bm{v}, \varepsilon, \sigma\right), \left(\varphi, \psi\right)\right\rangle\cdot e\left(\lambda\sigma\right)\cdot\textrm{d}\sigma\right).
\end{align}
Here, $c=c_P\left(\varphi, \psi, \bm{v}, \varepsilon\right)$ is at most equal to the integer $b_P\left(\varphi, \psi\right)$ present in the statement of Lemma~\ref{boundedvargelfler} and, for each $1\le j\le c=c_P\left(\varphi, \psi, \bm{v}, \varepsilon\right)$, the set $I_{j}=I^{(P)}_j\left(\varphi, \psi, \bm{v}, \varepsilon\right)$ is an interval with nonempty interior where the restriction of the Gel'fand--Leray function seen as a map in the variable $\sigma$ is continuously differentiable and monotonic. \\

Fixing the integer $j$ and denoting by $d_j=d_j^{(P)}\left(\varphi, \psi, \bm{v}, \varepsilon\right)$ and by $e_j=e^{(P)}_j\left(\varphi, \psi, \bm{v}, \varepsilon\right)$ the endpoints of the interval $I_{j}$, one obtains by integration by parts that 
\begin{align*}
&\left|\int_{I_{j}}\left\langle g_P \left(\bm{v}, \varepsilon, \sigma\right), \left(\varphi, \psi\right)\right\rangle\cdot e\left(\lambda\sigma\right)\cdot\textrm{d}\sigma\right|\\
&=\left|\frac{\left[\left\langle g_P \left(\bm{v}, \varepsilon, \sigma\right), \left(\varphi, \psi\right)\right\rangle\cdot e\left(\lambda\sigma\right)\right]_{\sigma=d_{j}}^{\sigma=e_{j}}}{2i\pi \lambda}-\frac{1}{2i\pi \lambda}\cdot\int_{d_{j}}^{e_{j}}\frac{\textrm{d}}{\textrm{d}\sigma}\left(\left\langle g_P \left(\bm{v}, \varepsilon, \sigma\right), \left(\varphi, \psi\right)\right\rangle\right)\cdot e\left(\lambda\sigma\right)\cdot\textrm{d}\sigma\right|\\
&\le\; \frac{\sup_{\sigma\in I_{j}} \left|\left\langle g_P \left(\bm{v}, \varepsilon, \sigma\right), \left(\varphi, \psi\right)\right\rangle\right|}{\lambda},
\end{align*}
where the last inequality follows from the Triangle Inequality and from the assumption of the monotonicity of the map $\sigma\mapsto \left\langle g_P \left(\bm{v}, \varepsilon, \sigma\right), \left(\varphi, \psi\right)\right\rangle$ restricted to the interval $I_{j}$. 
The definition of the quantity $M_P\left(\psi, \varphi, \varepsilon\right)$ and the decomposition~\eqref{decointergelfler} then yield that 
\begin{align*}
\left|\int_{\R}\left\langle g_P \left(\bm{v}, \varepsilon, \sigma\right), \left(\varphi, \psi\right)\right\rangle\cdot e\left(\lambda\sigma\right)\cdot\textrm{d}\sigma\right|\;&\le\;\frac{c_P\left(\varphi, \psi, \bm{v}, \varepsilon\right)\cdot M_P\left(\psi, \varphi, \varepsilon\right)}{\lambda}\\
&\le\; \frac{b_P\left(\varphi, \psi\right)\cdot M_P\left(\psi, \varphi, \varepsilon\right)}{\lambda},
\end{align*}
which completes the proof of the statement.
\end{proof}

\subsection{Counting Intervals of Monotonicity through O--Minimality}\label{domdefgl}

The goal in this section is to establish Lemma~\ref{boundedvargelfler}. To this end, recall that an \emph{o--minimal structure} is a collection of families of sets $\mathfrak{O}=\left\{\mathcal{O}_k\right\}_{k\ge 0}$ satisfying the following properties~:
\begin{itemize}
\item[(i)] for each $k\ge 0$, the family $\mathcal{O}_k$ is made of subsets of $\R^k$ stable under complementation and union and containing the empty set ;
\item[(ii)] for any integer $k\ge 0$, if $A\in\mathcal{O}_{k}$, then $A\times \R\in \mathcal{O}_{k+1}$ and $\R\times A\in \mathcal{O}_{k+1}$;
\item[(iii)] for any integer $k\ge 0$ and any indices $1\le i<j\le k$, the set $\left\{\left(x_1, \dots, x_k\right)\in\R^k\; :\; x_i=x_j\right\}$ lies in $\mathcal{O}_{k}$;
\item[(iv)] given an integer $k\ge 0$ and a set $A\in\mathcal{O}_{k+1}$, denoting by $\pi~: \R^{k+1}\rightarrow\R^k$ the canonical projection map, it holds that $\pi(A)\in\mathcal{O}_k$;
\item[(v)] Any $x\in\R$ is such that $\left\{x\right\}\in\mathcal{O}_1$; furthermore, $\left\{(x,y)\in\R^2\; :\; x<y\right\}\in\mathcal{O}_2$;
\item[(vi)] the family $\mathcal{O}_1$ consists precisely of the finite unions of points and open intervals in the real line.
\end{itemize}
A set is  \emph{definable} in the structure $\mathfrak{O}$ if it belongs to a family $\mathcal{O}_k$ for some $k\ge 0$, and a map is  definable if so is its graph.\\

The following generic process to generate o--minimal structures is described by Scanlon~\cite{scan}~: given $k\ge 0$, fix a class of so--called \emph{distinguished functions} $\mathcal{F}_k=\left\{f~:\R^k\rightarrow\R\right\}$. An \emph{atomic set} is a set of the form $$\left\{\bm{x}\in\R^k\; :\; a(\bm{x})<b(\bm{x})\right\}\qquad \textrm{or} \qquad \left\{\bm{x}\in\R^k\; :\; a(\bm{x})=b(\bm{x})\right\},$$where $a,b~:\R^k\rightarrow\R$ are any well--defined functions taken as the compositions of the coordinate maps, constant maps and distinguished maps. Letting
\begin{equation*}
\mathcal{F}\;=\; \bigcup_{k\ge 0}\mathcal{F}_k \,,
\end{equation*}
denote by $\mathcal{O}_{\mathcal{F}}$ the smallest collection of subsets of $\R^k$ (as $k\ge 0$ varies) containing the atomic sets and stable under the images of the canonical projections $\pi~: \R^{k+1}\rightarrow\R^k$ and by taking complements and finite unions. Then, $\mathcal{O}_{\mathcal{F}}$  defines an o--minimal structure. \\

The following examples are of particular importance in the proof of Lemma~\ref{boundedvargelfler}~: 
\begin{itemize}
\item[(a)] when $\mathcal{F}:=\mathcal{F}_{alg}$ is the set of all polynomials defined over $\R^k$ for some (varying) $k\ge 0$, the collection $\mathcal{O}_{\mathcal{F}}$ is the set of all \emph{semialgebraic sets}. That  the above condition (iv) should be satisfied is the content of the Tarski--Seidenberg Theorem~\cite[\S 2.10]{dreistame}, the other properties being easily verified.

\item[(b)]  when $\mathcal{F}:=\mathcal{F}_{an}$ is defined as the union of $\mathcal{F}_{alg}$ and of all restricted analytic functions (these are maps $f~:\R^k\rightarrow\R$ vanishing outside $[-1,1]^k$ which are the restriction to $[-1,1]^k$ of a function which is real analytic in a neighbourhood of $[-1, 1]^k$), one obtains the o--minimal structure of \emph{restricted analytic functions}. The proof in this case is due to Denef and van den Dries~\cite{driesdenef}. A detailed description of the sets in $\mathcal{F}_{an}$ is, furthermore, provided in~\cite{driesmel}~: these are precisely the \emph{globally subanalytic} ones; that is, sets  $V\subset\R^k$ such that for each $\bm{x}\in\R^k$, there exists a relatively compact semianalytic set $X$ contained in $\R^{k+l}$ for some integer $l\ge 1$ and a neighbourhood $U$ of $\bm{x}$ such that the intersection $V\cap U$ is the projection of $X$ onto the first $k$ coordinates (recall here that a set is \emph{semianalytic} if it is defined locally around any of its points as a finite sequence of  equations or inequalities involving analytic maps, or a finite union of such sets). 

\item[(c)] when $\mathcal{F}:=\mathcal{F}_{an, exp}$ is defined as the union of $\mathcal{F}_{an}$ and of the singleton consisting of the exponential function $\exp : \R\rightarrow\R$, one obtains the o--minimal structure of \emph{the expansian of restricted analytic functions by the exponential map}.  The proof in this case is due to van den Dries and Miller~\cite{driesmelbis}. It should be noted that the logarithm function $\log :\R_{>0}\rightarrow\R$ is definable in $\mathcal{F}_{an, exp}$ since its graph is the set $\left\{(x,y)\in\R^2\; :\; \left(x>0\right) \land (x=\exp(y))\right\}$.
\end{itemize}

The proof of Lemma~\ref{boundedvargelfler} requires to work with the parametric integral of a semialgebraic map (see the definition of the Gel'fand--Leray function in~\eqref{defgammamap0} and the equation~\eqref{defgammamap}). It turns out that such a parametric integral map need not remain in the structure of semialgebraic sets. To see this, it is enough to consider the semialgebraic function 
$$g(x,y)\;=\; 
\begin{cases}
1/y&\textrm{ if } \left(x>0\right) \land (1<y<x);\\
0 &\textrm{ otherwise.} 
\end{cases}
$$Integrating with respect to the variable $y$ then defines  a parametric integral expressed as the logarithm function, which is not semialgebraic. This is the reason explaining the introduction of the larger structure $\mathcal{F}_{an, exp}$ (which, in particular, comprises the logarithm function).

\begin{proof}[Proof of Lemma~\ref{boundedvargelfler}] Under the assumption that the weights $\psi$ and $\varphi$ are semialgebraic, the  Gel'fand--Leray map $
\left(\bm{v}, \varepsilon, \sigma\right)\in\Sph^{n-1}\times\R^2\;\mapsto\; \left\langle g_P \left(\bm{v}, \varepsilon, \sigma\right), \left(\varphi, \psi\right)\right\rangle$
defined in~\eqref{defgammamap0} is a parametric integral of a globally subanalytic function (since semialgebraic functions are globally subanalytic). The seminal work~\cite{LR2} by Lion and Rolin proves that a parametric integral of a globally subanalytic function, and therefore the Gel'fand--Leray map itself, is definable in the o--minimal structure $\mathcal{F}_{an, exp}$. (Kaiser exhibits in~\cite{kai} a smaller o-minimal structure containing all parametric integrals of semialgebraic functions. The so-called \emph{constructible functions} form a  class of functions containing all such integrals which is optimal  in a sense detailed in Section~\ref{preparationthmgl} below.)\\

The proof relies on the existence of a $\mathcal{C}^1$--cell decomposition of the domain of definition of the  Gel'fand--Leray map. For the purpose of the present considerations, the actual (technical) definition of a cell is irrelevant so that the reader is referred to~\cite[Chap.~3]{dreistame} for further details. It is enough to mention here that these are, in some suitable sense, the "building blocks" of an o--minimal structure (here $\mathcal{F}_{an, exp}$). They are definable and stable under the same operations as this structure (in particular under canonical projections such as in the above property (iv)) and can be used to partition any definable set into finitely many pieces. Furthermore, and this is the above mentioned existence claim, it follows from~\cite[Theorem~7.3.2]{dreistame} that the domain of definition of any definable function can be partitioned into finitely many cells where the function is $\mathcal{C}^1$~: in the case under consideration, this means that there exists a finite collection $\mathcal{D}$ of cells partitioning the set $\Sph^{n-1}\times\R^2$ such that the restriction of the  Gel'fand--Leray map to each cell $D$ in $\mathcal{D}$ is  continuously differentiable. In other words, for each $D\in\mathcal{D}$, there exists a definable open set $U_D$ and a $\mathcal{C}^1$ map $\Gamma_D~: U_D\rightarrow \R$ such that the restrictions to $D$ of $\Gamma_D$ and of the Gel'fand--Leray map $\left\langle g_P , \left(\varphi, \psi\right)\right\rangle$ coincide.\\

Theorem~7.3.2 in~\cite{dreistame} establishes that any definable set can be partitioned into finitely many cells which are $\mathcal{C}^1$ (this is imposing the regularity of the boundary of the cells, see \emph{loc.~cit.}~for further details). Consider then a partition $\mathcal{A}$ of the domain $\Sph^{n-1}\times\R^2$ into $\mathcal{C}^1$--cells inducing for each $D\in\mathcal{D}$  partitions of the sets $$A_D^{+}=\left\{\left(\bm{v}, \varepsilon, \sigma\right)\in D\; :\; \frac{\partial \Gamma_D}{\partial\sigma}\ge 0\right\}\qquad \textrm{and}\qquad A_D^{-}=\left\{\left(\bm{v}, \varepsilon, \sigma\right)\in D\; :\; \frac{\partial \Gamma_D}{\partial\sigma}< 0\right\}.$$These sets are definable as the partial derivatives of a definable function are definable on their domain of definition (see, e.g., \cite[Point~(3), \S 2.3]{wilkie} for a justification of this claim). Since the union over $D\in\mathcal{D}$ of the sets $A_D^{-}$ and $A_D^{+}$ covers the domain $\Sph^{n-1}\times\R^2$, each cell in the family $\mathcal{A}$ lies in some $A_D^{\pm}$.\\

Fix $\left(\bm{v}_0, \varepsilon_0\right)\in\Sph^{n-1}\times\R$ and denote by $\pi : \left(\bm{v}, \varepsilon, \sigma\right)\in\Sph^{n-1}\times\R^2\;\mapsto\; \left(\bm{v}, \varepsilon\right)\in\Sph^{n-1}\times\R$ the canonical projection. Consider the various  sets 
\begin{equation}\label{defAv0epsi0}
A_{\left(\bm{v}_0, \varepsilon_0\right)}\;=\;\left\{\sigma\in\R\; : \; \left(\bm{v}_0, \varepsilon_0, \sigma\right)\in A\right\}
\end{equation} 
obtained as the projection over the real line of the preimages under the projection map $\pi$ of all those cells $A\in\mathcal{A}$  such that $\left(\bm{v}_0, \varepsilon_0\right)\in \pi(A)$. From~\cite[Proposition~3.3.5]{dreistame}, they form  a cell decomposition of the real line denoted by $\mathcal{A}_{\left(\bm{v}_0, \varepsilon_0\right)}$. \\

Consider a set of the form~\eqref{defAv0epsi0} not reduced to a finite union of points and choose  $D\in\mathcal{D}$ such that $A\subset A_D^{\pm}$. Then, over each of the finitely many open intervals making its interior, the map $
\sigma\in\R\;\mapsto\; \left\langle g_P \left(\bm{v}_0, \varepsilon_0, \sigma\right), \left(\varphi, \psi\right)\right\rangle$ is $\mathcal{C}^1$ and monotonic as it coincides therein with the continuously differentiable map $\sigma\mapsto  \Gamma_D\left(\bm{v}_0, \varepsilon_0, \sigma\right)$ whose derivative keeps a constant sign. Lemma~\ref{boundedvargelfler} is then an immediate consequences of the following three observations~:
\begin{itemize}
\item the collection $\mathcal{A}_{\left(\bm{v}_0, \varepsilon_0\right)}$ is contained in $\mathcal{A}$;
\item the collection $\mathcal{A}$ has a finite cardinality;
\item given a set in $\mathcal{A}$, the number of connected components of its projection over $\R$ (with respect to the third coordinate $\sigma\in\R$ ) is finite from the above property (vi).
\end{itemize}
\end{proof}

\section{Completion of the Proof of the Main Theorem}\label{compleproof}

The first part (\S\ref{counbound})  of this antepenultimate section is devoted to the completion of the proof of the claim (Z1) in Theorem~\ref{thmgeneralsarnak} and the second one (\S\ref{smoothcasevariety}) to the proof of (Z2).

\subsection{The Counting Bound}\label{counbound}

The goal is here to establish  the counting bound~\eqref{estimvolsa}, which constitutes the substance of Claim (Z1) in Theorem~\ref{thmgeneralsarnak}. Recall that the discussion following the statement of the theorem makes it clear that it suffices to establish this claim in the case of a single polynomial $P(\bm{x})$.

\begin{proof}[Proof of Claim (Z1) in Theorem~\ref{thmgeneralsarnak}] Let 
\begin{equation}\label{ouestvarpesi}
T^{-\alpha}\;\le \;\varepsilon\; <\; T^{-1},
\end{equation} 
where, by assumption, $\alpha>1$. Let $\beta>0$ be a large parameter and let $\kappa\in (0,1)$ be a small parameter. Fix for the time being two nonnegative compactly supported semialgebraic maps $\varphi$ and $\psi$ which satisfy the inequalities~\eqref{bornesphi}. 
These maps are assumed to be $j\ge 1$ times continuously differentiable, where $j$ is any integer large enough (depending on $n$, $\beta$ and $\kappa$) so that the assumption~\eqref{boundjreg} of the Non--Stationary Phase Lemma~\ref{lemnonstatchap6} holds. Recall that the existence of such maps is guaranteed, e.g., by~\cite[Proposition~4.8]{kaw}.\\

From the Poisson summation estimate~\eqref{resultingsum} and from Lemma~\ref{lemnonstatchap6}, one has that
\begin{align*}
\mathcal{N}^{\dag}_{P}\left(K,T, \alpha\right)\; & \ll\; T^n\cdot \left(\int_{\R^n}\psi(\bm{x})\cdot\varphi\left(\varepsilon^{-1}\cdot P(\bm{x})\right)\cdot\textrm{d}\bm{x}\right)+ T^n \cdot \left(\varepsilon T\right)^\beta\nonumber\\
&+ T^n\cdot\left( \sum_{1\le \left\|\bm{k}\right\|\le  a_P(\psi, \kappa)\cdot \left(T\varepsilon\right)^{-1/(1-\kappa)}}\left|\int_{\R^n}\psi(\bm{x})\cdot\varphi\left(\varepsilon^{-1}\cdot P(\bm{x})\right)\cdot e\left(-T\bm{k\cdot x}\right)\cdot\textrm{d}\bm{x}\right|\right)
\end{align*}
for some implicit constant depending on $j, \beta, \kappa$ and on the weight functions $\psi$ and $\varphi$. In this inequality, Corollary~\ref{coromainestiunif} applied to the polar decomposition~\eqref{decompopolar} implies that for any given  $\bm{k}\in\Z^n\backslash\{\bm{0}\}$, 
\begin{align*}
\int_{\R^n}\psi(\bm{x})\cdot\varphi\left(\varepsilon^{-1}\cdot P(\bm{x})\right)\cdot e\left(-T\bm{k\cdot x}\right)\cdot\textrm{d}\bm{x}\;\ll\; \frac{M_P(\psi, \varphi, \varepsilon)}{T\cdot\left\|\bm{k}\right\|}
\end{align*}
with an implicit constant depending only on the polynomial $P(\bm{x})$ and on the weights $\psi$ and $\varphi$. As a consequence, 
\begin{align*}
\mathcal{N}^{\dag}_{P}\left(K,T, \alpha\right)\ll T^n\cdot &\left(\int_{\R^n}\psi(\bm{x})\cdot\varphi\left(\varepsilon^{-1}\cdot P(\bm{x})\right)\cdot\textrm{d}\bm{x}\right)+ T^n \cdot \left(\varepsilon T\right)^\beta\nonumber\\
& \qquad \qquad \qquad \qquad \qquad \quad \qquad + T^{n-1}\cdot\left(T\varepsilon\right)^{-(n-1)/(1-\kappa)}\cdot M_P(\psi, \varphi, \varepsilon)\nonumber
\end{align*}
in such a way that 
\begin{align}
\mathcal{N}^{\dag}_{P}\left(K,T, \alpha\right) \ll T^n\cdot &\left(\int_{\R^n}\psi(\bm{x})\cdot\varphi\left(\varepsilon^{-1}\cdot P(\bm{x})\right)\cdot\textrm{d}\bm{x}\right)+ T^n \cdot \left(\varepsilon T\right)^\beta\nonumber\\
& \qquad \qquad \qquad \qquad \qquad + T^{-\frac{\kappa\cdot (n-1)}{1-\kappa}}\cdot\varepsilon^{-\frac{n-1}{1-\kappa}}\cdot M_P(\psi, \varphi, \varepsilon).\label{npdagktalpha}
\end{align}
Recall that $K$ and $\mathcal{K}$ are two compact sets such that the former is contained in the interior of the latter and such that the latter is contained in the open set $U$ defined as part of the support restriction condition~\eqref{suppproper}. Given any real $c\in (0,1)$ and any smooth and compactly supported map $\widetilde{\psi}$, recall also the definition of the quantity $\mu_P\left(\widetilde{\psi}, c\right)$ in~\eqref{defmuvolweighted}. 
Choose then   $\widetilde{\psi}$ such that 
\begin{equation}\label{inegcpsitil2mathcakK}
\chi_K\;\le\;\psi\;\le\; \widetilde{\psi}\;\le\; \chi_{\mathcal{K}}.
\end{equation}  
Under the assumption $(\mathcal{H}_3)$ applied to the polynomial $P(\bm{x})$, and with the notation of Lemma~\ref{lemloir}, it holds that $\left(r_P\left(\widetilde{\psi}\right), m_P\left(\widetilde{\psi}\right)\right) = \left(r_P(K), m_P(K)\right)$. Here, under the assumption $(\mathcal{H}_2)$, the point $(X_1)$ in Proposition~\ref{proprlctgene} implies that $r_P(K)>0$.  Then, from Lemma~\ref{lemloir},
\begin{align*}
\int_{\R^n}\psi(\bm{x})\cdot\varphi\left(\varepsilon^{-1}\cdot P(\bm{x})\right)\cdot\textrm{d}\bm{x}\;&\underset{\eqref{bornesphi}}{\le}\;  \int_{\R^n}\widetilde{\psi}(\bm{x})\cdot\chi_{\left\{\left|P(\bm{x})\right|\le 2\varepsilon\right\}}\cdot\textrm{d}\bm{x}\;\underset{\eqref{defmuvolweighted}}{=}\; \mu_P\left(\widetilde{\psi}, 2\varepsilon\right)\\
&\underset{\eqref{asympmu0}}{\ll}\; \varepsilon^{r_P(K)}\cdot\left|\log \varepsilon\right|^{m_P(K)-1}.
\end{align*}
As a consequence, 
\begin{align*}
\mathcal{N}^{\dag}_{P}\left(K,T, \alpha\right) \; \underset{\eqref{npdagktalpha}}{\ll}\; T^n\cdot\varepsilon^{r_P(K)}\cdot\left|\log\varepsilon\right|^{m_P(K)-1}&+ T^n \cdot \left(\varepsilon T\right)^\beta + T^{-\frac{\kappa\cdot (n-1)}{1-\kappa}}\cdot\varepsilon^{-\frac{n-1}{1-\kappa}}\cdot M_P(\psi, \varphi, \varepsilon).\nonumber
\end{align*}
Let $\eta>0$. Provided that $\varepsilon>0$ is small enough, it follows from the definition of the measure of flatness $q_P(\mathcal{K})$ in~\eqref{conditionforall} that 
\begin{align}\label{vilespi}
\mathcal{N}^{\dag}_{P}\left(K,T, \alpha\right)\;  \ll\;  T^n\cdot\varepsilon^{r_P(K)}\cdot\left|\log\varepsilon\right|^{m_P(K)-1}&+ T^n \cdot \left(\varepsilon T\right)^\beta + \varepsilon^{-\frac{n-1}{1-\kappa}+q_P(\mathcal{K})-\eta}.
\end{align}
\\

The goal is first to prove the counting bound~\eqref{estimvolsa} in the generic case (i) of the statement of the theorem. To this end, note that the condition~\eqref{cndiestimvolsa1} implies that, provided that $\eta$  and $\kappa$ are small enough, 
\begin{equation}\label{recyclcondi1}
\rho_P(K, \mathcal{K}, \kappa, \eta)\; :=\; \frac{n}{\frac{n-1}{1-\kappa}-q_P(\mathcal{K})+r_P(K)+\eta}\;>\; 1.
\end{equation}
The counting bound is then established by  a  distinction of cases~:
\begin{itemize}
\item assume that 
\begin{equation}\label{recyclcondi1a}
\alpha\;\le \; \rho_P(K, \mathcal{K}, \kappa, \eta)
\end{equation}
and set $\varepsilon=T^{-\alpha}$. Then, the inequality~\eqref{vilespi} yields
\begin{align}\label{vilespibis}
\mathcal{N}^{\dag}_{P}\left(K,T, \alpha\right)\; \ll\; T^{n-r_P(K)\cdot\alpha}\cdot\left|\log T\right|^{m_P(K)-1}&+ T^{n+ \beta (1-\alpha)}.
\end{align}
Under the assumption that $\alpha>1$, it is enough to choose $\beta>\alpha\cdot r_P(K)/(\alpha-1)$ so that  this sum should be bounded by constant multiple of the first term which, from Case (2) in Theorem~\ref{volestim}, is the volume of the set $\mathcal{S}^{\dag}_{P}\left(K,T, \alpha\right)$. In other words,
\begin{align}\label{vilespibisbis}
\mathcal{N}^{\dag}_{P}\left(K,T, \alpha\right)\; \ll\; \V_n\left(\mathcal{S}^{\dag}_{P}\left(K,T, \alpha\right)\right).
\end{align}

\item Assume now that  the converse to the inequality~\eqref{recyclcondi1a} holds, namely that 
\begin{equation}\label{recyclcondi1b}
\alpha\; > \; \rho_P( K, \mathcal{K}, \kappa, \eta).
\end{equation}
In this case, choose 
\begin{equation}\label{secondchoicevarpesi}
\varepsilon=T^{-\rho_P( K, \mathcal{K}, \kappa, \eta)}
\end{equation} 
in such a way that the inequality~\eqref{vilespi} yields 
\begin{align*}
\mathcal{N}^{\dag}_{P}\left(K,T, \alpha\right)\; \ll\; T^{\widehat{\rho}_P( K, \mathcal{K}, \kappa, \eta) +\gamma}+ T^{n+ \beta(1-\alpha)}
\end{align*}
with $$\widehat{\rho}_P( K, \mathcal{K}, \kappa, \eta)\;:=\;n\cdot\left(1-\frac{r_P(K)}{\frac{n-1}{1-\kappa}-q_P(\mathcal{K})+r_P(K)+\eta}\right)$$for any $\gamma>0$ (this extra factor is to absorb the possible logarithmic contribution in~\eqref{vilespi}). The condition~\eqref{condiquP} ensures that one can choose $\delta_P(K, \mathcal{K})>0$ lying in the interval determined by the inequalities~\eqref{boundelat}. Provided that the parameters $\gamma, \eta$ and $\kappa$ are chosen small enough, it then holds that $$\widehat{\rho}_P(\mathcal{K}, K, \kappa, \eta)+\gamma<\tau_P(\mathcal{K})-\delta_P(K, \mathcal{K}).$$Under the assumption $\alpha>1$ and upon choosing the parameter $\beta$  large enough, this is saying that
\begin{align}\label{vilespiter}
\mathcal{N}^{\dag}_{P}\left(K,T, \alpha\right)\; \ll\; T^{\tau_P(\mathcal{K})-\delta_P(K,\mathcal{K})}.
\end{align}
\end{itemize}
Thus, under the assumptions~\eqref{cndiestimvolsa1} and~\eqref{condiquP}, the inequalities~\eqref{vilespibisbis} and~\eqref{vilespiter} establish the counting bound~\eqref{estimvolsa} for any $\alpha>1$ and for any exponent $\delta_P(K, \mathcal{K})>0$ in the range determined by~\eqref{boundelat}, namely, 
\begin{align}\label{vilespiterquater}
\mathcal{N}^{\dag}_{P}\left(K,T, \alpha\right)\; \ll\;  \V_n\left(\mathcal{S}^{\dag}_{P}\left(K,T, \alpha\right)\right)+T^{\tau_P(\mathcal{K})-\delta_P(K, \mathcal{K})}.
\end{align}
\\

It remains to prove this bound in the degenerate case (ii) of Claim (Z1) in Theorem~\ref{thmgeneralsarnak} when it is assumed that the condition~\eqref{cndiestimvolsa2} holds. The reason why the above argument is not applicable anymore is that the relation~\eqref{recyclcondi1} need not hold any longer when the parameters $\kappa$ and $ \eta$ are small enough. As a consequence, the condition~\eqref{recyclcondi1a} is neither met under the assumption $\alpha>1$ nor can one necessarily set $\varepsilon$ as in~\eqref{secondchoicevarpesi} when the complementary assumption~\eqref{recyclcondi1b} is satisfied  because of the constraints imposed in~\eqref{ouestvarpesi}.\\

To deal with this situation, set $\varepsilon=T^{-(1+\eta)}$ for the same $\eta>0$ as the one fixed above. Then, the inequa\-li\-ty~\eqref{vilespi} yields the upper bound
\begin{align*}
\mathcal{N}^{\dag}_{P}\left(K,T, \alpha\right)\; \ll\; T^{n-r_P(K)\cdot\left(1-\eta/2\right)}+T^{n-\eta\cdot\beta}+T^{\left(1+\eta\right)\cdot\left(\frac{n-1}{1-\kappa}-q_P(\mathcal{K})+\eta\right)}.
\end{align*}
If $\gamma>0$ is given and if $\eta$ and $\kappa$ are chosen suitably small, the assumption of de\-ge\-ne\-ra\-cy~\eqref{cndiestimvolsa2} implies that 
\begin{align*}
\mathcal{N}^{\dag}_{P}\left(K,T, \alpha\right)\; \ll\; T^{n-1-q_P(\mathcal{K})+\gamma}.
\end{align*}
Under the assumption~\eqref{cndiestimvolsa3}, one can fix a real $\delta_P(K, \mathcal{K})>0$ in the range determined by the inequalities~\eqref{boundelatbis}. It is then enough to choose $\gamma>0$ small enough so that $$n-1-q_P(\mathcal{K})+\gamma\;\le\; \tau_P(\mathcal{K})-\delta_P(K, \mathcal{K}).$$This establishes the bound~\eqref{vilespiterquater} in this degenerate case also and thus completes the proof of the claim (Z1) in Theorem~\ref{thmgeneralsarnak}.
\end{proof}

\subsection{The Case of Smooth Varieties}\label{smoothcasevariety}

The goal is now to establish Claim (Z2) of Theorem~\ref{thmgeneralsarnak}. To this end, fix $p\ge 1$ homogeneous forms $F_1(\bm{x}), \dots, F_p(\bm{x})\in\R[\bm{x}]$ such that the smoothness condition given by the nonvanishing of the map~\eqref{mapsmoothcominter} over $\mathcal{K}$ holds. In this case, the dimension of the algebraic variety $\mathcal{Z}_{\R}(\bm{F})$ over $\mathcal{K}$ equals $\tau_{\bm{F}}(\mathcal{K})=n-\widehat{\tau}_{\bm{F}}(\mathcal{K})=n-p$. Furthermore, from  Point (2) in Theorem~\ref{volestim}, one then has that $r_P(K)=p$ (recall here that $K$ is a compact set contained in the interior of $\mathcal{K}\subset U$ and that it contains the set $C$ in its interior, where $C$ and $U$ are defined in the support restriction condition~\eqref{suppproper}). Under these assumptions, the inequalities~\eqref{cndiestimvolsa1}  and~\eqref{condiquP} are the same~: showing that they hold amounts to establishing that
\begin{equation}\label{ineqsmoothcase}
q_{\bm{F}}(\mathcal{K})\;>\; \widehat{\tau}_{\bm{F}}(\mathcal{K})-1\; = \;p-1.\\
\end{equation}

\begin{proof}[Proof that the inequality~\eqref{ineqsmoothcase} holds under the assumptions of the claim (Z2) in Theorem~\ref{thmgeneralsarnak}.] Given $\rho>0$, denote by $\mathfrak{N}_{\bm{F}}\left(\mathcal{K}, \rho\right)$ the $\rho$--tubular neighbourhood of the variety $\mathcal{Z}_{\mathcal{K}}(\bm{F})$~: it is defined as the set of all those points in $\R^n$ lying at a distance less than $\rho$ from $\mathcal{Z}_{\mathcal{K}}(\bm{F})$.\\

Given open sets $\mathcal{O}_1$ and $\mathcal{O}_2$ contained in $\R^n$, the set $\mathcal{O}_1$ satisfies the \emph{unique nearest point property} with respect to $\mathcal{Z}_{\mathcal{O}_2}(\bm{F}):=\mathcal{Z}_{\R}(\bm{F})\cap\mathcal{O}_2$ if for every $\bm{y}\in \mathcal{O}_1$, there exists a unique point $\pi_{\bm{F}}(\bm{y})\in \mathcal{Z}_{\mathcal{O}_2}(\bm{F})$ minimising the distance from $\bm{y}$ to $\mathcal{Z}_{\mathcal{O}_2}(\bm{F})$. The \emph{reach} $\rho_{\bm{F}}(\bm{y})$ of a point $\bm{y}\in \mathcal{Z}_{\mathcal{O}_2}(\bm{F})$ is then defined as the supremum of all those real numbers $\rho>0$ such that the open ball with radius $\rho$ centered at $\bm{y}$ has the unique nearest point property with respect to $\mathcal{Z}_{\mathcal{O}_2}(\bm{F})$.\\ 

Under the assumption that the map~\eqref{mapsmoothcominter} does not vanish over $\mathcal{K}$, and therefore over a small enough neighbourhood $\mathcal{O}$ of it, a classical result by Federer (see~\cite[Theorem~4.12]{fed}) implies that $\rho_{\bm{F}}(\bm{y})>0$ at any point  $\bm{y}\in\mathcal{Z}_{\mathcal{O}}(\bm{F})$.   Since, from~\cite[Remark~4.2]{fed}, the reach map $\bm{y}\in \mathcal{Z}_{\mathcal{K}}(\bm{F})\mapsto \rho_{\bm{F}}(\bm{y})$ is continuous and since the set $\mathcal{K}$ is assumed to be compact, one infers the existence of a real $\rho_{\bm{F}}(\mathcal{K})>0$ such that the projection map
\begin{equation*}\label{projmapbis}
\pi_{\bm{F}} \; :\; \bm{y}\in\mathfrak{N}_{\bm{F}}\left(\mathcal{K}, \rho_{\bm{F}}(\mathcal{K})\right)\;\mapsto\;\pi_{\bm{F}}(\bm{y})\in \mathcal{Z}_\R(\mathcal{K})
\end{equation*}
is well--defined. From there on, establishing the inequality~\eqref{ineqsmoothcase} requires an elementary lemma~:

\begin{lem}[Tubular neigbourhoods are comparable to the set of small solutions to inequalities under the assumption of smoothness]\label{lemtubineg}
Let $\rho_{\bm{F}}(\mathcal{K})>0$ be as above. Define the constants
\begin{equation*}
\alpha_{\bm{F}}(\mathcal{K})= \min_{1\le i \le p}  \alpha_{i}(\mathcal{K})^{-1}\qquad \textrm{and}\qquad \beta_{\bm{F}}(\mathcal{K})\;=\; \max_{1\le i\le p}\left\{1, \beta_i(\mathcal{K})\right\},
\end{equation*}
where for $i=1, \dots, p$, 
\begin{equation*}
\alpha_{i}(\mathcal{K})= \max \left\{\left\|\nabla F_i (\bm{x})\right\|\; :\; \bm{x}\in \overline{\mathfrak{N}_{\bm{F}}\left(\mathcal{K}, \rho_{\bm{F}}(\mathcal{K})\right)}\right\}
\end{equation*}
and 
\begin{equation*}
\beta_i(\mathcal{K})= \max \left\{\vvvert\nabla^2 F_i (\bm{x})\vvvert\; :\; \bm{x}\in \overline{\mathfrak{N}_{\bm{F}}\left(\mathcal{K}, \rho_{\bm{F}}(\mathcal{K})\right)}\right\}.
\end{equation*}
Here, $\nabla^2 F_i (\bm{x})$ denotes the  $n\times n$  Hessian matrix associated to the homogeneous form $F_i$ at the point $\bm{x}$ and $\vvvert\,\cdot\,\vvvert$ denotes the matrix norm induced by the Euclidean norm. Set furthermore
\begin{equation*}
\theta_{\bm{F}}(\mathcal{K})\;=\;\max\left\{1, \left(\max \left\{\vvvert Gr_{\bm{F}}(\bm{x})^{-1}\vvvert \; :\; \bm{x}\in \overline{\mathfrak{N}_{\bm{F}}\left(\mathcal{K}, \rho_{\bm{F}}(\mathcal{K})\right)}\right\}\right)^{1/2}\right\},
\end{equation*}
where $Gr_{\bm{F}}(\bm{x})$ is the Gram $p\times p$ matrix whose $(i,j)^{\textrm{th}}$ entry is the scalar product \mbox{$\nabla F_i(\pi_{\bm{F}}(\bm{x}))\bm{\cdot}\nabla F_j(\pi_{\bm{F}}(\bm{x}))$}, and let
\begin{equation}\label{rhostar}
\rho^*_{\bm{F}}(\mathcal{K})\;=\; \min\left\{\rho_{\bm{F}}(\mathcal{K}), \frac{1}{2\sqrt{n}\cdot\theta_{\bm{F}}(\mathcal{K})\cdot \beta_{\bm{F}}(\mathcal{K})}\right\}\cdotp
\end{equation}

Then, all these constants are well--defined under the smoothness assumption that the map~\eqref{mapsmoothcominter} does not vanish. Moreover, whenever $0<\varepsilon<  \rho^*_{\bm{F}}(\mathcal{K})$, it holds that
\begin{equation}\label{inclutubneigh}
\mathfrak{N}_{\bm{F}}\left(\mathcal{K},\, \alpha_{\bm{F}}(\mathcal{K})\cdot\varepsilon\right) \;\subset\; \mathcal{M}_{\bm{F}}\left(\mathcal{K}, \varepsilon\right)\;\subset\; \mathfrak{N}_{\bm{F}}\left(\mathcal{K},\, 2\cdot\theta_{\bm{F}}(\mathcal{K})\cdot\varepsilon\right),
\end{equation}
where 
\begin{equation*}
\mathcal{M}_{\bm{F}}\left(\mathcal{K}, \varepsilon\right)\;=\; \left\{\bm{x}\in \mathfrak{N}_{\bm{F}}\left(\mathcal{K}, \rho^*_{\bm{F}}(\mathcal{K})\right) \; :\; \left|F_i(\bm{x})\right|\le\varepsilon\;\;\textrm{for all}\;\;1\le i\le p\right\}.
\end{equation*}
\end{lem}

\begin{proof}
Let  $\bm{x}\in \mathfrak{N}_{\bm{F}}\left(\mathcal{K}, \rho^*_{\bm{F}}(\mathcal{K})\right)$. Given $1\le i\le p$, Taylor--Lagrange's Theorem provides the existence of a point $\bm{y}_{\bm{x}}$ lying in the line segment $\left[\pi_{\bm{F}}(\bm{x}), \bm{x}\right]$ such that $$\left|F_i(\bm{x})\right|\;=\; \left|\nabla F_i(\bm{y}_{\bm{x}})\bm{\cdot}\left( \bm{x}-\pi_{\bm{F}}(\bm{x})\right)\right|\;\le\; \left\|\nabla F_i(\bm{y}_{\bm{x}})\right\|\cdot\left\|\bm{x}-\pi_{\bm{F}}(\bm{x})\right\|,$$where the last relation follows from the Cauchy--Schwarz inequality. This is easily seen to imply the first inclusion in~\eqref{inclutubneigh}.\\

As for the second one, fix again a point $\bm{x}\in \mathfrak{N}_{\bm{F}}\left(\mathcal{K}, \rho^*_{\bm{F}}(\mathcal{K})\right)$ and an index $1\le i\le p$. The Taylor--Lagrange Theorem at order two guarantees the existence of a point $\bm{y}_{\bm{x}}\in\left[\pi_{\bm{F}}(\bm{x}), \bm{x}\right]$ such that 
\begin{align}
\left|F_i(\bm{x})\right|\;&=\; \left|\nabla F_i(\pi_{\bm{F}}(\bm{x}))\bm{\cdot}\left( \bm{x}-\pi_{\bm{F}}(\bm{x}\right)\;+\;\frac{1}{2}\cdot\left( \bm{x}-\pi_{\bm{F}}(\bm{x})\right)^T\cdot \nabla^2F_i(\bm{y}_{\bm{x}})\cdot\left( \bm{x}-\pi_{\bm{F}}(\bm{x})\right)\right|\nonumber\\
&\ge\; \left|\nabla F_i(\pi_{\bm{F}}(\bm{x}))\bm{\cdot}\left( \bm{x}-\pi_{\bm{F}}(\bm{x})\right)\right|- \beta_i(\mathcal{K})\cdot \left\|\bm{x}-\pi_{\bm{F}}(\bm{x})\right\|^2,\label{normproj}
\end{align}
 where $\left( \bm{x}-\pi_{\bm{F}}(\bm{x})\right)^T$ denotes the transpose of the (column) vector $\bm{x}-\pi_{\bm{F}}(\bm{x})$. This vector lies in the subspace normal to the manifold $\mathcal{Z}_{\mathcal{K}}(\bm{F})$ at the point $\pi_{\bm{F}}(\bm{x})$. Under the smoothness assumption  that the map~\eqref{mapsmoothcominter} does not vanish, a basis for this normal space is the set of linearly independent vectors $\left(\nabla F_i(\pi_{\bm{F}}(\bm{x}))\right)_{1\le i\le p}$. The smoothness assumption also guarantees that the $p\times p$  Gram matrix $Gr_{\bm{F}}(\bm{x})$ formed from these vectors 
is invertible. \\

The norm of the vector $\bm{x}-\pi_{\bm{F}}(\bm{x})$ lying in the normal subspace under consideration then satisfies the relations
\begin{equation*}
\left\| \bm{x}-\pi_{\bm{F}}(\bm{x})\right\|^2\;=\; \bm{v}_{\bm{F}}(\bm{x})^T\cdot Gr_{\bm{F}}(\bm{x})^{-1} \cdot \bm{v}_{\bm{F}}(\bm{x})\;\le\; \theta_{\bm{F}}(\mathcal{K})^2\cdot \left\|\bm{v}_{\bm{F}}(\bm{x})\right\|^2,
\end{equation*}
where $\bm{v}_{\bm{F}}(\bm{x})$ is the $p$--dimensional (column) vector whose coordinates are the scalars \mbox{$\nabla F_i(\pi_{\bm{F}}(\bm{x}))\bm{\cdot}\left( \bm{x}-\pi_{\bm{F}}(\bm{x})\right)$} for $1\le i\le p$. If one assumes that $\left|F_i(\bm{x})\right|<\varepsilon$ for some  $\varepsilon>0$ and  for all $1\le i\le p$, the inequality~\eqref{normproj} then yields that 
\begin{equation*}
\varepsilon\;>\; \frac{\left\| \bm{x}-\pi_{\bm{F}}(\bm{x})\right\|}{ \theta_{\bm{F}}(\mathcal{K})}-\sqrt{n}\cdot\beta_{\bm{F}}(\mathcal{K})\cdot \left\| \bm{x}-\pi_{\bm{F}}(\bm{x})\right\|^2\;\underset{\eqref{rhostar}}{\ge}\;  \frac{\left\| \bm{x}-\pi_{\bm{F}}(\bm{x})\right\|}{2\cdot  \theta_{\bm{F}}(\mathcal{K})}\cdotp
\end{equation*}
This establishes the second inclusion in~\eqref{inclutubneigh} and thus completes the proof.
\end{proof}

To resume the proof of the inequality~\eqref{ineqsmoothcase}, note that the \L ojasiewicz Inequality (as stated, e.g., in~\cite[Theorem~4.1]{malideals}) guarantees the existence of constants $c>0$ and $\lambda>0$ depending on the compact set $\mathcal{K}$ and on the set of homogeneous forms $\bm{F}(\bm{x})=\left(F_1(\bm{x}), \dots, F_p(\bm{x})\right)$ such that $$\left\|\bm{F}(\bm{x})\right\|\;\ge\; c\cdot\textrm{dist}\left(\bm{x}, \mathcal{Z}_{\R}(\bm{F})\right)^\lambda.$$Here, the notation on the right--hand side stands for the Euclidean distance from a point $\bm{x}\in \mathcal{K}$ to the set $\mathcal{Z}_{\R}(\bm{F})$. This implies that for any $\rho>0$, there exists $\varepsilon_{\bm{F}}\left(\mathcal{K},\rho\right)>0$ such that for all $0<\varepsilon<\varepsilon_{\bm{F}}\left(\mathcal{K}, \rho\right)$, the  inequality $\left\|\bm{F}(\bm{x})\right\|<\varepsilon$ forces the point $\bm{x}\in \mathcal{K}$ to lie in the $\rho$--neighbourhood of the variety $\mathcal{Z}_{\R}(\bm{F})$. Combined with Lemma~\ref{lemtubineg}, one obtains the existence of real numbers $\widetilde{\varepsilon}_{\bm{F}}\left(\mathcal{K}\right)>0$ and $c_{\bm{F}}\left(\mathcal{K}\right)>0$ such that for any $0<\varepsilon<\widetilde{\varepsilon}_{\bm{F}}\left(\mathcal{K}\right)$ and any  $\bm{v}\in\Sph^{n-1}$ and $\sigma\in\R$, 
\begin{equation*}\label{incluslisse}
\left\{\bm{x}\in \mathcal{K}\cap\bm{v}^{\perp}(\sigma)\; :\; \left\|\bm{F}(\bm{x})\right\|<\varepsilon\right\}\;\subset\; \bm{v}^{\perp}(\sigma)\cap\mathcal{M}_{\bm{F}}\left(\mathcal{K},\varepsilon\right)     \;\subset\; \bm{v}^{\perp}(\sigma)\cap \mathfrak{N}_{\bm{F}}\left(\mathcal{K}, c_{\bm{F}}\left(\mathcal{K}\right)\cdot \varepsilon\right).
\end{equation*}
As a consequence, for $\varepsilon>0$ small enough, 
\begin{align}\label{volineqaffinter}
\V_{n-1}\left(\left\{\bm{x}\in \mathcal{K}\cap\bm{v}^{\perp}(\sigma)\; :\; \left\|\bm{F}(\bm{x})\right\|\le\varepsilon\right\}\right)\;&\le\; \V_{n-1}\left(\bm{v}^{\perp}(\sigma)\cap \mathfrak{N}_{\bm{F}}\left(\mathcal{K}, c_{\bm{F}}\left(\mathcal{K}\right)\cdot \varepsilon\right)\right).
\end{align}

To bound from above the volume on the right--hand side of this inequality, denote by $NM_{\bm{F}}(\mathcal{K})$ the normal bundle of the variety $\mathcal{Z}_{\mathcal{K}}(\bm{F})$~: it is the union over $\bm{x}\in\mathcal{Z}_{\mathcal{K}}(\bm{F})$ of the normal spaces $NM_{\bm{F}}(\mathcal{K}, \bm{x})$ where $\left(\bm{x}, \bm{y}\right)\in NM_{\bm{F}}(\mathcal{K}, \bm{x})$ if $\bm{x}\in\mathcal{Z}_{\mathcal{K}}(\bm{F})$ and if $\bm{y}$ lies in the subspace orthogonal to the tangent space to $\mathcal{Z}_{\mathcal{K}}(\bm{F})$ at the point $\bm{x}$. The \emph{Tubular Neighbourhood Theorem} states that  there exists a continous mapping $\bm{x}\in \mathcal{Z}_{\mathcal{K}}(\bm{F})\mapsto\delta(\bm{x})>0$ such that the map $ \left(\bm{x}, \bm{y}\right)\;\mapsto\; \bm{x}+ \bm{y} $ realises a diffeomorphism between the open subset of $NM_{\bm{F}}(\mathcal{K})$ defined as $\left\{\left(\bm{x},\bm{y}\right)\in NM_{\bm{F}}(\mathcal{K})\; :\; \left\|\bm{y}\right\|<\delta(\bm{x})\right\}$ and its image. The proof of this statement in~\cite[Theorem~6.24]{lee} makes it also clear that the mapping $\delta$ is uniformly bounded below by a strictly positive constant under the assumption that the reach of the manifold $\mathcal{Z}_{\mathcal{K}}(\bm{F})$ is strictly positive.\\

One thus obtains for some $\widetilde{\rho}=\widetilde{\rho}_{\bm{F}}(\mathcal{K})>0$ a diffeomorphism 
\begin{equation*}
 \left(\bm{x}, \bm{y}\right)\in NM_{\bm{F}}(\mathcal{K})\left[\widetilde{\rho}\right]\;\mapsto\; \bm{t}\;=\;\bm{x}+ \bm{y} \in \mathfrak{N}_{\bm{F}}\left(\mathcal{K}, \widetilde{\rho}\right)
\end{equation*}
with inverse
\begin{equation}\label{invdiffeo}
\bm{t}\in \mathfrak{N}_{\bm{F}}\left(\mathcal{K}, \widetilde{\rho}\right)\;\mapsto\;   \left(\bm{x}, \bm{y}\right)\;=\; \left(\pi_{\bm{F}}(\bm{t}), \bm{t}- \pi_{\bm{F}}(\bm{t})\right)\in NM_{\bm{F}}(\mathcal{K})\left[\widetilde{\rho}\right].
\end{equation}
Here, one has set $NM_{\bm{F}}(\mathcal{K})\left[\widetilde{\rho}\right]=\left\{\left(\bm{x},\bm{y}\right)\in NM_{\bm{F}}(\mathcal{K})\; :\; \bm{x}\in \mathcal{K}\;\; \textrm{and}\;\;  \left\|\bm{y}\right\|<\widetilde{\rho}\right\}$, where the vector $\bm{y}$ lies in a $p$--dimensional subspace (namely, the subspace normal to the $(n-p)$--dimensional tangent subspace to the variety $\mathcal{Z}_{\mathcal{K}}(\bm{F})$ at a given point). \\

Given $\varepsilon>0$ small enough, up to a multiplicative constant absorbing the sup--norm of the Jacobian of the change of variables~\eqref{invdiffeo}, the volume on the right--hand side of the inequality~\eqref{volineqaffinter} is thus at most the product of two nonnegative quantities~: on the one hand, the (surface) measure of the set of $\bm{x}\in \mathcal{Z}_{\mathcal{K}}(\bm{F})$ lying at a distance less than $c_{\bm{F}}\left(\mathcal{K}\right)\cdot \varepsilon$ from a point in the affine space $\bm{v}^{\perp}(\tau)$ and, on the other, the maximal volume of the intersection of a $p$--dimension ball of radius $c_{\bm{F}}\left(\mathcal{K}\right)\cdot \varepsilon$ with $\bm{v}^{\perp}(\sigma)$. In other words,   
\begin{equation}\label{volineqaffinterbis}
 \V_{n-1}\left(\bm{v}^{\perp}(\sigma)\cap \mathfrak{N}_{\bm{F}}\left(\mathcal{K}, c_{\bm{F}}\left(\mathcal{K}\right)\cdot \varepsilon\right)\right)\;\ll\; \tau_{\bm{F}}\left(\left\{\bm{x}\in \mathcal{K}\; :\; \left|\bm{v\cdot x}-\sigma\right|\le c_{\bm{F}}\left(\mathcal{K}\right)\cdot \varepsilon\right\}\right)\cdot \varepsilon^{p-1},
\end{equation}
where $\tau_{\bm{F}}$ denotes the surface measure of the restriction of the variety $\mathcal{Z}_{\R}(\bm{F})$ to a small enough neighbourhood of $\mathcal{K}$ where the gradient map~\eqref{mapsmoothcominter} does not vanish, and  where the implicit constant in~\eqref{volineqaffinterbis} is independent of $\bm{v}$ and $\sigma$. To complete the proof, it suffices to show that the surface measure of the set under consideration decays as a positive power in the variable $\varepsilon$.\\

From the compacity of $\mathcal{K}$, the problem can be reduced to finitely many similar local estimates which can be dealt with  the theory of $(C, \alpha)$--good functions along the lines of~\cite[Corollary~3]{interf}. As this analytic approach requires rather tedious calculations and a distinction of cases to obtain estimates uniform in the parameters $\bm{v}$ and $\sigma$, the claim is established hereafter with a more global approach.\\

To this end, let $\psi$ be in $\mathcal{C}_c^{\infty}(\R^n)$ such that it is constant equal to 1 over $\mathcal{K}$ and such that its support is contained in the interior of the neighbourhood of $\mathcal{K}$ used to define $\tau_{\bm{F}}$. Define then the volume element
\begin{equation}\label{ineqmeas}
\textrm{d} \mu_{\bm{F}}\;=\; \psi\cdot\textrm{d}\tau_{\bm{F}}.
\end{equation} 
In particular, the induced measure $\mu_{\bm{F}}$ is absolutely continuous with respect to the measure $\tau_{\bm{F}}$. \\

Under the assumption that the variety $\mathcal{Z}_{\mathcal{K}}(\bm{F})$ is not contained in a  linear subspace of dimension $n-p$, it follows from~\cite[p.351, Point (iii) \& Theorem~2]{stein} that there exists $\kappa>0$ such that the Fourier transform of the measure $\mu_{\bm{F}}$ satisfies the decay  property
\begin{equation}\label{decayfouriermeas}
\left|\widehat{\mu}_{\bm{F}}\left(\bm{\xi}\right)\right|\;\ll\; \left\|\bm{\xi}\right\|^{-\kappa}
\end{equation}
for all $\bm{\xi}\in\R^n\backslash\left\{\bm{0}\right\}$ and for some implicit constant depending on the choice of $\psi$. From the Esséen Concentration Inequality~\cite[\S 2.2.11]{taoess}, it then holds that for any $\bm{v}\in\Sph^{n-1}$ and $\sigma\in\R$, 
\begin{equation*}
\mu_{\bm{F}}\left(\left\{\bm{x}\in \mathcal{K}\; :\; \left|\bm{v\cdot x}-\sigma\right|\le c_{\bm{F}}\left(\mathcal{K}\right)\cdot \varepsilon\right\}\right)\;\ll\; \left(c_{\bm{F}}\left(\mathcal{K}\right)\cdot \varepsilon\right)\cdot\left(\int_{|t|\le \left(c_{\bm{F}}\left(\mathcal{K}\right)\cdot \varepsilon\right)^{-1}}\left|\widehat{\mu}_{\bm{F}}\left(t\bm{v}\right)\right|\cdot\textrm{d}t\right)
\end{equation*}  
for some implicit constant independent of the choice of $\bm{v}$, $\sigma$ and $\varepsilon$. Upon combining this relation with~\eqref{ineqmeas} and~\eqref{decayfouriermeas}, one obtains from the inequality~\eqref{volineqaffinterbis} that for some $\kappa>0$,
 \begin{equation*}
 \V_{n-1}\left(\bm{v}^{\perp}(\sigma)\cap \mathfrak{N}_{\bm{F}}\left(\mathcal{K}, c_{\bm{F}}\left(\mathcal{K}\right)\cdot \varepsilon\right)\right)\;\ll\; \varepsilon^{p-1+\kappa}.
\end{equation*}
As a consequence, the definition of the quantity $M_{\bm{F}}\left(\mathcal{K}, \varepsilon\right)$ in~\eqref{defmpueps} implies that $$M_{\bm{F}}\left(\mathcal{K}, \varepsilon\right)\;\ll\; \varepsilon^{p-1+\kappa}$$ in such a way that $$q_{\bm{F}}\left(\mathcal{K}\right)\; :=\; \liminf_{\varepsilon\rightarrow 0^+}\left(\frac{\log M_{\bm{F}}\left(\mathcal{K}, \varepsilon\right)}{\log \varepsilon}\right)\;\ge\; p-1+\kappa\; >\; p-1.$$ This completes the proof of the inequality~\eqref{ineqsmoothcase}.
\end{proof}

\section{Determining the Measure of Flatness}\label{preparationthmgl} 

The practical use of Theorem~\ref{thmgeneralsarnak} relies on the determination of the measure of flatness introduced in~\eqref{conditionforall}. This final section is devoted to the development of the theory ena\-bling one to determine effectivity the value of this quantity. To this end, following the discussion held after the statement of Theorem~\ref{thmgeneralsarnak}, consider, without loss of generality, the case of a single real homogeneous polynomial $P(\bm{x})$ of degree $q\ge 2$. One is thus interested in finding effectively the value of the real
\begin{equation}\label{conditionforallbis} 
q_{P}(\mathcal{K})\;=\; \liminf_{\varepsilon\rightarrow 0^+} \left(\frac{\log M_{P}\left(\mathcal{K}, \varepsilon\right)}{\log \varepsilon}\right), 
\end{equation}
where 
\begin{equation*}\label{defmpuepsbis} 
M_{P}\left(\mathcal{K}, \varepsilon\right)\;=\; \sup_{\bm{v}\in\Sph^{n-1}}\;\sup_{\sigma\in\R}\;\mu_{P,\mathcal{K}}\left(\bm{v}, \sigma, \varepsilon\right)
\end{equation*}
and
\begin{equation*}\label{defmupcvtebis}
\mu_{P, \mathcal{K}}\left(\bm{v}, \sigma, \varepsilon\right)\;=\; \V_{n-1}\left(\left\{\bm{x}\in \mathcal{K}\cap \bm{v}^{\perp}(\sigma)\; :\; \left|P(\bm{x})\right|\le\varepsilon\right\}\right).
\end{equation*}
Here, $\mathcal{K}$ is assumed to be any compact semialgebraic set with nonempty interior intersecting non trivially the algebraic variety $\mathcal{Z}_{\R}(P)$. Under these assumptions, the real $q_{P}(\mathcal{K})$ is well-defined.\\

Fix a real number $\varepsilon>0$, a vector $\bm{v}\in\Sph^{n-1}$ and a scalar $\sigma\in\R$. The first step is to make explicit the above-defined volume $\mu_{P, \mathcal{K}}\left(\bm{v},\sigma,\varepsilon\right)$ of the slice of the sublevel set $\left\{\bm{x}\in \mathcal{K}\; :\; \left|P(\bm{x})\right|\le \varepsilon\right\}$ with the affine space $\bm{v}^{\perp}(\sigma)$. To this end, denote by $R_{\bm{v}}$ the rotation matrix mapping the last element $\bm{e}_n\in\Sph^{n-1}$ of the canonical basis onto the unit vector $\bm{v}\in\Sph^{n-1}$  and leaving the orthogonal of the plane spanned by $\bm{e}_n$ and $\bm{v}$ unchanged (in the degenerate case when $\bm{v}=\bm{e}_n$, the matrix $R_{\bm{v}}$ is to be taken as the identity). Then, one can decompose the vector $\bm{x}\in \mathcal{K}\cap\bm{v}^{\perp}(\sigma)$ as $\bm{x}=R_{\bm{v}}\left(\bm{y}, \sigma\right)$, where $\bm{y}\in\R^{n-1}$  lies in the projection onto the first $n-1$ coordinates of the preimage of the semialgebraic set $\mathcal{K}$ by $R_{\bm{v}}^{-1}$. One is thus reduced to computing the $n-1$ volume of the set of vectors $\bm{y}$ such that 
\begin{equation}
\left|\left(P\circ R_{\bm{v}}\right)\left(\bm{y}, \sigma\right)\right|\;\le\;\varepsilon\qquad \textrm{with}\qquad \left(\bm{y}, \sigma\right)\in R_{\bm{v}}^{-1}\left(\mathcal{K}\right).
\end{equation}
This is achieved in the following statement.

\begin{prop}\label{propvol}
Keep the above notations and recall that the set $\mathcal{K}\subset\R^{n}$ is assumed to be bounded and semialgebraic. Then, there exists a semialgebraic set $$S_n(\mathcal{K})\subset \R^{n-2}\times\R\times \Sph^{n-1}\times (0,1)$$depending on $\mathcal{K}$ which  can be partitioned into $N\ge 1$ semialgebraic subsets $\mathcal{D}_j$ (where $1\le j\le N$) such that the following property holds~: over each of these subsets are defined $l(j)\ge 0$ analytic semialgebraic maps $\xi_1<\xi_2<\dots <\xi_{l(j)}\; :\; \mathcal{D}_j\rightarrow \R$ such that, defining for any given $\left(\sigma, \bm{v}, \varepsilon\right)\in \R \times \Sph^{n-1}\times \left(0, 1\right)$  the projected sets $$\mathcal{D}_j\left(\bm{v}, \sigma, \varepsilon\right)\;=\;\left\{\bm{t}\in\R^{n-2}\; :\; \left(\bm{t}, \bm{v}, \sigma, \varepsilon\right)\in \mathcal{D}_j\right\},$$ the volume $\mu_{P, \mathcal{K}}\left(\bm{v},\sigma,\varepsilon\right)$ is the sum of finitely many integrals of the form $$\int_{\mathcal{D}_j\left(\bm{v}, \sigma, \varepsilon\right)}\left(\xi_{i+1}\left(\bm{t}, \bm{v}, \sigma, \varepsilon\right)-\xi_i\left(\bm{t}, \bm{v}, \sigma, \varepsilon \right)\right)\cdot\textrm{d}\bm{t}.$$ Here, $1\le j\le N$ and $1\le i \le l(j)$ whenever $l(j)\ge 1$. 
\end{prop}

The definition of the semialgebraic set $S_n(\mathcal{K})$ is explicitly given in the proof of the proposition, which makes it also clear that it can be constructed effectively. This proof can be inferred from the following result  essentially  due to Coste~\cite{coste}~:

\begin{lem}\label{lemcoste}
Let $Q\left(\bm{x}\right)\in\R\left[\bm{x}\right]$ be a polynomial in $n\ge 1$ variables. Assume that $B\subset\R^{n-1}$ is a connected semialgebraic set and that $k$ and $d$ are integers such that $k\le d$ and such that for every point $\widehat{x}\in B$, the polynomial $Q\left(x_1, \widehat{x}\right)$ has degree $d$ in the variable $x_1$ and exactly $k$ distinct roots  in $\C$. Then, there exist $l\le k$ analytic semialgebraic maps $\zeta_1<\dots<\zeta_l~: B\rightarrow\R$ which can be effectively determined such that, for every $\widehat{x}\in B$, the set of real roots of the polynomial $Q\left(x_1, \widehat{x}\right)$ is exactly $\left\{\zeta_1\left(\widehat{x}\right), \dots, \zeta_l\left(\widehat{x}\right)\right\}$. Furthermore, for any $i=1, \dots, l$, the multiplicity of the root $\zeta_i\left(\widehat{x}\right)$ is constant for $\widehat{x}\in B$.
\end{lem}

\begin{proof}
This is~\cite[Proposition~2.6]{coste} with the exception that it is shown there that the maps $\xi_i$ are just continuous rather than analytic. This additional feature is obtained as a consequence of the continuity of the roots of a polynomial as functions of the coefficients. A stronger statement established in~\cite{roots} is true~: the distinct roots of a polynomial are analytic functions of the coefficients in the open set where the roots retain their multiplicities. This yields the above statement.
\end{proof}

\begin{proof}[Deduction of Proposition~\ref{propvol} from Lemma~\ref{lemcoste}] 
Consider the real polynomial $$Q(\bm{y}, \sigma, \bm{v}, \varepsilon)\;=\;\left(\left(P\circ R_{\bm{v}}\right)\left(\bm{y}, \sigma\right)\right)^2-\varepsilon^2.$$Even if it means reindexing the coordinates of the vector $\bm{y}=\left(y_1, \dots, y_{n-1}\right)\in\R^{n-1}$, assume without loss of generality that the variable $y_1$ appears therein explicitly. Then, denote by $S_n(\mathcal{K})$ the image of the set
\begin{equation}\label{semivol}
\left\{\left(\bm{y}, \sigma, \bm{v}, \varepsilon\right)\in \R^{n-1}\times\R\times \Sph^{n-1}\times (0,1)\; :\; \left(\bm{y}, \sigma\right)\in R_{\bm{v}}^{-1}\left(\mathcal{K}\right)\right\}
\end{equation}
under the projection map
\begin{equation*}\label{projmap}
\left(\bm{y}, \sigma, \bm{v}, \varepsilon\right)\in \R^{n-1}\times\R\times \Sph^{n-1}\times (0,1)\;\mapsto\; \left(\widehat{\bm{y}}, \sigma, \bm{v}, \varepsilon\right)\in \R^{n-2}\times\R\times \Sph^{n-1}\times (0,1),
\end{equation*}
where $\widehat{\bm{y}}=\left(y_2, \dots, y_{n-1}\right)\in\R^{n-2}$ when $\bm{y}=\left(y_1, y_2, \dots, y_{n-1}\right)\in\R^{n-1}$. From the Tarski-Seidenberg Theorem (see Point (a) in \S \ref{domdefgl}), the set $S_n(\mathcal{K})$ is also semialgebraic. \\

Partition $S_n(\mathcal{K})$   into finitely many semialgebraic subsets such that over each of them, the polynomial $Q$ seen as a function of the variable $y_1$ alone  with coefficients in $S_n(\mathcal{K})$  keeps a constant degree with its complex roots retaining their multiplicity. This is indeed possible since the degree of $Q$ in the variable $y_1$ is determined by mutually exclusive semialgebraic conditions (namely, the vanishing of a suitable set of coefficients) and so is the multiplicity of its set of roots (namely, the vanishing of discriminants polynomials whose variables are the set of coefficients). \\

From~\cite[Theorem~2.23]{coste}, each of the elements of the semialgebraic partition thus obtained has finitely many semialgebraic connected components. Lemma~\ref{lemcoste} then implies that over each such connected component, there exists a finite number (possibly zero) of semialgebraic analytic maps describing the roots of the polynomial $Q$ (in the variable $y_1$)  as functions of the coefficients $\left(\widehat{\bm{y}}, \sigma, \bm{v}, \varepsilon\right)$. Denote this finite number of maps by $\zeta_2, \dots, \zeta_{l-1}$ and define corresponding maps $\xi_2, \dots, \xi_{l-1}$ in the following way~: when $2\le i\le l-1$, $\xi_i$ is the restriction of $\zeta_i$ to the preimage $\zeta_i^{-1}\left(\pi_1(R_{v}^{-1}(\mathcal{K}))\right)$, where $\pi_1(R_{v}^{-1}(\mathcal{K}))$ is the semialgebraic image of $R_{v}^{-1}(\mathcal{K})$ under the projection map $\left(\bm{y}, \sigma\right)\mapsto y_1$. Since  from~\cite[Corollary~ 2.9]{coste}, the preimage of a semialgebraic set under a semialgebraic map is also semialgebraic, the maps $\xi_i$, where $1\le i\le l-1$, remain semialgebraic. Then, by definition, for all  $2\le i \le l-1$, the point $ \left(\left(\xi_i\left(\widehat{\bm{y}}, \sigma, \bm{v}, \varepsilon\right), \widehat{\bm{y}}\right), \sigma, \bm{v}, \varepsilon\right)$ lies in the set~\eqref{semivol}. Even if it means further  considering separately the restrictions of the maps $\xi_i$ ($2\le i \le l-1$)  to the (finite number of) connected components of their common domain of definition, they may be assumed to be defined over connected semialgebraic sets. \\

Let then $\xi_1$ and $\xi_l$ be the infimum and the supremum of the admissible values of $y_1$, respectively, when $\left(\widehat{\bm{y}}, \sigma, \bm{v}, \varepsilon\right)$ belongs to the connected semialgebraic set, say $\mathcal{D}$, defining the domain of definition of the maps $\xi_i$, where $2\le i \le l-1$ (these infimum and supremum are finite under the assumption that the set $\mathcal{K}$ is bounded). The Mean Value Theorem implies that if $Q(y_1, \widehat{\bm{y}}, \bm{v}, \sigma, \varepsilon)\le 0$ for some $(\left(y_1, \widehat{\bm{y}}\right), \sigma, \bm{v}, \varepsilon)$ lying in the set~\eqref{semivol}, there exists an index $1\le i\le l$ such that  $\xi_i\left(\widehat{\bm{y}}, \sigma, \bm{v}, \varepsilon\right)\le y_1\le \xi_{i+1}\left(\widehat{\bm{y}}, \sigma, \bm{v}, \varepsilon\right)$. \\

Fixing the components $\bm{v}, \sigma$ and $\varepsilon$, it follows from the Tarski--Seidenberg Theorem~\cite[\S 2.10]{dreistame} that the projected set $\mathcal{D}\left(\bm{v}, \sigma, \varepsilon\right) :=\left\{\widehat{\bm{y}}\in \R^{n-2}\; :\; \left(\widehat{\bm{y}}, \sigma,  \bm{v}, \varepsilon\right)\in \mathcal{D}\right\}$ is semialgebraic. Assuming that its $(n-2)$--dimensional Lebesgue measure is positive, the definition of the polynomial $Q$ yields that the measure of the set of points $\bm{y}\in\R^{n-1}$ such that $\left|\left(P\circ R_{\bm{v}}\right)\left(\bm{y}, \sigma\right)\right|\;\le\;\varepsilon$ with $\left(\bm{y}, \sigma\right)\in R_{\bm{v}}^{-1}\left(\mathcal{K}\right)$ is a finite sum of integrals of the form $$\int_{\mathcal{D}\left(\bm{v}, \sigma, \varepsilon\right)}\left(\xi_{i+1}\left(\widehat{\bm{y}}, \bm{v}, \sigma, \varepsilon\right)-\xi_i\left(\widehat{\bm{y}}, \bm{v}, \sigma, \varepsilon \right)\right)\cdot\textrm{d}\widehat{\bm{y}}.$$ The proof is complete upon indexing the (finitely many) semialgebraic sets $\mathcal{D}\left(\bm{v}, \sigma, \varepsilon\right)$ which have positive Lebesgue measure.
\end{proof}

The next step is to show that, given $\varepsilon>0$, the suprema defining the real $M_P(\mathcal{K}, \varepsilon)$ in~\eqref{defmpueps} are attained at some point $\left(\bm{v}, \sigma\right)=\left(\bm{v}(\varepsilon), \sigma(\varepsilon)\right)\in\Sph^{n-1}\times \R$. Combined with Proposition~\ref{propvol}, this implies that the value of $M_P(\mathcal{K}, \varepsilon)$ can be determined in a finite number of steps.

\begin{prop}\label{propcompact} 
Recall that the set $\mathcal{K}\subset\R^{n}$ is assumed to be compact. Fix $\varepsilon>0$. Then, there exists $\left(\bm{v}(\varepsilon), \sigma(\varepsilon)\right)\in\Sph^{n-1}\times \R$ such that $$M_P\left(\mathcal{K}, \varepsilon\right)\;=\; \mu_{P, \mathcal{K}}\left(\bm{v}(\varepsilon),\sigma(\varepsilon),\varepsilon\right).$$
\end{prop}

To prove the statement, it is convenient to introduce, given a value of $\varepsilon>0$, the map
\begin{equation}\label{mapvarpesmu}
\mu_{P, \mathcal{K}}\left(\; \cdot\;,\; \cdot\;,\varepsilon\right)\;:\; \left(\bm{v}, \sigma\right)\mapsto \mu_{P, \mathcal{K}}\left(\bm{v},\sigma,\varepsilon\right)\;=\; \int_{\R^{n-1}}\chi_{\left\{ \left|\left(P\circ R_{\bm{v}}\right)\left(\bm{y}, \sigma\right)\right|\;\le\;\varepsilon \right\}}\cdot \chi_{R_{\bm{v}}^{-1}\left(\mathcal{K}\right)}\left(\bm{y}, \sigma\right)\cdot \textrm{d}\bm{y}.
\end{equation}

\begin{proof}
Since the set $\mathcal{K}$ is bounded, the support of the map $\mu_{P, \mathcal{K}}\left(\; \cdot\;,\; \cdot\;,\varepsilon\right)$ is compact. To establish the claim, it is therefore enough to show that this map is upper semicontinuous (so that it then assumes a maximal value over its support).\\

Let then $\left(\bm{v}_k, \sigma_k\right)_{k\ge 1}$ be a sequence in the support of  $\mu_{P, \mathcal{K}}\left(\; \cdot\;,\; \cdot\;,\varepsilon\right)$ converging to a point $\left(\bm{v}_0, \sigma_0\right)$ in this support. The goal is to prove that 
\begin{equation}\label{upsimicionti}
\limsup_{k\rightarrow\infty } \mu_{P, \mathcal{K}}\left(\bm{v}_k, \sigma_k,\varepsilon\right)\;\le\; \mu_{P, \mathcal{K}}\left(\bm{v}_0, \sigma_0,\varepsilon\right).
\end{equation}
To this end,  one infers first from Fatou's Lemma (which is applicable from the assumption that the set $ \mathcal{K}$ is bounded) that $$\limsup_{k\rightarrow\infty } \mu_{P, \mathcal{K}}\left(\bm{v}_k, \sigma_k,\varepsilon\right)\;\le\; \int_{\R^{n-1}}\limsup_{k\rightarrow \infty}\left(\chi_{\left\{ \left|\left(P\circ R_{\bm{v}_k}\right)\left(\bm{y}, \sigma_k\right)\right|\;\le\;\varepsilon \right\}}\cdot \chi_{R_{\bm{v}_k}^{-1}\left(\mathcal{K}\right)}\left(\bm{y}, \sigma_k\right)\right)\cdot \textrm{d}\bm{y}.$$Fix $\bm{y}\in \R^{n-1}$ and let $\left(\bm{v}_{k_l}, \sigma_{k_l}\right)_{l\ge 1}$ be a subsequence (depending on $\bm{y}\in \R^{n-1}$) realising the upper limit of the sequence of maps inside the integral on the right--hand side. By the continuity of the polynomial map defined by $P(\bm{x})$ and from the closedness of the set $\mathcal{K}$, upon taking limits as the integer $l\ge 1$ tends to infinity in the relations $\left|\left(P\circ R_{\bm{v}_{k_l}}\right)\left(\bm{y}, \sigma_{k_l}\right)\right|\;\le\;\varepsilon $ and $R_{\bm{v}_{k_l}}\!\left(\bm{y}, \sigma_{k_l}\right)\in \mathcal{K}$, it holds that $\left|\left(P\circ R_{\bm{v}_{0}}\right)\left(\bm{y}, \sigma_{0}\right)\right|\;\le\;\varepsilon $ and $R_{\bm{v}_{0}}\left(\bm{y}, \sigma_0\right)\in \mathcal{K}$. As a consequence, $$\limsup_{k\rightarrow \infty}\left(\chi_{\left\{ \left|\left(P\circ R_{\bm{v}_k}\right)\left(\bm{y}, \sigma_k\right)\right|\;\le\;\varepsilon \right\}}\cdot \chi_{R_{\bm{v}_k}^{-1}\left(\mathcal{K}\right)}\left(\bm{y}, \sigma_k\right)\right)\;\le\; \chi_{\left\{ \left|\left(P\circ R_{\bm{v}_0}\right)\left(\bm{y}, \sigma_k\right)\right|\;\le\;\varepsilon \right\}}\cdot \chi_{R_{\bm{v}_0}^{-1}\left(\mathcal{K}\right)}\left(\bm{y}, \sigma_0\right)$$ in such a way that~\eqref{upsimicionti} is verified. This completes the proof.
\end{proof}

The conclusive step in the theory developed in this final section  is to show that a conjecture in o--minimality (communicated to the authors by Raf Cluckers) enables one to make the calculation of the measure of flatness $q_{P}(\mathcal{K})$ completely effective in the following sense~: the lower limit of the ratio $(\log M_{P}\left(\mathcal{K}, \varepsilon\right))/(\log \varepsilon)$ defining it in~\eqref{conditionforallbis} is an actual limit which differs from this ratio up to an explicit error term as $\varepsilon\rightarrow 0^+$.\\

To this end, one needs to consider a relevant class of functions to analyse the pro\-per\-ties of the map $\mu_{P, \mathcal{K}}\left(\; \cdot\;,\; \cdot\;,\varepsilon\right)$ introduced in~\eqref{mapvarpesmu} (given a fixed value of $\varepsilon>0$) and then of the map $\varepsilon>0\mapsto M_{P}\left(\mathcal{K}, \varepsilon\right)$. This is the class of \emph{constructible functions} over a given globally subanalytic set $X\subset\R^k$ defined in~\cite[\S2]{CCMRS} as the ring of real--valued functions over $X$ generated by functions which either are globally subanalytic or  the logarithms of positive globally subanalytic functions. In other words, a function $f$ is constructible over $X$ if it can be expressed as 
\begin{equation*}
f(\bm{x})\;=\; \sum_{i=1}^{a}f_j(\bm{x})\cdot\prod_{j=1}^{b}\log f_{ij}(\bm{x})
\end{equation*}
for some integers $a,b\ge 0$, where the $f_i$'s and the $f_{ij}$'s are  subanalytic with $f_{ij}>0$ for all $i$ and all $j$. More generally, a map is said to be constructible if it is constructible over some subanalytic set.\\

The main property satisfied by constructible functions, which is of fundamental importance in the present considerations, is the stability under parametric integration established in~\cite[Theorem~1.3]{CM}. This is saying that for any globally subanalytic subsets  $A\subset\R^n$ and $B\subset\R^m$  and for any map $f~: A\times B\rightarrow \R$ which is constructible (over $A\times B$), the map $\bm{y}\in A\mapsto\int_{B}f(\bm{y}, \bm{x})\cdot \textrm{d}\bm{x}$ is also constructible (over $A$). As a consequence, the map $\mu_{P, \mathcal{K}}\left(\; \cdot\;,\; \cdot\;,\varepsilon\right)$ is constuctible as a function defined over a semialgebraic set as the parametric integral over a semialgebraic domain of a semialgebraic function. \\

The function $$M_P\left(\mathcal{K}, \;\cdot\; \right)\; :\; \varepsilon>0\mapsto M_{P}\left(\mathcal{K}, \varepsilon\right)\;=\; \max_{\bm{v}\in \Sph^{n-1}}\;\max_{\sigma\in\R}\;  \mu_{P, \mathcal{K}}\left(\bm{v},\sigma,\varepsilon\right)$$thus turns out to be the parametric maximum of a constructible function. The class of constructible functions is nevertheless not stable under taking parametric maxima~: an explicit counterexemple to this claim is worked out in~\cite{kaiser}. A weaker form of this stability property is nevertheless expected to hold~:

\begin{conj}[Raf Cluckers]\label{conjcluck}
Given a constructible, non--negative function $F~: (\varepsilon, \bm{x})\in X\;\mapsto\; F(\varepsilon, \bm{x})$ defined over some subanalytic set $X\subset (0, \infty)\times \R^k$, there exist a constant $\delta>0$ and a constructible function $G~: (0, \infty)\rightarrow \R$ of the form 
\begin{equation}
G(\varepsilon)\;=\; c\cdot \left|\log \varepsilon\right|^l\cdot \varepsilon^a 
\end{equation}
for some real $c\ge 0$, some rational number $a$ and some integer $l\ge 0$ such that $$\delta\cdot G(\varepsilon)\;\le\; \sup_{\bm{x}\in\R^k} F(\varepsilon, \bm{x})\;\le G(\varepsilon)\qquad \textrm{for all} \qquad \varepsilon <\delta.$$
\end{conj}

A more general form of this conjecture would assert that, under suitable assumptions, the parametric supremum of a constructible function defined over any subanalytic set is "sandwiched" between constant multiples of a constructible function depending on the remaining variables. This is inspired by analogy to the $p$-adic case where the corresponding result holds (even uniformly in $p$) --- see~\cite[Theorem~B]{RGHbis} and~\cite[Theorem~2.1.3]{RGH}.\\

If Conjecture~\ref{conjcluck} holds, keeping the notation therein, the function $M_P(\mathcal{K}, \;\cdot\;)$ satisfies the inequalities  $$\delta\cdot c\cdot \left|\log \varepsilon\right|^l\cdot \varepsilon^a\;\le\; M_P(\mathcal{K}, \varepsilon)\;\le\; c\cdot \left|\log \varepsilon\right|^l\cdot \varepsilon^a $$ when $0<\varepsilon<\delta$. As a consequence, for such values of $\varepsilon$, $$ \frac{\log \left(c \cdot l\right)}{\log \varepsilon}+\frac{\log\left|\log \varepsilon\right|}{\log \varepsilon}\;\le\; \left|\frac{\log M_P(\mathcal{K}, \varepsilon)}{\log \varepsilon}-a\right| \;\le\; \frac{\log \left(c \cdot l\cdot \delta\right)}{\log \varepsilon}+\frac{\log\left|\log \varepsilon\right|}{\log \varepsilon},$$ implying in particular that $$\frac{\log M_P(\mathcal{K}, \varepsilon)}{\log \varepsilon}\;=\; q_P(\mathcal{K})+O\left(\frac{\log\left|\log \varepsilon\right|}{\log \varepsilon}\right), \qquad \textrm{where}\qquad q_P(\mathcal{K})=a.$$With the help of Propositions~\ref{propvol} and~\ref{propcompact}, this relation shows that the measure of flatness $q_P(\mathcal{K})$ can be determined effectively in a finite number of steps assuming the validity of Conjecture~\ref{conjcluck}.

\newpage
$\quad$
\newpage
\newpage
\bibliographystyle{unsrt}

\addcontentsline{toc}{chapter}{\protect\numberline{}References}


\end{document}